\documentclass[12pt]{article}
\usepackage{geometry}  \usepackage[french,british]{babel}
\usepackage{graphicx}
\usepackage{amssymb,amsthm,amsmath}
\usepackage{epstopdf}
\usepackage{fourier}
\usepackage[T1]{fontenc}
\usepackage{enumerate}
\usepackage[colorlinks=true,backref=true]{hyperref}
\theoremstyle{plain}
\newtheorem*{theosansnombre} {Theorem}
\newtheorem*{corsansnombre} {Corollary}

\newtheorem{theo}{Theorem} [section]
\newtheorem{prop} [theo]{Proposition}
\newtheorem{coro} [theo]{Corollary}
\newtheorem{lemm} [theo]{Lemma}

\theoremstyle{definition}
\newtheorem{defi}[theo]{Definition}
\newtheorem*{rema}{Remark}

\begin{document}
\newcommand{\supp}{\textrm{supp}}
\newcommand{\Aff}{\textrm{Aff}}
\newcommand{\End}{\textrm{End}}


\begin{center}
{\Large\textbf{
Spectral gap properties for linear random walks and Pareto's asymptotics for affine stochastic recursions
}}

Y. Guivarc'h, \'E. Le Page
\end{center}

\tableofcontents

\newpage
\selectlanguage{british}
{\small
\section*{Abstract}

  Let $V=\mathbb R^d$ be the Euclidean $d$-dimensional space, $\mu$ (resp $\lambda$) a probability measure on the linear (resp affine) group $G=G L (V)$ (resp $H= \Aff (V))$ and assume that $\mu$ is the projection of $\lambda$ on $G$. We study asymptotic properties of the iterated convolutions $\mu^n *\delta_{v}$ (resp $\lambda^n*\delta_{v})$ if $v\in V$, i.e asymptotics of the random walk on $V$ defined by $\mu$ (resp $\lambda$), if the subsemigroup $T\subset G$ (resp.\ $\Sigma \subset H$) generated by the support of $\mu$ (resp $\lambda$) is ``large''. We show spectral gap properties for the convolution operator defined by $\mu$ on spaces of homogeneous functions of degree $s\geq 0$ on $V$, which satisfy H\"older type conditions. As a consequence of our analysis we get precise asymptotics for the potential kernel
  $\Sigma_{0}^{\infty} \mu^k * \delta_{v}$, which imply its asymptotic homogeneity. Under natural conditions the $H$-space $V$ is a $\lambda$-boundary; then we use the above results and radial Fourier Analysis on $V\setminus \{0\}$ to show that the unique $\lambda$-stationary measure $\rho$ on $V$ is "homogeneous at infinity" with respect to dilations $v\rightarrow t v$  (for $t>0$), with a tail measure depending essentially of $\mu$ and $\Sigma$. Our proofs are  based on the simplicity of the dominant Lyapunov exponent for certain products of Markov-dependent random matrices, on the use of   renewal theorems for ``tame''  Markov walks, and on the dynamical properties of a conditional $\lambda$-boundary dual to $V$.

\selectlanguage{french}
\section*{Résumé}
 Soit $V$ l'espace Euclidien de dimension $d$, $\mu$ (resp.\ $\lambda$) une probabilité sur le groupe linéaire
  (resp.affine) $G=GL(V)$ (resp. $H=\Aff(V)$) et supposons que $\mu$ soit la projection de $\lambda$ sur
  $G$. Nous étudions certaines propriétés asymptotiques des convolutions itérées de $\mu$ (resp.\ $\lambda$) appliquées à un vecteur non nul $v\in V$, c'est à dire de la marche aléatoire sur $V$ définie par $\mu$ (resp.\ $\lambda$), si le semigroupe $T\subset G$ (resp.\ $\Sigma\subset H$) engendré par le support de $\mu$ (resp.\ $\lambda$) est \og grand\fg. Nous montrons des propriétés d'isolation spectrale pour l'opérateur de convolution défini par $\mu$ sur des espaces de fonctions homogènes de degré $s\geq 0$ sur $V$, qui satisfont des conditions du type de H\"older. Comme conséquence de notre analyse nous obtenons des asymptotiques précises pour le noyau potentiel $\Sigma_{0}^{\infty} \mu^k * \delta_{v}$, qui impliquent son homogénéité à l'infini.Sous des conditions naturelles, le $H$-espace $V$ est une $\lambda$-frontière; nous utilisons alors  les résultats précédents et l'analyse de Fourier radiale sur $V\setminus\{0\}$ afin de montrer que l'unique mesure $\lambda$-stationnaire est homogène à l'infini, par rapport aux dilatations $v\rightarrow t v$ ( pour $t>0$),avec une mesure de queue qui dépend essentiellement de $\mu$ et $\Sigma$. Nos preuves sont basÈes sur la simplicité de l'exposant de Lyapunov dominant de certains produits de matrices en dépendance markovienne, sur l'utilisation de théorèmes de renouvellement pour certaines marches markoviennes  et sur les propriétés dynamiques d'une $\lambda$-frontière duale de $V$.
  }
\selectlanguage{british}

\noindent
\textbf{Key words and phrases:} spectral gap, renewal theorem, Pareto asymptotics, random matrices, affine random recursions.

\section{Introduction, statement of results}

 We consider the $d$-dimensional Euclidean space $V=\mathbb R^d$, endowed with the natural scalar product $(x,y)\rightarrow  \langle x,y\rangle$, the associated norm $x\rightarrow |x|$, the linear  group $G=GL(V)$, and the affine group $H=\textrm{Aff}(V)$.  Let $\lambda$ be a probability measure on $H$ with projection $\mu$ on $G$,  such that $\supp\lambda$  has no fixed point in $V$.  We denote  by
\[T=[\supp\mu],  (\textrm{resp.}\
 \Sigma=[\supp\lambda])\]
 the closed subsemigroup of $G$ (resp $H$) generated by $\supp\mu$ (resp $\supp\lambda$).   Under natural conditions, including negativity of the dominant Lyapunov exponent $L_{\mu}$ corresponding to $\mu$,  for any $v\in V$, the sequence of iterated convolutions $\lambda^n * \delta_{v}$ converges weakly to $\rho$; the probability measure $\rho$ is the unique probability which solves the convolution equation $\lambda *\rho=\rho$, and $\supp \rho$ is unbounded if $[\supp\mu]$ contains an expanding element. Then, an important property of $\rho$ is the existence of $\alpha>0$ such that $\int |x|^{s} d\rho (x) <\infty$ for $s<\alpha$ and $\int |x|^{s} d\rho (x)=\infty$ for $s\geq \alpha$, if $\supp\lambda$ is compact. One of our main results below (Theorem C) gives the $\alpha$-homogeneity of $\rho$ at infinity, i.e.\ Pareto's asymptotics of $\rho$ (see \cite{47}, p. 74).

  In general, for the asymptotic behaviour of $\lambda^{n}*\delta_{v}$ and the ``shape at infinity'' of $\rho$  there are four cases:

\begin{enumerate}
\item The ``contractive'' case where the elements of $\supp \mu$ have norms less than 1, $\rho$ exists and is compactly supported.
\item
The ``expansive" case where $L_{\mu}>0$ and $\rho$ does not exist.
\item The `critical'' case where $L_{\mu}=0$ and $\rho$ does not exist.
\item The ``weakly contractive'' case where $L_{\mu}<0$ and $\rho$ exists with unbounded support.
\end{enumerate}

  Heuristically, cases 3, 4, mentioned above, can be considered as  transitions between the  cases 1, 2, which appear to be extreme cases. In this paper we are mainly interested in case 4 and in the shape at infinity of $\rho$; in the corresponding analysis  we use the approach of \cite{35}, based on the associated linear random walk, we develop methods and prove results which are of independent interest for products of random matrices. An important tool here is the ``Radon transform'' of $\rho$, i.e. the function on $V$ defined by
\[\widehat{\rho} (v)=\rho \{ x\in V;  \langle x,v\rangle  >1\},\] which allow us to transfer the shape problem for $\rho$ into an asymptotic problem for $\widehat{\rho}$.

  In \cite{35} the shape problem was connected to the study of a Poisson equation on the $G$-space $V\setminus\{0\}$ satisfied by $\widehat{\rho}$ and by a convolution operator associated with  $\mu$; the measure $\lambda$ was assumed to be supported on the positive matrices or to have a density on $H$. In the first case an important result was the validity of Pareto's asymptotics for the projection of $\rho$ on the positive directions.We observe that special cases of the above problem and various consequences of Pareto's asymptotics have been considered in the litterature (see for example \cite{10}, \cite{20}, \cite{38}), especially if $\mu$ has a density on $G$, a condition which implies spectral gap properties for the convolution operators associated to $\mu$ or $\lambda$ in suitable Hilbert spaces.  In contrast,our basic hypothesis which involves only $T$ and $L_{\mu}$, implies that the above operators satisfy Doeblin-Fortet inequalities (see \cite{33}), hence also spectral gap properties in spaces of H\"older functions. Here our main result, partly contained in Theorem C below, describes the general case and the homogeneity at infinity   stated in the theorem gives new results even for $d=1$ (see \cite{20}) or for the multidimensional situations considered in \cite{35},\cite{38}.

    We observe that,  more generally, the homogeneous behaviour at infinity of certain invariant measures is of interest for various questions in Probability Theory and Mathematical Physics (see \cite{10}, \cite{11}, \cite{12},  \cite{13},  \cite{35}, \cite{47}) but also in some geometrical questions such as  dynamical excursions of geodesic flow  and winding around cusps in hyperbolic manifolds (see \cite{1}, \cite{44}, \cite{49}),  or analysis of the $H$-space $(V,\rho)$ as a $\lambda$-boundary and its dynamical consequences (see \cite{3}, \cite{16}, \cite{17}, \cite{30}).

  Hence, following \cite{35}, we start with the linear situation i.e. the $G$-action on $V\setminus\{0\}$ and we consider a probability measure $\mu$ on $G$. As in \cite{17} we assume that      $T$  satisfies the so-called i-p condition (i-p for irreducibility and proximality), i.e. $T$ is strongly irreductible and contains at least one element with a unique simple dominant eigenvalue; if $d=1$, we assume furthermore that $T$ is non arithmetic, i.e. $T$ is not contained in a subgroup of $\mathbb R^{*}$ of the form $\{\pm a^n; n \in \mathbb Z\}$ for some $a>0$. We observe that for $d>1$ condition i-p is satisfied by $T$ if and only if it is satisfied by    the algebraic subgroup $Zc (T)$,which is the Zariski closure of $T$, hence condition i-p  is satisfied if $T$ is ``large''  (see \cite{22}, \cite {45}). On the other hand, the set of probability measures $\mu$ on $G$ such that the associated semigroup $T$ satisfies condition i-p is open and dense in the weak topology hence $\mu$ is '``generic''  if $d>1$ and $T$ satisfies i-p. Also, if $d>1$, an essential aperiodicity consequence of condition i-p is the density in  $\mathbb R_{+}^{*}$ of the multiplicative  subgroup generated by moduli of dominant eigenvalues of the elements of $T$ (see \cite{2}, \cite{25},  \cite{28}). We denote  by $|g|$ the norm of $g\in G$ and we write
 \[\gamma(g)=\sup (|g|, |g^{-1}|), \; \; I_{\mu}=\{s\geq 0: \ \int |g|^s d\mu (g)<\infty\},\]
    we denote  $]0,s_{\infty}[$ the interior of the interval $I_{\mu}$. For simplicity of exposition, and since linear maps commute with the symmetry $v\rightarrow -v$, it is  convenient to deal with the $G$-factor space $\breve{V}$ of $V\setminus \{0\}$ by symmetry, instead of $V\setminus \{0\}$ itself. We use the polar decomposition
\[\breve{V}=\mathbb P^{d-1}\times \mathbb R_{+}^{*},\]
  and the corresponding functional decompositions, where $\mathbb P^{d-1}$ is the projective space of $V$.

   We consider the convolution action of $\mu$ on continuous functions on $V\setminus\{0\}$ which are homogeneous of degree $s\geq0$, i.e.\ functions $f$ which satisfy $f(t v)=|t|^s f(v)$ $(t\in \mathbb R)$. This action reduces to the action of a certain positive operator $P^s$ on $C(\mathbb P^{d-1})$, the space of continuous functions on the projective space $\mathbb P^{d-1}$. More precisely, if $f(v)=|v|^s \varphi(\bar v)$ with $\varphi \in C (\mathbb P^{d-1})$, $\bar v \in \mathbb P^{d-1}$, then $P^s \varphi$ is given by
\[P^s \varphi (x)=\int |g x|^s \varphi (g\cdot x) d\mu (g),\]
  where $x\in \mathbb P^{d-1}$,  $x\rightarrow  g \cdot x$ denotes the projective action of $g$ on $x$, and $|g x|$ is the norm of any vector $gv$ with $|v|=1$ and $\bar v=x$. Also for $z=s+it \in \mathbb C$, with $s\in I_{\mu}$ and $t\in\mathbb{R}$,  we write $P^z \varphi(x)=\int |g x|^z \varphi ( g \cdot x) d\mu(g)$. By duality $P^{z}$ acts also on measures on $\mathbb P^{d-1}$ and for a measure $\nu$ we denote by $P^{z} \nu$ the new measure obtained from $\nu$. The space of endomorphisms of a Banach space $B$ will be denoted by End$B$. For $\varepsilon>0$ let $H_{\varepsilon} (\mathbb P^{d-1})$ be the space of $\varepsilon$-H\"older functions on $\mathbb P^{d-1}$, with respect to a certain natural distance. We denote by $\ell^s$ (resp $\ell$) the $s$-homogeneous (resp.\ Haar) measure on $\mathbb R_{+}^{*}$  and we write
\[\ell^{s}(dt)=\frac{dt}{t^{s+1}}, \ \ \ell (dt)=\frac{dt}{t}, \ \ h^{s}(v)=|v|^{s}.\]
 An s-homogeneous Radon measure $\eta$ on $\breve{V}= \mathbb P^{d-1} \times \mathbb R_{+}^{*}$ is written as $\eta=\pi \otimes \ell^s$ where $\pi$ is a bounded measure on $\mathbb P^{d-1}$. For $s\in I_{\mu}$ we define the function
\[k(s)=\lim_{n\rightarrow \infty} \left(\int |g|^s d\mu^n (g)\right)^{1/n},\]
  where $\mu^n$ is the $n^{-th}$ convolution power of $\mu$ on the group $G$ and we observe that $\log k(s)$ is a convex function on $I_{\mu}$. A key tool in our analysis for $d>1$ is the

\begin{theosansnombre}{\textbf{\textrm{A}.}}
\label{ThA}
  Assume $d>1$ and the subsemigroup $T\subset GL (V)$ generated by $\supp\mu$ satisfies condition i-p. Then, for any $s\in I_{\mu}$ there exists a unique probability measure $\nu^s$ on $\mathbb P^{d-1}$, a unique positive continuous function $e^s\in C(\mathbb P^{d-1})$ with $\nu^s (e^s)=1$ such that
  \[P^s \nu^s=k(s) \nu^s, P^s e^s=k(s) e^s.\]

  For  $s\in I_{\mu}$, if $\int |g|^s \gamma^{\tau} (g) d\mu (g)< \infty$ for some $\tau>0$  and  if $\varepsilon>0$ is sufficiently small, the action of $P^s$ on $H_{\varepsilon}(\mathbb P^{d-1})$ has a spectral gap:
\[P^s=k(s) (\nu^s \otimes e^s+U^{s}),\]
  where the operator $\nu^s\otimes e^s$ is the projection  on $\mathbb C e^s$ defined by $\nu^s, e^s$ and $U^{s}$  is an operator with spectral radius less than 1 which satisfies  $U^{s}(\nu^{s}\otimes e^{s})=(\nu^s \otimes e^s) U^{s}=0$.
Furthermore  the function $k(s)$ is analytic,  strictly convex  on $]0,s_{\infty}[$ and the function
$\nu^s \otimes e^s$ from $]0,s_{\infty}[$ to $\mathop{{End}}H_{\varepsilon} (\mathbb P^{d-1})$ is analytic.
The spectral radius of $P^z$ is less than $k(s)$ if $s=\textrm{Re} z \in [0,s_{\infty}[$ and $t= \textrm{Im} z\neq 0$.
 \end{theosansnombre}

   We observe that, since condition i-p is open, the last property is robust under perturbation of $\mu$ in the weak topology. If $d=1$, $k(s)$  is equal to $\int |x|^{s} d\mu (x)$, hence $k(s)$ is the Mellin transform of $\mu$ (see \cite{52}) and the above statements are also valid if $T$ is non arithmetic. However, the last property is not robust for $d=1$.

  If $s=0$, $P^s$ reduces to the convolution operator by $\mu$ on $\mathbb P^{d-1}$ and convergence to the unique $\mu$-stationary measure $\nu^0=\nu$ was studied in \cite{17} using proximality of the $T$-action on $\mathbb P^{d-1}$.
 In this case,   spectral gap properties  for $P^{z}$, if Rez$=s$ is small, were first proved in  \cite{40} using the simplicity of the dominant $\mu$-Lyapunov exponent (see  \cite{28}). Limit theorems of Probability Theory for the product $S_{n}=g_{n}\cdots g_{1}$ of the random i.i.d.\ matrices $g_{k}$, distributed according to $\mu$, are  consequences of this result and of  radial Fourier analysis  on $V\setminus \{0\}$ used in combination with boundary theory (see\cite{4}, \cite{6}, \cite{16}, \cite{21}, \cite{29}, \cite{40}). If $\mu$ has a density with compact support, Theorem A is valid for any $s\in \mathbb R$. In general and for $d>1$, it turns out that the function $k(s)$, as defined above, looses its analyticity at some $s_{1}<0$.
For a recent detailed study of the operators $P^{z}$ ($s=$Rez small) and their equicontinuous extensions in a geometrical setting which allows the algebraic group $Zc (T)$ to be reductive and defined over a local field of any characteristic, we refer to the forthcoming book \cite{4}. We observe that condition i-p used here and $s$ small imply that the  basic assumptions of (\cite{4}, chap 8) are satisfied. Also, for $s>0$, the properties described in the theorem  were considered in \cite{29}; they are basic ingredients for the study of precise large deviations of $S_{n}(\omega) v$ (\cite{41}).

  The Radon measure $\nu^s \otimes \ell^s$ on $\breve{V}$ satisfies the convolution equation
  \[\mu * (\nu^s \otimes \ell^s)=k(s) \nu^s \otimes \ell^s\]
  and the support of $\nu^s$ is the unique $T$-minimal subset of $\mathbb P^{d-1}$, the so-called limit set $\Lambda (T)$ of $T$ (see \cite{2}, \cite{4}, \cite{23}).
The function $e^s$ is an integral transform of the twisted $\mu$-eigenmeasure $^{*}\nu^s$.  For $s>0$ and $\sigma$ a  probability measure on $\mathbb P^{d-1}$ not concentrated on a proper subspace, $|g|^s$ is comparable to $\int|g x|^s d\sigma (x)$; the uniqueness properties of $e^s$ and $\nu^s$ are based on this geometrical fact.
The proof  of the spectral gap property depends on the simplicity of the dominant Lyapunov exponent for the product of random matrices $S_{n}=g_{n}\cdots g_{1}$ with respect to a  natural  shift-invariant Markov measure $\mathbb Q^s$ on $\Omega=G^{\mathbb N}$, which is locally equivalent to the product measure $\mathbb Q^0=\mu^{\otimes \mathbb N}$. A construction of a kernel-valued martingale (based on $^{*}\nu^s$) plays an essential role in the  proof of simplicity and in the comparison of $|S_{n}(\omega)|$ with $|S_{n}(\omega) v|$. For $s=0$ this  study corresponds to  \cite{28}.

  In order to develop probabilistic consequences of Theorem A we endow $\Omega=G^{\mathbb N}$   with the shift-invariant measure $\mathbb P=\mu^{\otimes \mathbb N}$ (resp $\mathbb Q^s$).  We know that if $\int \log \gamma (g) d\mu (g)$ (resp $\int |g|^s \log \gamma (g) d\mu (g))$ is finite the dominant Lyapunov exponent $L_{\mu}$  (resp $L_{\mu} (s)$) of $S_{n}=g_{n}\cdots g_{1}$ with respect to $\mathbb P$  (resp $\mathbb Q^s$) exists and
\[L_{\mu}=\lim_{n\rightarrow \infty} \frac{1}{n} \int \log |S_{n} (\omega)| d \mathbb P(\omega), \ \ L_{\mu} (s)=\lim_{n\rightarrow \infty} \frac{1}{n} \int \log |S_{n}(\omega)| d\mathbb Q^s (\omega).\]

   If $s\in  ]0,s_{\infty}[$, $k(s)$ has a continuous  derivative $k'(s)$ and $L_{\mu}(s)=\frac{k'(s)}{k(s)}$.   By strict convexity of $\log k(s)$, if $\displaystyle\mathop{\lim}_{s\rightarrow s_{\infty}} k(s)\geq 1$ and $s_{\infty}>0$, we can define $\alpha>0$ by $k(\alpha)=1$.

   We consider the potential kernel  $U$ on $\breve{V}$ defined by $U(v,\cdot)=\displaystyle\mathop{\Sigma}^{\infty}_{0} \mu^k *\delta_{v}$. Then we have the following multidimensional extensions of the classical renewal theorems  (see \cite{15}), which describes the asymptotic homogeneity  of $U(v,\cdot)$. We recall that for $v\in \breve{V}$ and $A\subset \breve{V}$, $U(v,A)$ is the mean number of visits of $S_{n}(\omega)v$ to $A$ for $n\geq 0$.

\begin{theosansnombre}{\textbf{\textrm{B}}.}
  Assume $T$ satisfies condition \textrm i-p, $\int \log \gamma (g) d\mu (g)<\infty$ and $L_{\mu}>0$. If $d=1$ assume furthermore that $\mu$ is non-arithmetic.
Then, for any $v\in \breve{V}$, $U(v, \cdot)$ is a Radon measure on $\breve{V}$ and we have the vague convergence
\[\displaystyle\mathop{\lim}_{t\rightarrow 0_{+}} U (t v,\cdot)=\frac{1}{L_{\mu}} \nu \otimes \ell,\]
  where $\nu$ is the unique $\mu$-stationary measure on $\mathbb P^{d-1}$.
\end{theosansnombre}

\begin{theosansnombre}{\textbf{\textrm{B}$^\alpha$}.}
Assume $T$ satisfies condition \textrm i-p,  $\int \log \gamma (g) d\mu (g)<\infty$, $L_{\mu}<0$, $s_{\infty}> 0$ and there exists $\alpha>0$ with $k(\alpha)=1$, $\int |g|^{\alpha} \log \gamma (g) d\mu (g)<\infty$. If $d=1$ assume furthermore that $\mu$ is non-arithmetic.
 Then  for any $v\in \mathbb P^{d-1}$ we have the vague convergence on $\breve{V}$

\[\lim_{t\rightarrow 0_{+}} t^{-\alpha} U (t v,\cdot)=\frac{e^{\alpha}(v)}{L_{\mu}(\alpha)} \nu^{\alpha}\otimes \ell^{\alpha}.\]

  Up to normalization the Radon measure $\nu^{\alpha}\otimes \ell^{\alpha}$ is the unique $\alpha$-homogeneous measure which satisfies the harmonicity equation $\mu * (\nu^{\alpha}\otimes \ell^{\alpha})=\nu^{\alpha}\otimes \ell^{\alpha}$.
  \end{theosansnombre}

  For $v\in \mathbb P^{d-1}$, we consider the random variable $M(v)=\sup\{|S_{n}v|\ ;\ n\in \mathbb N\}$. Then we have the following matricial version of   Cram\'er's  estimate   in  collective risk  theory (see \cite{15}).

\begin{corsansnombre}{}
  With the notations of Theorem B$^{\alpha}$, for any $u\in \mathbb P^{d-1}$, we have the convergence
\[\lim_{t\rightarrow \infty} t^{\alpha} \mathbb P \{M (u)>t\}=A e^{\alpha} (u)>0.\]
\end{corsansnombre}

    Theorems B, B$^{\alpha}$ 	are  consequences of the arguments used in the proof of Theorem A and of a renewal theorem for a class of Markov walks on $\mathbb R$ (see \cite{36}). An essential role is played by the law of large numbers for $\log |S_{n} v|$ under $\mathbb Q^s (s=0, \alpha)$; the comparison of $|S_{n}|$ and $|S_{n} v|$ follows from  the finiteness of the limit of $(\log |S_{n}|- \log |S_{n}v|)$. This is the essential property used in \cite{36}, in a more general framework. For the sake of brevity, we have formulated these theorems in the context of $\breve{V}$ instead of $V$.
Corresponding statements where $\mathbb P^{d-1}$ is replaced by the unit sphere $\mathbb S^{d-1}$ are given in section 4. Also the above  weak convergence  can be extended to a larger class of functions.

  In \cite{35},  renewal theorems as above were obtained for non negative matrices, the extension of these results to the general case was an open problem  and a partial solution was given in \cite{39}. Theorems B and B$^{\alpha}$   extend these results to a wider setting. In view of the interpretation of $U(v,\cdot)$ as a mean number  of visits, Theorem B is a strong reinforcement of the law of large numbers for $S_{n} (\omega) v$, hence it can be used in some problems of dynamics for group actions on $T$-spaces. In this respect we observe that  a specific version of the asymptotic homogeneity of $U$  stated in Theorem B has been of essential use in \cite{30} for the description of the $T$-minimal subsets of the action of a large subsemigroup $T\subset S L(d,\mathbb Z)$ of automorphisms of the torus $\mathbb T^d$.

  On the other hand Theorem B$^{\alpha}$ gives a description of the fluctuations of a linear random walk on $\breve{V}$ with $\mathbb P$-\textrm{a.e.} exponential convergence to zero,  under condition i-p and  the existence in $T$ of a matrix with spectral radius greater than one. These fluctuation properties are responsible for the homogeneity  at infinity of stationary measures for affine random walks on $V$, that we discuss now.

  Let $\lambda$ be a probability measure on the affine group $H$ of $V$, $\mu$ its projection on $G$ and $T$, $\Sigma$ as above.  We assume that  $T$  satisfies condition i-p, and   $\supp\lambda$ has no fixed point in $V$. If $d=1$ we assume that $T$ is non-arithmetic. We consider the affine stochastic recursion on $V$
\[X_{n+1}=A_{n+1} X_{n}+B_{n+1}, \quad\quad\quad (R)\]
 where $(A_{n}, B_{n})$ are $\lambda$-distributed i.i.d random variables. From a heuristic point of view, the corresponding affine random walk can be considered as a superposition of an additive random walk on $V$ governed by $B_{n}$ and a multiplicative random walk on $V\setminus\{0\}$ governed by $A_{n}$. Here, as it appears below in Theorem C, the non trivial multiplicative part $A_{n}$ plays a dominant role, while the additive part $B_{n}$ has a stabilizing effect since $\lambda$ has  a unique stationary probability $\rho$ on $V$ and $\rho$ is  not supported on a point.  If $\mathbb E (\log |A_{n}|)+\mathbb E(\log |B_{n}|)<\infty$ and the dominant Lyapunov exponent $L_{\mu}$ for the product $A_{1}\cdots A_{n}$ is negative, then $R_{n}= \Sigma^{n-1}_{0} A_{1}\cdots A_{k} B_{k+1}$ converges $\lambda^{\otimes \mathbb N}-\textrm {a.e.}$ to $R$, the law $\rho$ of $R$  is  the unique $\lambda$-stationary measure on $V,(V,\rho)$ is a $\lambda$-boundary  (see \cite {17})and $\supp \rho=\Lambda_{a} (\Sigma)$ is the unique $\Sigma$-minimal subset of $V$. If $T$ contains at least one matrix with an expanding direction, then $\Lambda_{a} (\Sigma)$ is unbounded and if $I_{\mu}=[0,\infty[$  there exists $\alpha>0$ with $k(\alpha)=1$, hence we can inquire about the ``shape at infinity'' of $\rho$.  According to a conjecture of F. Spitzer   the measure  $\rho$ should belong to the domain of attraction of a stable law with index $\alpha$ if $\alpha \in [0,2[$ or a Gaussian law if $\alpha \geq 2$. Here we prove a multidimensional precise form of this conjecture, and more generally the $\alpha$-homogeneity at infinity of $\rho$, where we assume that $\mu$ satisfies the conditions of Theorem B$^{\alpha}$,  $\lambda$ satisfies moment conditions and $\supp\lambda$ has no fixed point in $V$. We denote by $\widetilde{\Lambda} (T)$ the inverse image of the limit set $\Lambda(T)$ in $\mathbb S^{d-1}$ and by $\widetilde{\nu}^{\alpha}$  the symmetric lifting of $\nu^{\alpha}$ to $\mathbb S^{d-1}$. Then our main result implies the following

\begin{theosansnombre}{\textbf{\textrm{C}}.}
  With the above notation we assume that $T$ satisfies condition \textrm i-p, $\supp \lambda$ has no fixed point in $V$, $s_{\infty}>0$, $L_{\mu}<0$ and $\alpha \in ]0,s_{\infty}[$ satisfies $k(\alpha)=1$. If $d=1$ we assume also $\mu$ is non-arithmetic. Then, if $\mathbb E (|B|^{\alpha+\tau})<\infty$ and $\mathbb E(|A|^{\alpha} \gamma^{\tau}(A))<\infty$ for some $\tau>0$, the unique $\lambda$-stationary measure $\rho$ on $V$ satisfies the following vague convergence on $V\setminus \{0\}$
\[\lim_{t\rightarrow 0_{+}} t^{-\alpha} (t\cdot\rho)=C \sigma^{\alpha}\otimes \ell^{\alpha},\]
where $C>0$, $\sigma^{\alpha}$ is a probability measure  on $\widetilde{\Lambda} (T)$ and $\sigma^{\alpha} \otimes \ell^{\alpha}$ is a $\mu$-harmonic Radon measure supported on $\mathbb R^{*} \widetilde{\Lambda} (T)$.
 If $T$ has no proper convex invariant cone in $V$, we have $\sigma^{\alpha}=\widetilde{\nu}^{\alpha}$.
The above convergence is also valid on any Borel function $f$ such that  the set of discontinuities of $f$ is $\sigma^{\alpha} \otimes \ell^{\alpha}$-negligible and such that for some $\varepsilon >0$ the function  $ |v|^{-\alpha} |\log |v||^{1+\varepsilon}|f(v)|$ is bounded.
\end{theosansnombre}

Briefly, we say that $\rho$ satisfies Pareto's asymptotics of index $\alpha$ (see \cite{47}, page 74).
The convergence in Theorem C can be considered as a Cram\'er type  estimate for the random variable R and  was stated in \cite{26}.This statement  gives the homogeneity at infinity of  $\rho$, hence the measure $C \sigma^{\alpha}\otimes \ell^{\alpha}$  defined by the theorem can be interpreted as the "tail measure" of $\rho$. In the context of extreme value theory for the process $X_{n}$, the convergence stated in the theorem implies that $\rho$ has ``multivariate regular variation'' and this property plays an essential role in the theory (see \cite{19}, \cite{47}). If $T$ has a proper convex invariant cone, then $C\sigma^{\alpha}$ can be decomposed as
\[C\sigma^{\alpha}=C_{+} \nu_{+}^{\alpha}+ C_{-} \nu_{-}^{\alpha},\]
 where $C_{+}$, $C_{-}\geq 0$ and $\nu_{+}^{\alpha} \otimes \ell^{\alpha}$, $\nu_{-}^{\alpha} \otimes \ell^{\alpha}$ are $\mu$-harmonic extremal measures on $V\setminus\{0\}$. In  section 5 the discussion of positivity for $C$, $C_{+}$, $C_{-}$ in terms of $\Lambda(T)$ and $\Lambda_{a}(\Sigma)$ lead us, via Radon transforms, to consider an associated linear random walk on the vector space $V\times \mathbb R$, which plays a dual role to the original $\lambda$-random walk $X_{n}$ on $V$.  The proof of positivity for $C$ depends on the use of Kac's recurrence theorem for this dual random walk. On the other hand, the proof of $\alpha$-homogeneity of $\rho$ at infinity follows from a Choquet-Deny type property for the linear $\mu$-random walk on $V\setminus\{0\}$.
 The spectral gap property, stated in theorem A, plays an essential role in this study.

   For  $d=1$, positivity of $C=C_{+}+C_{-}$ was proved in \cite{20} using Levy's symmetrisation
 argument,  positivity of $C_{+}$ and $C_{-}$ was  tackled in \cite{26} by a complex analytic method introduced in \cite{11}. For $d=1$ we have $\nu_{+}^{\alpha}=\delta_{1}$, $\nu_{-}^{\alpha}=\delta_{-1}$ and  the precise form of Theorem C gives that the condition $C_{+}=0$ is equivalent to $\supp\rho \subset ]-\infty,c]$ with $c\in \mathbb R$. In the Appendix we give an approach to part of Theorem C using tools familiar in Analytic Number Theory like Wiener-Ikehara's  theorem and a lemma of E. Landau but also results for Radon transforms of positive measures which are only valid for $\alpha \notin \mathbb N$ (see \cite{5}, \cite{51}). However the discussion of positivity for $C_{+}$, $C_{-}$ seems to be not possible using only these analytical tools.

  A natural question is the speed of convergence in Theorem C. For $d=1$ see \cite{20}, if $\lambda$ has a density. For $d>1$ and under condition i-p, this question is connected with the possible uniform spectral gap for the operator $P^{z}$ of Theorem A,if $z=\alpha+it$.

   To go further we observe that Theorem C  gives a natural construction for a large class of probability measures in the domain of attraction of a stable law. Using also spectral gaps and  weak dependence properties of the process $X_{n}$, Theorem C allow us   to prove  convergence to stable laws for  normalized Birkhoff sums along the affine  $\lambda$-walk  on $V$ (see  \cite{18}) Furthermore if $d=1$, and   conditionally on regularity assumptions usual in extreme value theory, such convergences were shown in \cite{10}; hence the results of \cite{18} improve and extends the results of \cite{10} in the case of GARCH processes $(d\geq1)$. If $d>1$  the above convergence is robust under perturbation of $\lambda$ in the  weak topology. These  convergences to stable laws are connected with the study of random  walk in a random medium on the line or the strip (see \cite{12}) if $\alpha<2$.  On the other hand the study of the extremal value behaviour of the process $X_{n}$ can be fully developed on the basis of Theorem C and on the above weak dependence properties of $X_{n}$; in particular the asymptotics of the extremes of   $|X_{n}|$ are given by Fr\'echet type laws with index $\alpha$ (\cite{31}), a result which extends the main result of \cite{38} to the generic case. In a geometrical context,  as observed in \cite{44} for excursions of  geodesic flow around the cusps of the modular surface, the famous Sullivan's logarithm law is a  simple consequence of Fr\'echet's law for the continuous fraction expansion of a real number uniformly distributed in $[0,1]$. Here also a logarithm law is valid for the random walk $X_{n}$ (see \cite{31}). The arguments developed in the proof of homogeneity at infinity for $\rho$ can also be used in the study of certain quasi-linear  equations which occur in  various domains related to branching random walks in particular (see \cite{20}), \cite{42}). For example the description of the shape at infinity of the fixed points of the multidimensional version of the ``smoothing transformation'' considered in \cite{13} in the context of infinite particle systems in interaction depends on such arguments (see \cite{9}). In an econometrical context, the stochastic recursion $(R)$ can be interpreted as a mechanism which, in the long run, produces debt or wealth accumulation with a specific homogeneous structure  at large values; hence,in the natural setting of affine stochastic recursions, this mechanism ``explains''  the remarkable power law asymptotic shape of wealth distribution empirically discovered by the economist V. Pareto (\cite{43}).

  For  information on the role of spectral gap properties in limit theorems for Probability theory and Ergodic theory we refer to \cite{1}, \cite{4}, \cite{6}, \cite{18}, \cite{21},  \cite{27}, \cite{28}, \cite{40}. For information on products of random matrices we refer to \cite{4}, \cite{6},  \cite{16}, \cite{24}. Theorem A (resp B, B$^{\alpha}$ and C) is proved in sections 2,3 (resp 4 and 5).

  We thank Ch.M. Goldie,  I. Melbourne and D. Petritis for useful informations on stochastic recursions and extreme value theory.We thanks also the referees for careful reading and very useful sugestions.

\section{Ergodic properties of transfer operators on projective spaces}

  In this section we study the qualitative properties of transfer operators on $\mathbb P^{d-1}$ or $\mathbb S^{d-1}$. As mentioned in the introduction one can find in \cite{4}  a detailed study of a general class of transfer operators on flag manifolds with equicontinuity properties. Here, our transfer operators depend on a complex parameter $z$ which is typically large with Rez$>0$. Hence, for self-containment reasons in particular, we develop from scratch our study on $\mathbb P^{d-1}$ or $\mathbb S^{d-1}$ for Rez$=s\geq 0$. A first step is to reduce these transfer operators to Markov operators (Theorems 2.6, \ref{thm:2.17}) with equicontinuity properties.
\subsection{Notation and preliminary results}

  Let  $V=\mathbb R^d$ be the Euclidean space endowed with the scalar product $\langle x , y \rangle=\sum_{1}^{d} x_{i} y_{i}$ and the norm $|x|=\left(\sum_{1}^{d}|x_{i}|^{2}\right)^{1/2}$, $\breve{V}$ the factor space of $V\setminus\{0\}$ by the finite group $\{\pm Id\}$. We denote by $\mathbb P^{d-1}$ (resp $\mathbb  S^{d-1})$, the projective space (resp unit sphere) of $V$ and by $\bar v$ (resp $\tilde v$) the projection of $v\in V$ on $\mathbb P^{d-1}$ (resp $\mathbb  S^{d-1}$). The linear group $G=GL(V)$ acts on $V$, $\breve{V}$ by $(g,v) \rightarrow g v$. If $v\in V\setminus \{0\}$, we write $g\cdot
  v=\frac{gv}{|gv|}$ and we observe that $G$ acts on $\mathbb S^{d-1}$ by $(g,x) \rightarrow  g \cdot x$.  We will also write the action of $g\in G$ on $x\in \mathbb P^{d-1}$ by $ g \cdot x$; we define $|g x|$ as $|g \tilde x|$ if $x\in \mathbb P^{d-1}$ and $\tilde x \in\mathbb S^{d-1}$ has projection $x\in \mathbb P^{d-1}$. Also, if $x,y \in \mathbb P^{d-1}, |\langle x , y \rangle|$ is defined as $| \langle \widetilde x, \widetilde y\rangle |$ where $\widetilde x, \widetilde y\in \mathbb S^{d-1}$ have projections $x,y$. Corresponding  notations will be  taken when convenient.  For a subset $A\subset\mathbb S^{d-1}$ the convex envelope $\textrm{Co}(A)$ of $A$ is defined as the intersection with $\mathbb S^{d-1}$ of the closed convex cone generated by $A$ in $V$. We denote by $O(V)$ the orthogonal group of $V$ and by $m$ the $O(V)$-invariant measure on $\mathbb P^{d-1}$. A positive measure $\eta$ on $\mathbb P^{d-1}$ will be said to be proper if $\eta(U)=0$ for every proper projective subspace $U\neq \mathbb P^ {d-1}$.We denote by $\End V$ the space of endomorphisms of the vector space $V$.

  Let $P$ be a positive kernel on a Polish space $E$ and let $e$ be a  positive function on $E$ which satisfies $Pe=ke$ for some $k>0$. Then we can define a Markov kernel $Q_{e}$ on $E$ by Doob's relativisation procedure: $Q_{e} \varphi=\frac{1}{ke} P(\varphi e)$. This procedure will be used frequently here.  For a Polish $G$-space $E$ we denote by $M^1(E)$ the space of probability measures on $E$. If $\nu \in$ $M^1 (E)$, and $P$ is as above, $\nu$ will  be said to be $P$-stationary if $P\nu=\nu$, i.e for any Borel function $\varphi$, $ \nu (P\varphi )=\nu (\varphi)$. We will write $C(E)$ (resp $C_{b} (E)$ for the space of continuous (resp bounded continuous) functions on $E$. If $E$ is a locally compact $G$-space, $\mu \in$ $M^1 (G)$, and $\rho$ is a Radon measure on $E$, we recall that the convolution $\mu * \rho$ is defined as a Radon measure by $\mu * \rho=\int \delta_{g x} d\mu (g) d\rho (x)$, where $\delta_{y}$ is the Dirac measure at $y\in E$.  A $\mu$-stationary measure on $E$ will be a probability measure $\rho \in$ $M^1 (E)$ such that $\mu * \rho=\rho$. In particular, if $E=V$ or $\breve{V}$ and $\mu \in$ $M^1 (G)$ we will consider the Markov kernel $P$ on $V$ (resp $\breve{P}$ on $\breve{V}$) defined by $P(v, \cdot)=\mu * \delta_{v}$, (resp $\breve{P} (v,\cdot)=\mu * \delta_{v})$. On $\mathbb P^{d-1}$ (resp $\mathbb  S^{d-1})$ we will write $\bar P (x,\cdot)=\mu * \delta_{x}$ (resp $\tilde P(x,\cdot)=\mu * \delta_{x})$.

  If u is an endomorphism of $V$, we denote $u^*$
   its adjoint map, i.e.\ $\langle u^*x, y\rangle =\langle x, u y\rangle $ if $x, y \in V$.  If $\mu \in M^1 (G)$ we will write $\mu^*$ for its push forward by the map  $g\rightarrow g^*$ and we define the kernel ${^{*} P}$ on $V$ by ${^{*} P(v,.)}=\mu^{*}* \delta_{v}$.
 For $s\geq 0$ we denote $\ell^s$ (resp $h^s$) the $s$-homogeneous measure (resp function) on $\mathbb R_{+}^*=\{t\in \mathbb R \ ; \ t>0\}$ given by
 \[\ell^s (dt)=\frac{dt}{t^{s+1}}\  (\hbox{\rm resp}\  h^s(t)=t^s).\]
  For $s=0$ we write $\ell (dt)=\frac{dt}{t}$.

  Using the polar decomposition $V\setminus\{0\}=\mathbb S^{d-1} \times \mathbb R_{+}^*$ and the corresponding functional decompositions on $V\setminus\{0\}$, every $s$-homogeneous measure $\eta$ (resp function $\psi$) on $V\setminus \{0\}$ can be written as
\[\eta=\pi \otimes \ell^s \qquad (\textrm{resp.} \ \psi=\varphi \otimes h^s),\]
where $\pi$ (resp.\ $\varphi$) is a measure (resp.\  function) on $\mathbb S^{d-1}$. Similar decompositions are valid on $\breve{V}=\mathbb P^{d-1}\times \mathbb R_{+}^{*}$. If $g\in G$,  and $\eta=\pi \otimes \ell^{s}$ (resp $\psi=\varphi \otimes h^s)$ the directional component of $g\eta$ (resp.\ $\psi \circ g$) is given by
\[\rho^s(g) (\eta)=\int |g x|^s \delta_{ g \cdot x} d\eta (x),  \qquad   \rho_{s} (g) (\psi) (x)=|g x|^s \psi ( g \cdot x).\]
  The representations $\rho^s$ and $\rho_{s}$ extend to measures on $G$ by the formulae
\[\rho^{s} (\mu) (\eta)=\int |gx|^{s} \delta_{ g \cdot x} d\mu (g) d\pi(x),  \qquad
\rho_{s}(\mu) (\psi) (x)=\int |g x|^{s} \psi ( g \cdot x) d\mu (g).\]

  We will write, for $\varphi \in C (\mathbb P^{d-1})$ (resp.\  $\psi \in C(\mathbb S^{d-1})$)
 \[P^s \varphi=\rho_{s}(\mu) (\varphi), \ \ (\textrm{resp.} \tilde P^s \psi=\rho_{s}(\mu) (\psi)),\ \  ^{*}P^s \varphi=\rho_{s} (\mu^*)
(\varphi),\ \ (\textrm{resp.} {^{*}\tilde{P}}^s \psi=\rho_{s}(\mu^*) (\psi)).\]
  We endow $\mathbb S^{d-1}$ (resp.\ $\mathbb P^{d-1}$) with the distance $\tilde{\delta}$ (resp.\ $\delta$) defined by
\[\tilde{\delta}(x,y)=|x-y|\  (\hbox{\rm resp.\ }\  \delta (\bar x, \bar y)= \inf \{|x-y| ; |x|=|y|=1\}).\]

  For $\varepsilon>0$, $\varphi\in C(\mathbb P^{d-1})$ (resp.\  $\psi \in C (\mathbb S^{d-1}))$, we denote
\begin{align*}
[\varphi]_{\varepsilon}=\displaystyle\sup_{x\neq y} \frac{|\varphi(x)-\varphi(y)|}{\delta^{\varepsilon}(x,y)}\qquad  & \mbox{}(\textrm{resp.}\  [\psi]_{\varepsilon}=\displaystyle\sup_{x\neq y} \frac{|\psi (x)-\psi (y)|}{\tilde{\delta}^{\varepsilon}(x,y)}),\\
|\varphi|=\displaystyle\sup\{|\varphi (x)| ; x\in \mathbb P^{d-1}\} \qquad &\mbox{}
 (\textrm{resp.}\  |\psi |=\displaystyle\sup \{|\psi (x)| ; x\in \mathbb S^{d-1}\}),
 \end{align*}
 and we write
 \[H_{\varepsilon}(\mathbb P^{d-1})=\{\varphi \in C(\mathbb P^{d-1}) ; [\varphi]_{\varepsilon}<\infty\},  (\textrm{resp.}\ H_{\varepsilon}(\mathbb S^{d-1})=\{\psi \in C(\mathbb S^{d-1}) ; [\psi]_{\varepsilon}<\infty\}.\]

   The set of positive integers will be denoted by $\mathbb N$. We denote by $\mu^{n}$ the $n^{th}$ convolution power of $\mu$, i.e for $\psi \in C_{b}(G)$, $\mu^{n}(\psi)=\int \psi (g_{n}\cdots g_{1}) d\mu^{\otimes n}(g_{1},\cdots, g_{n})$.

\begin{defi}{}
\label{D2.1}
  If $s\in [0, \infty[$ and $\mu\in M^{1} (G)$, we denote
\[k(s)=k_{\mu}(s)=\lim_{n\rightarrow \infty} \left(\int |g|^{s} d\mu^{n} (g)\right)^{1/n},
\ \ I_{\mu}= \{s\geq 0 ; k_{\mu}(s) <+\infty\}.\]
\end{defi}

  We observe that the above limit exists, since by subadditivity of $g\rightarrow \log |g|$, the quantity $u_{n}(s)=\int |g|^{s} d\mu^{n} (g)$ satisfies $u_{m+n}(s) \leq u_{m}(s) u_{n}(s)$. Also $k_{\mu}(s)=\inf_{n\in \mathbb N} (u_{n}(s))^{1/n}$, which implies $I_{\mu}=\{s\geq 0 ; \int |g|^{s} d\mu (g) < +\infty\}$. Furthermore, H\"older inequality implies that $I_{\mu}$ is an  interval of the form $[0, s_{\infty}[$ or $[0, s_{\infty}]$, and $\log k_{\mu}(s)$ is  convex on $I_{\mu}$. Also $k_{\mu^{*}}=k_{\mu}$ since $|g|=|g^{*}|$. If $\mu$ and $\mu'$ commute and c$\in [0,1]$, $\mu''= \mu+(1-$c) $\mu'$ then $k_{\mu''}(s)=ck_{\mu}(s) +(1-$c) $k_{\mu'}(s)$, if $s\in I_{\mu}\cap I_{\mu'}$.
\begin{defi}{}
\label{D2.2}
\begin{enumerate}
\item
An element $g\in \textrm{End} V$ is said to be proximal if $g$ has a unique  eigenvalue $\lambda_{g} \in \mathbb R$ of maximum modulus and $\lambda_{g}$ is simple.
\item A semigroup $T \subset G$ is said to be strongly irreducible if no finite union of proper subspaces is $T$-invariant.
\end{enumerate}
\end{defi}
  Proximality of $g$ means that we can write $V=\mathbb R v_{g}\oplus V_{g}^{<}$ with $g v_{g}=\lambda_{g} v_{g}, \ \  g V_{g}^{<}\subset V_{g}^{<}$ and the restriction of $g$ to $V_{g}^{<}$ has spectral radius less than $|\lambda_{g}|$.  In this case $\displaystyle\mathop{\lim}_{n\rightarrow+\infty} g^{n}\cdot \bar x=\bar v_{g}$ if $x \notin V_{g}^{<}$ and we say that $\lambda_{g}$ is the dominant eigenvalue of $g$. If $E\subset G$ we denote by $E^{prox}$ the set of proximal elements of $E$.The closed subsemigroup (resp.\ group) generated by $E$ will be denoted $[E]$ (resp.\ $\langle E\rangle $). In particular we will consider below the case $E=\supp \mu$ where $\supp \mu$ is the support of $\mu \in M^1 (G)$.

\begin{defi}{}
\label{D2.3}
  A semigroup $T \subset G$ is said to satisfy condition i-p if $T$ is strongly irreducible and $T^{\textrm{prox}}\neq \emptyset$.
\end{defi}
  As shown in (\cite{22}, \cite{45}) this property of $T$ is satisfied  if it is satisfied by $Z c (T)$, the  Zariski closure of  $T$. It can be proved that condition i-p is valid if and only if the connected component of the closed subgroup $Z c (T)$ is locally  the product of a similarity group and a semi-simple real Lie group without compact factor which acts proximally and irreducibly on $\mathbb P^{d-1}$. In this sense $T$ is ``large''. For example, if $T$ is a countable subgroup of $G$ which satisfies condition i-p then $T$ contains a free subgroup with two generators.

  We recall that in $\mathbb C^{n}$, the Zariski closure of $E\subset \mathbb C^{n}$ is the set of zeros of the set of polynomials which  vanish on $E$. The group $G=GL (V)$ can be considered as a Zariski-closed subset of $\mathbb R^{d^2+1}$. If $T$ is a semigroup, then $Zc (T)$ is a closed subgroup of $G$ with a finite number of connected components. If $d=1$, condition i-p is always satisfied. Hence, when using condition i-p, $d>1$ will be understood.
\begin{rema}{}
The above definitions will be used below in the analysis of laws of large numbers and renewal theorems. A corresponding analysis has been developed in \cite{35} for the case of non-negative matrices. We observe that proximality of an element in $G$  is closely related to the Perron-Frobenius property for a positive matrix. If $T$ is not irreducible, we can consider the subspace $V^{+}(T)$ generated by the dominant eigenvectors of the elements of $T^{ \textrm{prox} }$; then $V^{+}(T)$ is $T$-invariant (see \cite{21} p 120) and, if $T^{+}$  is the  restriction of $T$ to $V^{+} (T)$, then $T^{+}$ satisfies condition i-p. In that way our results below could be used in reducible situations (see for example \cite{9}).
\end{rema}
\begin{defi}{}
\label{D2.4}
  Assume $T$ is a subsemigroup of $G$ which  satisfies condition i-p. Then the closure of the set $\{\bar v_{g} ; g\in T^{\textrm{prox}}\}$ will be called the limit set of $T$ and will be denoted $\Lambda(T)$.
\end{defi}
We recall that for a semigroup $T$ acting on a topological space $E$, a subset $X\subset E$ is said to be $T$-minimal if any orbit $T y (y\in X)$ is contained in $X$ and dense in $X$. For the minimality of $\Lambda(T)\subset \mathbb P^{d-1}$, see \cite{2}, \cite{23}.

With these definitions we have the
\begin{prop}
\label{P2.5}
Assume $T \subset G$ is a subsemigroup which satisfies condition i-p and $S\subset T$ generates $T$. Then $T \Lambda(T)=\Lambda(T)$ and $\Lambda (T)$ is the unique $T$-minimal subset of $\mathbb P^{d-1}$. If $\mu \in M^{1} (G)$ is such that $T=[\supp \mu]$ satisfies \textrm {i-p}, there exists a unique $\mu$-stationary measure $\nu$ on $\mathbb P^{d-1}$. Also $\supp \nu=\Lambda (T)$ and $\nu$ is proper. Furthermore, if $d>1$,  the subgroup of $\mathbb R_{+}^{*}$ generated by the set $\{|\lambda_{g}| ; g\in T^{\textrm{prox}} \}$ is dense  in $\mathbb R_{+}^{*}$.
 In particular, if $\varphi \in C(\Lambda (T))$ satisfies for some $t\in \mathbb R$, $|e^{i\theta}|=1$ : $\varphi ( g \cdot x)\  |g x|^{it}=e^{i\theta} \varphi (x)$
for any $g\in S$, $x\in \Lambda (T)$
 then $t=0$, $e^{i\theta}=1$, $\varphi=$ constant.
 \end{prop}

\begin{rema}
  The first part of the above statement is essentially due to H. Furstenberg (\cite{17}, Propositions \ref{pro:4.8}, 7.4) . The second part is proved in (\cite{28}. Proposition 3). For another proof and extensions of this property see \cite{2}, \cite{25}. This property  plays an essential role in the renewal theorems of section 4 as well as in section 5 for $d>1$. In the context of non negative matrices a modified form is also valid (see \cite{9}); in (\cite{35}), Theorem A) under weaker conditions on $T$, its conclusion is assumed as an hypothesis.  If $d=1$, we will need to assume it, i.e we will assume that $T$ is non-arithmetic ; if $T=[\supp\mu]$ satisfies this condition, we say that $\mu$ is non-arithmetic.
  \end{rema}

\subsection{Uniqueness of eigenfunctions and eigenmeasures on $\mathbb P^{d-1}$}

  Here we consider $s\in I_{\mu}$ and the operator $P^{s}$ (resp ${{^{*}P^{s}}}$) on $C(\mathbb P^{d-1})$ defined by
\[P^{s} \varphi (x)=\int |g x|^{s} \varphi ( g \cdot x) d\mu (g), \qquad (\textrm{resp.}\ {^{*}P^{s}} \varphi (x)=\int |g x|^{s} \varphi ( g \cdot x) d\mu^* (g)).\]

  For a measure $\nu$ on $\mathbb P^{d-1}$, $P^{s}\nu$ is defined by duality against $C(\mathbb P^{d-1})$.
For $z=s+it \in \mathbb C$ we will also write $P^z \varphi (x)=\int |g x|^z \varphi ( g \cdot x) d\mu (g)$.
Below  we study existence and uniqueness for eigenfunctions or eigenmeasures of $P^{s}$. We show equicontinuity properties of the normalized iterates of $P^{s}$ and $\tilde P^{s}$.
\begin{theo}
\label{T2.6}
Assume $\mu \in M^{1} (G)$ is such that the semigroup $[\supp \mu]$ satisfies (i-p) and let  $s\in I_{\mu}$. If $d=1$ we assume that $\mu$ is non arithmetic. Then the equation $P^{s} \varphi=k (s) \varphi$ has a unique continuous  solution $\varphi=e^{s}$, up to normalization. The function $e^{s}$ is positive and $\bar s$ -H\"older with $\bar s=\inf (1,s)$.

  Furthermore there exists a unique $\nu^{s} \in M^{1} (\mathbb P^{d-1})$ such that $P^{s} \nu^{s}$ is proportional  to $\nu^{s}$. One has $P^{s}\nu^{s}=k(s) \nu^{s}$ and $\supp \nu^{s}=\Lambda ([\supp \mu])$. If ${^{*}\nu^s} \in M^{1} (\mathbb P^{d-1})$ satisfies ${^{*}P^{s}} ({^{*}\nu^s)}=k(s) ^{*}\nu^s$, and $e^s$ is normalized by $\nu^s (e^s)=1$, one has $ p(s) e^s (x)=\int |\langle x , y \rangle|^s \ d {^{*}\nu^s} (y)$ where $p(s)=\int |\langle x , y \rangle|^s d\nu^s (x) d {^{*}\nu^s} (y)$.

  The map $s\rightarrow \nu^{s}$ (resp $s\rightarrow e^{s})$ is continuous in the weak topology (resp uniform topology) and the function $s\rightarrow \log k (s)$ is strictly convex.

  The Markov operator $Q^{s}$ on $\mathbb P^{d-1}$ defined by $Q^{s}\varphi=\frac{1}{k(s) e^{s}} P^{s} (\varphi e^{s})$ has a unique stationary measure $\pi^{s}$ given by $\pi^{s}=e^{s} \nu^{s}$ and we have for any $\varphi \in C(\mathbb P^{d-1})$  the uniform convergence of $(Q^s)^n \varphi$ towards $\pi^s(\varphi)$.  If $z=s+it$, $t\in \mathbb R$ and $Q^z$ is defined by $Q^z \varphi=\frac{1}{k(s) e^s} P^z (e^s \varphi)$,then the equation $Q^z \varphi=e^{i\theta} \varphi$ with $\varphi \in C (\mathbb P^{d-1})$, $\varphi\neq 0$ implies $e^{i\theta}=1$, $t=0, \varphi=$ constant.
    \end{theo}

\begin{rema}{}
\begin{enumerate}
\item If $s=0$, then $e^{s}=1$, and $\nu^s=\nu$ is the unique $\mu$-stationary measure \cite {17}.
 The fact that $\nu$ is proper is of essential  use in \cite{40}, \cite{6} and \cite{28}, for the study of limit theorems.
\item
 In section 3  we will also construct a suitable kernel-valued martingale which allows to prove that $\nu^{s}$ is proper (see Theorem 3.2), if $s\in I_{\mu}$. We  note that analyticity of $k(s)$ is proved in Corollary 3.20 below. Continuity of the derivative of $k$ will be essential  in sections 3,4 and is proved in Theorem 3.10. If $s=0$ the corresponding martingale construction was done in \cite{17}.
 \item
 By definition of $P^{s}$ and $\breve{P}$, the function (resp.\ measure)  $e^s \otimes h^s$ (resp.\ $\nu^s\otimes \ell^s$) satisfies the equation
 \[\breve{P}(e^s \otimes h^s)=k(s)(e^s \otimes h^s), (\text{resp. } \breve{P}(\nu^s \otimes \ell^s)=k(s)\nu^s \otimes\ell^s)\]

 \end{enumerate}
 \end{rema}
The proof of the theorem depends of a proposition and the following lemmas  improving corresponding results for positive matrices in \cite{35}.

\begin{lemm}
\label{L2.7}
Assume $\sigma \in M^{1}(\mathbb P^{d-1})$ is not supported  by a hyperplane. Then, there exists a constant $c_{s}(\sigma)>0$ such that, for any $u$ in  $\End V$

\[ \int |u x|^{s} d\sigma (x) \geq c_{s} (\sigma) |u|^s.\]
\end{lemm}

\begin{proof}
  Clearly it suffices to show the above inequality if $|u|=1$. The fonction $u \rightarrow \int |u x|^{s} d\sigma(x)$ is continuous on $\textrm{End} V$, hence its attains its infimum $c_{s}(\sigma)$ on the compact subset of $\textrm{End} V$ defined  by $|u|=1$. If $c_{s} (\sigma)=0$, then for some $u\in \textrm{End} V$ with $|u|=1$, we have $\int |u x|^{s} d\sigma (x)=0$ hence,  $u x=0$, $\sigma - a.e$. In other words,  $\supp \sigma \subset \textrm{Ker}(u)$, which contradicts the hypothesis on $\sigma$. Hence $c_{s}(\sigma) >0$.
  \end{proof}
\begin{lemm}
\label{L2.8}
  If $s\in I_{\mu}$ there exists $\sigma \in M^{1} (\mathbb P^{d-1})$ such that $P^{s} \sigma =k\sigma$ for some $k>0$. For any such $\sigma$, we have $k=k(s)$ and $\sigma$ is not supported on a hyperplane.
Furthermore for every $n\in \mathbb N$
\[\int |g|^{s} d\mu^{n} (g) \geq k^{n} (s) \geq c_{s} (\sigma) \int |g|^{s} d\mu^{n}(g) .\]
\end{lemm}
\begin{proof}
  We consider the non-linear operator $\widehat{P}^{s}$  on $M^{1}(\mathbb P^{d-1})$ defined by $\widehat{P}^{s}\sigma=\frac{P^{s} \sigma}{(P^{s}\sigma) (1)}$.
Since $\int |g|^{s}d\mu (g) <+\infty$, this operator is continuous in the weak topology. Since $M^{1} (\mathbb P^{d-1})$ is compact and convex, Schauder-Tychonov theorem implies the existence of $k>0$ and $\sigma \in M^{1} (\mathbb P^{d-1})$ with $P^{s} \sigma=k\sigma$, hence $k=(P^{s}\sigma) (1)$.
For  such a $\sigma$, the equation
\[k \sigma (\varphi) = \int \varphi (g\cdot x) |g x|^{s} d\mu (g) d\sigma (x)\]
  implies that if $x\in \supp \sigma$, then $ g \cdot x \in \supp \sigma$, $\mu$\textrm{-a.e}.

  Then for any $g\in \supp \mu$ we have  $g\cdot \supp \sigma \subset \supp \sigma$. In particular the projective subspace $H$ generated by $\supp \sigma$ satisfies $[\supp \mu] \cdot H=H$. Since $[\supp \mu]$ satisfies i-p, we have $H=\mathbb P^{d-1}$. Then Lemma 2.7 gives,  for any $g$,
\[\int |g x|^{s} d\sigma (x) \geq c_{s}(\sigma) |g|^{s}.\]
  The relation $(P^{s})^{n}\sigma=k^{n} \sigma$ implies $k^{n}= \int |g x|^{s} d\mu^{n} (g) d\sigma(x)$; hence, using  Lemma 2.7 we get  $c_{s} (\sigma) \int |g|^{s} d\mu^{n}(g) \leq k^{n} \leq \int |g|^{s}d\mu^{n}(g)$. It follows that
  \[k=\displaystyle\mathop{\lim}_{n\rightarrow +\infty} \left(\int |g|^{s}d\mu^{n} (g)\right)^{1/n}=k(s).\]
  \end{proof}

  Assume $e \in C (\mathbb P^{d-1})$ is positive and satisfies $P^{s} e=k(s) e$. Then we can define the Markov kernel $Q_{e}^{s}$ and the cocycle $\theta_{e}^{s}$ by
\[Q_{e}^{s} \varphi (x) = \frac{1}{k(s)} \int \varphi ( g \cdot x) \frac{e( g \cdot x)}{e(x)} |g x|^{s} d\mu (g),\ \ \  \theta_{e}^{s}(x,g)=|gx|^{s} \frac{e( g \cdot x)}{e(x)}.\]
 In view of the cocycle property of $\theta^{s}_{e} (x,g)$ we can calculate the iterate $(Q_{e}^{s})^{n}$ by the formula $(Q_{e}^{s})^{n} \varphi (x) =\int \varphi ( g \cdot x) q^{s}_{e,n} (x,g) d\mu^{n} (g)$
 with $q_{e,n}^{s} (x,g)=\frac{1}{k^{n}(s)} \ \frac{e( g \cdot x)}{e(x)} |g x|^{s}$,
 and $\int q_{e,n}^{s} (x)   (g) =1$. For $s=0$ we have $e=1$, $Q_{e}^{s}=P$. Also we write $q_{e,1}^{s}=q_{e}^{s}$.

\begin{lemm}
\label{L2.9}
Assume $e$ is as above, $f\in C (\mathbb P^{d-1})$ is real valued and satisfies $Q_{e}^{s} f \leq f$. Then, on $\Lambda ([\supp \mu])$, $f$ is constant and equal to its infimum on $\mathbb P^{d-1}$.
\end{lemm}

\begin{proof}
Let $M^{-}=\{x \in \mathbb P^{d-1}$; $f(x)=\textrm{inf }\{f(y); y\in \mathbb P^{d-1}\}$. The relation $f(x) \geq \int q_{e}^{s} (x,g) f( g \cdot x) d\mu (g)$ implies that if $x\in M^{-}$ then $ g \cdot x \in M^{-}, \mu \textrm{-a.e}$. Hence $[\supp \mu]\cdot M^{-} \subset M^{-}$. Since $\Lambda([\supp \mu])$ is the unique $[\supp \mu]$ -minimal subset of $\mathbb P^{d-1}$, we get $\Lambda ([\supp \mu])\subset M^-$, i.e.\
\[f(x)=\inf \{f(y); y\in \mathbb P^{d-1}\},\ \  \hbox{\rm if}\ x\in \Lambda ([\supp\mu]).\]
\end{proof}
 Using Lemma 2.8,   the existence of $e\in C (\mathbb P^{d-1})$ with $P^{s} e=k(s) e$ is obtained by the following
\begin{lemm}
\label{L2.10}
  Assume $\sigma \in M^{1} (\mathbb P^{d-1})$ and $k>0$ satisfy $^{*}P^{s} \sigma=k\sigma$. Then the function $\widehat{\sigma}^{s}$ on $\mathbb P^{d-1}$ defined by,
\[\widehat{\sigma}^{s} (x) =\int |\langle x, y \rangle |^{s}d\sigma (y)\]
satisfies $P^{s}\widehat{\sigma}^{s}=k \widehat{\sigma}^{s}$. Furthermore
$\widehat{\sigma
}^{s}$ is  positive and H\"older of order $\bar s=\inf (1,s)$.
\end{lemm}

\begin{proof}
  We have
$|g x|^{s} \widehat{\sigma}^{s} ( g \cdot x)=\int |\langle x,g^*\cdot y\rangle|^{s} |g^{*} y|^{s} d\sigma (y)$ and $^{*}P^s \sigma=k\sigma$; hence
\[P^{s} \widehat{\sigma}^{s} (x) =k \int |\langle x,z \rangle|^{s} d\sigma (z)=k \widehat{\sigma}^{s}(x).\]
 If $\widehat{\sigma}^{s} (x)=0$ for some $x$, then $|\langle x, y \rangle |=0$, $\sigma$\textrm{-a.e}.; hence
\[\supp \sigma \subset \{ y\in \mathbb P^{d-1}; \langle x , y \rangle=0\}.\]
This contradicts Lemma 2.8, since $[\supp \mu]$ satisfies i-p. Hence $\widehat{\sigma}^{s}$ is positive.

  In order to show the H\"older property of $\widehat{\sigma}^{s}$, we use the inequality
  $|a^{s}-b^{s}| \leq \widehat s |a-b|^{\bar s}$  where $a,b \in [0,1] \ \ \widehat s=\sup (s,1)$, $\bar s=\inf (s,1)$. Then
\[||\langle x , y \rangle|^{s}-|\langle x', y\rangle |^s| \leq \widehat s |x-x'|^{\bar s} |y|^{\bar s}, |\widehat{\sigma}^{s} (x)- \widehat{\sigma}^{s}(x')| \leq \widehat s \delta^{\bar s} (x,x').\]
\end{proof}
\begin{lemm}
\label{L2.11}
Let $e$ be a positive and $\bar s$-H\"older function on $\mathbb P^{d-1}$  with $\bar s=\inf (s,1)$, $s> 0$. There exists a constant $b_{s}>0$ such that for any $(x,y)\in \mathbb S^{d-1} \times \mathbb S^{d-1}, g\in G$
\begin{align*}
||g x|^{s} - |g y|^{s}| &\leq (s+1)  |g|^{s} \widetilde{\delta}^{\bar s} (x,y),\\
\widetilde{\delta} ( g \cdot x,  g \cdot y)&\leq 2 \frac{|g|}{|g x|} \widetilde{\delta} (x,y),\\
|\theta_{e}^{s} (x,g) -\theta_{e}^{s} (y,g) &\leq b_{s} |g|^{s}  \widetilde{\delta}^{\bar s} (x,y).
\end{align*}
\end{lemm}
\begin{proof}
We use the inequality
$|a^{s}-b^{s} | \leq |a-b|^{s}$ if $a, b \geq 0$, $s\leq 1$ to get
\[||g x|^{s}-|g y|^{s} |\leq ||g x|-|gy||^{s} \leq |g (x-y)|^{s}\leq |g|^{s} |x-y|^{s}.\]
 Hence $||gx|^{s}-|gy|^{s}|\leq |g|^{s} \widetilde{\delta}^{\bar s} (x,y)$.
If $s>1$, we use $\frac{1}{s} |a^s-b^s|\leq \sup (a,b)^{s-1} |a-b|$ if $a,b\geq 0$. We get
\[| |gx|^{s}- |g y |^{s} | \leq s |g|^{s-1} | | g x|-|g y|| \leq s |g|^{s-1} | g( x- y)| \leq s |g|^{s} | x- y|.\]
  Hence the first inequality follows. Furthermore
\[\widetilde{\delta}( g \cdot x, g. y)=\left |\frac{g  x}{|g x|}-\frac{g  y}{|gy|}\right | \leq \frac{|g( x- y)|}{|g x|}+ |g y| \left |\frac{1}{|gx|}-\frac{1}{|gy|}\right   |\leq 2 \frac{|g| | x- y|}{|g x|} .\]
  Hence $\widetilde{\delta}( g \cdot x,  g \cdot y)\leq 2 \frac{|g|}{|g x|} \widetilde{\delta} (x,y)$.
We write
\[|\theta_{e}^{s} (x,g)-\theta_{e}^{s} (y,g) | \leq ||g x|^{s}-|gy|^{s}| \frac{e( g \cdot x)}{e(x)}+|g y|^{s} \frac{|e( g \cdot x)-e( g \cdot y)|}{e(x)}+ |g y|^{s} e ( g \cdot y) \left | \frac{1}{e(x)}-\frac{1}{e(y)}\right |.\]
  In view of the first inequality the first term satisfies the required bound. The last term also satisfies it, since $\frac{1}{e(x)}$ is $\bar s$- H\"older and $|g y|^{s}\leq |g|^{s}$. For the second term we write :
$$|e ( g \cdot x)-e(g\cdot y)| \leq [e]_{\bar s} \widetilde{\delta}^{\bar s} ( g \cdot x , g\cdot y)\leq 2^{\bar s} [e]_{\bar s} \left(\frac{|g|}{|g y|}\right)^{\bar s} \widetilde{\delta}^{\bar s} (x,y).$$
  Hence this term is bounded by :
$$|e^{-1}| 2^{\bar s} [e]_{\bar s} \widetilde{\delta}^{\bar s} (x,y) |g y|^{s-\bar s} |g|^{\bar s}\leq 2^{\bar s} |e^{-1}|  [e]_{\bar s} |g|^{s} \widetilde{\delta}^{\bar s} (x,y).$$
  The  above inequalities imply the lemma.
  \end{proof}

\vskip 2mm
  The following is well known (see for example \cite{48} Theorem 6)

\begin{lemm}
\label{L2.12}
  Let $X$ be a compact metric space, $Q$ a Markov operator on $X$, which preserves $C(X)$. Assume that all the $Q$-invariant continuous functions are constant and for any $\varphi \in C(X)$, the sequence $Q^{n} \varphi$ is equicontinuous. Then $Q$ has a unique stationary measure $\pi$, hence for any $\varphi \in C(X)$ the sequence $\frac{1}{n} \displaystyle\mathop{\Sigma}_{0}^{n-1} Q^{k} \varphi$ converges uniformly to $\pi(\varphi)$. Furthermore if the equation $Q\psi=e^{i\theta}\psi$, $\psi \in C(X)$ implies $e^{i\theta}=1$, then, for any $\varphi \in C(X)$, $Q^n \varphi$ converges uniformly to $\pi (\varphi)$.
  \end{lemm}

\begin{prop}
\label{P2.13}
Let $\mu \in M^{1} (G)$ and assume that the semigroup $[\supp \mu]$ satisfies  $\ i-p.$
Let $s\in I_{\mu}$, $s>0\  \varepsilon \in ]0,\bar s]$ with $\bar s=\inf (1,s),\  e\in C (\mathbb P^{d-1})$ is positive, $\bar s$-H\"older  with $P^{s} e=k(s) e$.
Then there exists $a_{s}\geq 0$ such that for any $n\in \mathbb N$ and any $\varepsilon$-H\"older function $\varphi$ on $\mathbb P^{d-1}$
\[[(Q_{e}^{s})^{n} \varphi ]_{\varepsilon} \leq a_{s}|\varphi|+ \rho_{n,s}(\varepsilon) [\varphi]_{\varepsilon},\]
  where
  $\rho_{n,s}(\varepsilon) = \sup_{x,y}\int q_{e,n}^{s} (x,g) \frac{\delta^{\varepsilon}( g \cdot x,  g \cdot y)}{\delta^{\varepsilon} (x,y)} d\mu^{n} (g)$ is bounded independently of $n$.
 In particular  for any $\psi \in C (\mathbb P^{d-1})$ the sequence $(Q_{e}^s)^n \psi$ is equicontinuous.
 Furthermore any continuous $Q_{e}^{s}$-invariant function is constant and $Q_{e}^{s}$ has a unique stationary measure $\pi_{e}^{s}$. If $s=0$ and $\varphi \in C (\mathbb P^{d-1})$ we have the uniform convergence $\lim_{n\rightarrow \infty} \overline{P}^n \varphi=\nu (\varphi)$.
  \end{prop}

\begin{proof}
  The definition of $Q_{e}^{s}$ \ \ gives for any $\varepsilon$-H\"older function $\varphi$, and $q_{e,n}$ as before Lemma 2.9,
\[|(Q_{e}^{s})^{n} \varphi (x) -(Q_{e}^{s})^{n} \varphi (y) | \leq |\varphi| \int |q_{e,n}^{s} (x,g)- q_{e,n}^{s} (y,g) | d\mu^{n} (g)+[\varphi]_{\varepsilon} \int q_{e,n}^{s} (y,g) \delta^{\varepsilon} ( g \cdot x, g\cdot y) d \mu^{n} (g).\]
  Lemma 2.11 shows that the first integral is dominated by $\frac{b_{s}}{k^{n} (s)} \delta^{\bar s} (x,y) \int |g|^{s} d\mu^{n} (g)$ and Lemma 2.8 gives $k^{n}(s) \geq c_{s} \int |g|^{s} d\mu^{n} (g)$.  Hence, $[(Q_{e}^{s})^{n}\varphi]_{\varepsilon}\leq a_{s} |\varphi|+\rho_{n,s} (\varepsilon) [\varphi]_{\varepsilon}$ with $a_{s} =b_{s}/c_{s}$. Lemma 2.11 allows to bound $\rho_{n,s}(\varepsilon)$ as follows
\[\delta^{\varepsilon} ( g \cdot x , g\cdot y)\leq 2^{\varepsilon} \frac{|g|^{\varepsilon}}{|g x|^{\varepsilon}} \delta^{\varepsilon} (x,y), \ \ \rho_{n,s} (\varepsilon) \leq \frac{2^{\varepsilon}}{k^{n}(s)} \displaystyle\sup_{x} \int \frac{e(g\cdot x)}{e(x)} |g x|^{s-\varepsilon} |g|^{\varepsilon} d\mu^{n} (g).\]
  We denote $c=\displaystyle\sup_{g,x}\frac{e( g \cdot x)}{e(x)}<\infty$, hence using $s\geq \varepsilon$, $|g x|^{s-\varepsilon} \leq|g|^{s-\varepsilon}$ we get
\[\frac{e(g\cdot x)}{e(x)} |g x|^{s-\varepsilon} |g|^{\varepsilon} \leq c|g|^{s}, \qquad \rho_{n,s} (\varepsilon) \leq \frac{c 2^{\varepsilon}}{k^{n}(s)} \int |g|^{s} d\mu^{n} (g) \leq c \frac{2^{\varepsilon}}{c_{s}}.\]

  Assume that $\varphi \in C(\mathbb P^{d-1})$ satisfies $Q_{e}^{s} \varphi=\varphi$ and denote
\[M^{+}=\{x\in \mathbb P^{d-1}; \varphi (x) = \displaystyle\sup_{y \in \mathbb P^{d-1}} \varphi (y)\}, \
M^{-}=\{x\in \mathbb P^{d-1}; \varphi (x)=\inf_{y \in \mathbb P^{d-1}} \varphi (y)\}.\]
Then, as in the proof of Lemma 2.9, $\supp\mu\cdot  M^{+}\subset M^{+}$, $\supp\mu\cdot M^{-} \subset M^{-}$, hence by minimality of $\Lambda([\supp\mu])$ we have $\Lambda ([\supp\mu]) \subset M^{+}\cap M^{-}$.
It follows $M^{+} \cap M^{-}\neq \phi$, $\varphi=$constant on $\mathbb P^{d-1}$.
If $\psi \in C(\mathbb P^{d-1})$ is $\varepsilon$-H\"older the above inequality gives for any $x,y\in \mathbb P^{d-1}$
\[|(Q_{e}^s)^n \varphi (x)-(Q_{e}^s)^n \varphi (y)|\leq (a_{s}|\varphi|+\rho_{n,s}(\varepsilon)) [\varphi]_{\varepsilon}) \delta^{\varepsilon} (x,y).\]
Since $\rho_{n,s}(\varepsilon)$ is bounded this shows that the sequence $(Q_{e}^{s})^{n} \psi$ is equicontinuous. By density this remains valid for any $\psi\in C(\mathbb P^{d-1})$.

  Hence we can apply Lemma 2.12 to $Q=Q_{e}^{s}$ : there is a unique $Q_{e}^{s}$-stationary measure. If $s=0$, we have $\varepsilon=0$, hence the above inequality does not show the equicontinuity of $\overline{P}^n \varphi$. In this case the equicontinuity follows from Theorem 3.2 in the next section; we have for $\varepsilon>0$ the convergence  $\lim_{n\rightarrow \infty} \int \delta^{\varepsilon} ( g \cdot x, g \cdot y) d\mu^n (g)=0$, which implies for $\varphi \in H_{\varepsilon} (\mathbb P^{d-1})$,
\begin{align*}
|\overline{P}^n \varphi) (x)- \overline{P}^n \varphi (y)|&\leq [\varphi]_{\varepsilon} \int \delta^{\varepsilon} ( g \cdot x,  g \cdot y) d\mu^n (g),\\
|\overline{P}^n \varphi (x)-\nu (\varphi)|&\leq [\varphi]_{\varepsilon} \int \delta^{\varepsilon} ( g \cdot x,  g \cdot y) d\nu (y) d\mu^n (g), \\
\lim_{n\rightarrow \infty} |\overline{P}^n \varphi (x)-\nu (\varphi)|&=0.
\end{align*}
\end{proof}

\begin{rema}
\begin{enumerate}
\item   If for some $\tau>0$ and $s\in[0, s_\infty[$,  $\int |g|^{s} \gamma^{\tau} (g) d\mu (g) <\infty$ with $\gamma(g)=\sup (|g|, |g^{-1}|)$ then it is proved in section 3, Corollary 3.18 that  $\lim_{n\rightarrow \infty} \rho_{n,s} (\varepsilon)=0$, hence $\rho_{n,s} (\varepsilon) <1$ for some $n=n_{0}$. Hence $(Q_{e}^{s})^{n_{0}}$ satisfies a  Doeblin-Fortet inequality (see \cite{27}).
\item
 Let $\tilde Q_{e}^{s}$ be the Markov kernel on $\mathbb S^{d-1}$ defined by $\tilde Q_{e}^{s} \varphi=\frac{1}{k(s) e} \tilde P^{s}(\varphi e)$ where $e$ still denotes the function on $\mathbb S^{d-1}$ corresponding to $e\in C(\mathbb P^{d-1})$.  Then the inequality and its proof remain valid for $\tilde Q_{e}^{s}$ instead of $Q_{e}^{s}$. In particular for any $\psi \in C(\mathbb S^{d-1})$ the sequence $(\tilde Q_{e}^{s})^{n} \psi$ is equicontinuous. This fact will be used in the next paragraph.
 \end{enumerate}
 \end{rema}

\noindent
\textit{Proof of Theorem 2.6:}
 As in the proof  of Lemma 2.8, we consider the non linear operator ${^{*}\widehat{P}^{s}}$ on $M^{1} (\mathbb P^{d-1})$ defined by  ${^{*}\widehat{P}^{s}} \sigma=\frac{^{*}P^{s}\sigma}{(^{*}P^{s}\sigma) (1)}$. The same argument gives the existence of $k$ and $\sigma \in M^{1} (\mathbb P^{d-1})$ such that $^{*}P^{s}\sigma=k\sigma$, with $k=(^{*}P^{s} \sigma) (1) >0$. We consider only the case $d>1$.

  Since $[\supp \mu^{*}]=[\supp \mu]^{*}$ satisfies i-p, Lemma 2.8 applied to $\mu^{*}$ gives $k=k(s)$ and $\sigma$ is not supported by a hyperplane. Then Lemma 2.10 implies that $\widehat{\sigma}^{s}(x)=\int |\langle x, y \rangle |^{s} d\sigma (y)$  satisfies $P^{s} \widehat{\sigma}^{s}=k(s) \widehat{\sigma}^{s}$ and is positive, H\"older  continuous of order $\bar s=\inf (1,s)$. Hence we can apply Proposition 2.13 with $e=\widehat{\sigma}^{s}$; then we get existence and uniqueness of $e^{s} \nu^{s}\in M^{1}(\mathbb P^{d-1})$ with $P^{s} e^{s}=k(s) e^{s}$, $P^{s}\nu^{s}=k(s)\nu^{s}, \nu^{s}(e^{s})=1$ and $e^{s}$ satisfies $p(s) e^s (x)=\int |\langle x , y \rangle|^s d$ $\sigma (y)$ where $p(s)=\nu^{s}(\widehat{\sigma}^{s})$. Also $Q^{s}=Q_{e^{s}}^{s}$ has a unique stationary measure $\pi^{s}$. The uniqueness of $\nu^{s} \in M^{1} (\mathbb P^{d-1})$ with $P^{s} \nu^{s}=k(s) \nu^{s}$ follows. Also $\sigma={^{*}\nu^s}$ by the same proof.

   Lemma 2.8 implies that if some $\eta \in M^{1} (\mathbb P^{d-1})$ satisfies $P^{s} \eta=k \eta$, then $k=k(s)$. Since $\supp \nu^{s}$ is $[\supp \mu]$-invariant and $\Lambda ([\supp \mu])$ is minimal we get $\supp \nu^{s} \supset \Lambda ([\supp \mu])$. We can again use Schauder-Tychonoff theorem in order to construct $\sigma' \in M^{1} (\mathbb P^{d-1})$ with $\supp \sigma' \subset \Lambda ([\supp \mu]), P^{s} \sigma'=k\sigma'$. Since $\sigma'=\nu^{s}$, we get finally $\supp \nu^{s}=\Lambda ([\supp \mu])$.

  In order to show the continuity of $s\rightarrow \nu^{s}$, $s\rightarrow e^{s}$ we observe that, from the above argument, $\nu^{s}$ is uniquely defined by  $P^{s} \nu^{s}=k(s) \nu^{s}$, $\nu^{s} \in M^{1} (\mathbb P^{d-1})$. Also, by convexity,  $k(s)$ is continuous. On the other hand, the uniform continuity of $(x,s)\rightarrow |g x|^{s}$ and the fact that $|g x|^s\leq |g|^s$ is bounded by the $\mu$-integrable function $\sup(|g|^{s_{1}}, |g|^{s_{2}})$ on $[s_{1} , s_{2}] \subset I_{\mu}$ implies the uniform continuity of $P^{s} \varphi$ if $\varphi$ is fixed. Then we consider a sequence $s_{n} \in I_{\mu}$, $s_{0} \in I_{\mu}$ with $\displaystyle\mathop{\lim}_{s_{n}\rightarrow s_{0}} \nu^{s_{n}}=\eta \in M^{1} (\mathbb P^{d-1})$. We have

\centerline {$P^{s_{n}} \nu^{s_{n}} (\varphi)=\nu^{s_{n}} (P^{s_{n}} \varphi)$, \ \ $\displaystyle\mathop{\lim}_{s_{n}\rightarrow s_{0}} P^{s_{n}}\nu^{s_{n}}(\varphi) = \displaystyle\mathop{\lim}_{s_{n}\rightarrow s_{0}} k(s_{n}) \nu^{s_{n}}
 (\varphi)=k(s_{0}) \eta (\varphi)$. }

   Then the uniform continuity in $(s,x)$ of $P^s \varphi (x)$ implies $P^{s_{0}} \eta=k(s_{0}) \eta$. The uniqueness of $\nu^{s_{0}}$ implies $\nu^{s_{0}}=\eta$, and the arbitrariness of $s_{n}$ gives the continuity of $s\rightarrow \nu^{s}$ at $s_{0}$. The same property is true for the operator $^{*}P^{s}$ and the measure $^{*}\nu^{s}$ defined by $^{*}P^{s}\  ({^{*}\nu^{s}})=k(s)\  ({^{*}\nu^{s}})$,\ \  ${^{*}\nu^{s}} \in M^{1}(\mathbb P^{d-1})$.
 Lemma 2.10, implies $p(s) e^s (x)=\int |\langle x , y \rangle|^{s} d {^{*}\nu^{s}} (y)$, and since the set of functions  $x\rightarrow |\langle x,y\rangle |^{s}$ $(y\in \mathbb P^{d-1}, s\in I_{\mu})$ is locally equicontinuous  we have  $\displaystyle\mathop{\lim}_{s\rightarrow s_{0}} |e^{s}- e^{s_{0}}|=0$.

  In order to show the strict convexity of $\log k(s)$ we take $s, t\in I_{\mu}$, $p\in (0,1)$ and we observe that from H\"older inequality,
$P^{ps+(1-p)t} [(e^{s})^{p} (e^{t})^{1-p}] \leq k^{p} (s) k^{1-p} (t) (e^{s})^{p} (e^{t})^{1-p}$.
We denote $f=(e^{s})^{p} (e^{t})^{1-p}$ and assume $k(p s+(1$-$p)t)=k^p (s) k^{1-p}(t)$ for some $s\neq t$.
Then Lemma 2.9 can be used with $e=e^{ps+(1-p)t}$ and $Q^{p q+(1-p)t}_{e} \varphi=\frac{1}{k(p s+(1-p) t)e} P^{p s+(1-p)t}
 (\varphi e)$. It gives on $\Lambda ([\supp \mu]) : f=c e^{ps+(1-p)t}$ for some constant $c>0$.

   Hence, on $\Lambda([\supp\mu])$ we have
\[P^{p s+(1-p)t} [(e^{s})^{p} (e^{t})^{1-p}]=k^{p} (s) k^{(1-p)} (t) (e^{s})^{p} (e^{t})^{1-p}.\]
   This means that there is equality in the above H\"older inequality. It follows that, for some positive function $c(x)$ and any $x$ in $\Lambda([\supp \mu])$, $g\in \supp\mu$
\[ |g x|^{s} \frac{e^{s} ( g \cdot x)}{e^{s}(x)}= c(x) |g x|^{t} \frac{e^{t}( g \cdot x)}{e^{t} (x)}.\]
   Integration with respect to $\mu$ gives: $c(x)=\frac{k(s)}{k(t)}$. Since $s\neq t$, we get, for some constant $c>0$ and $\varphi \in C(\mathbb P^{d-1})$ positive, for any $(x,g)$ as above
 $ |g x|=c \frac{\varphi ( g \cdot x)}{\varphi (x)}.$
   It follows, if $g\in (\supp \mu)^{n}$ and $x\in \Lambda ([\supp \mu])$, $|gx| = c^{n} \frac{\varphi ( g \cdot x)}{\varphi (x)}$. If $g \in [\supp \mu]^{prox}$, we get $|\lambda_{g}| \in c^{\mathbb N}$. This contradicts Proposition 2.5.

  In order to show the  convergence of $(Q^s)^n\varphi$, since by Proposition 2.13 the family  $(Q^{s})^{n} \varphi$ is equicontinuous,  it suffices to show in view of  Lemma 2.12  that the  relation $Q^{s} \varphi=e^{i\theta} \varphi$ with $\varphi \in C(\mathbb P^{d-1})$, $|e^{i\theta}|=1$ implies $e^{i\theta}=1$, $\varphi=$constant.  Taking absolute values we get $|\varphi| \leq Q^{s}|\varphi|$. As in Lemma 2.9, we get that for any $x$ in $\Lambda ([\supp \mu])$,
 $ |\varphi(x)|= \sup \{ |\varphi (y)|; y\in \mathbb P^{d-1}\}$.
  Hence we can assume  $|\varphi (x)|=1$ on $\Lambda ([\supp \mu])$.  Now we can use the equation
\[e^{i\theta} \varphi (x)=\int q^{s} (x,g) \varphi ( g \cdot x) d\mu (g),\]
   where $q^{s} (x,g)=\frac{1}{k(s)} \ \frac{e^{s}( g \cdot x)}{e^{s}(x)} |g x|^{s}$, hence $\int q^s (x,g) d\mu (g)=1$.  Strict convexity  yields the equality $e^{i\theta} \varphi (x)=\varphi ( g \cdot x)$,
for any $x \in \Lambda ([\supp \mu]) , g\in \ \supp \mu$.

  We know, from Proposition 2.13 that $\bar P^{n} \varphi$ converges uniformly to $\nu (\varphi)$ where $\nu$ is the unique $\bar P$-stationary measure on $\mathbb P^{d-1}$. Furthermore, on $\Lambda([\supp \mu])$ we have $\overline P^{n} \varphi=e^{i n \theta} \varphi$. The above convergence gives $e^{i\theta}=1$, since $\varphi \neq 0$ on $\Lambda ([\supp \mu])$. The fact that $\varphi$ is constant follows from Proposition 2.13.

  In order to show the last assertion in case $d>1$ we write $Q^z \varphi (x)$ as
\[Q^z \varphi (x)=\int |g x|^{it} q^s (x,g) \varphi ( g \cdot x) d\mu (g).\]
   We observe that the absolute value of the function $Q^z \varphi$  is bounded by the function  $Q^s |\varphi|$. Hence, from above, the equation $Q^z \varphi=e^{i\theta} \varphi$ gives $Q^s |\varphi| \geq |\varphi|$, hence $Q^s |\varphi|=|\varphi|$ and $|\varphi|=$cte. Then the equation $Q^z \varphi=e^{i \theta} \varphi$ gives for any $x$ and $g \in \supp\mu$, since $\int q^s (x,g) d\mu (g)=1$  we have
$|g x|^{it} \frac{\varphi ( g \cdot x)}{\varphi (x)}=e^{i\theta}$, $\mu$\textrm{-a.e}.
 This contradicts Proposition 2.5 if $t\neq 0$.
 If $t=0$,from above we have $e^{i \theta}=1$, $\varphi=$constant .
 \mbox{}\hfill  $\square$

\subsection{Eigenfunctions, limit sets, and eigenmeasures on $\mathbb S^{d-1}$}

  Here we study the operator  $\tilde P^s$ on $\mathbb S^{d-1}$ defined by $\tilde P^{s} \varphi (x)=\int \varphi ( g \cdot x) |g x|^{s} \mu (dg)$. We show that there are 2 cases, depending on the existence of a $[\supp \mu]$-invariant proper convex cone in $V$ or not. We still denote by $e^{s}$ the function on $\mathbb S^{d-1}$ lifted from $e^{s}\in C(\mathbb P^{d-1)}$. We denote $\tilde Q^{s}$ the operator on $\mathbb S^{d-1}$ defined by  $\tilde Q^{s} \varphi=\frac{1}{k(s) e^{s}} \tilde P^{s} (\varphi e^{s})$.

  We already know, using the remark which follows Proposition 2.13, that for $s>0$ and any given $\varphi \in C(\mathbb S^{d-1})$, the sequence $(\tilde Q^{s})^{n} \varphi$ is equicontinuous. For any subsemigroup $T$ of $G$ satisfying condition i-p, we denote by $\tilde{\Lambda}(T)$ the inverse image of $\Lambda(T)$ in $\mathbb S^{d-1}$. We  begin by considering the dynamics of $T$ on $\mathbb S^{d-1}$. For analogous results in more general situations see \cite{24}. We recall that a convex cone in $V$ is said to be proper if it does not contain a line.

\begin{prop}
\label{pro:2.15}
  Assume $T \subset G$ is a subsemigroup which satisfies condition$\ i-p $. If $d=1$, we assume that $T$ is non-arithmetic. Then the action of $T$ on $\mathbb S^{d-1}
$ has one or two minimal sets whose union is $\tilde{\Lambda} (T)$ :
\begin{enumerate}[I]
\item   There is no $T$-invariant proper convex cone in $V$ and
in that case,  $\tilde{\Lambda} (T)$ is the unique $T$-minimal subset of $\mathbb S^{d-1}$.
\item
$T$ preserves a closed proper convex cone $C\subset V$ and
then the action of $T$ on $\mathbb S^{d-1}$ has two and only two minimal subsets $\Lambda_{+} (T)$, $\Lambda_{-}(T)$ with $\Lambda_{-}(T)=-\Lambda_{+} (T), \Lambda_{+} (T) \subset \mathbb S^{d-1} \cap\  C$. The convex cone generated by $\Lambda_{+}(T)$ is proper and $T$-invariant.
  \end{enumerate}

\end{prop}

  The proof depends of the following lemma.

\begin{lemm}
\label{lem:2.16}
  Let $V_{i} (1\leq i \leq r)$ be vector subspaces of $V$. If condition i-p is valid, then there exists $g\in T^{\textrm{prox}} $ such that the hyperplane $V_{g}^{<}$ does not contain any $V_{i}$  $(1\leq i\leq r)$.
  \end{lemm}

\begin{proof}

 The dual semigroup $T^{*}$ of $T$ satisfies also condition i-p hence we can also consider its limit set $\Lambda (T^{*}) \subset \mathbb P(V^{*})$. Let $\bar v (g^{*})$ be the point of $\mathbb P(V^{*})$ corresponding to a dominant eigenvector of $g^{*}$. Observe that the condition that an hyperplane contains $V_{i}$ defines a subspace of $V^{*}$. If for any $g^{*} \in (T^{*})^{prox}$ the hyperplane $\bar v(g^{*})$ contains some $V_{i}$ then by density any $x\in \Lambda (T^{*})$ contains some $V_{i}$. Then the $T^{*}$-invariance of $\Lambda(T^{*})$ implies that $T^{*}$ leaves invariant a finite union of subspaces of $\mathbb P(V^{*})$, which contradicts condition i-p.
 \end{proof}

 \textit{Proof of  Proposition \ref{pro:2.15}:}
 Let $x \in \tilde{\Lambda} (T)$ and $S=\overline{T \cdot y}$. We observe that if $y\in \mathbb S^{d-1}$, then $\overline{T \cdot y}$ contains $x$ or $-x$, since the projection of $\overline{T\cdot x}$ in $\mathbb P^{d-1}$ contains $\Lambda (T)$. Assume first $-x \notin  \overline{T \cdot x}$. If $y\in \overline{T\cdot x}$, then $\overline{T \cdot y} \subset S$, hence $x \in \overline{T\cdot y}$. This shows the $T$-minimality of $S$. The same argument shows that $-y\notin S$, hence $S\cap-S=\phi$.
 Since the projection of $S$ in $\mathbb P^{d-1}$ is $\Lambda(T)$, we see that the projection of $\mathbb S^{d-1}$ on $\mathbb P^{d-1}$ gives a $T$-equivariant homeomorphism of $S$ on $\Lambda(T)$. Since $-x\notin S$, there are two $T$-minimal sets, $S$ and $-S$. Since for any $y\in \mathbb S^{d-1}$, $\overline{T\cdot y}$ contains $S$ or $-S$, these sets are the unique minimal sets.

   Assume now $-x\in \overline{T \cdot x}$, hence $S=-S$. Since the projection of $S$ in $\mathbb P^{d-1}$ is $\Lambda (T)$, we see that $S=\tilde{\Lambda}(T)$.

  Assume now that $C$ is a $T$-invariant closed proper convex cone. Then $C\cap \mathbb S^{d-1}$ is $T$-invariant and closed, hence $C\cap \mathbb S^{d-1} \supset \Lambda_{+}(T)$ or $\Lambda_{-}(T)$ in the first situation, $(-x\notin \overline{T \cdot x})$.
 In the second situation $C$ cannot exists, since $C\cap \mathbb S^{d-1}$ would contain $\tilde{\Lambda}(T)$, which is symmetric.

  It remains to show that, in the first situation, there exists a $T$-invariant  closed proper convex cone. Let $C$ be the convex cone generated by $\Lambda_{+} (T)$ and let us show $C\cap -C=\{0\}$. Assume $C \cap-C\neq \{0\}$; then we can find $y_{1}, \ldots , y_{p}\in C$,  $z_{1}, \ldots , z_{q} \in -C$ and convex combinations $y=\displaystyle\mathop{\Sigma}_{1}^{p} \alpha_{i} y_{i}, z=\displaystyle\mathop{\Sigma}_{1}^{q} \beta_{j} z_{j}$ with $y=z$.  Lemma \ref{lem:2.16} shows that there exists $g\in T^{\textrm{prox}} $ such that $y_{i} (1\leq i\leq p)$ and $z_{j}(1\leq j \leq q)$ do not belong to $V_{g}^{<}$. Hence, with $n\in 2\mathbb N$ :
\[\lim_{n\rightarrow+\infty} \frac{g^{n}y}{|g^{n}|}=\lim_{n\rightarrow \infty} \Sigma_{1}^{p} \alpha_{i} \frac{g^{n} y_{i}}{|g^{n} y_{i}|} \ \frac{|g^{n} y_{i}|}{|g^{n}|}=\left(\Sigma_{1}^{p} \alpha_{i} u_{i}\right) v_{g},\]
  where  $u_{i}= \displaystyle\mathop{\lim}_{n\rightarrow \infty} \frac{|g^{n} y_{i}|}{|g^{n}|}>0$ and $v_{g} \in \Lambda_{+} (T)$ is the unique dominant eigenvector of $g$ in $\Lambda_{+} (T)$.
In the same way :
$$\displaystyle\mathop{\lim}_{n\rightarrow \infty}  \frac{g^{n} z}{|g^n|}=-\left(\displaystyle\mathop{\Sigma}_{1}^{q} \beta_{j} u_{j}'\right) v_{g},$$
  with $u'_{j}>0$. Since $y=z$ we have a contradiction. Hence we have the required dichotomy. The last assertion follows. \hfill $\square$

  We denote by $\tilde{\nu}^s$ the  symmetric measure on $\mathbb S^{d-1}$ with projection $\nu^s$ on $\mathbb P^{d-1}$. In case II, we denote by $\nu_{+}^s$ (resp $\nu_{-}^s)$ the normalized restrictions of $\tilde{\nu}^s$ to $\Lambda_{+}(T)$ (resp $\Lambda_{-}(T))$. For a subset $X\subset \mathbb S^{d-1}$ we recall that $\textrm{Co}(X)$ is the convex envelope of $X$ in $\mathbb S^{d-1}$.

  \begin{theo}
  \label{thm:2.17}
  Let $\mu \in M^{1} (G), \ s\in I_{\mu}$ and assume $T=[\supp\mu]$ satisfies$\ i-p$.  If $d=1$ we assume that $\mu$ is non-arithmetic. Then for any $\varphi \in C (\mathbb S^{d-1})$, $x\in \mathbb S^{d-1}$, we have the uniform convergence
\[\lim_{n\rightarrow \infty} \frac{1}{n} \Sigma_{1}^{n} (\tilde Q^{s})^{n} \varphi (x)=\widetilde{\pi}^{s} (x) (\varphi),\]
  where,  $\tilde{\pi}^{s} (x)  \in M^{1} (\mathbb S^{d-1})$ is supported on $\tilde{\Lambda} (T)$ and is $\tilde Q^{s}$-stationary.

  Furthermore there are 2 cases given by Proposition \ref{pro:2.15}.

 \begin{enumerate}[I]
 \item
  $\tilde Q^{s}$ has a unique stationary measure $\tilde{\pi}^{s}$ with $\supp \tilde{\pi}^{s}= \tilde{\Lambda}(T)$ and $\tilde{\pi}^{s}(x)=\tilde{\pi}^{s}$ for any $x\in \mathbb S^{d-1}$. The $\widetilde Q^s$-invariant functions are constant.
We have $\tilde{\pi}^s=e^s \tilde{\nu}^s$ and $\tilde{P}^s \tilde{\nu}^s=k(s) \tilde{\nu}^s$.
\item $\tilde Q^{s}$ has two and only two extremal stationary measures $\pi_{+}^{s}$, $\pi_{-}^{s}$. We have $\supp \pi_{+}^{s}=\Lambda_{+}(T)$ and $\pi_{-}^{s}$ is symmetric of $\pi_{+}^{s}$. If $\pi_{+}^s=e^s \nu_{+}^s$, then $\widetilde P^s \nu_{+}^s=k(s) \nu_{+}^s$.
Also, there are 2 minimal $\tilde Q^{s}$-invariant continuous functions $p_{+}^{s}$, $p_{-}^{s}$ and we have
\[\tilde{\pi}^{s} (x)=p_{+}^{s} (x) \pi_{+}^{s} + p^{s}_{-} (x) \pi_{-}^{s}.\]

  Furthermore $p_{+}^{s} (x)$ is equal to the entrance probability in the convex envelope $\textrm{Co}(\Lambda_{+} (T))$ for the Markov chain defined by  $\tilde Q^{s}$. In particular $p_{+}^{s} (x)=1$ (resp $p^{s}_{+}(x)=0$) if $x\in \Lambda_{+} (T)$ (resp $\Lambda_{-}(T)$.

  If $^{*}\nu_{+}^s \in M^1 (\Lambda_{+}(T^{*})$ satisfies $^{*}\tilde{P}^s$ $^{*}\nu_{+}^s=k(s) ^{*}\nu_{+}^s$, we have for $u\in \mathbb S^{d-1}$ and $p(s)$ as in Theorem 2.6,
$p(s) e^s_{+} (u)=\int \langle u,u'\rangle_{+}^s d$  $^{*}\nu_{+}^s (u')$ with $ \langle u,u'\rangle_{+}=\sup (0, \langle u,u'\rangle)$ and $e_{+}^s=p^{s}_{+}e^s$,  (resp $e^{s}_{-}=p^{s}_{-} e^{s}$) satisfies $\tilde{P}^s e^s_{+}=k(s) e^s_{+}$ (resp $\widetilde{P}^{s} e^{s}_{-}=k(s) e^{s}_{-}$).

   The space of  continuous $\widetilde Q^s$-invariant  (resp $\tilde{P}^s)$-eigenfunctions is generated by $p^s_{+}$ and $p^s_{-}$ (resp $e^s_{+}$ and $e^s_{-}$).
   \end{enumerate}
\end{theo}

 For $s=0$ we will need the following lemma, which uses results of section 3.

 \begin{lemm}
 \label{lem:2.18}
   For $u\in \mathbb S^{d-1}$,$t>0$ we denote $\Delta_{u}^t=\{y\in \mathbb P^{d-1}\ ;\ |\langle u,y\rangle | <t\}$. Then, for any $\varepsilon, t>0$, $x,y\in \mathbb S^{d-1}$
 \[\limsup_{n\rightarrow \infty} \int \widetilde{\delta}^{\varepsilon} ( g \cdot x,  g \cdot y) d\mu^n (g)\leq \frac{2^{\varepsilon}}{t^{\varepsilon}} \widetilde{\delta}^{\varepsilon} (x,y)+2^{\varepsilon}\nu (\Delta_{x}^{t}).\]
 In particular, for any $\varphi \in H_{\varepsilon} (\mathbb S^{d-1})$ the sequence $\widetilde{P}^n \varphi$ is equicontinuous.
 \end{lemm}
\begin{proof}
 We write
\begin{align*}
\int \widetilde{\delta}^{\varepsilon} ( g \cdot x,  g \cdot y) d\mu^n (g)&=\int 1_{[1/t,\infty[} (\frac{|S_{n}|}{|S_{n}x|})  \widetilde{\delta}^{\varepsilon} (S_{n}\cdot x, S_{n} \cdot y) d\mathbb P(\omega)\\
& \ \ +\int 1_{]0,1/t[} (\frac{|S_{n}|}{|S_{n} x|})  \widetilde{\delta}^{\varepsilon} (S_{n} \cdot  x, S_{n} \cdot y) d\mathbb P(\omega).
\end{align*}
 Using Lemma 2.11 we have
 $ \widetilde{\delta}^{\varepsilon}(S_{n} \cdot  x, S_{n} \cdot y)\leq (2 \frac{|S_{n}|}{|S_{n}x|}  \widetilde{\delta} (x,y))^{\varepsilon}$.

   On the other hand, using Theorem 3.2, we know that  $\displaystyle\mathop{\lim}_{n\rightarrow \infty} \frac{|S_{n}|}{|S_{n}x|}=\frac{1}{|\langle z^{*}(\omega),x\rangle |}$ where $z^{*}(\omega) \in \mathbb P^{d-1}$ has law $\nu$. Hence
 \begin{align*}
 \limsup_{n\rightarrow \infty} \int  \widetilde{\delta}^{\varepsilon}( g \cdot x, g \cdot y) d\mu^n (g) &\leq \frac{2^{\varepsilon}}{t^{\varepsilon}}  \widetilde{\delta}^{\varepsilon} (x,y)+ 2^{\varepsilon} \mathbb P \{|\langle z^{*}(\omega), x\rangle| <t\}\\
 &=\frac{2^{\varepsilon}}{t^{\varepsilon}}  \widetilde{\delta}^{\varepsilon} (x,y)+2^{\varepsilon}\nu (\Delta_{x}^{t}).
 \end{align*}

  We have $|\widetilde{P}^n \varphi (x)- \widetilde{P}^n \varphi (y)|\leq [\varphi]_{\varepsilon} \int  \widetilde{\delta}^{\varepsilon}( g \cdot x,  g \cdot y)d\mu^n (g)$. From Theorem 3.2, we know that $\nu$ is proper, hence $\displaystyle\mathop{\lim}_{t\rightarrow 0} \nu (\Delta_{x}^{t})=0$. Then, for $x$ fixed we use the above estimation of  $\displaystyle\mathop{\limsup}_{n\rightarrow \infty} \int  \widetilde{\delta}^{\varepsilon} ( g \cdot x,  g \cdot y) d\mu^n (g)$  to choose $t$ sufficiently small in order to  get the continuity  of $\displaystyle\limsup_{n\rightarrow \infty} \int  \widetilde{\delta}^{\varepsilon} ( g \cdot x,  g \cdot y) d\mu^n (g)$.This gives that  $|\widetilde{P}^n \varphi (x)-\widetilde{P}^n \varphi (y)|$ depends continuously of $y$, hence the equicontinuity of the sequence $\widetilde{P}^n \varphi$.
  \end{proof}

\textit{Proof of  Theorem \ref{thm:2.17}}.
As observed in remark 2 after Proposition 2.13, if $s>0$ for any $\varphi \in C (\mathbb S^{d-1})$ the set  of functions $\{(\widetilde Q^{s})^{n} \varphi$; $n\in \mathbb N\}$ is equicontinuous. In view of Lemma \ref{lem:2.18}, this is also valid if $s=0$. Hence we can use here Lemma 2.12 . This gives the first convergence.
 Since $\widetilde{\pi}^{s} (x)$ is $\widetilde Q^{s}$-stationary, its projection on $\mathbb P^{d-1}$ is equal to the unique $Q^{s}$-stationary measure $\pi^{s}$, hence $\supp \widetilde{\pi}^{s} (x) \subset \widetilde{\Lambda}(T)$. On the other hand $\supp \widetilde{\pi}^{s} (x)$ is closed and $T$-invariant, hence contains a $T$-minimal set.

  In case I, $\widetilde{\Lambda}(T)$ is the unique minimal set, hence $\supp \widetilde{\pi}^{s} (x) =\widetilde{\Lambda} (T)$.
  Furthermore, if $\varphi \in C (\mathbb S^{d-1})$ is  $\widetilde Q^{s}$-invariant, the sets
\[
M^{-}=\{x\  ; \ \varphi(x)=\inf \{\varphi (y) ; y\in \mathbb P^{d-1}\}\},\ \
M^{+}=\{x\  ; \ \varphi(x)=\sup \{ \varphi (y) ; y\in \mathbb P^{d-1}\}\},\]
  are closed and $T$-invariant, hence they contain minimal sets. Since $\tilde{\Lambda} (T)$ is the unique minimal set, $M^{+} \cap\ M^{-} \supset \tilde{\Lambda} (T)\neq \phi$, hence $\varphi$ is constant.

  Then, using Proposition 2.13 and Lemma 2.12, we get that there exists a unique stationary measure $\tilde{\pi}^{s}$. It follows that $\tilde{\pi}^{s}$ is symmetric with projection $\pi^{s}$ on $\mathbb P^{d-1}$ and $\tilde{\pi}^s=e^s \tilde{\nu}^s$.

  In case II, the restriction to  the convex envelope $\textrm{Co} (\Lambda_{+} (T))=\Phi$ of the projection on $\mathbb P^{d-1}$ is a $T$-equivariant  homeomorphism. If we denote by $i_{+}$ its inverse, we get that $i_{+}(\pi^{s})$ is the unique $\tilde Q^{s}$-stationary measure supported in $\Phi$. Hence $i_{+}(\pi^{s})=\pi_{+}^{s}$. Then $\pi_{+}^{s}$ and $\pi_{-}^{s}$ are extremal $\tilde Q^{s}$-stationary measures.

  Since the projection of $\widetilde{\pi}^{s} (x)$ on $\mathbb P^{d-1}$ is $\pi^{s}$, we can write
\[\widetilde{\pi}^{s} (x)=\int (p_{+}^{s} (x,y) \delta_{y}+p_{-}^{s} (x,y) \delta_{-y}) d\pi^{s}(y)=p_{+}^{s} (x) \pi_{+}^{s}+p_{-}^{s}(x) \pi_{-}^{s},\]
  where $p_{+}^{s} (x)=p_{+}^s (x,.)$ and $p_{-}^{s} (x)=p_{-}^s (x,.)$ are Borel functions of $y \in \tilde{\Lambda} (T)$ such that $p_{-}^{s} (x)+ p_{-}^{s} (x)=1$. Then $p_{+}^{s} (x) \pi_{+}^{s}$ is the restriction of $\widetilde{\pi}^{s} (x)$ to $\Lambda_{+} (T)$, hence is a $\tilde Q^{s}$- invariant measure. In view of the uniqueness of the stationary measure of $\tilde Q^{s}$ restricted to $\Lambda_{+}(T)$, we get that $p_{+}^{s} (x) \pi_{+}^{s}$ is proportional to $\pi_{+}^{s}$, $\textrm i.e. $\ \  $p_{+}^{s} (x)$ is independent of $y$, $\pi_{+}^{s}$\textrm{-a.e.} .

  Hence, the first assertion of the theorem implies that the only extremal $\tilde Q^{s}$-stationary measures are $\pi_{+}^{s}$ and $\pi_{-}^{s}$. The corresponding facts for $\nu_{+}^s$ and $\nu_{-}^s$ follow.
From the mean ergodic theorem in $C(\mathbb S^{d-1})$ and the equicontinuity property of $(Q^{s})^{n}$ we know that the operator defined by $\displaystyle\mathop{\lim}_{n\rightarrow \infty} \frac{1}{n} \displaystyle\mathop{\Sigma}_{0}^{n-1} (\tilde Q^{s})^{k}$ is the projection on the  space of $\widetilde Q^s$-invariant functions and is equal to  $p_{+}^{s}(x) \pi_{+}^{s} +p_{-}^{s} (x) \pi_{-}^{s}$. The continuity and the extremality of the $\widetilde Q^{s}$-invariant functions $p_{+}^{s} (x)$ and $p_{-}^{s} (x)$ follows. The corresponding facts for $e_{+}^s$ and $e_{-}^s$ follow as in the proof of Theorem 2.6.

  If we restrict the convergence of $\frac{1}{n} \displaystyle\mathop{\Sigma}_{0}^{n-1} (\tilde Q^{s})^{k} (\delta_{x})$  to $x\in\Phi$, in view of the fact that the restriction to $\Phi$ of the projection on $\mathbb P^{d-1}$ is a homeomorphism onto its image, we get  $p_{+}^{s} (x)=1$, $p^{s}_{-} (x)=0$ if $x\in \Phi$.

  Let us denote by $\tau$ the entrance time of $S_{n}(\omega)\cdot x$ in $\Phi \cup-\Phi$ and by $^{a}\mathbb E_{x}^{s}$ the expectation symbol associated with the Markov chain $S_{n} (\omega) \cdot x$ defined by $\widetilde Q^{s}$.
Using theorem 2.6 we get $^{a}\mathbb E_{x}^{s} (1_{\Phi \cup-\Phi} (S_{\tau}\cdot x))=1$. Since $p_{+}^{s} (x)$
is a $\tilde Q^{s}$-invariant function $p_{+}^{s} (S_{n} \cdot x)$ is a martingale, hence $p_{+}^{s} (x)=$ $^{a}\mathbb E_{x}^{s} (p_{+}^{s}(S_{\tau}\cdot x))$.
Since $p_{+}^{s} (x)=1$ on $\Phi$ and $p_{+}^{s} (x) =0$ on $-\Phi$ we get $p_{+}^{s} (x)=$ $^{a}\mathbb E_{x}^{s} (1_{\Phi} (S_{\tau}\cdot x))$, hence the  stated interpretation for $p_{+}^s (x)$.

  As in Lemma 2.10, we verify  that the function $\varphi (u)=p(s) \int \langle u,u'\rangle_{+}^s d{^{*}\nu_{+}^s} (u')$ on $\mathbb S^{d-1}$ satisfies $\widetilde P^s \varphi=k(s) \varphi$, hence the function $\frac{\varphi}{e^s}$ satisfies $\widetilde Q^s(\frac{\varphi}{e^s})=\frac{\varphi}{e^s}$. By duality, cases II for $\mu$ and $\mu^{*}$ are the same, hence there are two minimal $T^{*}$-invariant subsets
  $\Lambda_{+} (T^{*})$ and $\Lambda_{-}(T^{*})=-\Lambda_{+}(T^{*}).$
   On the other hand, the set
\[\Lambda_{+} (T)=\{u \in \mathbb S^{d-1}:
\langle u,u'\rangle\geq0 \ \textrm{  for any } \ u'\in   \Lambda_{+} (T^{*}) \}\]
  is non trivial, closed, $T$-invariant and has non zero interior, hence $\widehat{\Lambda}_{+} (T)$ contains either $\Lambda_{+} (T)$ or $\Lambda_{-} (T)$ and has trivial intersection with one of then. We can assume $\widehat{\Lambda}_{+}(T) \supset \Lambda_{+} (T)$. Then, for $u\in \Lambda_{+} (T)$ and any $u' \in \Lambda_{+} (T^{*})$, we have $\langle u,u'\rangle_{+}=\langle u,u'\rangle$, hence,
\[\varphi (u)=p(s)\int |\langle u,u'\rangle|^s d^{*}\nu^s_{+} (u')=e^s (u),\]
i.e.,$\frac{\varphi}{e^s}=1$ on $\Lambda_{+}(T)$. Also we have $\langle u,u'\rangle_{+}=0$ for $u\in \Lambda_{-}(T)$, $u'\in \Lambda_{+}(T^{*})$. Since $\widetilde Q (\frac{\varphi}{e^s})=\frac{\varphi}{e^s}  \geq 0$ and $\frac{\varphi}{e^s}=1$ on $\Lambda_{+}(T)$, we conclude from above that $p_{+}^s=\frac{\varphi}{e^s}$, hence we get the last formula and last assertions with $e^{s}_{-}(u)=p(s) \int \langle u,x\rangle ^{s}_{+} d{^{*}\nu^{s}_{-}}(x)=p^{s}_{-} (u) e^{s}(u)$. \hfill $\square$

   From above we know that if $s\geq 0$ and
   $\varphi \in C(\mathbb S^{d-1})$ the sequence
   $(\widetilde Q^s)^n\varphi$ is equicontinuous. Lemma 2.12 reduces the discussion of the behaviour of $(\widetilde Q^s)^n \varphi$ to the existence of eigenvalues $z$ of $\widetilde Q^s$ with $|z|=1$. In this direction we have the following

 \begin{coro}

 \label{cor:2.19}
 For $s\in I_{\mu}$, the equation $\widetilde Q^s \varphi=e^{i\theta} \varphi$ with $e^{i\theta}\neq 1$, $\varphi \in C (\mathbb S^{d-1})$ has a non trivial solution only in case I. In that case $e^{i\theta}=-1$, $\varphi$ is antisymmetric, satisfies $\varphi (g \cdot x)=-\varphi (x)$ on $\supp\mu \times \widetilde{\Lambda} (T)$  and is uniquely defined up to a coefficient.
   \end{coro}

 \begin{proof}
    We observe that, since $\varphi$ satisfies $\widetilde Q^s \varphi=e^{i\theta}\varphi$, the function $\varphi'$ defined by $\varphi' (x)=\varphi(-x)$ satisfies also $\widetilde Q^s \varphi' = e^{i\theta} \varphi'$. Then $\varphi+\varphi'$ is symmetric and defines a function $\overline{\varphi}$ in $C(\mathbb P^{d-1})$ with $Q^s \displaystyle\mathop{\varphi}^{-}=e^{i\theta}\displaystyle\mathop{\varphi}^{-}$. If $e^{i\theta}\neq 1$, Theorem 2.6 gives $\displaystyle\mathop{\varphi}^{-}=0$, i.e $\varphi$ is antisymmetric.

   Furthermore, in case II, the restriction of $\varphi$ to $\widetilde{\Lambda}_{+} (T)$ satisfies the same equation and the projection of $\widetilde{\Lambda}_{+} (T)$ on $\Lambda (T)$ is an equivariant homeomorphism. Then Theorem 2.6 gives a contradiction. Hence if $\varphi \in C(\mathbb S^{d-1})$ satisfies $\widetilde Q^s \varphi=e^{i\theta} \varphi$, then we are in case I. Also, passing to absolute values as in the proof of Theorem 2.6 we get $\widetilde Q^s |\varphi|=|\varphi|$, $|\varphi|=$cte. Furthermore by strict convexity we have on $\supp \mu\times \widetilde{\Lambda} (T)$, $\varphi ( g \cdot x)=e^{i\theta}\varphi(x)$, hence $\varphi^2 ( g \cdot x)=e^{2 i\theta}\varphi^2 (x)$. Since $\varphi^2$ is symmetric and satisfies $\widetilde Q^s \varphi^2=e^{2i\theta} \varphi^2$, we get $e^{2i\theta}=1$, i.e $e^{i\theta}=-1$; in particular $\varphi ( g \cdot x)=-\varphi(x)$ on $\supp\mu \times \widetilde{\Lambda} (T)$. If $\varphi'\in C (\mathbb S^{d-1})$ satisfies also $\widetilde Q^s \varphi'=-\varphi'$, we get from above $\widetilde Q^s \frac{\varphi'}{\varphi}=\frac{\varphi'}{\varphi}$, $\frac{\varphi'}{\varphi}(-x)=\frac{\varphi'}{\varphi} (x)$, hence $\varphi'$ is proportional to $\varphi$.
   \end{proof}

\section{Laws of large numbers and spectral gaps}

  Here we develop section 2 in a quantitative direction. A martingale construction plays an essential role in this study. The renewal theorems of section 4 are based on Theorem 3.2 below. Also the Doeblin-Fortet inequality in  Corollary 3.21 is used in section 5 to show the homogeneity at infinity of the stationary measure for an affine random walk on $V$. Theorem A in the introduction follows directly from Corollaries 3.19, 3.20.

\subsection{Notation}

  As in section 2, we assume that condition i-p is valid for the semigroup $T=[\supp\mu]$. If $d=1$ we assume that $T=[\supp\mu]$ is non arithmetic. For $s\in I_{\mu}$ we consider the functions $q^{s}$ and $q_{n}^{s}$ $(n >0)$ on $\mathbb P^{d-1} \times G$, defined by
\[q^{s} (x,g)=\frac{1}{k(s)} \ \frac{e^s ( g \cdot x)}{e^s (x)} |g x|^{s}=q_{1}^{s} (x,g), \ \ q_{n}^{s} (x,g)=\frac{1}{k^n (s)}\ \frac{e^s ( g \cdot x)}{e^s (x)} |g x|^{s},\]
hence by the definition of  $e^s$ we have $\int q_{n}^s (x,g)$ $d\mu^n (g)=1$.

  We denote by $\mathbb P^{d-1}_{2}$ the flag manifold of planes and by $\mathbb P^{d-1}_{1,2}$ the manifold of contact elements on $\mathbb P^{d-1}$. Such a plane is defined up to normalisation by a 2-vector $x \wedge y\in \wedge^2 V$ and we can assume $|x\wedge y| =1$. We write $g(x\wedge y)=g x \wedge g y$. Also a contact element $\xi$ is defined by its origin $x\in \mathbb P^{d-1}$ and  a line through $x$. Hence we can write $\xi=(x, x\wedge y)$ where $|x|=|x \wedge y|=1$. The following additive cocycles of the actions of $G$ on $\mathbb P^{d-1}$, $\mathbb P^{d-1}_{2}, \mathbb P^{d-1}_{1,2}$ will play an essential role:
\[\sigma_{1} (g,x)= \log |g x|, \ \sigma_{2} (g,x \wedge y)=\log |g (x \wedge y)|,\ \sigma(g,\xi)=\log |g(x\wedge y)|-2 \log |g x|.\]

  In addition to the norm of $g$ we will need to use the quantity
\[\gamma(g)=\sup (|g|, |g^{-1}|)\geq 1.\]
 Clearly, for any $x\in V$, with $|x|=1$, we have $-\log\  \gamma(g) \leq \log |g x| \leq \log\  \gamma (g)$.

  For a finite sequence $\omega=(g_{1}, g_{2}, \cdots, g_{n})$ we  write
\[S_{n} (\omega)=g_{n}\cdots  g_{1}\in G, \ \ q_{n}^{s} (x,\omega)=\displaystyle\mathop{\Pi}_{k=1}^{n} q^s (S_{k-1}\cdot x, g_{k}).\]
 We denote by $\Omega_{n}$ the space of finite sequences $\omega=(g_{1}, g_{2},\cdots,g_{n})$ and we write $\Omega=G^{\mathbb N}$.

   We observe that $\theta^s(x,g)=|g x|^s \frac{e^s( g \cdot x)}{e^s (x)}$ satisfies  the cocycle relation $\theta^s(x, gg')=\theta^s (g'\cdot x,g) \theta^s (x,g')$, hence $q_{n}^s (x,\omega)=\frac{1}{k^n (s)} \theta^s (x, S_{n}(\omega))$.

\begin{defi}{}
\label{def:3.1}
   We denote $\mathbb{Q}_{x}^s \in M^1 (\Omega)$ the limit of the projective system of probability measures $q_{k}^s (x,.) \mu^{\otimes k}$ on $\Omega_{k}$. We write $\mathbb{Q}^s=\int \mathbb{Q}_{x}^s d\pi^s (x)$ where $\pi^s$ is the unique $Q^s$-stationary measure on $\mathbb P^{d-1}$.
\end{defi}

  We recall that Theorem 2.6 implies that $\pi^{s}$ is not supported by an hyperplane.

The corresponding expectation symbol will be written $\mathbb E_{x}^{s}$ and the shift on $\Omega$ will be denoted by $\theta$. We write also $\mathbb E^s (\varphi)=\int \mathbb E_{x}^s (\varphi) d \pi^s (x)$.
   The path space of the Markov chain defined by $Q^{s}$ is a factor space of $^a\Omega=\mathbb P^{d-1} \times \Omega$,
and the corresponding shift on $^a\Omega$ will be written $^a\theta$ with $^a\theta (x,\omega) = (g_{1} (\omega)\cdot x, \theta \omega)$. Hence ($\mathbb P^{d-1}\times \Omega$,  $^a\theta)$ is a skew product over $(\Omega, \theta)$. The projection on $\Omega$ of the Markov measure $^a\mathbb{Q}_{x}^{s}=\delta_{x}\otimes \mathbb{Q}_{x}^{s}$ is $\mathbb{Q}_{x}^{s}$, hence $^a\mathbb{Q}^{s}=\int \delta_{x} \otimes \mathbb{Q}_{x}^{s} d\pi^{s}(x)$ projects on $\mathbb{Q}^{s}$. The uniqueness of the $Q^{s}$-stationary measure $\pi^{s}$ implies the ergodicity of the $^a\theta$-invariant measure $^a\mathbb{Q}^{s}$, hence $\mathbb{Q}^{s}$ is also $\theta$-invariant and ergodic.

  If $s=0$, the random variables $g_{k}(\omega)$ are i.i.d with law $\mu$ and $\mathbb{Q}^{\circ}=\mathbb P=\mu^{\otimes \mathbb N}$. Here, under condition i-p, we extend the results of \cite{28} to the case $s\geq 0$, in particular we construct a suitable kernel-valued martingale with contraction properties as in \cite{17}. This will allow us to prove strong forms of the law of large numbers for $S_{n} (\omega)$ and to compare the measures $\mathbb Q_{x}^s$ when $x$ varies. Then we can deduce the simplicity of the dominant Lyapunov exponent of $S_{n}(\omega)$ under the $\theta$-invariant probability $\mathbb Q^s$ for $s\geq 0$. Spectral gap properties for twisted convolution operators on the projective space and on the unit sphere will follow.

\subsection{A martingale and the equivalence  of $\mathbb{Q}_{x}^{s}$ to $\mathbb{Q}^{s}$}

  When convenient we identify $x\in \mathbb P^{d-1}$ with one of its representatives $\widetilde{x}$ in $\mathbb S^{d-1}$. We recall that the Markov kernel $^{*}Q^{s}$ is defined by $^{*}Q^{s} \varphi=\frac{1}{k(s) {^{*}e^{s}}} {^{*}P^{s}} (\varphi {^{*}e^{s}})$ where $^{*}P^{s} \varphi (x)=\int \varphi ( g \cdot x) |g x|^s d\mu^* (g)$, $^{*}P^{s} ({^{*}e^{s}})=k(s) ({^{*}e^{s}})$ and $^{*}Q^{s}$ has a unique stationary measure $^{*}\pi^{s}$. Furthermore we have $^{*}\pi^{s}={^{*}e^{s}} {^{*}\nu^{s}}$ where $^{*}\nu^{s} \in M^{1} (\mathbb P^{d-1})$ is the  unique solution of $^{*}P^{s} ({^{*}\nu^{s}})=k(s) ({^{*}\nu^{s}})$. We denote by $m$ the unique rotation invariant probability measure on $\mathbb P^{d-1}$.

\begin{theo}
\label{thm:3.2}
 Let $\Omega'\subset \Omega$ be the (shift-invariant) Borel subset of elements $\omega \in \Omega$ such that $S_{n}^{*} (\omega).m$  converges to a Dirac measure $\delta_{z^*} (\omega)$. Then $g_{1}^{*} \cdot z^*(\theta \omega)=z^*(\omega)$, $\mathbb{Q}^{s}(\Omega')=1$, the law of $z^*(\omega)$ under $\mathbb{Q}^{s}$ is $^{*}\pi^{s}$ and $^{*}\pi^{s}$ is proper.

   In particular if $\omega\in \Omega'$ and  $|\langle x,z^*(\omega) \rangle \langle y, z^* (\omega)\rangle |\neq 0$,  then
\[\displaystyle\mathop{\lim}_{n\rightarrow \infty} \delta (S_{n} (\omega)\cdot x, S_{n}(\omega)\cdot y)=0.\]
If $\omega \in \Omega'$ and $\xi=(x,x \wedge y)\in \mathbb P_{1,2}^{d-1}$

\[\lim_{n\rightarrow \infty} \frac{|S_{n}(\omega) x|}{|S_{n}(\omega)|}=|\langle z^* (\omega), x \rangle|,\ \displaystyle\mathop{\lim}_{n\rightarrow \infty} S_{n}^{*} \cdot m=\delta_{z^{*}(\omega)}.\]
In particular, if $\langle z^{*}(\omega), x \rangle \neq 0$ then $\displaystyle\mathop{\lim}_{n\rightarrow \infty} \sigma (S_{n},\xi)=-\infty$.

  Also, for any $x\in \mathbb P^{d-1}$  $\mathbb{Q}_{x}^{s}$ is equivalent to $\mathbb{Q}^{s}$ and $\frac{d \mathbb{Q}_{x}^{s}}{d \mathbb{Q}_{y}^{s}} (\omega)=\left | \frac{\langle z^*(\omega), x\rangle }{\langle  z^* (\omega), y\rangle }\right |^{s} \frac{e^{s} (y)} {e^{s}(x)}$.
  \end{theo}

  The proof of Theorem 3.2 is based on the following lemmas, in particular on the study of a kernel-valued martingale.

\begin{lemm}
\label{lem:3.3}
  Assume $z^{*}\in \mathbb P^{d-1}$ and $u_{n}\in G$ is a sequence such that  $\lim_{n\rightarrow \infty} u_{n}^*
 \cdot m=\delta_{z^*}$. Then, for any $x,y \in \mathbb P^{d-1}$ with $|\langle z^*,x\rangle  \langle z^*,y\rangle|\neq 0$ we have  $\lim_{n\rightarrow \infty} \delta (u_{n}\cdot x, u_{n}\cdot y)=0$. If $\xi=(x,x \wedge y)\in \mathbb P_{1,2}^{d-1}$ and $\langle z^{*}, x\rangle \neq 0$ then
\[\lim_{n\rightarrow \infty} \frac{|u_{n} x|}{|u_{n}|}=|\langle z^*,x\rangle|, \  \lim_{n\rightarrow \infty} \sigma (u_{n}, \xi)=-\infty.\]

  These convergences are uniform on any compact subset  on which $\langle z^*,x\rangle $ do not vanish.
  \end{lemm}

\begin{proof}
  We denote by $e_{i} (1\leq i\leq d)$ an  orthonormal basis of $V$, by $\bar e_{i}$ the projection of $e_{i}$ in $\mathbb P^{d-1}$,
 by $\bar{A}^+$ the set of diagonal matrices  $a= \textrm{diag}(a^1,a^2,\cdots,a^d)$ and $a^1\geq a^2\geq \cdots\geq a^d>0$. We write $u_{n}=k_{n} a_{n} k'_{n}$ with $a_{n} \in \bar{A}^+, \ \ k_{n}, \ k'_{n} \in O(d)$. Then, for $x\in \mathbb P^{d-1}$,
\[|u_{n} x|^2=|a_{n}k'_{n}x|^2=\Sigma^{d}_{1}(a_{n}^i)^2 |\langle k'_{n} x, e_{i}\rangle |^2.\]
  Also, $u^*_{n}\cdot m=(k'_{n})^{-1} a_{n}\cdot m$ converges to $z^* \in \mathbb P^{d-1}$, which implies
\[\lim_{n\rightarrow \infty} a_{n}\cdot m=\delta_{\overline{e_{1}}},\ \  \lim_{n\rightarrow  \infty} (k'_{n})^{-1}.\bar e_{1}=z^*.\]
  In particular, if $i>1$, we have $a^i_{n}=o(a^1_{n})$ and
\[\lim_{n\rightarrow \infty} |\langle k'_{n} x,\ e_{1}\rangle |=|\langle z^*, x\rangle |\neq 0.\]
  It follows that $|u_{n}x| \sim a^1_{n} |\langle z^*, x \rangle|$. Since $|u_{n}|=a^1_{n}$, we get $\displaystyle\mathop{\lim}_{n\rightarrow \infty} \frac{|u_{n}x|}{|u_{n}|}=|\langle z^*, x\rangle |$, as asserted. We  get also, if $|\langle  y, z^*\rangle|\neq 0$, $|u_{n} y| \sim a^1_{n} |\langle z^*, y\rangle|$.

  On the exterior  product space $\wedge^2 V$ there exists an
  $O(d)$-invariant scalar product such that on any decomposable 2-vector  $x\wedge y$  :
\[|x \wedge y|^2 = |x|^2 |y|^2 - |\langle x , y \rangle|^2.\]

  For $x,y\in \mathbb P^{d-1}$ and corresponding $\tilde x, \tilde y \in \mathbb S^{d-1}$ we write $|x \wedge y|=|\tilde x \wedge \tilde y |$.
  Then on $\mathbb P^{d-1}$, there is an associated distance $\delta_{1}$ given by  $\delta_{1} (x,y)=| x \wedge  y|$ and we have $\frac{1}{2} \delta \leq \delta_{1}\leq \delta$. We observe that $\delta_{1} (u_{n}\cdot x, u_{n} \cdot  y)=\frac{|u_{n} x \wedge u_{n}y|}{|u_{n} x| \ |u_{n} y|}$.
 Also
 \[|u_{n} x \wedge u_{n} y|^2=\displaystyle\mathop{\Sigma}_{i>j} (a^i_{n} \ a^j_{n})^2\  |\langle k_{n} (\widetilde x\wedge \widetilde y), e_{i} \wedge  e_{j}\rangle|^2 \leq \frac{d (d-1)}{2} (a^1_{n} a^2_{n} |\widetilde x \wedge \widetilde y|)^2.\]

  It follows
\[\delta_{1} (u_{n}\cdot x,\ u_{n}\cdot y) \leq \left(\frac{d (d-1)}{2}\right)^{1/2}\  \frac{a^1_{n} a^2_{n}}{|u_{n}x|\ |u_{n}y|} \ |\widetilde x\wedge \widetilde y|.\]
  Since $|u_{n}x| \sim a^1_{n} |\langle z^*, x\rangle|, \ |u_{n}y| \sim a^1_{n} |\langle z^*, y\rangle|$ and $a^2_{n}=o (a^1_{n}), \langle z^{*}, x\rangle \langle z^{*},y\rangle\neq 0$, we get
\[\lim_{n\rightarrow \infty} \delta_{1} (u_{n}\cdot x,\ u_{n}\cdot  y)=0.\]
  It follows, for any $x,y \in \mathbb P^{d-1}$, that  $\lim_{n\rightarrow \infty} \delta (u_{n}\cdot  x,\ u_{n}\cdot y)=0$.
Also, since $a_{n}^2=o (a_{n}^1)$, and $\langle z^{*},x\rangle\neq 0$ we get $\lim_{n\rightarrow \infty} \sigma (u_{n}, \xi)=-\infty$.

  The above calculations imply the uniformities in the convergences. \end{proof}

\begin{lemm}
\label{lem:3.4}
  Assume $\nu_{n} \in M^1 (\mathbb P^{d-1})$ is a sequence such that $\nu_{n}$ is relatively compact in variation, and each $\nu_{n}$ is proper. Let $u_{n}\in G$ be a sequence such that $u_{n}^{*}\cdot \nu_{n}$ converges weakly to $\delta_{z*}$ $(z^{*}\in \mathbb P^{d-1})$. Then for any proper $\rho \in M^1(\mathbb P^{d-1})$, $u_{n}^{*}\cdot\rho$ converges weakly to $\delta_{ z^*}$.
  \end{lemm}

\begin{proof}
We can assume, in variation, $\lim_{n\rightarrow \infty} \nu_{n}=\nu_{0}$ where $\nu_{0}$ is proper. Also we can assume, going to subsequences, that $u_{n}^{*}$ converges to a quasi-projective map of the form $u^*$, defined and continuous outside a projective subspace $H\subset \mathbb P^{d-1}$. Let $\varphi \in C(\mathbb P^{d-1})$ and denote
\[I_{n}=(u_{n}^*\cdot \nu_{n}) (\varphi)-(u^* \cdot  \nu_{0}) (\varphi)=(\nu_{n}-\nu_{0}) (\varphi \circ u_{n}^*)+\nu_{0}(\varphi \circ u_{n}^*-\varphi \circ u^*) .\]
  The first term is bounded by $|\varphi|\  \|\nu_{n}-\nu_{0}\|$, hence it converges to zero. Since $\nu_{0} (H)=0$ and $\varphi \circ u_{n}^{*}$ converges to $\varphi\circ u^*$ outside $H$,  we can use dominated convergence for the second term,
$\displaystyle\mathop{\lim}_{n\rightarrow \infty} \nu_{0} (\varphi \circ u_{n}^*-\varphi\circ u^*)=0$; hence $\displaystyle\mathop{\lim}_{n\rightarrow \infty} I_{n}=0$. Then $u_{n}^{*} \cdot  \nu_{n}$ converges to $u^{*} \cdot  \nu_{0}$ weakly. In particular $u^{*} \cdot  \nu_{0}= \delta_{z^{*}}$, hence $u^{*} \cdot  y=z^{*}$, ${\nu_{0}\textrm{-a.e}}$. Since $\nu_{0}(H)=0$, we have $u^{*} \cdot  y=z^*$ on $\mathbb P^{d-1} \setminus H$.

  Since $\rho$ is proper $u^* \cdot  \rho=\delta_{z^*}$, hence $\displaystyle\mathop{\lim}_{n\rightarrow \infty} u_{n}^* \cdot  \rho=\delta_{z^*}$.
  \end{proof}

\begin{lemm}
\label{lem:3.5}
  For $x,y\in \mathbb P^{d-1}$ the total variation measure of $\mathbb Q_{x}^s-\mathbb Q_{y}^s$ is bounded by $B \delta^{\bar s} (x,y) \mathbb Q^s$.
Furthermore, there exists $c(s)>0$ such that, for any $x\in \mathbb P^{d-1}$ we have  $ \mathbb Q_{x}^{s} \leq c(s)\  \mathbb Q^{s}$.
\end{lemm}

\begin{proof}
  We write $q^s_{n}(g)=\int q^s_{n} (x,g) d\pi^s (x)$ and we observe that for any measurable $\varphi$ depending on the first $n$ coordinates,
$$\int \varphi (\omega) d\mathbb Q^s (\omega)=\int q^s_{n} (S_{n} (\omega)) \varphi (\omega) d\mu^{\otimes n} (\omega).$$
  Also, $|(\mathbb Q^s_{x}-\mathbb Q^s_{y}) (\varphi)|\leq \int |q_{n}^s (x, S_{n})-q_{n}^s (y, S_{n}) | \ |\varphi (\omega)| d\mu^{\otimes n}$. Using Lemma 2.11 we have for any $g\in G$
$$|q_{n}^s (x,g)-q_{n}^s (y,g)|\leq b_{s} \frac{|g|^s}{k^n (s)} \delta^{\bar s} (x,y).$$
   Since $\pi^s$ is not supported by an hyperplane we can use Lemma 2.7, hence for some $b>0$
\[ q_{n}^s (g)\geq b \frac{|g|^s}{k^n (s)} \ \textrm{ and }\  |q_{n}^s (x,g)-q_{n}^s (y,g)|\leq \frac{b_{s}}{b} q_{n}^s (g) \delta^{\bar s} (x,y).\]

  It follows:
$$|(\mathbb Q_{x}^s-\mathbb Q_{y}^s) (\varphi)| \leq \frac{b_{s}}{b} \delta^{\bar s} (x,y) \int |\varphi (\omega)| d\mathbb Q^s (\omega),$$
  hence the first conclusion with $B=\frac{b_{s}}{b}$.

  Integrating with respect to $\pi^s$ we get, since $\delta (x,y)\leq \sqrt 2$, $\mathbb Q_{x}^s \leq (1+B(\sqrt 2)^{\bar s}) \mathbb Q^{s}$ hence the second formula with $c(s)=1+B(\sqrt 2)^{\bar s}$.
  \end{proof}

\begin{lemm}
\label{lem:3.6}
  We consider the positive kernel $\kappa_{x}^s$ from $\mathbb P^{d-1}$ to $\mathbb P^{d-1}$ given by
  $\kappa_{x}^s=\frac{|\langle x,\cdot\rangle  |^s }{e^s (x)}{^{*}\nu^s}$. Then:
\[\int g^*\cdot\kappa_{ g \cdot x}^s q^s(x,g) d\mu (g)=\kappa_{x}^s, \quad \kappa_{x}^s (1)=\frac{1}{e^s (x)} \int |\langle x,y\rangle |^s d{^{*}\nu^s} (y)=p(s)\in ]0,1]\]
and $x\rightarrow \kappa_{x}^s$ is continuous in variation.

  In particular $S_{n}^{*} \cdot \kappa^s_{_{S_{_{n}}\cdot x}}$ is a bounded martingale with respect to $\mathbb Q_{x}^s$ and the natural filtration.
  \end{lemm}

\begin{proof}
  We consider the $s$-homogeneous measure $\lambda^s$ on $\breve{V}$ defined by $\lambda^s={^{*}\nu^s} \otimes \ell^s$. By definition of $^{*}\nu^{s}, \int g^{*} \lambda^{s} d\mu (g)=k(s) \lambda^s$. Then the Radon measure $\lambda_{v}^s$ defined by $\lambda_{v}^s=|\langle v,\cdot \rangle|^s \lambda^s$  satisfies $\int g^* \lambda_{g v}^s d\mu (g)= k(s) \lambda_{v}^s$. This can be written, by definition of $\kappa_{x}^s$ and $q^s (x,g)$,
$$\int g^* \cdot \kappa_{ g \cdot x}^s q^s (x,g) d\mu (g)=\kappa_{x}^s.$$
  The martingale property of $S_{n}^{*}\cdot \kappa^s_{_{S_{n} \cdot  x}}$ follows.

  Furthermore since $p(s) e^s (x)$ is equal to $\int |\langle x, y \rangle |^s d\  {^{*}\nu^s} (y)$, Lemma 2.10 gives
$$\kappa_{x}^s (1) =\frac{1}{e^s (x)} \int |\langle  x,y\rangle |^s d\ {^{*}\nu^s} (y)=p(s)\in ]0,1].$$
  The continuity in variation of $x\rightarrow \kappa_{x}^s$ follows from the definition.
  \end{proof}
\begin{lemm}
\label{lem:3.7}
  Let $\rho\in M^1 (\mathbb P^{d-1})$, $\mathcal H$ be the set of projective subspaces $H$ of minimal dimension such that $\rho(H)>0$. Then the subset of elements $H\in \mathcal H$  such that $\rho (H)=\sup\{\rho (L) \ ;\ L\in  \mathcal H\}$ is finite and non void. Furthermore, there exists $\varepsilon_{\rho}>0$ such that  for any $H\in \mathcal H$ :
$$ \rho (H)=c_{\rho} \ \ \hbox{\rm or}\ \ \rho(H) \leq c_{\rho}-\varepsilon_{\rho},$$
  where $c_{\rho}= \sup \{ \rho (L)\   ; \  L\in \mathcal H\}$.
  \end{lemm}
\begin{proof}
   If $H$, $H'\in \mathcal H$, $H\neq H'$, then $\dim H\cap H'< \dim H$, hence $\rho (H \cap H')=0$. Then, for any $\beta>0$, the cardinality of the set of elements $H\in \mathcal H$ with $\rho (H) \geq \beta$ is bounded by $\frac{1}{\beta}$. The first assertion follows. Assume the second assertion is false. Then there exists a sequence $H_{n} \in \mathcal H$ with $\frac{c_{\rho}}{2}< \rho (H_{n})< c_{\rho}$, $\displaystyle\mathop{\lim}_{n\rightarrow \infty}  \rho (H_{n})=c_{\rho}$, and $\rho (H_{n})\neq \rho (H_{m})$ if $n\neq m$. This contradicts the fact that the cardinality of the sequence $H_{n}$ is at most $\frac{2}{c_{\rho}}$.
   \end{proof}

\begin{lemm}
\label{lem:3.8}
  Assume that the Markovian kernel $x\rightarrow \nu_{x} \in M^1 (\mathbb P^{d-1})$ is continuous in variation and satisfies
$$\nu_{x}=\int q^s (x,g) g^* \cdot  \nu_{ g \cdot x} d\mu (g)\ .$$
  Let $\mathcal H_{p,r}$ the set of finite unions of $r$ distinct subspaces of dimension $p$ and  let $h$ be the function $h(x)=\sup \{\nu_{x} (W)$; $W \in \mathcal H_{p,r}\}$. Then $h$ is continuous and the set
\[\{x: h(x)=\sup\{h(y), y\in\mathbb P^{d-1}\}\},\]
  is closed and $[\supp \mu]$-invariant.
  \end{lemm}
\begin{proof}
  If $W\in \mathcal H_{p,r}$ is fixed the function $x\rightarrow \nu_{x} (W)$ is continuous since $|\nu_{x} (W)-\nu_{y}(W)|\leq \|\nu_{x}-\nu_{y}\|$. This implies $|h(x)-h(y)|\leq \|\nu_{x}-\nu_{y}\|$, hence the continuity of $h$.

  We have for any $W\in \mathcal H_{p,r}$
\[\nu_{x}(W)=\int q^s (x,g) \nu_{ g \cdot x} ((g^*)^{-1} W) d\mu (g).\]
Hence: $h(x)\leq \int q^s (x,g) h( g \cdot x) d\mu (g)$. Then, as in Lemma 2.9, $X$ is $[\supp \mu]$-invariant and closed.
  \end{proof}

\begin{lemm}
\label{lem:3.9}
  Let $\nu_{x}$ be as in Lemma 3.8. Then for any $x\in \mathbb P^{d-1}$, $\nu_{x}$ is proper.
  \end{lemm}
\begin{proof}
  We write $\pi_{x}=\frac{\nu_{x}}{\nu_{x} (1)}$, denote by $\mathcal H_{k}$   the set of projective subspaces  of dimension $k$  and
  \[\mathcal H=\cup_{k\geq 0} \mathcal H_{k}, \ \
d(x)=\inf \{\textrm{dim} H ; H\in \mathcal H, \pi_{x}(H)>0\}, \]
  $m(x)= \sup \{\pi_{x} (H) ; H\in \mathcal H, dim H=d(x)\}$,
 $\mathcal W(x)=\{H\in \mathcal H ; \pi_{x} (H)=m (x)\}$.

  Lemma 3.7 implies  that the set $\mathcal W(x)$ has finite cardinality $n(x)>0$. Also we denote $p=\inf \{d(x) ; x\in \mathbb P^{d-1}\}$, $h_{p}(x)=\sup \{\pi_{x} (H) ; H\in \mathcal H_{p}\}$.

  Lemma 3.8 shows that $h_{p}(x)$ reaches its maximum $\beta$ on a closed $[\supp \mu]$-invariant subset $X\subset \mathbb P^{d-1}$. Hence on $\Lambda ([\supp \mu])$ we have $h_{p} (x)=\beta=m(x)$.
It follows $d(x)=p$ on $\Lambda ([\supp \mu])$. The relation $n(x) m(x) \leq 1$ implies $n(x) \leq \frac{1}{\beta}$ on $\Lambda ([\supp \mu])$.

  Let $r= \sup \{n(x) ; x\in \Lambda ([\supp \mu])\}$ and denote $h_{p,r} (x)=\sup \{\pi_{x}(W) ; W\in \mathcal H_{p,r}\}$. Then Lemma 3.8 implies  $h_{p,r} (x)=r \beta$ on $\Lambda ([\supp \mu])$. Since $m(x)=\beta$, this relation implies $n(x)=r$ on $\Lambda ([\supp \mu])$. Let $W(x)=\cup \{H ; H\in \mathcal{W}(x)\}$ and let us show the local constancy of the function $ W(x)$. Using Lemma 3.7 we get
$$\beta (x)=\sup \{\pi_{x}(H) ; H\in \mathcal H_{p}, H \notin \mathcal W (x)\}<\beta.$$
  Let $x\in \Lambda([\supp \mu])$, $U_{x}=\{y ; \|\pi_{y} - \pi_{x}\| < \beta-\beta (x)\}$ and $H_{y} \in \mathcal H_{p}$ with $\pi_{y} (H_{y})=\beta$.
Then,
\[\beta-\pi_{x} (H_{y})=\pi_{y}(H_{y})- \pi_{x} (H_{y}) \leq \|\pi_{y}-\pi_{x}\|<\beta-\beta(x).\]
  Hence $\pi_{x} (H_{y})> \beta (x)$ and, by definition of $\beta(x)$ we get $H_{y}\in \mathcal W(x)$ for any $y\in U_{x}$. Since  $\pi_{x}$ is continuous in variation, $U_{x}$ is a neighbourbood of $x$, hence $W(x)$ is locally constant. Since $\Lambda([\supp \mu])$ is compact, $W=\displaystyle\mathop{\cup W (x)}_{x\in \Lambda([\supp \mu])}$ is a finite union of subspaces.

  On the other hand, the relations
\[r \beta=\pi_{x} (W(x))=\int q^s (x,g) g^*\cdot  \pi_{ g \cdot x} (W (x)) d\mu (g)\ \  , \ \  r\beta \geq (g^*\cdot \pi_{ g \cdot x}) (W(x))\]
  imply that, for any $x\in \Lambda([\supp \mu])$,
$ r \beta = g^*\cdot \pi_{ g \cdot x} (W(x))\  \mu$\textrm{-a.e}.
By definition of $W( g \cdot x)$, we get, $(g^*)^{-1} W(x)=W ( g \cdot x)$  $\mu$\textrm{-a.e}.
Hence, for any $g\in \supp \mu$,
$ (g^*)^{-1} (W (x))=W( g \cdot x).$

  The relation $(g^*)^{-1} (W)=\displaystyle\mathop{\cup (g^*)^{-1}}_{x\in \Lambda ([\supp \mu])}  (W (x))= \displaystyle\mathop{\cup W ( g \cdot x)}_{x\in\Lambda ([\supp \mu])}$ shows that $W$ is $[\supp \mu^*]$-invariant. Then condition i-p implies $W=\mathbb P^{d-1},\ \  r=1, p=d-1, \ \  d(x)=d-1, m(x)=1$ for any $x\in \mathbb P^{d-1}$, hence the Lemma.
  \end{proof}

\textit{Proof of  Theorem 3.2}.
  We use the Markov kernel $\pi^s_{x}$ defined by $\pi^s_{x}=\frac{\kappa_{x}^s}{\nu_{x}^s (1)}$ with $\kappa_{x}^s
$ given in Lemma 3.6.

  Then we have the harmonicity equation,
\[\pi_{x}^s=\int q^s (x,g) g^* \cdot \pi^s_{ g \cdot x} d\mu (g),\]
  and the continuity in variation of $\pi_{x}^s$. The above equation implies that the sequence of kernels $S_{n}^* (\omega)\cdot \pi^s_{S_{n}(\omega)\cdot x}$ is a $\mathbb{Q}^s_{x}$-martingale with respect to the natural filtration on $\Omega$. Since for $x\in \mathbb P^{d-1}$, $\pi_{x}^s \in M^1 (\mathbb P^{d-1})$ we can apply the martingale convergence theorem.

  Since, by Lemma 3.9, $\pi^s_{x}$ is proper and, by definition, $x\rightarrow \pi^s_{x}$ is continuous in variation, we can use the same method as in  \cite{28} ; because of i-p condition, the martingale $S_{n}^* (\omega). \pi^s_{S_{n}(\omega)\cdot x}$ converges $\mathbb{Q}_{x}^s$\textrm{-a.e} to a Dirac measure. Then, using Lemma 3.4, $S_{n}^*\cdot m$ converges $\mathbb{Q}_{x}^s$\textrm{-a.e}.\ to  a Dirac measure $\delta_{z^*(\omega)}$. Then, from above, for any $x\in \mathbb P^{d-1}$,
\[\mathbb{Q}^s_{x} (\Omega')=1\ \  , \ \  z^* (\mathbb{Q}^s_{x})=\pi_{x}^s.\]
  Hence the law of $z^*(\omega)$ under $\mathbb{Q}^s_{x}$ is $\pi^s_{x}$. It follows, by integration,
\[\mathbb{Q}^s (\Omega')=1\ \ , \  \  z^*(\mathbb{Q}^s)=\int z^* (\mathbb{Q}^s_{x}) d\pi^s (x)=\int \pi_{x}^s d\pi^s (x).\]
  In view of the formulae  for $\nu_{x}^s$, $\pi_{x}^s$ and the relation $\pi^s=\frac{e^s \nu^s}{\nu^s (e^s)}$, we define ${^{*}\pi}^{s}$ by $z^*(\mathbb{Q}^s)={^{*}\pi^s}$. Then Lemma 3.9 and the definition of $\pi^s_{x}$ give that ${^{*}\pi^s}$ is proper. The relations:

  $\displaystyle\mathop{\lim}_{n\rightarrow \infty} \delta (S_{n}(\omega)\cdot x,\  S_{n} (\omega) \cdot  y)=0,  \ \  \displaystyle\mathop{\lim}_{n\rightarrow \infty} \frac{|S_{n}(\omega) x|}{|S_{n}(\omega)|}=|\langle z^*(\omega), x\rangle |,\ \ \displaystyle\mathop{\lim}_{n\rightarrow \infty} \sigma (S_{n},\xi)=-\infty,
$
  follow from the geometrical Lemma 3.3, since $S^*_{n}\cdot m$ converges to $\delta_{z^* (\omega)}$.

  Using Lemma 3.5 we know that $\mathbb{Q}^s_{x}$ is absolutely continuous with respect to $\mathbb{Q}^s$. We calculate $\frac{d \mathbb{Q}^s_{x}}{d \mathbb{Q}^s} (\omega)$ as follows.

  By definition of $\mathbb{Q}^s_{x}$ and $\mathbb{Q}^s$,
$$\mathbb E_{x}^s (\frac{d \mathbb{Q}^s_{x}}{d \mathbb{Q}^s} (\omega) |g_{1}, \cdots , g_{n})=\frac{q_{n}^s (x, S_{n}(\omega))}{\int q_{n}^s (y,S_{n} (\omega)) d\pi^s (y)} .$$
  Furthermore
$$\frac{q^s_{n} (x, S_{n} (\omega))}{q_{n}^s (y, S_{n}(\omega))}=\frac{|S_{n}(\omega)x|^s}{|S_{n}(\omega) y|^s} \ \frac{e^s (S_{n}(\omega)\cdot x)}{e^s (S_{n} (\omega) \cdot  y)}\ \frac{e^s (y)}{e^s (x)}.$$
  The martingale convergence theorem gives
$$\frac{d\mathbb{Q}^s_{x}}{d\mathbb{Q}^s} (\omega)=\displaystyle\mathop{\lim}_{n\rightarrow \infty} \frac{q^s_{n} (x,S_{n}(\omega))}{\int q^s_{n} (y ,  S_{n} (\omega)) d\pi^s (y)}.$$
  Using the relation $\displaystyle\mathop{\lim}_{n\rightarrow \infty} \frac{|S_{n} (\omega) x|}{|S_{n}(\omega)|}=|\langle z(\omega),x\rangle|$,\ \  if  \ \ $\omega\in \Omega'$, we get
\[\lim_{n\rightarrow \infty} \frac{q_{n}^s (x,S_{n} (\omega))}{q_{n}^s (y, S_{n}(\omega))}=\left |\frac{\langle z^*(\omega),\rangle}{\langle z^*(\omega),y\rangle}\right |^s \frac{e^s (y)}{e^s (x)}.\]
  Hence $\frac{d\mathbb{Q}^s_{x}}{d\mathbb{Q}^s} (\omega)=\frac{|\langle z^*(\omega),x\rangle|^s}{e^s (x)}\left [ \int \frac{|\langle z^*(\omega),y\rangle|^s}{e^s (y)} d\pi^s (y) \right ]^{-1}$. Since,  from above $\pi^s$ is proper and the $\mathbb{Q}^s$-law of $z^*(\omega)$ is $\pi^s$, we have for any $x\in \mathbb P^{d-1}$,
$|\langle z^*(\omega), x\rangle| >0$, ${ \mathbb{Q}^s}$\textrm{-a.e}

We conclude  that  $\mathbb{Q}^s_{x}$ is equivalent to $\mathbb{Q}^s$ because
  $\frac{d\mathbb{Q}^s_{x}}{d\mathbb{Q}^s} (\omega) > 0$,  $\mathbb{Q}^s$-a.e.  Also, using the formulae above we have
  \[\frac{d\mathbb{Q}^s_{x}}{d\mathbb{Q}^s_{y}} (\omega)= \left |\frac{\langle z^*(\omega),x\rangle}{\langle z^* (\omega), y\rangle}\right |^s \frac{e^s (y)}{e^s (x)}.\]
  \hfill $\square$

\subsection{The law of large numbers for $\log |S_{n}(\omega) x|$ with respect to $\mathbb{Q}^s$}

  Here, by derivatives of a function $\varphi$ at the boundaries of an interval $[a,b]$ we will mean finite half derivatives i.e. we write

\centerline{$\varphi' (a)=\varphi' (a_{+}) \in \mathbb R \ \ , \ \ \varphi'(b)=\varphi' (b_{-}) \in \mathbb R$ .}
\begin{theo}
\label{thm:3.10}
  Let $\mu\in M^1 (G)$, $s\in I_{\mu}$.  Assume $[\supp \mu]$ satisfies condition i-p, and $\log \gamma (g)$ is $\mu$-integrable. Assume also that $|g|^s \log\  \gamma (g)$ is $\mu$-integrable and write
\[L_{\mu}(s)=\int \log |g x| q^s (x,g) d\pi^s (x)\ d\mu (g).\]
Then, for any $x\in \mathbb P^{d-1}$, $\mathbb Q^s\textrm{-a.e}$, we have
\[\lim_{n\rightarrow \infty} \frac{1}{n} \log |S_{n} (\omega) x|=\lim_{n\rightarrow \infty} \frac{1}{n} \log |S_{n} (\omega)|=L_{\mu}(s).\]
  This convergence is valid in $\mathbb L^1 (\mathbb{Q}^s)$ and in $\mathbb L^1 (\mathbb{Q}^s_{x})$ for any $x\in \mathbb P^{d-1}$. Furthermore $k(t)$ has a continuous derivative on $[0,s]$ and  if $t\in [0,s]$, $x\in \mathbb P^{d-1}$,
\[L_{\mu} (t)=\frac{k' (t)}{k(t)}=\displaystyle\mathop{\lim}_{n\rightarrow \infty} \frac{1}{n k^n (t)} \int |g x|^t \log |g x| d\mu^n (g)=\displaystyle\mathop{\lim}_{n\rightarrow \infty} \frac{1}{n}\frac{\int |g|^t \log |g| d\mu^n (g)}{\int |g|^td\mu^n(g)}.\]
  In particular if $\alpha >0$ satisfies $k(\alpha)=1$, then $k'(\alpha)>0$.
  \end{theo}

\begin{proof}
  We consider the function $f(x,\omega)$ on $^a\Omega$ defined by $f(x,\omega)=\log |g_{1}(\omega)x|$. If $|x|=1$,
we have $-\log|g^{-1}|\leq \log |g x|\leq \log |g|$,\ \  hence $f(x,\omega)$ is $^a Q^s$-integrable.
Moreover $\int f(x,\omega) d^a\mathbb{Q}(x,\omega)=\int q^s (x,g)$  $\log |g x| d\nu^s (x) d\mu (g)=L_{\mu} (s)$, and
\[\sum^{n-1}_{0} (f\circ^a\theta^k) (x,\omega)=\log |S_{n}(\omega) x|.\]
As mentioned above $^a\mathbb{Q}^s$ is $^a\theta$-ergodic, hence we get using Birkhoff's theorem,
\[\lim_{n\rightarrow \infty} \frac{1}{n} \log |S_{n}(\omega) x|=L_{\mu}(s), \ ^a\mathbb{Q}^s\textrm{-a.e.}\]
  On the other hand, we can apply the subadditive ergodic theorem to the sequence
$\log |S_{n} (\omega)|$ and  to the ergodic system $(\Omega, \theta, \mathbb{Q}^s)$.

  This implies that there exists $L(s)\in \mathbb R$ such that, $\mathbb{Q}^s$\textrm{-a.e.} and in $\mathbb L^1(\mathbb Q^s)$, the sequence $\frac{1}{n} \log |S_{n} (\omega)|$ converges to $L(s)$. We know, using Theorem 3.2, that for fixed $x$ and $\mathbb{Q}^s$\textrm{-a.e.},
\[\lim_{n\rightarrow \infty} \frac{|S_{n} (\omega) x|}{|S_{n}(\omega)|}=|\langle z^{*}(\omega), x\rangle |,\]
  and furthermore the law of $z^{*}(\omega)$ under $\mathbb{Q}^s$ is proper. Hence, for fixed $x$ we have
$$|\langle z^{*}(\omega), x\rangle|>0,\ \  \mathbb{Q}^s\textrm{-a.e}.$$
 Then for fixed $x\in \mathbb P^{d-1}$ and $\mathbb{Q}^s\textrm{-a.e.}$
$$\displaystyle\mathop{\lim}_{n\rightarrow \infty} \frac{1}{n} \log |S_{n}(\omega) x|=\displaystyle\mathop{\lim}_{n\rightarrow \infty} \frac{1}{n} \log |S_{n} (\omega)|=L_{\mu}(s).$$
  Using Lemma 3.6, since  $\mathbb{Q}^s_{x}\leq c(s)\  \mathbb{Q}^s$, this convergence is also valid $\mathbb{Q}^s_{x}$\textrm{-a.e}. Hence by definition of $^a\mathbb{Q}^s$,  we have  $L(s)=L_{\mu}(s)$. The first assertion follows.

  In order to get the $\mathbb L^1$-convergences, we observe that Fatou's Lemma gives
$$\displaystyle\mathop{\liminf}_{n\rightarrow \infty} \frac{1}{n} \int \log |S_{n}(\omega) x| d\mathbb{Q}^s (\omega)\geq L_{\mu} (s).$$
  On the other hand, the subadditive ergodic theorem gives
$$\displaystyle\mathop{\lim}_{n\rightarrow\infty} \frac{1}{n} \int \log |S_{n} (\omega)| \ \ d\mathbb{Q}^s (\omega)=L(s)=L_{\mu} (s).$$
  Since $|S_{n} (\omega) x| \leq |S_{n} (\omega)|$ if $|x|=1$, these two relations imply, for every $x\in \mathbb P^{d-1}$,
$$\displaystyle\mathop{\lim}_{n\rightarrow\infty} \frac{1}{n} \int \log |S_{n} (\omega) x| d\mathbb{Q}^s (\omega)=L_{\mu} (s).$$
  Now we write
$$\frac{1}{n} |\log | S_{n} (\omega) x| - L(s)| \leq \frac{1}{n} (\log |S_{n} (\omega)| -\log |S_{n}(\omega) x|)+\frac{1}{n}  |\log |S_{n} (\omega)|-L(s) |.$$
  From  the above calculation, the integral of the first term converges to zero. The subadditive ergodic theorem implies the same for the second term.

  Hence $\displaystyle\mathop{\lim}_{n\rightarrow+\infty} \int |\log| S_{n}(\omega) x| - L(s)| d\mathbb{Q}^s(\omega)=0$. Since $\mathbb{Q}_{x}^s \leq c(s) \mathbb{Q}^s$, this convergence is also valid in $\mathbb L^1(\mathbb{Q}_{x}^s)$ for any fixed $x$. This gives the second assertion, in particular
$$L_{\mu} (s)=\displaystyle\mathop{\lim}_{n\rightarrow \infty} \frac{1}{n} \int \log |S_{n}(\omega) x| d\mathbb{Q}_{x}^s (\omega)=\displaystyle\mathop{\lim}_{n\rightarrow +\infty} \frac{1}{n} \mathbb E_{x}^s (\log |S_{n} (\omega) x|).$$
  The above limit can be expressed as follows.

  Let $\varphi$ be a continuous function on $\mathbb P^{d-1}$. Then Theorem 2.6 implies
$$\displaystyle\mathop{\lim}_{n\rightarrow \infty} \mathbb E^s_{x}(\varphi (S_{n} (\omega)\cdot x))=\displaystyle\mathop{\lim}_{n\rightarrow \infty} (Q^s)^n \varphi (x)=\pi^s(\varphi),$$
  uniformly in $x\in \mathbb P^{d-1}$.

  Hence $L_{\mu}(s) \pi^{s} (\varphi)=\displaystyle\mathop{\lim}_{n\rightarrow \infty} \frac{1}{n} \mathbb E^s_{x} (\varphi (S_{n}(\omega)\cdot x) \log |S_{n}(\omega) x|)$.

  In particular with $\varphi=\frac{1}{e^s}$ and any $x$ we have

\centerline{$ e^{s}(x) L_{\mu} (s) =\displaystyle\mathop{\lim}_{n\rightarrow \infty} \frac{1}{nk^n(s)} \int |g x|^s \log |g x| d\mu^n (g)$. }

  We denote $v_{n}(s)=\int  |g x|^s d\mu^n (g)$ and we observe that $v'_{n}(s)=\int |g x|^s \log |g x|$  $d\mu^n (g)$.

  Using Theorem 2.6, we get
$\displaystyle\mathop{\lim}_{n\rightarrow \infty} \frac{v_{n}(s)}{k^n (s) e^{s}(x)}=\pi^s \left(\frac{1}{e^s}\right)=1.$

  Then the above formula for $L_{\mu}(s)$ reduces to  $L_{\mu}(s)=\displaystyle\mathop{\lim}_{n\rightarrow \infty}\frac{1}{n} \ \frac{v'_{n}(s)}{v_{n}(s)}$.

  On the other-hand,
$\frac{1}{n} \log v_{n} (s) = \frac{1}{n} \int_{0}^{s} \frac{v'_{n}(t)}{v_{n}(t)} dt,\ \  \  \displaystyle\mathop{\lim}_{n\rightarrow \infty} \frac{1}{n} \log v_{n}(s)=\log k(s).$

  The convexity of $\log v_{n}$ on $[0,s]$ gives

\centerline{$ \frac{v'_{n} (0)}{v_{n} (0)} \leq \frac{v'_{n}(t)}{v_{n}(t)} \leq   \frac{v'_{n}(s)}{v_{n}(s)}.$}

  Then the convergence of  the sequences $\frac{1}{n} \ \frac{v'_{n}(0)}{v_{n}(0)}$ and $\frac{1}{n} \ \frac{v'_{n}(s)}{v_{n}(s)}$ implies that the sequence $\frac{1}{n} \ \frac{v'_{n}(t)}{v_{n}(t)}$ is uniformly bounded on $[0,s]$. On the other hand, H\"older inequality implies the $\mu$-integrability of $|g|^t \log |g|$ if $t\in [0,s]$, and the bound:
$$\int |g|^t |\log|g||d\mu(g) \leq \left(\int |\log|g|| d\mu (g)\right)^{\frac{s-t}{s}} \left(\int |g|^s |\log|g|| d\mu (g)\right)^{t/s}.$$
  Hence, as above, we have the convergence of $\frac{1}{n}\ \frac{v'_{n}(t)}{v_{n}(t)}$, to $L(t)$. Then  dominated convergence gives the convergence of $\frac{1}{n} \log v_{n}(s)$ to $\log k(s)=\int_{0}^{s} L(t) dt$.

  We  have $L(t)=\int \log |gx| q^t  (x,g) d\mu (g) d\pi^t (x)$, and the continuity of $q^s$, $\pi^s$ given by Theorem 2.6. Then the bound of $\int |g|^t |\log| g||d\mu(g)$ given above implies the continuity of $L(t)$ on $[0,s]$.
The  integral expression of $\log k(s)$ in terms of $L(t)$ implies that $k$ has a derivative and $\frac{k'(t)}{k(t)}=L(t)$ if $t\in [0,s]$. This gives the first part of the last relation in the theorem.  In order to get the rest, we consider $u_{n} (t)=\int |g|^t d\mu^n (g)$ and write for $t\in [0,s]$
$$\frac{u'_{n}(t)}{u_{n}(t)}=\frac{\int |g|^t \log |g| d\mu^n (g)}{\int |g|^t d\mu^n (g)}.$$
  The convergence of $\frac{1}{n} \log u_{n} (t)$ to $\log k(t)$ and the convexity of the functions $\log u_{n}(t)$,  $\log k(t)$ give for $t\in [0,s]$,
$$\frac{k'(t_{-})}{k(t)}\leq \displaystyle\mathop{\liminf}_{n\rightarrow \infty} \frac{1}{n} \ \frac{u'_{n}(t)}{u_{n}(t)} \leq \displaystyle\mathop{\limsup}_{n\rightarrow \infty}    \frac{u'_{n}(t)}{u_{n}(t)}\leq \frac{k'(t_{+})}{k(t)}.$$
  Since if $t\in ]0,s[$, we have $k'(t_{-})=k'(t_{+})=k'(t)$, we get $\displaystyle\mathop{\lim}_{n\rightarrow \infty} \frac{1}{n} \  \frac{u'_{n} (t)}{u_{n}(t)}=\frac{k'(t)}{k(t)}$ if $t<s$.

  Furthermore, by continuity we have $\displaystyle\mathop{\lim}_{n\rightarrow \infty} \frac{1}{n} \  \frac{u'_{n} (s)}{u_{n}(s)}=\frac{k'(s_{-})}{k(s)}$. Now the rest of the formula follows from the expression of $\frac{u'_{n}(t)}{u_{n}(t)}$ given above.
The relation $L_{\mu} (\alpha)>0$ follows from the formula $L_{\mu}(t) =\frac{k'(t)}{k(t)}$ and the strict convexity of  $\log \ k(t)$.
\end{proof}

\subsection{Lyapunov exponents and spectral gaps}

  We begin with a more general situation than above. As special cases, it contains the Markov chains  on $\mathbb P^{d-1}$ considered in section 2. In particular, simplicity of the dominant Lyapunov exponent given by Theorem 3.17 below will be a simple consequence of their special properties and of the general Proposition 3.11 below. For $s=0$, this result was shown in \cite{28} under condition i-p. For the use of the Zariski closure as a tool to show condition i-p see \cite{22}, \cite{45}. We give corresponding notation.

  Let $X$ be a compact metric space, $\mathcal C (X,X)$ the semigroup of continuous maps of $X$ into itself which is  endowed with   the topology of uniform convergence. We denote by $ g \cdot x$ the action of $g\in \mathcal C(X,X)$ on $x\in X$ and we consider a closed subsemigroup $\Sigma$ of $\mathcal C(X,X)$. Let $\mu$ be a probability measure on $\Sigma$ and let $q(x,g)$ be a continuous non negative function on $X \times \supp \mu$ such that $\int q(x,g)d\mu(g)=1$. We will denote by $(X, q\otimes \mu, \Sigma)$ this set of data and we will say that $(X, q\otimes \mu, \Sigma)$ is a Markov system on $(X, \Sigma)$. We write $\Omega=\Sigma^{\mathbb N}$, we denote by $\Omega_{n}$ the set of finite sequences of length $n$ on $\Sigma$ and for $\omega=(g_{1},g_{2},\cdots,g_{n})$ in $\Omega_{n}$, $x\in X$, we write $q_{n} (x,\omega)=\displaystyle\mathop{\Pi}_{1}^{n} q(S_{k-1}\cdot x, g_{k})$ where $S_{n}=g_{n}\cdots g_{1}$, $S_{0}=Id$.

  We define a probability measure $\mathbb Q_{x}^{n}$ on $\Omega_{n}$ by $\mathbb Q_{x}^{n}=q_{n}(x,\cdot ) \mu^{\otimes n}$ and we denote  by $\mathbb Q_{x}$ the probability measure on $\Omega$ which is the projective limit of this system. If $\nu$ is a probability measure on $X$ we set $\mathbb Q_{\nu}=\int \mathbb Q_{x} d\nu(x)$. We will consider the extended shift $^{a}\theta$ on ${^{a}\Omega}=X\times \Omega$ which is defined by $^{a}\theta (x,\omega)=(g_{1}\cdot x, \theta \omega)$, and also the Markov chain on $X$ with kernel $Q$ defined by $Q\varphi (x)=\int \varphi (g\cdot x) q(x,g) d\mu (g)$. Clearly, when endowed with the corresponding shift, the space of paths of this Markov chain is a factor system of $(X\times \Omega, ^{a}\theta)$. If $\pi$ is a $Q$-stationary measure on $X$, the measure $\mathbb Q_{\pi}$ on $\Omega$ is shift-invariant and $^{a}\mathbb Q_{\pi}=\int \delta_{x} \otimes \mathbb Q_{x} d\pi (x)$ is $^{a}\theta$-invariant. In this situation we will say that $(X, q \otimes \mu, \Sigma, \pi)$ is a stationary Markov system. If $\pi$ is  $Q$-ergodic, then $^{a}\mathbb Q_{\pi}$ is $^{a}\theta$-ergodic and $\mathbb Q_{\pi}$ is $\theta$-ergodic. We will denote by $\mathbb E_{x}, \mathbb E_{\pi}$ the corresponding expectations symbols.

  In particular we will consider below Markov systems of the form $(X,q\otimes \mu, T)$ where $\Sigma=T \subset GL(d,\mathbb R)$ $(d>1)$, and also extensions of them. We can extend the action of $g\in T$ to $X \times \mathbb P^{d-1}$ by $g(x,v)=( g \cdot x, g\cdot v)$ and define a new Markov chain with kernel $Q_{1}$ by $Q_{1}\varphi (x,v)=\int \varphi ( g \cdot x, g\cdot v) q (x,g) d\mu (g)$. Given a $Q$-stationary probability measure $\pi$,  we will denote by $C_{1}$ the compact convex set of probability measures on $X\times \mathbb P^{d-1}$ which have projection $\pi$ on $X$. The same considerations apply if $\mathbb P^{d-1}$ is replaced by $\mathbb P^{d-1}_{2}$, the manifold of 2-planes or $\mathbb P^{d-1}_{1,2}$, the manifold of contact elements in $\mathbb P^{d-1}$. Then we define similarly the kernels $Q_{2}, Q_{1,2}$ and  the convex sets $C_{2}, C_{1,2}$.

  Since $S_{n}=g_{n} \cdots g_{1}$ and the $g_{k}$ are $\mathbb Q_{\pi}$-stationary random variables where $\mathbb Q_{\pi}$ is $\theta$-invariant and ergodic, the Lyapunov exponents of the product $S_{n}$ are well defined as soon as
\[\int \log |g_{1} (\omega)| d \mathbb Q_{\pi} (\omega)\ \ \textrm{ and }\ \ \int \log |g_{1}^{-1}(\omega)| d\mathbb Q_{\pi}(\omega)\]
are finite (see \cite{46}). In particular the two largest ones $\gamma_{1}$ and $\gamma_{2}$ are given by
\[\gamma_{1}=\displaystyle\mathop{\lim}_{n\rightarrow \infty} \frac{1}{n} \int \log |S_{n}(\omega)| d\mathbb Q_{\pi}(\omega), \ \
\gamma_{1}+\gamma_{2}=\displaystyle\mathop{\lim}_{n\rightarrow \infty} \frac{1}{n} \int \log |\wedge^2 S_{n} (\omega)| d\mathbb Q_{\pi}(\omega).\]

   In order to study the values of $\gamma_{1}, \gamma_{2}$ we need to consider the above Markov operators $Q_{1}, Q_{2}, Q_{1,2}$ and the convex sets $C_{1}, C_{2}, C_{1,2}$ of corresponding stationary measures. We denote $q(g)=\displaystyle\sup_{x\in \mathbb P^{d-1}} |q(x,g)|$ and we assume $\int \log \gamma (g) q(g) d\mu (g)<\infty$. For $\eta_{1}\in C_{1}$, we will write $I_{1}(\eta_{1})=\int \sigma_{1} (g,v) d\eta_{1} (x,v) d\mu (g)$, and similarly with $\eta_{2}\in C_{2}, \eta_{1,2}\in C_{1,2}$ we define $I_{2}(\eta_{2}), I_{1,2}(\eta_{1,2})$ where we use the cocycles $\sigma_{1}, \sigma_{2}, \sigma$ defined at the beginming of this section. The following result will be a basic tool in this subsection.

\begin{prop}
\label{prop:3.11}
With the above notation, let $T$ be a closed subsemigroup of $GL(d,\mathbb R)$,  $(X, q\otimes \mu, T, \pi)$ a stationary and uniquely ergodic Markov system. Assume that $S_{n}^{*}\cdot m$ converges $\mathbb Q_{\pi}\textrm{-a.e}$ to a Dirac measure $\delta_{z^*(\omega)}$ such that for any $v\in \mathbb P^{d-1},\langle v,z^* (\omega)\rangle\neq 0$  $\mathbb Q_{\pi}\textrm{-a.e}$. Assume that $\int \log \gamma (g) q(g) d\mu (g)$ is finite. Then we have $\gamma_{2}-\gamma_{1}=\sup\{I_{1,2}(\eta) ; \eta \in C_{1,2}\} <0$, and the sequence $\frac{1}{n} \displaystyle\sup_{x,v,v'} \mathbb E_{x} (\log \frac{\delta (S_{n} (\omega)\cdot v, S_{n}(\omega)\cdot v')}{\delta (v,v')})$ converges to $\gamma_{2}-\gamma_{1}<0$.
\end{prop}

  The proof  uses the same arguments as in  \cite{29}. The condition $\displaystyle\mathop{\lim}_{n\rightarrow \infty} S_{n}^*\cdot{m}=\delta_{z^*(\omega)}$ is satisfied in the examples of subsection 3.2 (see Theorem 3.2).

\begin{lemm}
\label{lem:3.12}
  Let $m_{p}$ be the natural $SO(d)$-invariant probability measure on the submanifold of $p$-decomposable unit multivectors $x=v_{1}\wedge v_{2}\wedge \cdots \wedge v_{p}$. Then there exists $c>0$ such that for any $u \in \textrm{End} V$ : $0<\log|\wedge^p
u|-\int \log |(\wedge^p u) (x)| dm_{p}(x) \leq c$.
\end{lemm}
\begin{proof}

  We write $u$ in polar form $u=kak'$ with $k, k' \in SO(d), a=\textrm{diag}(a_{1},a_{2}\cdots, a_{d})$ and $a_{1}\geq a_{2}\geq \cdots \geq a_{d}>0$. We write also $x=k'' \varepsilon^p$ with $k''\in SO(d), \varepsilon^p=e_{1}\wedge e_{2}\wedge \cdots \wedge e_{p}$, hence,
\[\int \log |\wedge^p u x| dm_{p} (x) = \int \log |\wedge^p a k' k'' \varepsilon^p|d\widetilde{m} (k'')=\int \log |\wedge^p  a k \varepsilon^p| d\widetilde{m} (k),\]
where $\widetilde{m}$ is the normalized Haar measure on $SO(d)$. Furthermore
\begin{align*}
|\wedge^p a k \varepsilon^p| &\geq \langle\wedge^p a k \varepsilon^p, \varepsilon^p\rangle| = |\wedge^p a| |\langle k\varepsilon^p, \varepsilon^p \rangle|,\\
\int \log |\wedge^p u x |dm_{p} (x) &\geq  \log |\wedge^p u|+\int \log |\langle k \varepsilon^p, \varepsilon^p\rangle| d\widetilde{m} (k)\\
&= \log |\wedge^p u|+\int \log |\langle x, \varepsilon^p\rangle| dm_{p} (x).
\end{align*}

  Hence it suffices to show the finiteness of the integral in the right hand side. But the set of unit decomposable $p$-vectors is an algebraic submanifold of the unit sphere of $\wedge^p V$ and $m_{p}$ is its natural Riemannian measure. Since the map $x\rightarrow  \langle x, \varepsilon^p \rangle^2$ is polynomial, there exists $q\in \mathbb N$, $c>0$ such that :
$c t^q\leq m_{p}\{x ;  \langle x ,\varepsilon^p\rangle ^2 \leq t\}\leq 1$.

  The push forward $\sigma$ of $m_{p}$ by this map is an absolutely continuous probability measure on [0,1] which satisfies  $\sigma(0,t)\geq c t^{q/2}$. Then
\[\int \log | \langle x, \varepsilon^p  \rangle | dm_{p} (x)=\int^1_{0} \log \ t \ d\sigma (t) >-\infty,\]
  since the integral $\int^1_{0} t^{q/2-1}  d t$ is finite for $q>0$.
  \end{proof}

  We recall from \cite{37} the following
\begin{lemm}
\label{lem:3.13}
  Let $(E, \theta, \nu)$ be a dynamical system where $\nu$ is a $\theta$-invariant probability measure, $f$ a $\nu$-integrable function. If  $\displaystyle\mathop{\Sigma}^{n}_{1} f \circ \theta^k$ converges $\nu\textrm{-a.e}$ to $-\infty$, then one has $\nu (f)<0$.
  \end{lemm}

\begin{lemm}
\label{lem:3.14}
Let $E$ be a compact metric space, $P$ a Markov kernel on $E$ which preserves the space of continuous functions on $E$, $C(P)$ the compact convex set of $P$-stationary measures. Then, for every continuous function $f$ on $E$, the sequence $\displaystyle\sup_{x\in E} \frac{1}{n} \displaystyle\mathop{\Sigma}^{n}_{1} P^k f(x)$ converges to $\sup\{\eta (f) ; \eta \in C(P)\}$. In particular, if for all $\eta, \eta' \in C(P)$ we have $\eta (f)=\eta'(f)$, then we have the uniform convergence, $\displaystyle\mathop{\lim}_{n\rightarrow \infty} \frac{1}{n} \displaystyle\mathop{\Sigma}^{n}_{1} P^k f(x)=\eta(f)$.
 \end{lemm}

\begin{proof}
  Let $J\subset \mathbb R$ be the set of cluster values of the sequences $\frac{1}{n} \displaystyle\mathop{\Sigma}_{0}^{n-1} (P^k f) (x_{n})$ with $x_{n}\in E$. We will show that the convex envelope of $J$ is equal to $\{\eta (f)\  ; \  \eta\in C(P)\}$. If the sequence $\frac{1}{n_{k}} \displaystyle\mathop{\Sigma}^{n_{k}}_{0} (P^i f) (x_{n_{k}})$ converges to $c\in \mathbb R$, we can consider the sequence of probability measures $\eta_{k}=\frac{1}{n_{k}} \displaystyle\mathop{\Sigma}^{n_{k}-1}_{0} P^{i} \delta_{x_{n_{k}}}$ and extract  a convergent subsequence with limit $\eta \in C(P)$. Then, since $f$ is continuous we have
\[\eta(f)=\lim_{k\rightarrow \infty} \frac{1}{n_{k}} \displaystyle\mathop{\Sigma}^{n_{k}-1}_{0} (P^{i} f) (x_{n_{k}})=c.\]
Conversely, if $\eta \in C(P)$ is ergodic, Birkhoff's theorem applied to the sequence $\frac{1}{n} \displaystyle\mathop{\Sigma}^{n-1}_{0} (P^{i}f) (x)$ gives $\eta $\textrm{-a.e},
\[\lim_{n\rightarrow \infty} \frac{1}{n} \displaystyle\mathop{\Sigma}^{n-1}_{0} (P^{i} f) (x)=\eta(f),\]
hence there exists $x\in E$ such that $\eta  (f)$ is the limit of $\frac{1}{n} \displaystyle\mathop{\Sigma}^{n-1}_{0} (P^{i} f) (x)$. If $\eta$ is not ergodic, $\eta$ is a barycenter of ergodic measures, hence $\eta (f)$ belongs to the convex envelope of $J$. In view of the above, this shows the first claim. Since $J$  is closed, the convex envelope of $J$ is a closed interval $[a,b]$, hence $b=\displaystyle\mathop{\lim}_{n\rightarrow \infty} \frac{1}{n} \displaystyle\sup_{x\in E} \displaystyle\mathop{\Sigma}^{n-1}_{0} (P^{i} f) (x)=\displaystyle\sup_{\eta\in C(P)} \eta (f)$.
  \end{proof}

\begin{lemm}
\label{lem:3.15}
  We have the formulae
\begin{align*}
\gamma_{1} &= \displaystyle\mathop{\lim}_{n\rightarrow \infty} \frac{1}{n} \displaystyle\sup_{x,v} \int \log |S_{n}(\omega) v| d\mathbb Q_{x}(\omega)=\displaystyle\sup_{\eta\in C_{1}} I_{1}(\eta),\\
\gamma_{1}+\gamma_{2}&=\displaystyle\mathop{\lim}_{n\rightarrow \infty} \frac{1}{n} \displaystyle\sup_{x,v,v'} \int \log |S_{n}(\omega) v \wedge S_{n}(\omega) v'| d\mathbb Q_{x} (\omega)=\displaystyle\sup_{\eta\in C_{2}} I_{2} (\eta).
\end{align*}
\end{lemm}

\begin{proof}
  We consider the Markov chain on $X\times \mathbb P^{d-1}$ with kernel $Q_{1}$ defined by
\[Q_{1}\varphi (x,v)=\int \varphi (g\cdot x, g \cdot v) q (x,g) d\mu (g),\]
and the function
$\psi (x,v)=\int \sigma_{1} (g,v) q(x,g) d\mu (g).$
We observe that
\[\int \sigma_{1}(S_{n} (\omega), v) d\mathbb Q_{x} (\omega)=\displaystyle\mathop{\Sigma}^{n-1}_{0} Q_{1}^k \psi (x,v),\]
and $\psi$ is continuous since $\int \log \gamma (g) q(g) d\mu(g)<\infty$. Also, since $\pi$ is the unique $Q$-stationary measure, any $Q_{1}$-stationary measure has projection $\pi$ on $X$. Then, using Lemma 3.14,we have
\[\sup_{\eta \in C_{1}} I_{1} (\eta)=\displaystyle\mathop{\lim}_{n\rightarrow \infty} \frac{1}{n} \displaystyle\sup_{x,v} \int \sigma_{1} (S_{n} (\omega), v) d\mathbb Q_{x} (\omega),\]
 which gives the second part of the first formula. In order to show the first part we consider $\eta \in C_{1}$ which is $Q_{1}
$-ergodic, the extended shift $\widetilde{\theta}$ on $X\times \mathbb P^{d-1} \times \Omega$ and the function $f(x,v,\omega)=\sigma_{1}(g_{1}(\omega),v)$. Then
\[\widetilde{\theta} (x,v,\omega)=(g_{1}\cdot x, g_{1}\cdot v, \theta \omega) \textrm{ and }  \sigma_{1} (S_{n} (\omega),v)=\displaystyle\mathop{\Sigma}_{0}^{n-1} f \circ \widetilde{\theta}^k (x,v,\omega).\]

  Also, $f(x,v,\omega)$ is $\widehat{\mathbb Q}_{\eta}$-integrable where $\widetilde{\mathbb Q}_{\eta}=\int \delta_{(x,v)} \otimes \mathbb Q_{x}d\eta (x,v)$. Using the subadditive ergodic theorem, we have
\[I_{1}(\eta)=\frac{1}{n}\int \sigma_{1} (S_{n}(\omega), v) d \mathbb Q_{x} (\omega) d\eta (x,v)\leq \displaystyle\mathop{\lim}_{n\rightarrow \infty} \frac{1}{n} \int \log |S_{n} (\omega)| d\mathbb Q_{\pi} (\omega)=\gamma_{1}.\]
We show now that for some $\eta \in C_{1}$ we have $\gamma_{1}=I_{1}(\eta)$.
Using Lemma 3.12, we know that
\[0\leq \log |S_{n}(\omega)|- \int \log |S_{n}(\omega) v| dm (v)\leq c,\]
hence, integrating with respect to $\mathbb Q_{\pi}$, we have
\[0\leq \int dm (v) \int (\log |S_{n}(\omega)|-\log |S_{n} (\omega) v| d\mathbb Q_{\pi} (\omega) \leq c.\]
  Then the sequence of non negative functions $h_{n}(v)$ on $\mathbb P^{d-1}$ given by
$$h_{n} (v)=\frac{1}{n} \int (\log |S_{n} (\omega)|- \log |S_{n} (\omega) v|) d\mathbb Q_{\pi} (\omega)$$
  satisfies $0\leq h_{n}(v) \leq \frac{c}{n}$ with $c$ given by Lemma 3.12, $\displaystyle\mathop{\lim}_{n\rightarrow \infty} \int h_{n}(v) dm (v)=0$. It follows that  we can find a subsequence $h_{n_{j}}$ such that $h_{n_{j}}(v)$ converges to zero $m$\textrm{-a.e}, hence
\[\gamma_{1}=\displaystyle\mathop{\lim}_{j\rightarrow \infty} \frac{1}{n_{j}} \int \sigma_{1} (S_{n_{j}}(\omega), v) d\mathbb Q_{\pi} (\omega), m-\textrm{a.e.}, \]
  in particular  this convergence is valid for some $v\in \mathbb P^{d-1}$. The sequence of probability measures $\frac{1}{n_{j}} \displaystyle\mathop{\Sigma}_{1}^{n_{j}} Q_{1}^k
(\pi \otimes \delta_{v})$ has a weakly convergent subsequence $\eta_{j}$ with a $Q_{1}$-invariant limit $\eta$. Furthermore, the function $\psi$ considered above is continuous, hence
\[\eta(\psi) =\displaystyle\mathop{\lim}_{j\rightarrow \infty} \frac{1}{n_{j}} \int \sigma_{1} (S_{n_{j}} (\omega) v) d\mathbb Q_{\pi} (\omega)=\gamma_{1},\  \gamma_{1}=I_{1}(\eta) =\displaystyle\sup_{\eta_{1}\in C_{1}} I_{1}(\eta_{1}).\]

  The same argument is valid for $\log |S_{n} (\omega) (v\wedge v')|$ with $m$ replaced by $m_{2}$, hence the second formula.
  \end{proof}

\begin{lemm}
\label{lem:3.16}
  For any $\eta \in C_{1}$, we have $\gamma_{1}=I_{1} (\eta)$ .
\end{lemm}
\begin{proof}
  As in the proof of Lemma 3.15 we have $I_{1}(\eta)=\displaystyle\mathop{\lim}_{n\rightarrow \infty} \frac{1}{n} \sigma_{1} (S_{n} (\omega), v)$, $\mathbb Q_{\eta}$\textrm{-a.e.}, hence the existence of $v\in \mathbb P^{d-1}$ such that  $\mathbb Q_{\pi} $ \textrm{-a.e.},
$$I_{1} (\eta) =\displaystyle\mathop{\lim}_{n\rightarrow \infty} \frac{1}{n} \log |S_{n}(\omega) v| .$$
  Then, using  Theorem 3.2 and Lemma 3.15 we get, $\mathbb Q_{\pi} $\textrm{-a.e.}\ \,
\[\lim_{n\rightarrow \infty} \frac{1}{n} \log \frac{|S_{n}(\omega) v|}{|S_{n}(\omega)|}=\lim_{n\rightarrow \infty} \frac{1}{n} \log |\langle z^* (\omega), v\rangle| =0,\]
  since $\langle z^*(\omega), v\rangle \neq 0$, $\mathbb Q_{\pi}$\textrm{-a.e}.
Hence $I_{1}(\eta)=\displaystyle\mathop{\lim}_{n\rightarrow \infty} \frac{1}{n} \int \log |S_{n}(\omega)| d\mathbb Q_{\pi} (\omega)=\gamma_{1}$.
  \end{proof}

\textit{Proof of proposition 3.11}
  We have $\gamma_{2}-\gamma_{1}=(\gamma_{1}+\gamma_{2})-2\gamma_{1}$, $\gamma_{1}=I_{1} (\eta_{1})$ for any $\eta_{1} \in C_{1}$ and $\gamma_{1}+\gamma_{2}=\displaystyle\sup_{\eta_{2}\in C_{2}} I_{2} (\eta_{2})$. Using the theorem of Markov-Kakutani for the inverse image of $\eta_{2} \in C_{2}$ in $C_{1,2}$ we know that any $\eta_{2}\in C_{2}$ is the projection of some $\eta_{1,2} \in C_{1,2}$, hence $\gamma_{1}+\gamma_{2}=\displaystyle\sup_{\eta_{1,2}\in C_{1,2}} I_{2} (\eta_{1,2})$. If $\eta'_{1}$ is the projection of $\eta_{1,2}$ on $\mathbb P^{d-1}$, we have $I_{1,2} (\eta_{1,2})=I_{2} (\eta_{1,2})-2 I_{1}(\eta'_{1})$ and from Lemma 3.16, $\gamma_{1}=I_{1} (\eta'_{1})$. It follows $\gamma_{2}-\gamma_{1}=\displaystyle\sup_{\eta_{1,2}\in C_{1,2}} I_{1,2} (\eta_{1,2})$.

  Since $I_{1,2} (\eta_{1,2})$ depends continuously of $\eta_{1,2}$ and $C_{1,2}$ is compact, in order to show that $\gamma_{2}-\gamma_{1}$ is negative it suffices to prove $I_{1,2} (\eta)<0$,  for any $\eta \in C_{1,2}$. We consider the extended shift $\widetilde{\theta}$ on $X\times \mathbb P^{d-1}_{1,2} \times \Omega$ defined by $\widetilde{\theta} (x, \xi, \omega)=(g_{1}\cdot x, g_{1}.\xi, \theta \omega)$, the function $f(\xi,\omega)=\sigma (g_{1}, \xi)$ and the $\widetilde{\theta}$-invariant measure $\widetilde{\mathbb Q}_{\eta}=\int \delta_{(x,\xi)} \otimes \mathbb Q_{x} d\eta (x,\xi)$. Since $S^*_{n}\cdot m$ converges $\mathbb Q_{\pi}$\textrm{-a.e.} to $\delta_{z^*(\omega)}$, Lemma 3.3 implies $\displaystyle\mathop{\lim}_{n\rightarrow \infty} \sigma (S_{n} (\omega), \xi)=-\infty$, $\mathbb Q_{\pi}$\textrm{-a.e.} if the origine $v$ of $\xi$ satisfies $\langle v, z^* (\omega)\rangle \neq 0$. By hypothesis, this condition is satisfied for any $\xi$ and $\mathbb Q_{\pi}$\textrm{-a.e.}, hence we have $\displaystyle\mathop{\lim}_{n\rightarrow \infty} \displaystyle\mathop{\Sigma}_{1}^{n} f \circ \theta^k=-\infty\  \mathbb Q_{\pi}$\textrm{-a.e.} for any $\xi$. It follows that this convergence is valid $\widetilde{\mathbb Q}_{\eta}$\textrm{-a.e.}, hence Lemma 3.13 gives
$\eta (f)=I_{1,2} (\eta)<0.$
  We consider $I_{n}=\frac{1}{n} \mathbb E_{x} \left(\log \ \frac{\delta(S_{n}(\omega)\cdot v,S_{n}(\omega)\cdot v')}{\delta (v,v')}\right)$.   With $|v|=|v'|=1$, $\delta_{1} (v, v')=|v \wedge v'|$ we have
$$ \log \frac{\delta_{1}(S_{n}(\omega)\cdot v, S_{n}(\omega)\cdot v')}{\delta_{1} (v,v')}= \log \frac{|S_{n}(\omega) v \wedge S_{n} (\omega) v'|}
{|v\wedge v'|} - \log |S_{n}(\omega) v|-\log|S_{n}(\omega) v'|.$$
  By Lemma 3.15, we have also
$$\gamma_{1}+\gamma_{2}=\displaystyle\mathop{\lim}_{n\rightarrow \infty} \frac{1}{n} \displaystyle\sup_{x,v,v'} \mathbb E_{x} \left(\log \frac{| (S_{n}(\omega) v \wedge S_{n}(\omega) v'|}{|v\wedge v'|}\right).$$
  Furthermore, by Lemmas 3.14 and 3.16, we have the convergence of $\displaystyle\sup_{x,v} \frac{1}{n} \mathbb E_{x} (\log |S_{n}(\omega) v|)$ and $\displaystyle\mathop{\inf}_{x,v} \frac{1}{n} \mathbb E_{x} (\log |S_{n}(\omega) v|)$  to $I_{1}(\eta)=\gamma_{1}$.
 The uniform convergence of $\frac{1}{n} \mathbb E_{x} (\log |S_{n}(\omega) v|)$ to $\gamma_{1}$ follows. Then the equivalence of $\delta_{1}$, $\delta$ implies that $\displaystyle\sup_{x,v,v'} \frac{1}{n} I_{n}$ converges to $\gamma_{2}-\gamma_{1}$. \hfill  $\square$

  With the notations of paragraph 3 above we have the following corollaries, for products of random matrices.
\begin{theo}
\label{the:3.17}
  Assume $d>1$, the closed subsemigroup $T\subset GL(d, \mathbb R)$ satisfies condition i-p, $s\in I_{\mu}$ and $\int |g|^s \log \gamma (g) d\mu (g)$ is finite. Then the dominant Lyapunov exponent of $S_{n}(\omega)$ is simple and
\[\lim_{n\rightarrow \infty} \frac{1}{n} \displaystyle\sup_{x,v,v'} \mathbb E_{x}^{s} \left(\log \ \frac{\delta (S_{n}(\omega)\cdot v, S_{n}(\omega)\cdot v')}{\delta (v,v')}\right)=L_{\mu,2} (s)-L_{\mu,1}(s)<0, \]
where $L_{\mu,1}(s)$, $L_{\mu,2}(s)$ are the two highest Lyapunov exponents of $S_{n}(\omega)$ with respect to $\mathbb Q^s$. In particular,
\[\lim_{n\rightarrow \infty} \frac{1}{n} \displaystyle\sup_{v,v'} \mathbb E^{s} \left(\log \ \frac{\delta (S_{n}(\omega)\cdot v, S_{n}(\omega)\cdot v')}{\delta (v,v')}\right)\leq L_{\mu,2} (s) -L_{\mu,1} (s)<0.\]
  \end{theo}

\begin{proof}
  In view of Theorems 3.2, 2.6, the conditions of Proposition 3.11 are satisfied by
 $X=\mathbb P^{d-1}$, $q \otimes \mu=q^s \otimes \mu$, and $\pi=\pi^s$. On the other hand we have $\mathbb Q^s=\int \mathbb Q_{x}^s d\pi (x)$, hence the second formula.
 \end{proof}

   We will use Theorem 3.17 to establish certain functional inequalities for the operators $Q^s,\   \widetilde Q^s$ on $\mathbb P^{d-1}, \mathbb S^{d-1}$ defined in section 2 and acting on $H_{\varepsilon} (\mathbb P^{d-1})$ or $H_{\varepsilon} (\mathbb S^{d-1})$. Using \cite{34}, spectral gap properties will follow. We will say that $Q$ satisfies a ``Doeblin-Fortet  inequality'' on $H_{\varepsilon} (X)$, where $X$ is a compact metric space if we have for any $\varphi \in H_{\varepsilon} (X)$,  $[Q^{n_{0}} \varphi]_{\varepsilon}\leq \rho[\varphi]_{\varepsilon}+D|\varphi|$ for some $n_{0}\in \mathbb N$ where $\rho<1$, $D\geq 0$.

\begin{coro}
\label{corol:3.18}
  For $\varepsilon$ sufficiently small and $s\in [0,s_{\infty}[$, if $\int |g|^s \gamma^{\tau} (g) d\mu (g)<\infty$ for some $\tau>0$, then
\[\lim_{n\rightarrow \infty} ( \displaystyle\sup_{x,y} \mathbb E^s (\frac{\delta^{\varepsilon} (S_{n} \cdot  x, S_{n} \cdot y)}{\delta^{\varepsilon} (x,y)}))^{1/n}=\rho (\varepsilon) <1.\]
 If $k'(s)>0$, then
 $\displaystyle\mathop{\lim}_{n\rightarrow \infty}  (\displaystyle\sup_{x} \mathbb E^{s} (\frac{1}{|S_{n} x|^{\varepsilon}}))^{1/n} <1$.
\end{coro}

\begin{proof}
  The proof of the first formula is based on Theorem 3.17 and is given below. The proof of the second formula  follows from a similar argument (see also \cite{39}, Theorem 1, for $s=0$).

  We denote $\alpha_{n}(x,y,\omega)=\log \frac{\delta (S_{n}(\omega)\cdot x, S_{n} (\omega) \cdot  y)}{\delta (x,y)}$ and we observe that
$$e^{\varepsilon\alpha_{n}}\leq 1+\varepsilon \alpha_{n}+\varepsilon^2 \alpha^2_{n} e^{\varepsilon |\alpha_{n}|}, \ \ |\alpha_{n}|\leq 2 \log \gamma (S_{n}).$$
  Since $t^2 e^{|t|} \leq e^{3|t|}$, there exists $\varepsilon_{0}>0$ such that for $\varepsilon \leq \varepsilon_{0}$
$$\alpha^2_{n} \ e^{\varepsilon |\alpha_{n}|}\leq \frac{1}{\varepsilon_{0}^2} \ e^{3 \varepsilon_{0}|\alpha_n|} \leq \frac{1}{\varepsilon_{0}^2}\  (\gamma(S_{n}))^{6 \varepsilon_{0}}.$$
  We observe that $I_{n}=\frac{1}{\varepsilon_{0}^2}\ \mathbb E^s (\gamma^{6\varepsilon_{0}} (S_{n}))$ is finite for $s<s_{\infty}$ and $\varepsilon_{0}$ sufficiently small (see the proof of Corollary 3.20. below).
 It follows
\begin{eqnarray*}
\mathbb E^s (e^{\varepsilon \alpha_{n} (x,y,\omega)}) &\leq &1+\varepsilon \ \mathbb E^s (\alpha_{n} (x,y,\omega))+\varepsilon^2 I_{n},\\
 \displaystyle\sup_{x,y} \mathbb E^s \left(\frac{\delta^{\varepsilon}(S_{n} \cdot  x, S_{n} \cdot y)}{\delta^{\varepsilon}(x,y)}\right) &\leq &1 +\varepsilon \displaystyle\sup_{x,y} \mathbb E^s\left(\log \frac{\delta (S_{n} \cdot  x, S_{n} \cdot y)}{\delta^{\varepsilon}(x,y)}\right)  + \varepsilon^2 I_{n}.
 \end{eqnarray*}

  Also the quantity $\rho_{n}(\varepsilon)=\displaystyle\sup_{x,y} \mathbb E^s \left(\frac{\delta^{\varepsilon} (S_{n} \cdot  x, S_{n} \cdot y)}{\delta^{\varepsilon} (x,y)}\right)$ satisfies $\rho_{m+n} (\varepsilon) \leq \rho_{m} (\varepsilon) \rho_{n} (\varepsilon)$, hence we have $\rho (\varepsilon)=\displaystyle\mathop{\lim}_{n\rightarrow \infty} \rho_{n}(\varepsilon)^{1/n}=\displaystyle\mathop{\inf}_{n\in \mathbb N} (\rho_{n} (\varepsilon))^{1/n}$. It follows that, in order to show $\rho (\varepsilon)<1$, for $\varepsilon$ small it suffices to show $\rho_{n_{0}} (\varepsilon)<1$ for some $n_{0}$. To do that, we choose $n_{0}$ such that $\displaystyle\sup_{x,y} \mathbb E^s \left(\log \frac{\delta (S_{n_{0}}\cdot x, S_{n_{0}} \cdot  y)}{\delta (x,y)}\right)=c<0$ which is possible using  Theorem 3.17, and we take $\varepsilon$ sufficiently small so that $\varepsilon^2 I_{n_{0}}+\varepsilon c<0$. Then we get
$\rho_{n_{0}}(\varepsilon) \leq 1 + \varepsilon^2 I_{n_{0}}+\varepsilon c <1$.
\end{proof}

\begin{coro}
\label{coro:3.19}
   Let $H_{\varepsilon} (\mathbb P^{d-1})$ be the Banach space of $\varepsilon$-H\"older functions  on $\mathbb P^{d-1}$ with the norm $\varphi\rightarrow [\varphi]_{\varepsilon}+|\varphi|$ and assume $0\leq s< s_{\infty}$, $\int|g|^s \gamma^{\tau} (g) d\mu (g)$ for some $\tau>0$.   Then    for $\varepsilon$ sufficiently small  the operator  $Q^s$ (defined in Theorem 2.6) on $H_{\varepsilon} (\mathbb P^{d-1})$ satisfies the following Doeblin-Fortet inequality
   $[(Q^s)^{n_{0}} \varphi]_{\varepsilon}\leq \rho' (\varepsilon) [\varphi]_{\varepsilon}+B |\varphi|$,
       where $ B\geq 0\ n_{0}\in \mathbb N$,   $\rho^{n_{0}} (\varepsilon) < \rho' (\varepsilon) <1$.

  In particular the operator $P^s$ admits the following spectral decomposition in $H_{\varepsilon}(\mathbb P^{d-1})$
\[P^s=k(s) (\nu^s \otimes e^s+U^{s}),\]
  where $U^{s}$   has spectral radius less than 1,  and satisfies  $U^s(\nu^s \otimes e^s)=(\nu^s\otimes e^s) U^s=0$.
  \end{coro}

\begin{proof}
From Lemma 3.5, we know that $\mathbb Q_{x}^{s}\leq c(s) \mathbb Q^{s}$, hence, using Corollary 3.18,  for $n\geq n_{0}$ sufficiently large and with $\rho' (\varepsilon) \in ]\rho^{n_{0}} (\varepsilon),1[$, we have $\displaystyle\sup_{x,y} \mathbb E_{x}^{s} \left(\frac{\delta^{\varepsilon}(S_{n} \cdot  x, S_{n} \cdot y)}{\delta^{\varepsilon} (x,y)}\right) \leq \rho' (\varepsilon)$.  We can write
\begin{align*}
(Q^{s})^n \varphi (x)-(Q^{s})^n \varphi (y)&=\mathbb E_{x}^{s} (\varphi(S_{n} \cdot  x))-\mathbb E_{y}^{s}(\varphi (S_{n} \cdot y))\\
&=\mathbb E_{x}^{s} (\varphi (S_{n} \cdot  x))-\varphi (S_{n} \cdot y))+
(\mathbb E_{x}^{s}-\mathbb E_{y}^{s}) (\varphi (S_{n} \cdot y)).
\end{align*}

  The first term in the right hand side is bounded by $[\varphi]_{\varepsilon} \delta^{\varepsilon}(x,y) \displaystyle\sup_{x,y} \mathbb E_{x}^{s} \left(\frac{\delta^{\varepsilon}(S_{n} \cdot  x, S_{n} \cdot y)}{\delta^{\varepsilon}(x,y)}\right)$ i.e. by $[\varphi]_{\varepsilon} \delta^{\varepsilon} (x,y) \rho'(\varepsilon)$.
 Using Lemma 3.5, we  know that the second term is bounded by $ B |\varphi| \delta^{\bar s} (x,y)$. Hence with $\varepsilon<\bar s$ we get the required inequality.

  From Theorem 2.6 we know that $Q^s$ has a unique stationary measure $\pi^s$ and 1 is the unique eigenvalue with modulus one. Then the above Doeblin-Fortet inequality implies (see \cite{33}) the relation $Q^s=\pi^s \otimes 1+V^{s}$ where $V^{s}$ commutes with the projection $\pi^s \otimes 1$  has spectral radius less then one and satisfies $V^s(\pi^s\otimes 1)=0$, hence the required formula for $P^s$.
  \end{proof}

\begin{coro}
\label{cor:3.20}
  With the notation and hypothesis of Corollary 3.19, the following Doeblin-Fortet inequality is valid, if $z=s+it$,  $0\leq s<s_{\infty}$ :
\[[(Q^z)^{n_{0}}\varphi]_{\varepsilon} \leq \rho' (\varepsilon) [\varphi]_{\varepsilon}+(B+A_{n_{0}} (\varepsilon) |t|^{\varepsilon}) |\varphi|,\]
where $0\leq A_{n_{0}} (\varepsilon)<\infty$. For $t\neq 0$, the spectral radius of $Q^z$ is less than 1.
Furthermore $k(s)$ and the projection $\nu^{s}\otimes e^{s}$ are analytic on $]0,s_{\infty}[$, and 1 is a simple eigenvalue of $Q^s$.
\end{coro}

\begin{proof}
  By definition of $Q^z=Q^{s+it}$ we have $(Q^{s+it})^n \varphi (x)=\mathbb E_{x}^s (|S_{n}x|^{it} \varphi (S_{n} \cdot  x))$, hence
$|(Q^z)^n \varphi (x)-(Q^z)^n \varphi (y)|$ is bounded by the expression
\[|(\mathbb E_{x}^s-\mathbb E_{y}^s) (|S_{n}x|^{it} \varphi (S_{n} \cdot  x))+|\mathbb E_{y}^s (|S_{n}x|^{it} \varphi (S_{n} \cdot  x)-|S_{n}y|^{it} \varphi (S_{n} \cdot y))|.\]
Using Lemma 3.5 the first term is bounded by $B \delta^{\varepsilon} (x,y) \mathbb E^s (|\varphi|)$ i.e by $B|\varphi| \delta^{\varepsilon} (x,y)$. The second term is  dominated by
$\mathbb E_{y}^s (|S_{n} x|^{it}-|S_{n} y|^{it}|) |\varphi|+ \mathbb E_{y}^s (|\varphi (S_{n} \cdot  x)-\varphi (S_{n} \cdot y)|)$.

  As in the proof of Corollary 3.19, for $n\geq n_{0}$ we write
\begin{align*}
\mathbb E_{y}^s (|\varphi (S_{n} \cdot  x)-\varphi (S_{n} \cdot y)|)&\leq [\varphi]_{\varepsilon} \mathbb E_{y}^s (\delta^{\varepsilon} (S_{n} \cdot  x, S_{n} \cdot y))\\
& \leq [\varphi]_{\varepsilon} \delta^{\varepsilon} (x,y) \displaystyle\sup_{x,y} \mathbb E_{y}^s (\frac{\delta^{\varepsilon}(S_{n} \cdot  x,S_{n} \cdot y)}{\delta^{\varepsilon} (x,y)})\leq [\varphi]_{\varepsilon} \delta^{\varepsilon} (x,y) \rho' (\varepsilon).
\end{align*}

  On the other hand, using the relation
$||u|^{it}-|v|^{it}|\leq 2|t|^{\varepsilon} |\log |u|-\log |v||^{\varepsilon}\leq 2|t|^{\varepsilon} \sup (\frac{1}{|u|}, \frac{1}{|v|})^{\varepsilon} ||u|-|v||^{\varepsilon}$,
we get
\[|(|S_{n}x|^{it}-|S_{n}y|^{it}||\leq 2|t|^{\varepsilon} \displaystyle\sup_{|v|=1} \frac{1}{|S_{n}v|^{\varepsilon}} |S_{n} (x-y)|^{\varepsilon}\leq 2|t|^{\varepsilon} \displaystyle\sup_{|v|=1} \frac{|S_{n}|^{\varepsilon}}{|S_{n}v|^{\varepsilon}} \delta^{\varepsilon} (x,y).\]
Since $|S_{n}v|\geq |S_{n}^{-1}|^{-1}$ we get
$\mathbb E_{y}^s (|S_{n} x|^{it}-|S_{n} y|^{it})|) \leq 2 c(s) |t|^{\varepsilon} \delta^{\varepsilon} (x,y) \mathbb E^s (\gamma^{2\varepsilon} (S_{n}))$.
 Since $\gamma (S_{m+n})\leq \gamma (S_{m}) \gamma (S_{n}\circ \theta^m)$ and $\mathbb Q^s$ is shift-invariant,
 $\mathbb E^s (\gamma^{2\varepsilon} (S_{n}))\leq (\mathbb E^s (\gamma^{2\varepsilon}  (S_{1}))^n<\infty$.

  Then for $n$ fixed and $\varepsilon$ sufficiently small the hypothesis implies that $\mathbb E_{y}^s (|S_{n} x|^{it}-|S_{n}y|^{it}))$ is bounded by $A_{n}(\varepsilon) |t|^{\varepsilon} \delta^{\varepsilon} (x,y)$. Finally, for $n=n_{0}$,
  \[[(Q^z)^{n_{0}} \varphi]_{\varepsilon} \leq \rho' (\varepsilon) |\varphi]_{\varepsilon}+(B+A_{n_{0}} (\varepsilon) |t|^{\varepsilon}) |\varphi|.\]

  Then, using \cite{33}, one gets that the possible unimodular spectral values of $Q^z$ are eigenvalues. Using Theorem 2.6, if $t\neq 0$, one get that no such eigenvalue exists, hence the spectral radius of $Q^z$ is less than 1.

  In order to show the analyticity of $k(s)$ and $\nu^{s}\otimes e^{s}$ on $]0,s_{\infty}[$, we consider the operator $P^z$ for $z\in \mathbb C$ close to $s$. We begin by showing the holomorphy of $P^z$ for $Re \ z\in ]0,s_{\infty}[$. Let $\gamma$ be a loop contained in the strip $Re\ z\in  {\displaystyle\mathop{I}^{\circ}}_{\mu}$ and $\varphi \in H_{\varepsilon} (\mathbb P^{d-1})$. Then, since $z\rightarrow |g x|^z$ is holomorphic
\[\int_{\gamma} P^z \varphi (x) dz=\int_{G\times \gamma} \varphi ( g \cdot x) |g x|^z d\mu (g) dz=\int_{G}\varphi ( g \cdot x) d\mu(g) \int_{\gamma} |g x|^z dz=0,\]

  On the other hand,  the spectral gap property of the operator $P^s$ implies that $k(s)$ is a simple pole of the function $\zeta\rightarrow(\zeta I-P^s)^{-1}$, hence by functional calculus if $\gamma$ is a small circle of center
  $k(s)\in\mathbb C$
$$k(s) \nu^s \otimes e^s=\frac{1}{2 i \pi} \int_{\gamma} (\zeta I- P^s)^{-1} d\zeta.$$
  Since $P^z$ depends continuously of $z$, the function $(\zeta I-P^z)^{-1}$ has a pole inside the small disk defined by $\gamma$, if $z$ is close to $s$. Then by perturbation theory $P^z$ has  an isolated spectral value $k(z)$ close to $k(s)$. The corresponding projection $\nu^z \otimes e^z$ satisfies
$$k(z) \nu^z \otimes e^z=\frac{1}{2 i\pi} \int_{\gamma}(\zeta I-P^z)^{-1} d\zeta.$$
  This formula and the holomorphy of $P^z$ shows that $k(z)$ and $\nu^{z}\otimes e^{z}$ are holomorphic in a neighbourhood of $s$. The analyticity of $k(s)$ and $\nu^{s}\otimes e^{s}$ follow. The fact that 1 is a simple eigenvalue of $Q^s$ follows from Theorem 2.6.
  \end{proof}

  \begin{coro}
    \label{cor:3.21}
  Assume $\int|g|^s \gamma^{\tau} (g) d\mu(g)<\infty$ for some $\tau>0$. Then given $\varepsilon>0$ sufficiently small, for any $\varepsilon_{0}>0$ there exists $\delta_{0}=\delta_{0}(\varepsilon_{0})$, $n_{0}=n_{0} (\varepsilon_{0})$ such that if $x,y\in \mathbb S^{d-1}$ satisfy $\widetilde{\delta} (x,y) \leq \delta_{0}$, then $\mathbb E^s(\widetilde{\delta}^{\varepsilon} (S_{n_{0}}\cdot x, S_{n_{0}} \cdot  y))\leq \varepsilon_{0} \widetilde{\delta}^{\varepsilon} (x,y)$.

  One has the following Doeblin-Fortet inequality on $H_{\varepsilon}(\mathbb S^{d-1}) $ with $D\geq 0$ and $\rho_{0}=\varepsilon_{0} c(s)<1$,
$[(\widetilde Q^s)^{n_{0}} \varphi]_{\varepsilon} \leq\rho_{0} [\varphi]_{\varepsilon}+D|\varphi|$, where $c(s)$ satisfies $\mathbb Q_{x}^s \leq c(s) \mathbb Q^s.$
In particular the spectral value 1 is isolated and the corresponding finite dimensional projector depends analytically on $s \in ]0,s_{\infty}[$.

  In case I, the $\widetilde Q^s$-invariant functions are constant. In case II, the space of $\widetilde Q^s$-invariant functions is generated by $p^s_{+}$ and $p^s_{-}$.
If $t\neq 0$ the spectral radius of $\widetilde Q^z$ is less than 1.

  Furthermore, 1 is the unique unimodular eigenvalue of $\widetilde Q^s$ except in case I, where $-1$ is the unique non trivial possibility.
\end{coro}

\begin{proof}

  Assume $\varepsilon$ is as in Corollary 3.18. We will use for any $n\in \mathbb N$ and $t>0$, the relation
\[\mathbb E^s(\widetilde{\delta}^{\varepsilon} (S_{n} \cdot  x, S_{n} \cdot y))=\mathbb E^s(\widetilde{\delta}^{\varepsilon} (S_{n} \cdot  x, S_{n} \cdot y) 1_{\{\gamma (S_{n})>t\}})+\mathbb E^s(\widetilde{\delta}^{\varepsilon} (S_{n} \cdot  x, S_{n} \cdot y) 1_{\{\gamma (S_{n})\leq t\}}).\]

  In view of Corollary 3.18 we have for some $n_{0}$, any $\bar x, \bar y \in \mathbb P^{d-1}$ and given $\varepsilon_{0}>0$, we have
$\mathbb E^s(\delta^{\varepsilon} (S_{n_{0}}\cdot \bar x, S_{n_{0}}\cdot\bar y))\leq \frac{\varepsilon_{0}}{2} \delta^{\varepsilon} (\bar x, \bar y)$.

  Using Lemma 2.11 we have, for $x,y \in \mathbb S^{d-1}$,  $\widetilde{\delta} (S_{n_{0}}\cdot x, S_{n_{0}} \cdot  y)\leq 2 \gamma^2 (S_{n_{0}}) \widetilde{\delta} (x,y)$, hence
\[\mathbb E^s(\widetilde{\delta}^{\varepsilon} (S_{n_{0}}\cdot x, S_{n_{0}} \cdot  y) 1_{\{\gamma (S_{n_{0}})>t\}})\leq 2 \mathbb E^s (\gamma^{2\varepsilon} (S_{n_{0}}) 1_{\{\gamma (S_{n_{0}})>t\}}).\]
Since, as in the proof of Corollary 3.20, we have  if $\varepsilon$ is sufficiently small, $\mathbb E^s(\gamma^{2\varepsilon} (S_{n_{0}}))<\infty$, we can choose $t_{0}>0$ so that $\mathbb E^s (\widetilde{\delta}^{\varepsilon}(S_{n_{0}}\cdot x, S_{n_{0}} \cdot  y) 1_{\{\gamma (S_{n_{0}})>t_{0}\}})\leq \frac{\varepsilon_{0}}{2} \widetilde{\delta}^{\varepsilon} (x,y)$. Then, on the set $\{\gamma (S_{n_{0}})\leq t_{0}\}$ we have
$\widetilde{\delta} (S_{n_{0}}\cdot x, S_{n_{0}} \cdot  y)\leq 2 \gamma^2 (S_{n_{0}}) \widetilde{\delta} (x,y)\leq 2\  t_{0}^2 \ \widetilde{\delta} (x,y)$.

  We observe that, if $\widetilde{\delta}(u,v)\leq \sqrt2$ with $u, v \in \mathbb S^{d-1}$, then $\delta (\bar u, \bar v)=\widetilde{\delta} (u,v)$. Hence, if $2\  t_{0}^2\  \widetilde{\delta} (x,y)\leq \sqrt2$, we get  $\widetilde{\delta} (S_{n_{0}}\cdot x, S_{n_{0}} \cdot  y) =\delta (S_{n_{0}}\cdot \bar x,S_{n_{0}}\cdot \bar y)$ on the set $\{\gamma (S_{n_{0}}) \leq t_{0}\}$.
It follows, if $\widetilde{\delta}(x,y)\leq \frac{\sqrt2}{2 t_{0}^2}=\delta_{0}$,
$\mathbb E^s(\widetilde{\delta}^{\varepsilon} (S_{n_{0}}\cdot x, S_{n_{0}} \cdot  y) 1_{\{\gamma (S_{n_{0}})\leq t_{0}\}})\leq \frac{\varepsilon_{0}}{2} \widetilde{\delta}^{\varepsilon}(x,y)$.
Hence we get, if $\widetilde{\delta} (x,y)\leq \delta_{0}$,
 $\mathbb E^s(\widetilde{\delta}^{\varepsilon} (S_{n_{0}}\cdot x, S_{n_{0}} \cdot  y)) \leq  \varepsilon_{0}\  \widetilde{\delta}^{\varepsilon} (x,y)$.

  Using $\mathbb Q_{x}^s < c(s) \mathbb Q^s$ we obtain
 $\displaystyle\sup_{\widetilde{\delta} (x,y) \leq \delta_{0}} \mathbb E^s_{x} (\frac{\widetilde{\delta}^{\varepsilon} (S_{n_{0}}\cdot x,S_{n_{0}} \cdot  y)}{\widetilde{\delta}^{\varepsilon} (x,y)})\leq c(s) \varepsilon_{0}$.

  On the other hand, for $\varphi \in H_{\varepsilon} (\mathbb S^{d-1})$:
\[(\widetilde Q^s)^n \varphi (x)- (\widetilde Q^s)^n \varphi (y)=\mathbb E^s_{x} (\varphi (S_{n} \cdot  x)-\varphi (S_{n} \cdot y))+(\mathbb E_{x}^{s}-E_{y}^{s}) (\varphi (S_{n}\cdot  y)).\]

   In view of Lemma 3.5, the second term is bounded by $B\ \widetilde{\delta}^{\bar s} (x,y) |\varphi|$. Then, for $\widetilde{\delta} (x,y|\leq \delta_{0}$ we obtain, since $\varepsilon \leq \bar s$,
\[|(\widetilde Q^s)^{n_{0}} \varphi (x)-(\widetilde Q^s)^{n_{0}} \varphi (y) |\leq c(s) \varepsilon_{0} [\varphi]_{\varepsilon} \widetilde{\delta}^{\varepsilon} (x,y)+ B|\varphi| \widetilde{\delta}^{\varepsilon} (x,y).\]

  If $\widetilde{\delta} (x,y)\geq \delta_{0}$ we have  trivially $\widetilde{\mathbb E}_{x} (|\varphi (S_{n_{0}}\cdot x)-\varphi (S_{n_{0}} \cdot  y)|)\leq 2\ c(s) \frac{\widetilde{\delta}^{\varepsilon}(x,y)}{\delta^{\varepsilon}_{0}} |\varphi|$ .

  Finally, on $\mathbb S^{d-1}$,  $[(\widetilde Q^s)^{n_{0}} \varphi]_{\varepsilon} \leq c(s) \varepsilon_{0} [\varphi]_{\varepsilon}+(B+2 \frac{c(s)}{\delta^{\varepsilon}_{0}}) |\varphi|$,
hence the result with $D=B+2 \frac{c(s)}{\delta^{\varepsilon}_{0}}$. The structure of the space of $\widetilde Q^s$-invariant functions is given by Theorem \ref{thm:2.17}. It follows from \cite{33} that 1 is an isolated spectral value of $\widetilde{Q}^{s}$ and the corresponding projection has finite rank. The same argument as in the proof of Corollary 3.20 gives the analyticity of this projection. Doeblin-Fortet inequality implies that the possible unimodular spectral values of $\widetilde Q^z$ are eigenvalues. Then, as in the end of proof of Theorem 2.7, one would have for some $\varphi \in H_{\varepsilon} (\mathbb S^{d-1})$, $e^{i\theta} \in \mathbb C$, and any $g\in \supp\mu$,   $|g x|^{it} \varphi ( g \cdot x)=e^{i\theta} \varphi (x)$. This would  contradicts Proposition 2.5 if $t\neq 0$.

  The last assertion is a direct consequence of Corollary 2.19.
  \end{proof}

\textit{Proof of Theorem A.}
  The spectral decomposition $P^s=k(s) (\nu^s \otimes e^s+U^s)$ is part of Corollary 3.19. The analyticity of $k(s)$ and $\nu^{s}\otimes e^{s}$ on $]0, s_{\infty}[$ is stated in Corollary 3.20. The strict convexity of $\log k(s)$ is stated in Theorem 2.6. The fact that the spectral radius of $P^z$ is less than $k(s)$ follows from the corresponding assertion for $Q^z$ in Corollary 3.20. \hfill $\square$

\section{Renewal theorems for products of random matrices}

  It is well known that the potential theory for a random walk on $\mathbb R$ with positive drift is closely related to renewal theory (see \cite{15}). In this context,  one basic result gives the homogeneous  behaviour at infinity of the potential measure; another basic result gives the convergence $(t\rightarrow \infty)$  of the entrance measure of the random walk into $]t,\infty[$ towards a certain probability which has a density with respect to Lebesgue measure, with a simple expression in terms of the associated ladder random walk. In this section we extend these two results to linear random walks on the $G$-spaces $V\setminus \{0\}=\mathbb S^{d-1}\times \mathbb R_{+}^{*}$ and $\breve{V}=\mathbb P^{d-1} \times \mathbb R_{+}^{*}$. We denote by $\mu\in M^{1}(GL(V))$ the law of our random walk and we identify $\mathbb R$ to $\mathbb R_{+}^{*}$ via the exponential map. For the proofs we use the results of \cite{36} which give renewal theorems for a class of Markov walks on $\mathbb R$, which satisfy the tameness conditions explained below. An important observation of \cite{35} is that, if $\mu$ is supported on the positive matrices, these tameness conditions are satisfied. Here we assume that the semigroup $[\supp\mu]$ satisfies condition i-p and we use the results of sections 2 and 3 to show that the tameness conditions of \cite{36} are still valid in our generic situation. Hence we extend the results of \cite{35}, \cite{39} to the general case. This extension will play an essential role in section 5.

\subsection{The renewal theorem for a class of fibered Markov chains}

  We begin by summarizing, with a few changes and comments, the basic notations of  \cite{36}. Let $(S,\delta)$ be a complete separable metric space, $P$ (resp $\overline P$) a Markov kernel on $S\times \mathbb R$ (resp $S$) which preserves $C_{b} (S\times \mathbb R)$ (resp $C_{b} (S)$). We assume that $P$ commutes with the translations $(x,t)\rightarrow (x, t+a)$ on $S\times \mathbb R$, and  $\overline P$ is the factor kernel of $P$ on $S$.If $\pi$ is a $\overline P$-invariant probability measure, and $\ell$ is Lebesgue measure on $\mathbb R$, we note that the measure $\pi\otimes\ell$ is $P$-invariant.  We will say that $P$ is
  a fibered Markov kernel  and defines a "fibered Markov chain over $S$". More generally,if $P$ is a measurable kernel on $S\times \mathbb R$ which satisfies the above commutation we will say that $P$ is a measurable fibered kernel over $S$.In this situation,if $\pi\otimes\ell$ is a $P$-invariant measure,  we will say that P is a measurable fibered Markov kernel if P is positive and $P1=1$, $\pi\otimes\ell$-a.e. These measurable fibered kernels will play an important role in subsection 4.3 below .In this section, from now on we denote by P a Markov fibered kernel.In our applications,for Markov fibered kernels, we will have $S$ compact and $S=\mathbb P^{d-1}$ or $S\subset\mathbb S^{d-1}$, hence $S\times \mathbb R$ will be identified with a cone in $\breve{V}$ or $V\setminus\{0\}$. Here we consider paths on $S\times \mathbb R$ starting from $(x,0)\in S\times \{0\}$. Such a  path can be written as $(x_{n}, V_{n})_{n\in \mathbb N}$ with $V_{0}=0$, $x_{0}=x$, $V_{n} -V_{n-1}=U_{n}$ $(n\geq1)$.

  The corresponding space  of  paths for a fibered Markov kernel $P$ will be denoted ${^a \Omega}=S \times \displaystyle\mathop{\Pi}_{1}^{\infty} (S\times \mathbb R)$, the  Markov measure on ${^{a}\Omega}$ associated with $P$ and starting from $x\in S$ will be denoted by $^{a}\mathbb P_{x}$ and  the expectation symbol will be written ${^{a}\mathbb E}_{x}$. The space of bounded measurable functions on a measurable space $Y$ will be denoted $\mathcal B(Y)$. We observe that the Markov kernel $P$ on $S\times \mathbb R$ is completely defined by the family of measures $F(du |  x,y)$ $(x,y\in S)$ where $F(du | x,y)$ is the conditional law of $V_{1}$ given $x_{0}=x$, $x_{1}=y$.  Given a fixed $\overline{P}$-stationary probability $\pi$ on $S$,
the number  $\int u F(du |x,y) \overline P (x , dy)  d\pi (x)$\  with be called the mean of $P$,\  if the corresponding integral $\int |u| F(du | x,y) \overline P (x,dy) d\pi (x)$ is finite. In that case, we say that $P$ has a 1-moment.

  If $t\in \mathbb R_{+}$ we define the hitting time $N(t)$ of  the interval $]t, \infty[$ by :
\[N(t)=\inf \{n\geq 1 \ ;\ V_{n} > t\};  (=+\infty \  \textrm { if no such $n$ exists}). \]

  On the event $N(t)<+\infty$  we take $W(t)=V_{N(t)}-t \ ,  Z(t)=x_{N(t)}$. If $V_{1}$ has a lifetime interpretation then $W(t)$ is the residual waiting time of the interval $]t, \infty[$ (see \cite{15}). In general the law of $(Z(t), W(t))$ under $\mathbb P_{x}$ is the hitting measure of $S\times ]t, \infty[$ starting from $(x,0)$. In particular, in analogy with \cite{15}, we define the "ladder kernel" $H$ of $P$,  starting from $(x,t)\in S\times \mathbb R$, to be the law (possibly defective) of $(x_{N(0)}, t+ V_{N(0)})$ under ${^{a}\mathbb P}_{x}$. In analogy with \cite{15}, $N(0)$ (resp $V_{N(0)})$ is the first "ladder index" (resp "ladder height") of $P$. By definition $H$ is a  measurable fibered kernel over $S$; it plays an important role in the expression of various asymptotic quantities for $P$.

   One needs some technical definitions  concerning direct Riemann integrability, aperiodicity of the Markov chain defined by $P$ on $S\times \mathbb R$, and the possibility of comparing ${^{a}\mathbb P}_{x}, {^{a}\mathbb P}_{y}$ in a weak sense for different points $x, y$ in $S$. We add some comments as follows.

   Given a fibered Markov chain on
 $S\times \mathbb R$, we denote
 \[ C_{0}=\phi, C_{k}=\{x\in S\ ;\ ^{a}{\mathbb P}_{x} \{ \frac{V_{m}}{m}> \frac{1}{k}\ \  \forall m\geq k\}\geq \frac{1}{2}\}
 \textrm{ for } k\geq 1.\]

\begin{defi}{}
\label{def:4.1}
  A Borel function $\varphi\in \mathcal B(S\times \mathbb R)$ is said to be $d.R.i$ (for directly Riemann integrable) if
$$\sum_{0}^{\infty} \ \displaystyle\mathop{\sum}_{\ell=-\infty}^{+\infty} (k+1) \sup \{|\varphi(x,t)| \ ; \ x\in C_{k+1}\setminus C_{k}, \ \ell\leq t \leq \ell+1\}<+\infty\ ,$$
  and for every fixed $x\in S$ and  any $\beta >0$, the function $t\rightarrow \varphi(x,t)$ is Riemann integrable on $[-\beta,\beta]$.
\end{defi}
  In our setting below we will have $C_{k}=S$ for some $k>0$ and for some $\varepsilon>0$, any $x\in S$ and $m$ sufficiently large ${^{a}\mathbb P}_{x}(\frac{V_{m}}{m} \geq \varepsilon )\geq \frac{1}{2}$. Then the following stronger form of the above definition will be used.

\begin{defi}{}
\label{def:4.1'}
 The function   $\varphi \in \mathcal B (S\times \mathbb R)$ is said to be boundedly Riemann integrable $(b.R.i)$ if the following holds:
$\displaystyle\mathop{\Sigma}_{\ell=-\infty}^{\ell=\infty} \sup \{|\varphi (x,t)| ; x\in S, t\in [\ell, \ell+1[\}<\infty$,
and for any fixed $x\in S$, any $\beta>0$, the function $t\rightarrow \varphi (x,t)$ is Riemann integrable on $[-\beta,\beta]$.
\end{defi}
\begin{rema}{}

  Definition  \ref{def:4.1'} corresponds to $\sup\{ |\varphi (x,t)| ; x\in S\}$ directly Riemann integrable in the sense of \cite{15}. If $C_{k}=S$ for some $k\in \mathbb N$, then clearly condition $b.R.i$ implies condition $d.R.i$
\end{rema}
  The following definitions will help us to express the appropriate tameness conditions for $(P, \pi)$, where $\pi$ is a $\overline{P}$-stationary probability.
\begin{defi}{}
\label{def:4.2}
  The kernel $P$, the space $(S,\delta)$ and the measure $\pi \in M^1 (S)$ being as above we consider a point $(\zeta,\lambda,y)\in \mathbb R \times [0,1] \times S$ and we say that $(P,\pi)$ has distortion $(\zeta,\lambda)$ at $y$ if for any $\varepsilon >0$, there exists $A\in \mathcal B (S)$ with $\pi (A)>0$ and $m_{1}, m_{2} \in \mathbb N$, such that for any $x\in A$,
$^{a}{\mathbb P}_{x}\{ \delta (x_{m_{1}},y)+\delta(x_{m_{2}},y) <\varepsilon, |V_{m_{1}}- V_{m_{2}}-\zeta|< \lambda\}>0$.
\end{defi}

  For any $f\in B({^{a}\Omega})$, $\varepsilon >0$ we write
\[f^{\varepsilon}(x_{0}, x_{1}, \cdots,v_{1}\cdots)=\displaystyle\mathop{\lim\sup}_{n\rightarrow \infty} \{f(y_{0}, y_{1}, \cdots, w_{1}, \cdots) \ ; \delta (x_{i}, y_{i})+|v_{i}-w_{i}| < \varepsilon \textrm{ if } i\leq n\}.\]
\begin{defi}{}
\label{def:4.3}
  We will say that the kernel $P$ on $S\times \mathbb R$ is non-expanding  if for each fixed  $x\in S$, $\varepsilon >0$, there exists $r_{0}=r_{0} (x,\varepsilon)$ such that for all real valued $f\in \mathcal B ({^{a}\Omega})$ and for all $y$ with $\delta(x,y)<r_{0}$, one has
\[^{a}{\mathbb P}_{y} (f) \leq {^{a}{\mathbb P}}_{x} (f^{\varepsilon})+ \varepsilon |f|, \
 ^{a}{\mathbb P}_{x} (f) \leq {^{a}{\mathbb P}}_{y} (f^{\varepsilon})+ \varepsilon |f|.\]
 \end{defi}
 This condition of non expansion says that, in probability, if $x,y \in S$ are close then the paths along the fibered Markov chain starting from x,y and defined by $P$ remain close.

   One can see that, if $\varphi \in C_{b} (S\times \mathbb R)$ is uniformly continuous, then the condition of non expansion for $P$ implies that the set of functions $\{P^{n} \varphi\ ;\ n\in \mathbb N\}$ is uniformly equicontinuous.

  We denote by I.1--I.4 the following conditions, where $\pi$ denotes a given $\overline P$-stationary probability on $S$,
  \begin{itemize}
  \item{I.1.} For every open set $O$ in $S$ with $\pi(O)>0$, and   $^{a}{\mathbb P}_{x}$\textrm{-a.e}., for each $x\in S$, we have
\[^{a}{\mathbb P}_{x} \{x_{n}\in O \textrm{ for some } n\}=1.\]
  \item{I.2.} $P$ has a 1-moment and for all $x\in S$, $^{a}{\mathbb P}_{x}$\textrm{-a.e.} we have
\[\lim_{n\rightarrow \infty} \frac{V_{n}}{n}=L=\int u F (du |x,y) \overline P (x,dy) d\pi (x)>0.\]
\item{I.3.} There exists a sequence $(\zeta_{i})_{i\geq 1}\subset \mathbb R$ such that the group generated by $\zeta_{i}$ is dense in $\mathbb R$ and such that for any $i \geq1$ and $\lambda \in [0,1]$, there exists $y=y(i,\lambda)\in S$ such that $(P,\pi)$ has distortion $(\zeta_{i},\lambda)$ at $y$.
\item{I.4.} The kernel $P$ on $S\times \mathbb R$ is non-expanding.
  \end{itemize}

  Condition I.1 means that for any $x\in S$ the trajectories of $\overline{P}$ visit any non $\pi$-negligible open set with probability one. It implies that continuous $\overline{P}$-invariant functions are constant on $S$. Condition I.3 guarantees that there is no $r>0$ such that, with probability one,  the values of $V_{n}-V_{m} (n,m \in \mathbb N)$ are  of the form $k(n,m)r$ with $k(n,m)\in \mathbb Z$.
\vskip 2mm
  One can see that conditions I.1 and I.4 imply that $\overline{P}$ admits the unique invariant probability $\pi$ and that for $\varphi \in C_{b} (S)$,  the set $\{\overline{P}^{n} \varphi\ ;\ n\in \mathbb N\}$ is equicontinuous. We will denote $\ell_{+}=1_{[0,\infty[}\ell$ where $\ell$ is Lebesgue measure on $\mathbb R$. For a bounded measure $\theta$ on $S\times \mathbb R$, we will write $\overline{\theta}$ its projection on $S$. We recall that the Radon transform of a bounded measure on $\mathbb R$ is defined by $\widehat{\theta}(t)=\theta (]t,\infty[)$ if $t\in \mathbb R$.

  The ladder kernel $H$ can be written as
\[H ((y,t),\cdot)=\int (\delta_{z} \otimes H_{y}^{z} * \delta_{t}) \overline{H} (y,dz),\]
where $H_{y}^{z}$ is the conditional law of $V_{N(0)}$ given $x_{N(0)}=z$, starting from $y\in S$. Condition I.2 implies that for any $y\in S$, $H(((y,t),\cdot)$ is a probability.

  Then, the following  extension of the classical renewal theorem  is proved in \cite{36}.
\begin{theo}
\label{the:4.4}
  Assume conditions I.1-I.4 are satisfied for  the fibered Markov kernel $P$ . Then there exists a positive measure $\chi$ on $S\times \mathbb R_{+}$ absolutely continuous with respect to $\pi \otimes \ell_{+}$ such that for any $x \in S$ and $\varphi \in C_{b} (S\times ]0, \infty[)$,
\[\lim_{t\rightarrow \infty} {^{a}{\mathbb E}}_{x} \varphi (Z(t), W(t))=\frac{1}{L}\int \varphi (z,s) 1_{[s,\infty[}(t) {^{a}{\mathbb P}}_{y} \{x_{N(0)} \in dz, V_{N(0)} \in dt\} d \overline{\chi} (y) ds=\chi(\varphi),\]
 i.e. $\chi= \frac{1}{L}\int(\delta_{z} \otimes \widehat{H}_{y}^{z}  \ell_{+}) \overline{H}(y, dz) d\overline{\chi} (y)$.

  Moreover, if $\varphi \in C_{b} (S\times \mathbb R)$ is $d.R.i$, then, for any $x\in S$,
  \begin{align*}
\lim_{t\rightarrow \infty} {^{a}{\mathbb E}}_{x} (\displaystyle\mathop{\sum}_{0}^{\infty} \varphi (x_{n}, V_{n}-t))&=\frac{1}{L}\int \varphi (y,s)d\pi (y)ds,\\
\lim_{t\rightarrow -\infty} \displaystyle\mathop{\sum}_{0}^{\infty} P^{k}
\varphi (x,t)&=\frac{1}{L} (\pi\otimes \ell) (\varphi).
\end{align*}

  Furthermore $\overline{\chi}$ is an invariant measure for the measurable Markov chain on $S$ with kernel $\overline{H} (y,dz)= {^{a}\mathbb P}_{y} (x_{N(0)} \in dz)$ and
$\int \mathbb E_{y} (N(0))d\overline{\chi} (y)<\infty$,
$\int \mathbb E_{y} (V_{N(0)}) d\overline{\chi} (y)<\infty$, where $\overline{\chi}$ is absolutely continuous with respect to $\pi$.
\end{theo}
\begin{rema}
\begin{enumerate}
\item  If $S$ is compact, condition I.1 is a consequence of uniqueness of the $\overline P$-stationary measure $\pi$. This follows from the law of large numbers for Markov chains with a unique stationary measure \cite{8}: for any continuous function with $0\leq f \leq 1$, $\pi (f)>0$ we have ${^{a}{\mathbb P}}_{x}$\textrm{-a.e}.\ for all $x\in S$, $\displaystyle\mathop{\lim}_{n\rightarrow \infty} \frac{1}{n} \displaystyle\mathop{\sum}_{0}^{n-1} f(X_{k})=\pi (f) > 0$. This implies condition I.1.
\item  The construction of $\overline{\chi}$ in \cite{36} is based on Kac's recurrence theorem and implies the absolute continuity of $\overline{\chi}$ with respect to $\pi$, hence the measure $\chi$ is independent on $x$ and absolutely continuous with respect to $\pi \otimes \ell$.
\end{enumerate}
  \end{rema}

\subsection{Tameness conditions I.1--I.4 are valid for linear random walks}

  We verify conditions I.1--I.4 in four related situations. As in (\cite{35}, Proposition 1) the basic property used for the validity of I.1, I.2, I.4 is the fact that for a product $S_{n}$ of random matrices, under condition i-p the lengths of colum vectors of $S_{n}$ are comparable to the norm of $S_{n}$ (see Theorem 3.2). We observe that this property plays also an important role in the general context of \cite{4} (see Theorem 6.9). The aperiodicity property in condition I.3 is verified below using the properties of dominant eigenvalues of elements of $T$ (see Proposition 2.5). Here $\mathbb R$ is identified with $\mathbb R_{+}^*$ via the map $t\rightarrow e^t$. If $d>1$ we use condition i-p. If $d=1$ we use  non arithmeticity of $\mu$.

  The first (and simpler) situation corresponds to $S=\mathbb P^{d-1}$, $S\times \mathbb R_{+}^{*}=\breve{V}$, $P(v,\cdot )=\mu * \delta_{v}$  where $P$ is the operator on $\breve{V}$ denoted $\breve{P}$   in section 2. Also we write on $\mathbb P^{d-1} ,\overline{P}(x,\cdot )=\mu * \delta_{x}$ if $x\in \mathbb P^{d-1}$.  We will begin the verifications by this case and show how to modify the arguments in the other cases corresponding to $s=\alpha$ or $S\subset \mathbb S^{d-1} \subset V\setminus \{0\}$.
 In the  case,  of $V \setminus \{0\}$, $S$ is a compact subset of $\mathbb S^{d-1}$ and $P$ (resp $\overline{P}$) will be the restriction to $S\times \mathbb R_{+}^{*}$ (resp $S$) of the kernel already denoted $P$ (resp $\tilde P$) in section 2. Since, for any $t\in \mathbb R_{+}^*$ and $g\in G$, we have $g(t v)=t g (v)$, the kernels $P$ and $\breve{P}$ define fibered Markov chains on $S\times \mathbb R_{+}^{*}$. As shown at the end of section 2, (Theorem \ref{thm:2.17}) we need to consider  two cases for $\tilde P$, depending of the fact that $P$ preserves a proper convex cone (case II) or not (case I). In case I (resp II) we will have $S=\mathbb S^{d-1}$ (resp $S=\textrm{Co} (\Lambda_{+} ([\supp \mu])$. With these choices, there exists a unique $\tilde P$-stationary measure on $S$, as follows from Theorem \ref{thm:2.17}. In paragraph 5 below we state the detailed results for $V\setminus \{0\}$. We denote by $\alpha \in I_{\mu}$ the  positive number (if it exists) such that $k(\alpha)=1$, where $k(s)=\displaystyle\mathop{\lim}_{n\rightarrow \infty} (\int |g|^s d\mu^n (g))^{1/n}$.

  We know from section 2, that for any $s\in I_{\mu}$, there are two  Markov kernels $Q^s$ on $\mathbb P^{d-1}$ and $\tilde Q^s$ on $\mathbb S^{d-1}$, naturally associated with the operator $P^s$ considered in section 2. We are here mainly interested in the cases $s=\alpha$ and $s=0$, with $Q^0=\overline P$ and $\tilde Q^0=\tilde P$, but we observe that our considerations are valid for any $s\in I_{\mu}$.

  We denote by $^{a}{\mathbb Q}^s_{x}$ the natural Markov measure on the path space $^{a}{\Omega}$ associated with  $\mu$ and $s\geq 0$. If $s=0$  we use the notation $^{a}\mathbb P_{x}$. In our linear situations we have ${^{a}\mathbb P}_{x}=\delta_{x} \otimes \mathbb P$, ${^{a}\mathbb Q}_{x}^{s}=\delta_{x} \otimes \mathbb Q_{x}^{s}$ where $\mathbb Q_{x}^{s}$ is defined in section 3. We write $V_{k}=\log |S_{k} x|,\  x_{k}=S_{k}\cdot x$ and we denote by $\Delta_{x}$ the map from $\Omega=G^{\mathbb N}$ to $^{a}{\Omega}$ given by $(g_{1}, g_{2}, \cdots) \rightarrow (x, x_{1}, V_{1}, x_{2},V_{2}, \cdots)$. Clearly ${^{a}{\mathbb Q}}_{x}^{\alpha}$ (resp $^{a}{\mathbb P}_{x}$) is the push-forward of $\mathbb Q_{x}^{\alpha}$ (resp $\mathbb P$) by $\Delta_{x}$, hence we can translate the results of section 3 in the new setting.

  The validity of condition I.1, in all cases, is a direct consequence of the remark following   Theorem \ref{the:4.4},  since by Theorem 2.6 and Theorem \ref{thm:2.17} the kernels $\bar P$, $Q^{\alpha}$, $\tilde P$, $\tilde Q^{\alpha}$ have unique stationary measures on  $S$.

  In order to verify I.2 we begin by $S=\mathbb P^{d-1}$, $S\times \mathbb R_{+}^{*}=\breve{V}$,
$ \overline P(x,.)=\mu * \delta_{x}$, $\breve{P} (v,\ .)=\mu * \delta_{v}$. Then $F(A | x,y)=\mu \{g\in G, \log |g x| \in A,  g \cdot x=y\}$ for any Borel set $A\subset \mathbb R$, and $\pi=\nu$ with $\mu * \nu=\nu$. We observe that   $|\log |g x|| \leq \log \gamma (g)$ if $|x|=1$ and $\gamma(g)= \sup (|g|, |g^{-1}|)$. The finiteness of $\int |u | F(d u| x,y) \bar P (x, dy) d\pi (x)$ follows since $\log \gamma (g)$ is $\mu$-integrable. Also  $\int u F (d u | x,y) \bar P (x, dy) d\pi (x)=\int \log |g x| d\mu (g) d\nu (x)=L_{\mu}$. Then the relation $^{a}{\mathbb P}_{x}=\Delta_{x} (\mathbb P)$ and Theorem 3.10 imply, for every $x\in \mathbb P^{d-1}$  (case $s=0$), and $\mathbb P-$a.e.:
\[L_{\mu}=\displaystyle\mathop{\lim}_{n\rightarrow \infty} \frac{1}{n} \log |S_{n} x|=\int \log |g x| d\mu (g) d\nu (x).\]

  Except for $L_{\mu}>0$, this is condition I.2 in the first case. If $S \subset \mathbb S^{d-1}$ is as above, the result is the same, since the involved quantities depend only on $|g x|$ with $x\in \mathbb P^{d-1}$, and $\tilde P$ has a unique stationary measure on $S$.

  In the cases of $ Q^{\alpha}$ and $\widetilde{Q}^{\alpha}$ it suffices also to consider the case $S=\mathbb P^{d-1}$.
  The 1-moment condition in I.2 follows from $\int |g|^{\alpha} |\log \gamma (g)| d\mu(g)<+\infty$. The convergence part follows from Theorem 3.10  with
  \[L_{\mu}(\alpha)=\int q^{\alpha} (x,g) \log |g x| d\pi^{\alpha}(x) d\mu (g)=\frac{k'(\alpha)}{k(\alpha)}>0.\]

  We show I.3 as follows.
If $d>1$, since the semigroup $T=\displaystyle\mathop{\cup}_{n\geq 0} (\supp \mu)^n$ satisfies (i-p), we know using Proposition 2.5 that the set
$\Delta=\left\{\log |\lambda_{h}| \ ; \ h\in T^{\textrm{prox}} \right\}$
 is dense in $\mathbb R$.
The same is true of $2\Delta=\{\log \lambda_{h^2}$; $h\in T^{\textrm{prox}} \}$.

  If $d=1$, the same properties follow from the non arithmeticity of $\mu$.

  We take for $\zeta_{i}$ $(i\in \mathbb N)$ a dense countable subset of $2\Delta$. Let $\zeta_{i}=\log \lambda_{g} \in 2\Delta$, with $\lambda_{g}>0$, $g=h^2$, $h=u_{1}\cdots u_{n}$, $u_{i}\in \supp \mu$ $(1\leq i \leq n)$ and $y=y(\zeta_{i}, \lambda)=\bar v_{g}\in \mathbb P^{d-1}=S$. We observe that, if $\varepsilon$ is sufficiently small and $B_{\varepsilon}=\{x\in \mathbb P^{d-1}; \delta (x, \bar v_{g})\leq \varepsilon\}$, then $g\cdot B_{\varepsilon} \subset B_{\varepsilon'}$, with $\varepsilon'<\varepsilon$ and $g$ as above,
 Also, $\lambda>0$ being fixed, and $\varepsilon$ sufficiently small, we have $|\log \lambda_{g}-\log |g x||<\lambda$  if $x\in  B_{\varepsilon}$. These relations remain valid for $g'$ instead of $g$ if $g'$ is sufficiently close to $g$.
Then we have for $x\in  B_{\varepsilon}$, and $S_{2n}=g \in (\supp \mu)^{2n}$ as above,
$${^{a}\mathbb P}_{x} \{\delta(S_{2n}. x, \bar v_{g}) <\varepsilon,\ \  |\log|S_{2n} x|-\log \lambda_{g}|<\lambda\}>0.$$
  With $\zeta_{i}=\log \lambda_{g}, y=\bar v_{g}, A=B_{\varepsilon},\  \tau=0,\  m_{1}=0,\ m_{2}=2n$, this implies condition I.3 for the probability $\mbox{}^{a}{\mathbb P}_{x}=\Delta_{x} (\mathbb P)$.

  The definition of $\mathbb Q_{x}^{\alpha}$ shows its equivalence to $\mathbb P$ on the $\sigma$-algebra of the sets depending of the first $n$ coordinates. Then the relation $^{a}{\mathbb Q}_{x}^{\alpha} =\Delta_{x} (\mathbb Q_{x}^{\alpha})$ implies with $g$ as above
\[^{a}{\mathbb Q}_{x}^{\alpha}\{ \delta (S_{2n}. x, \bar v_{g})<\varepsilon,\ \  |\log| S_{2n}x|-\log \lambda_{g} |<\lambda\}>0.\]
  Hence condition I.3 is valid for $\mbox{}^{a}\mathbb Q_{x}^{\alpha}$ also.
If we consider $\mathbb S^{d-1}$ instead of $\mathbb P^{d-1}$, i.e.\ $S=\mathbb S^{d-1}$ or $S=\textsf{Co} (\Lambda_{+} ([\supp\mu]))$, and the metric $\tilde{\delta}$ on $S$, the above geometrical argument remains valid with $g=h^{2}$, $y=\tilde v_{g} \in \Lambda_{+}\ ([\supp \mu])$ in the second case,  $\lambda_{g}>0$ and $\varepsilon$ sufficient small.  This shows I.3 in this setting.

   Condition I.4 follows from the proof of Proposition 1 of \cite{35}. The proof of the corresponding part of this proposition  is a consequence of the condition,
$${^{a}{\mathbb P}}_{x} \{\exists C>0\ \hbox{\rm with}\ |S_{n} x|\geq C|S_{n}| \ \hbox{\rm for\ all}\  n\}=1,$$
  for all $x\in S$, which implies that $|S_{n} x|$ and $|S_{n} y|$ are comparable if $x$ and $y$ are close.

  For the proof of the above condition, we observe that if $x\in S$ and $s\in I_{\mu}$, in particular if $s=0$ or $\alpha$, this condition has been proved in the stronger form
\[\lim_{n\rightarrow \infty} \frac{|S_{n}x|}{|S_{n}|}=|\langle z^*(\omega),x\rangle|>0 \ \  \mathbb Q^s_{x}\textrm{-a.e}., \]
in Theorem 3.2,  hence condition I.4 is valid in all the  cases under consideration.

\subsection{Direct Riemann integrability}

  In case of the spaces $S=\mathbb P^{d-1}$ or $S\subset \mathbb S^{d-1}$ considered above,  under condition i-p for $[\supp \mu]$, the $d.R.i$ condition takes the  simple form given by Lemma \ref{lem:4.4}  below, in multiplicative notation.

  We  assume now that the hypothesis  of Theorem 3.10 is satisfied, we use the corresponding notations, and $C_{k}$ is as in Definition \ref{def:4.1} .

\begin{lemm}
\label{lem:4.5}

  For $k$ large we have $C_{k}=S$.
\end{lemm}

\begin{proof}
We consider  the case  $S=\mathbb P^{d-1}$. Using Theorem 3.10 for $s=0$, we get for any fixed $x$
\[\lim_{n\rightarrow +\infty} \frac{\log |S_{n} (\omega) x|}{n} = L_{\mu},\ \  {\mathbb P-\textrm{a.e.} .}\]
We observe that for any $x$,  $y\in \mathbb S^{d-1}$, $\|S_{n} y|-|S_{n} x\| \leq |S_{n}| \widetilde{\delta} (x,y) \leq 2|S_{n}|$.  It follows
$$\left |\frac{|S_{n} y|}{|S_{n}x|}-1\right|\leq 2 \frac{|S_{n}|}{|S_{n}x|},\ \ \ \  \left | \log |S_{n} y|-\log |S_{n} x|\right|\leq 2 \frac{|S_{n}|}{|S_{n}x|} .$$
  Using Theorem 3.2, we get that the sequence $\frac{|S_{n}|}{|S_{n} x|}$ converges ${\mathbb P-\textrm{a.e.}}$ to
  $\frac{1}{|\langle z^{*}(\omega), x\rangle |} <\infty$, hence the sequence $\frac{1} {n}$ $\frac{|S_{n}|}{|S_{n} x|}$ converges ${\mathbb P-\textrm{a.e.}}$ to zero.
It follows that $\frac{1}{n} \log |S_{n} x|-\frac{2}{n}\ \frac{|S_{n}|}{|S_{n} x|}$ converges $\mathbb P$-\textrm{a.e}.\ to $L_{\mu}$.

   Hence there exists $m_{0}>0$ such that
 $\mathbb P\{\frac{1}{n} \log |S_{n} x|-\frac{2}{n} \ \frac{|S_{n}|}{|S_{n} x|} >\frac{1}{2} L_{\mu}$ for all $n\geq m_{0}\}\geq 1/2$.
  In view of the inequality $\frac{1}{n} \log |S_{n} y| \geq\frac{1}{n} \log |S_{n} x|-\frac{2}{n}\ \frac{|S_{n}|}{|S_{n}x|}$, we have for any $y\in \mathbb S^{d-1}$,
$\mathbb P\{\frac{1}{n} \log |S_{n} y|>L_{\mu}/2$ for all $n\geq m_{0}\}> \frac{1}{2}$.
This implies $C_{k}=\mathbb P^{d-1}, C_{k+1}\setminus C_{k}=\phi$ if $\frac{1}{k} \leq \inf \left(\frac{1}{m_{0}},  \ \frac{L_{\mu}}{2}\right)$.

  If $s=\alpha$, the argument is the same with $\mathbb P$ replaced by $\mathbb{Q}^{\alpha}$ and the relation $\mathbb{Q}_{x}^{\alpha} \leq c(\alpha) \mathbb{Q}^{\alpha}$ is used as follows.
$$\mathbb{Q}^{\alpha}\{\frac{1}{n} \log |S_{n} y| > L_{\mu} (\alpha)/2\  \hbox{\rm for \ all}\  n\geq m_{0}\}>1-\frac{1}{2 c(\alpha)}.$$
  Since $\mathbb{Q}_{y}^{\alpha} \leq c(\alpha) \mathbb{Q}^{\alpha}$, this gives for any $y\in \mathbb P^{d-1}$,
\[\mathbb{Q}_{y}^{\alpha}\{\frac{1}{n} \log |S_{n} y| > L_{\mu} (\alpha)/2\  \hbox{\rm for \ all}\  n\geq m_{0}\}>\frac{1}{2}.\]
  Then also $C_{k}=\mathbb P^{d-1}$ for $\frac{1}{k} \leq \inf \left(\frac{1}{m_{0}}, \ \frac{L_{\mu} (\alpha)}{2}\right)$. Hence we conclude as above.
  \end{proof}

\subsection{The renewal theorems for linear random walks}

  We consider $\breve{V}=\mathbb P^{d-1} \times \mathbb R^*_{+}$ $V\setminus \{0\}=\mathbb S^{d-1} \times \mathbb R_{+}^*$, and we study the asymptotics of the potential kernels of the corresponding random walks defined by $\mu$. We denote
$$\displaystyle\mathop{\breve{V}_{1}}=\{ v\in \breve{V}\ ;\ |v|>1\},\ \ V_{1}=\{v\in V\ ;\ |v|>1\},$$
  and we consider also the entrance measures $ \breve{H}(v, \ .)$ or $H(v,\ .)$ of $S_{n} v$ in $\breve{V}_{1}$ or $V_{1}$, starting from $v$. Since conditions $I$ are valid, their behaviour for $v$ small are given by Theorem \ref{the:4.4},  and we will state them below. We denote by $\breve{\Lambda} ([\supp \mu])$ the inverse image of $\Lambda ([\supp \mu])$ in $\breve{V}$. Also we denote
\begin{align*}
\breve{\Lambda}_{1} ([\supp \mu])&=\{v\in\breve{V}\ ;\ \bar v \in \Lambda ([\supp \mu]),\ |v| \geq 1\},\\
\Lambda_{1} ([\supp \mu])&=\{ v\in V\  ; \  \bar v \in \Lambda ([\supp \mu]),\ |v| \geq 1\}.
\end{align*}
As shown below these closed sets support the limits $(v\rightarrow 0)$ of $\breve{H} (v,.)$ and $H (v,\cdot )$.
The function $e^{s}$ (resp measure $\nu^{s}$) defined in Theorem 2.2 plays an essential role for $s=0$, $\alpha$ in the following theorems.

  The results will take two forms according as $L_{\mu}>0$ or $L_{\mu}<0$.

\begin{theo}
\label{the:4.6}
  Assume $\mu \in M^1 (G)$ is such that the semigroup $[\supp \mu]$ satisfies condition i-p if $d>1$ or $\mu$ is non arithmetic if $d=1$, $\log\gamma(g)$ is $\mu$-integrable and
  \[L_{\mu}= \displaystyle\mathop{\lim}_{n\rightarrow \infty} \frac{1}{n} \int \log |g| d\mu^n(g)>0.\]
Then  if $v\in \breve{V}$, $\displaystyle\mathop{\sum}^{\infty}_{0} \mu^k * \delta_{v}$ is a Radon measure on $\breve{V}$ such that on $C_{c} (\breve{V})$  we have the vague convergence,
\[\lim_{v\rightarrow 0} \displaystyle\mathop{\sum}^{\infty}_{0} \mu^k *\delta_{v}=\frac{1}{L_{\mu}} \nu \otimes \ell,\]
  where $\nu \in M^1 (\Lambda ([\supp \mu]))$ is the unique $\bar P$-invariant measure on $\mathbb P^{d-1}$.
This convergence is valid on any bounded continuous function $f$ which satisfy on $\breve{V}$,
\[\sum^{+\infty}_{-\infty} \sup \{|f(v)| \ ; \ 2^{\ell} \leq |v| \leq 2^{\ell+1} \} <\infty.\]

  Furthermore  the ladder kernel $\breve{H}$ satisfies the  following   weak convergence,
\[\lim_{v\rightarrow 0} \breve{H} (v, \cdot )=\breve{\chi} \in M^1 (\breve{\Lambda}_{1} ([\supp \mu])),\]
where $\breve{\chi}$ is defined by this convergence and is absolutely continuous with respect to $\nu \otimes \ell$
 \end{theo}

\begin{proof}
In view of the verifications of conditions I in subsections 2, 3 and of Lemma 4.6,  this is a direct consequence of Theorem 4.5  applied in the case $S=\mathbb P^{d-1}$, $S\times \mathbb R_{+}^{*}=\breve{V} \setminus \{0\}$.
 \end{proof}

  If $L_{\mu}<0$ the following gives the asymptotic behaviour $(v\rightarrow 0)$ of the potential measure $\displaystyle\mathop{\sum}^{\infty}_{0} \mu^k *\delta_{v}$; this asymptotics allow us to obtain a Cram\'er  estimate for the random variable
$$M(v)=\sup \{|S_{n}v| \ ;  n\geq 1\}.$$

\begin{theo}
 \label{the:4.7}
Assume that $\mu \in M^1 (G)$ is such that $[\supp \mu]$ satisfies i-p, if $d>1$ or $\mu$ is non arithmetic if $d=1$.
 Assume $L_{\mu}<0, \alpha>0$ exists with  $k(\alpha)=1$, $\int  |g|^{\alpha} \log \gamma (g) d\mu (g)<\infty$ and write
\[L_{\mu} (\alpha)=\displaystyle\mathop{\lim}_{n\rightarrow \infty} \frac{1}{n} \int |g|^{\alpha} \log |g| d\mu^n (g)=\frac{k'(\alpha)}{k(\alpha)}.\]
Then $L_{\mu} (\alpha)>0$ and for any $u\in \mathbb P^{d-1}$,  we have the vague convergence in $\breve{V}$,
\[\lim_{t\rightarrow 0}  t^{- \alpha} \displaystyle\mathop{\Sigma}^{\infty}_{0} \mu^k * \delta_{tu} =\frac{e^{\alpha}(u)}{L_{\mu} (\alpha)}  \nu^{\alpha}\otimes \ell^{\alpha},\]
 where $\nu^{\alpha} \in M^1 (\mathbb P^{d-1})$ (resp $e^{\alpha} \in C(\mathbb P^{d-1})$,  $\nu^{\alpha}(e^{\alpha})=1)$ is the unique solution of the equation $P^{\alpha} \nu^{\alpha}=\nu^{\alpha}$ (resp $P^{\alpha} e^{\alpha}=e^{\alpha})$ and $\nu^{\alpha}$ has support $\Lambda([\supp \mu])$.

  Furthermore, on $C_{b} (\breve{V}_{1})$ and for any $u\in \mathbb P^{d-1} \subset \breve{V}$,  the ladder kernel $\breve{H}(tu,\cdot)$ satisfies the vague convergence,
\[\lim_{t\rightarrow 0} t^{-\alpha} \breve{H} (t u, \cdot )= e^{\alpha} (u) \breve{\chi}^{\alpha},\]
where $ \breve{\chi}^{\alpha}$, defined by this convergence,  is a positive measure supported on $\breve{\Lambda}_{1} ([\supp \mu])$ and is absolutely continuous with respect to $\nu^{\alpha}\otimes \ell^{\alpha}$.

  In particular, for  $A=\breve{\chi}^{\alpha} (\breve{V_{1}})>0$ and any $u\in \mathbb P^{d-1}$,  $\displaystyle\mathop{\lim}_{t\rightarrow \infty} t^{\alpha} \mathbb P \{M (u)>t\}=A e^{\alpha} (u)$.

  The first convergence is valid on  any continuous function $f$ which satisfies
\[\sum^{+\infty}_{-\infty} 2^{-\ell \alpha}  \sup \{ |f (v)| \  ;\  2^{\ell} \leq |v| \leq 2^{\ell+1}\}<\infty.\]
 \end{theo}
 \begin{proof}
  We observe that the function $e^{\alpha}\otimes h^{\alpha}$ on $\breve{V}$ satisfies $\breve{P} (e^{\alpha} \otimes h^{\alpha})=e^{\alpha} \otimes h^{\alpha}$, hence we can consider the associated Markov operator $\breve{Q}_{\alpha}$  on $\breve{V}$ defined by
\[\breve{Q}_{\alpha} (f)=\frac{1}{e^{\alpha} \otimes h^{\alpha}} \breve{P} (f e^{\alpha} \otimes h^{\alpha}).\]
  Then the potential kernel of ${\breve{Q}}_{\alpha}$ is given by
\[\sum^{\infty}_{0} (\breve{Q}_{\alpha})^k  (f)=\frac{1}{e^{\alpha}\otimes h^{\alpha}}\displaystyle\mathop{\sum}^{\infty}_{0} \breve{P}^k (f e^{\alpha} \otimes h^{\alpha}).\]
Clearly $\breve{Q}_{\alpha}$ commutes with dilations, hence defines a fibered Markov kernel on $\breve{V}$.

  Also the mean of $\breve{Q}_{\alpha}$ is $L_{\mu} (\alpha)>0$. Then, taking $f=\frac{\varphi}{e^{\alpha}\otimes h^{\alpha}}$, since conditions I are valid, the result follows from Theorem \ref{the:4.4}. Cram\'er's estimation for $\mathbb P\{\sup |S_{n} u|>t\ ;\ n\in \mathbb N\}$ follows with $A=\breve{\chi}^{\alpha} (\breve{V}_{1})>0$.
  \end{proof}

\textit{Proofs of Theorems B, B$^{\alpha}$ and Corollary.}
These results are simple consequences of  Theorems \ref{the:4.6}, \ref{the:4.7}.
\hfill $\square$

\subsection{Renewal theorems on $V\setminus\{0\}$}

  Here we extend Theorems \ref{the:4.6} et \ref{the:4.7} to the natural setting of $V\setminus \{0\}$. We note that for $d=1$, we have $V\setminus\{0\}=\mathbb R^{*}$, then the results reduce to the classical renewal theorems (\cite{15}).
\vskip 3mm
\begin{theo}
\label{the:4.6'}
 Assume $\mu$ is as in Theorem \ref{the:4.6}. Then there are 2 cases  as in Theorem \ref{thm:2.17},
\begin{itemize}
\item{Case I:} No proper convex cone in $V$ is $\supp \mu$-invariant.
Then, in vague topology,
\[\lim_{v\rightarrow 0} \displaystyle\mathop{\sum}^{\infty}_{0} \mu^k * \delta_{v}=\frac{1}{L_{\mu}} \widetilde{\nu}\otimes \ell,\]
where $\widetilde{\nu}$ is the unique $\mu$-stationary measure on $\mathbb S^{d-1}$.
  \item{Case II:}  Some proper convex cone in $V$ is $\supp \mu$-invariant.
  Then, for any $u\in \mathbb S^{d-1}$, in vague topology,
\[\lim_{t\rightarrow 0_{+}} \displaystyle\mathop{\sum}^{\infty}_{0} \mu^k * \delta_{t u}=\frac{1}{L_{\mu}} (p_{+}(u) \nu_{+}\otimes \ell + p_{-}(u) \nu_{-}\otimes \ell),\]
where $\nu_{+}$ is the unique $\mu$-stationary measure on $\Lambda_{+}([\supp \mu])$, $\nu_{-}$ is symmetric of $\nu_{+}, \ \ p_{+}(u)$ is the entrance probability of  $S_{n}\cdot u$ in $\textrm{Co} (\Lambda_{+}([\supp \mu]))$, $p_{-}(u)=1 -p_{+} (u)$.
\end{itemize}
  In the two cases these convergences are also valid on any bounded continuous functions $f$ on $V \setminus \{0\}$  such that $\displaystyle\mathop{\sum}_{-\infty}^{+\infty} \sup \{|f(v)| \ ;\ 2^{\ell} \leq |v| <2^{\ell+1}\}<\infty$.

   In addition, for any $u\in \mathbb S^{d-1}$, in weak topology,
\[\lim_{t\rightarrow 0_{+}} H (t u, \cdot )=\chi_{u} \in M^1 (\Lambda_{1} ([\supp\mu]),\]
  where $\chi_{u}$ is defined by this convergence and is absolutely continuous with respect to $\widetilde{\nu}\otimes \ell$.
\begin{itemize}
\item
  In case I, $\chi_{u}=\chi$ is independent on $u$.
\item
  In case II, with $\Lambda_{1,+} ([\supp \mu])=\{v\in V \ ; \ \tilde v \in \Lambda_{+} ([\supp \mu]), \ |v| \geq 1\}$ we have
\[\chi_{u}=p_{+}(u) \chi_{+} +p_{-}(u) \chi_{-},\]
  where $\chi_{1,+}\in M^1 (\Lambda_{+} (\supp \mu))$, and $\chi_{-}$ is symmetric of $\chi_+$.
  \end{itemize}
  \end{theo}

  \begin{proof}
  In case I, the proof is the same as for Theorem \ref{the:4.6} with $S=\mathbb S^{d-1}$ instead of $\mathbb P^{d-1}$.
  In case II, we take $S=\textrm{Co} (\Lambda_{+} ([\supp \mu])$ and we observe that $S\times \mathbb R_{+}^*$ is a $\supp \mu$-invariant convex cone with non zero interior to which,  as in the proof of Theorem \ref{the:4.6}, we can apply Theorem \ref{the:4.4} .

  If $u\in S$ (resp $u\in-S$) we have
\[\lim_{t\rightarrow 0_{+}} \displaystyle\mathop{\sum}^{\infty}_{0} \mu^k * \delta_{t u} =\frac{1}{L_{\mu}} (\nu_{+}\otimes \ell), \ \ (\textrm{resp.}
\ \ \lim_{t\rightarrow 0_{+}} \displaystyle\mathop{\sum}^{\infty}_{0} \mu^k * \delta_{t u}=\frac{1}{L_{\mu}}(\nu_{-}\otimes \ell)).\]

  If $u\in V \setminus \{0\}$, we denote by $p_{+} (u, dv)$ (resp $p_{-}(u, dv))$ the entrance measure of $S_{n}\cdot u$ in the cone $\Phi=S\times \mathbb R_{+}^*$ (resp $-\Phi)$. Clearly the mass of $p_{+} (u, dv)$ is $p_{+}(u)$, and $p_{+}(u, dv)$ is supported on $\Phi$.

  We denote also by  $U (v,\cdot)=\displaystyle\mathop{\sum}^{\infty}_{0} \mu^k * \delta_{v}$ the potential kernel of the  linear random walk $S_{n} (\omega) v$ starting from $v$ in $V\setminus\{0\}$. Then, for any $\varphi \in C_{c} (\Phi \cup-\Phi)$,
\[U(t u, \varphi)=\int U (v, \varphi) (p_{+} (tu, dv)+p_{-} (t u, dv)).\]
  Clearly the kernel $p_{+}(x, dv)$ commutes with the scaling $x\rightarrow t x$ $(t>0)$. Then it follows from above that, on $C_{c} (\Phi\cup-\Phi)$ :
\[\lim_{t\rightarrow 0_{+}}\  U (t u,.)=\frac{1}{L_{\mu}} (p_{+}(u) \nu_{+} \otimes \ell+p_{-}(u) \nu_{-}\otimes \ell).\]

  If $\varphi \in C_{c} (V \setminus \{0\})$ vanishes on $\Phi\cup-\Phi$, Theorem \ref{the:4.6} implies  $\displaystyle\mathop{\lim}_{t\rightarrow 0_{+}} U (t u, \varphi)=0$.

   Finally we have  $\displaystyle\mathop{\lim}_{t\rightarrow 0_{+}} \displaystyle\mathop{\sum}^{\infty}_{0} \mu^k * \delta_{t u}=\frac{1}{L_{\mu}} (p_{+}(u) \nu_{+} \otimes \ell + p_{-}(u) \nu_{-}\otimes \ell)$. The existence of $\chi_{u}$ follows from the first formula in Theorem \ref{the:4.4} . In particular the right hand side of this formula is independent on $u\in S$. Hence, in case I, $\chi_{u}$ is  independent on $u$. In case II, we use $S=\textrm{Co}(\Lambda_{+}([\supp\mu]))$ and we argue as above in order to obtain the formula $\chi_{u}=p_{+}(u) \chi_{+}+p_{-} (u) \chi_{-}$ where $\chi_{+}=\chi_{u}$ for $u\in \textrm{Co}(\Lambda_{+}([\supp\mu])$ and $\chi_{-}=\chi_{u}$ for $u\in \textrm{Co} (\Lambda_{-} ([\supp\mu]))$.
\end{proof}
  We have also the following analogue of Theorem \ref{the:4.7}. The proof is a combination of the arguments in the proofs of Theorems \ref{the:4.7} and \ref{the:4.6}.

\begin{theo}
\label{the:4.7'}
  Assume $\mu$ and  $\alpha$ are as in Theorem \ref{the:4.7}. Then for any $u\in \mathbb S^{d-1}$ we have the vague convergence : $\displaystyle\mathop{\lim}_{t\rightarrow 0_{+}}t^{-\alpha} \displaystyle\mathop{\sum}^{\infty}_{0} \mu^k * \delta_{t u}=\frac{e^{\alpha}(u)}
{L_{\mu} (\alpha)} \widetilde{\nu}_{u}^{\alpha}\otimes \ell^{\alpha}$, where $\widetilde{\nu}^{\alpha}_{u}\in M^1 (\widetilde{\Lambda}(T))$ is $\widetilde P^{\alpha}$-invariant. There are 2 cases like in Theorem \ref{the:4.6'}.
\begin{itemize}
\item
  Case I: $\widetilde{\nu}_{u}^{\alpha}=\widetilde{\nu}^{\alpha}$ has support $\widetilde{\Lambda}(T)$,
\item
  Case II: $\widetilde{\nu}_{u}^{\alpha}= p_{+}^{\alpha} (u) \nu_{+}^{\alpha} +  p_{-}^{\alpha}(u) \nu_{-}^{\alpha}$,
  where $p_{+}^{\alpha} (u)$ (resp $p_{-}^{\alpha} (u)$) denotes the entrance probability under $\mathbb Q_{u}^{\alpha}$  of $S_{n}\cdot u$ in the convex envelope of $\Lambda_{+}(T)$ (resp $\Lambda_{-} (T)$).
  \end{itemize}
  The above convergences are valid on any continuous function $f$ which satisfies
\[\sum^{\ell=+\infty}_{\ell=-\infty} 2^{-\ell \alpha} \sup \{|f(v)| \ ;\ 2^{\ell} \leq |v| \leq 2^{\ell+1}\}<\infty.\]

  Furthermore, on $C_{b} (V_{1})$ for any $u\in \mathbb S^{d-1}$, we have the vague convergence,
$$\lim_{t\rightarrow 0_{+}} t^{-\alpha} H (t u,.)= e^{\alpha} (u) \ (p_{+}^{\alpha} (u) \chi^{\alpha}_{-} + p^{\alpha}_{-}(u) \chi_{+}^{\alpha}),$$
  where $\chi_{+}^{\alpha}, \chi_{-}^{\alpha}$ are defined by this convergence and are absolutely continuous with respect to $\widetilde{\nu}^{\alpha}\otimes \ell$.

  In case I, $\chi_{+}^{\alpha}=\chi_{-}^{\alpha}=\chi^{\alpha}$ is a positive measure supported on $\Lambda_{1}([\supp \mu])$. In case II,  $\chi^{\alpha}_{+}$ is a positive measure supported on $\Lambda_{+} ([\supp \mu])$ and $\chi_{-}^{\alpha}$ is symmetric of $\chi^{\alpha}_{+}$.
  \end{theo}

\subsection{On the asymptotics of $k(s)$ $(s\rightarrow \infty)$}

  For the existence of $\alpha>0$ such that $k(\alpha)=1$, we have the following sufficient condition where  we denote by $r(g)$ the spectral radius of $g\in G$.

\begin{prop}
\label{pro:4.8}
  Let $\mu\in M^1 (G)$ and assume that $k(s)={\displaystyle\mathop{\lim}_{n\rightarrow \infty}} (\int |g|^s d\mu^n(g))^{1/n}$ is finite for any $s>0$. For any $p\in \mathbb N$ and $g\in (\supp\mu)^p$ we have

\centerline{$\displaystyle\mathop{\lim}_{s\rightarrow \infty} \frac{\log k(s)}{s} \geq \ p^{-1} \ \log \ r(g).$}

  In particular if some $g\in [\supp\mu]$ satisfies $r(g)>1$, then $k(s)>1$ for $s$ sufficiently large.
\end{prop}
  The proof is based on the following elementary lemma which we state without proof.
\begin{lemm}
\label{lemm:4.9}

  Let $g\in G$. Then   for any $\varepsilon>0$  there exists $c(\varepsilon)>0$ and a neighbourhood $V(\varepsilon)$ of $g$ such for any sequence $g_{k}\in V(\varepsilon)$ one has $|g_{n}\cdots g_{1}|\geq c(\varepsilon)\ r^n (g) (1-\varepsilon)^n$.
\end{lemm}
\textit{Proof of  Proposition \ref{pro:4.8}.}
  The convexity of $\log\ k(s)$ implies that ${\displaystyle\mathop{\lim}_{s\rightarrow \infty} } \frac{\log\ k(s)}{s}$ exists. Let $g\in \supp\mu$, hence given $\varepsilon>0$ these exists a neighbourhood $V(\varepsilon)$ of $g$ as in Lemma \ref{lemm:4.9} such that $\mu (V(\varepsilon))=C(\varepsilon)>0$. From the lemma we have,
  \begin{align*}
  k(s)&=\displaystyle\mathop{\lim}_{n\rightarrow \infty}  (\int |g_{n}\cdots g_{1}|^s d\mathbb P(\omega))^{1/n}\\
& \geq \displaystyle\mathop{\lim}_{n\rightarrow \infty}  (c^s(\varepsilon) r^{ns} (g) (1-\varepsilon)^{ns} C^n (\varepsilon))^{1/n}=r^s (g) (1-\varepsilon)^s C(\varepsilon).
\end{align*}
Hence
$\frac{\log\ k(s)}{s}\geq \log(1-\varepsilon) +\log\ r(g)+\frac{\log\ C(\varepsilon)}{s}$, therefore ${\displaystyle\mathop{\lim}_{s\rightarrow \infty}}  \frac{\log \ k(s)}{s}\geq  \log \ r (g)$.

   We observe that if $\mu$ is replaced by  $\mu^p$, then $k(s)$ is replaced by $k^p(s)$. Hence for $g\in (\supp\mu)^p$ we have from above the required inequality.
If $g\in[\supp\mu]$ then we can assume $g\in (\supp\mu)^p$ for some $p\in \mathbb N$;  since $r(g)>1$, we have $\log\ r(g)>0$, hence ${\displaystyle\mathop{\lim}_{s\rightarrow \infty}}  \frac{\log\ k(s)}{s} >0$. \hfill $\square$

\section{The tails of an affine stochastic recursion}

\subsection{Notation and main result}

  Let $H$ be the affine group of the $d$-dimensional Euclidean space, i.e. the set of maps $f$ of $V$ into itself of the form $f (x)=g x+b$ where $g\in GL (V)=G, b\in V$. Let $\lambda$ be a probability measure on $H$, $\mu$ its projection on $G$. We denote by $\Sigma$ (resp $T$) the closed subsemigroup of $H$ (resp $G$) generated by $\supp\lambda$ (resp $\supp\mu$).  We consider the  affine random walk on $V=\mathbb R^d$ defined by $\lambda$, i.e. the Markov chain on $V$ described by the stochastic recursion,
\[X^x_{n+1}=A_{n+1}X^x_{n}+B_{n+1},   X^x_{0}=x\in V,\]
   where  $(A_{n}, B_{n})$ are  $H$-valued i.i.d random variables with law $\lambda$. We denote $\widehat{\Omega}=H^{\mathbb N}$ and we endow $\widehat{\Omega}$ with the shift $\widehat{\theta}$ and the product measure $\widehat{\mathbb P}=\lambda^{\otimes \mathbb N}$; by abuse of notation the expectation symbol with respect to $\widehat{\mathbb P}$ will be denoted by $\mathbb E$. We have
\[X^x_{n}  =A_{n}\cdots A_{1} x+\displaystyle\mathop{\sum}^n_{1} A_{n}\cdots A_{k+1} B_{k}.\]

  We are interested in the case where $R_{n}=\displaystyle\mathop{\sum}^n_{1} A_{1}\cdots A_{k-1} B_{k}$ converges $\hat{\mathbb P}$-\textrm{a.e.} to a random variable $R$ and $X^x_{n}$ converges in law to $R$. We observe that $X^x_{n}  -A_{n}\cdots A_{1} x$ and $R_{n}$ have the same law. In that case we have
\[R=\displaystyle\mathop{\sum}^{\infty}_{0} A_{1}\cdots A_{k} B_{k+1},\]
  hence the random variable $R$ satisfies the  equation
\[ R=A R\circ \widehat{\theta}+B, \quad (S)\]
  and the law $\rho$ of $R$ satisfies the convolution equation $\rho=\lambda * \rho=\int h \rho d\lambda (h)$. Also, if $R$ is unbounded, we will be interested in the tail of $R$  in direction $u$, i.e. the asymptotics ${t\rightarrow \infty}$ of $\hat{\mathbb{P}}\{\langle R,u \rangle>t\}$ (resp $\hat{\mathbb P} \{|\langle  R, u \rangle|> t\} $), where $u\in \mathbb S^{d-1}$ (resp $u\in \mathbb P^{d-1}$). We are mainly interested in the "shape at infinity" of $\rho$ i.e.\ the asymptotics $(t\rightarrow 0_{+})$ of the measure $t \cdot  \rho$ where $t \cdot  \rho$ is the push-forward of $\rho$ by the dilation $v\rightarrow t v$ in $V (t>0)$. It turns out that this "shape at infinity" depends essentially of the semigroups $T$ and $\Sigma$ defined above. A basic role will be played by the top Lyapunov exponent $L_{\mu}$ of the product of random matrices $S_{n}=A_{n}\cdots A_{1}$, and $\mu$ will be assumed to satisfy
\[\int \log \gamma (g) d\mu (g) <\infty \textrm{ where } \gamma (g)=\sup (|g|, |g^{-1}|).\]
  The main hypothesis will be on $\mu$, which is  always assumed to satisfy $L_{\mu}<0$ and condition  i-p of section 2  if $d>1$, or $\mu$ non arithmetic if $d=1$. We recall that the function $k(s)$ is defined on the interval $I_{\mu}\subset [0,\infty[$ by $k(s)=\displaystyle\mathop{\lim}_{n\rightarrow \infty} (\int |g|^s d\mu^n (g))^{1/n}$ and $\log k(s)$ is strictly convex (see Theorem 2.6). It is natural to assume that $\supp \lambda$ has no fixed point in $V$ since otherwise the affine recursion reduces to a  linear recursion.   We denote by $\Delta_{a} (\Sigma)$ the set of fixed attractive points of the elements of $\Sigma$, i.e. fixed points $h^+ \in V$ of elements $h=(g,b) \in \Sigma$ such that $\displaystyle\mathop{\lim}_{n\rightarrow \infty} |g^n|^{1/n}<1$. For $v\in V\setminus \{0\}$ we  denote
\[H_{v}^+=\{x \in V; \langle v, x \rangle >1\},\]
  and for a bounded measure $\xi$ on $V$ we  consider its Radon transform $\widehat{\xi}$, i.e. the function on $V\setminus \{0\}$  defined by
\[\widehat{\xi}(v)=\xi(H_{v}^+) \textrm{ with } H_{v}^{+}=\{x \in V\ ;\langle x,v \rangle >1\}.\]

  We  also write $u=t v$ with $u\in \mathbb S^{d-1}$, $t>0$ and $\widehat{\xi}(u,t)=\widehat{\xi}(\frac{u}{t})$. In particular,  the directional tails of $\xi$ are described by the function $\widehat{\xi}(v)$ $(v\rightarrow 0)$. We start with the basic

\begin{prop}
\label{pro:5.1}
  Assume $L_{\mu} <0$ and $\mathbb E(\log |B|)<\infty$. Then $R_{n}$ converges $\widehat{\mathbb P}-\textrm{a.e}$ to $R=\displaystyle\mathop{\sum}^{\infty}_{1} A_{1} \cdots A_{k-1} B_{k}$, and for any $x\in V$, $X^x_{n}$ converges in law to $R$. For all  $\beta \in I_{\mu}$ with $k(\beta)<1$ and $\mathbb E(|B|^{\beta})<\infty$, we have $\mathbb E (|R|^{\beta})<\infty$.

  The law $\rho$ of $R$ is the unique $\lambda$-stationary measure on $V$. The closure $\overline{\Delta_{a}(\Sigma)}=\Lambda_{a}(\Sigma )$ is the unique $\Sigma$-minimal subset in $V$ and is equal to $\supp\rho$. If the semigroup $T$ contains an element $g$ with $\displaystyle\mathop{\lim}_{n\rightarrow \infty} |g^n|^{1/n}>1$ and $T$ has no fixed point   then $\supp\rho$ is unbounded.

  If $T$ satisfies condition i-p and $\supp \lambda$ has no fixed point in $V$, then $\rho$ $(W)=0$ for any affine subspace $W$.
\end{prop}

\begin{proof}

  Under the conditions $L_{\mu}<0$ and $\mathbb E(\log |B|)<\infty$ the  $\widehat{\mathbb P}$-\textrm{a.e.} convergence of $R_{n}$ to $R$ is well known as well as the moment condition $\mathbb E (|R|^{\beta})<\infty$ if $k(\beta)<1$ and $\beta\in I_{\mu}$ (see for example \cite{7}). We complete the argument by observing that, since $L_{\mu}<0$, we have $\displaystyle\mathop{\lim}_{n\rightarrow \infty} |A_{n} \cdots A_{1} x|=0$ hence,  since $X^x_{n}-A_{n}\cdots A_{1} x$ has the same law as $R_{n}$, the convergence in law of $X^x_{n}$ to $R$  for any $x$ follows. In particular, if $x \in V$ is distributed according to $\xi \in M^1 (V)$, the law of $X^x_{n}$ is $\lambda^n*\xi=\int \lambda^n *\delta_{x} d\xi(x)$, hence has limit $\rho$ at $n=\infty$. If $\xi$ is $\lambda$-stationary, we have $\lambda^n*\xi=\xi$, hence $\xi=\rho$.

  Since $L_{\mu}<0$, there exists $h=(g,b) \in \Sigma$, such that $|g|<1$, hence $\displaystyle\mathop{\lim}_{n\rightarrow \infty} |g^n|^{1/n}<1$.
 If $h=(g,b)\in \Sigma$ satisfies $\displaystyle\mathop{\lim}_{n\rightarrow \infty} |g^n|^{1/n}<1$, then $I-g$ is invertible, hence the unique fixed point $h^+$ of $h$ satisfies $(I-g) h^+=b$, and for any $x\in V$, $h^n x-h^+=g^n(x-h^+)$,  hence $\displaystyle\mathop{\lim}_{n\rightarrow \infty} h^n x=h^+$. Taking $x$ in $\supp\rho$ we get $h^+ \in \supp\rho$, since $\supp\rho$ is $h$-invariant. Furthermore, for any $x\in V$ and $h' \in \Sigma$ we have $\displaystyle\mathop{\lim}_{n\rightarrow \infty} h' h^n x=h' (h^+)$ and $h' h^n \in \Sigma$ satisfies $\displaystyle\mathop{\lim}_{n\rightarrow \infty} |g' g^n|=0$, hence the unique fixed point $x_{n}$ of $h' h^n
$ satisfies $\displaystyle\mathop{\lim}_{n\rightarrow \infty} x_{n}=h' (h^+)$. Then $\overline{\Delta_{a}(\Sigma)}=\Lambda_{a}(\Sigma)$ is a $\Sigma$-invariant non trivial, closed subset of $\supp\rho$.

  On the other hand, for $x\in \Delta_{a} (\Sigma)$ we have
$\displaystyle\mathop{\lim}_{n\rightarrow \infty} \lambda^n *\delta_{x}=\rho,\ \ (\lambda^n *\delta_{x}) (\Lambda_{a}(\Sigma))=1$ for all $n$
 hence $\rho$ $(\Lambda_{a} (\Sigma)=1$, i.e $\Lambda_{a} (\Sigma)=\supp\rho$. The $\Sigma$-minimality of $\Lambda_{a} (\Sigma)$ follows from the fact that, for any $x\in V$ and $h=(g,b)$ with $|g|<1$, one has $\displaystyle\mathop{\lim}_{t\rightarrow \infty} h^n x=h^+ \in \Lambda_{a} (\Sigma)$ hence $\overline{\Sigma x} \supset \Lambda_{a}(\Sigma)$. This implies also the uniqueness of the $\Sigma$-minimal set.

  Observe that, if $\supp\rho$ is bounded, then the convex envelope $\textrm{Co}(\supp\rho)$ is a compact subset of $V$ . Also any $h\in \Sigma$ preserves $\supp\rho$ and $\textsf{Co}(\supp\rho)$. Then  Markov-Kakutani theorem implies that the affine map $h$ has a fixed point $h^0$ in $\textrm{Co} (\supp\rho)$. If $h=(g,b)\in \Sigma$ satisfies $\displaystyle\mathop{\lim}_{n\rightarrow \infty} |g^n|^{1/n}>1$ we have :
\[(I-g) h^0=b\ \ \hbox{\rm and}\ \ h^n x-h^0=g^n (x-h^0),\]
  hence if $x\neq h^{0}$, we have $\displaystyle\mathop{\lim}_{n\rightarrow \infty} |h^{n}x|=\infty$.  Then $\supp\rho$ is unbounded since if $x\in \supp\rho \neq \delta_{h_{0}}$ the point $h^n x$ belongs to $\supp \rho$.

   Let $\mathcal W=\{W_{i}\  ; \ i \in I\}$ be the set of affine subspaces of minimal dimension  with $\rho (W_{i})>0$. Since $\textrm{dim}(W_{i} \cap W_{j})<\textrm{dim} W_{i}$ if $i\neq j$, we have $\rho(W_{i} \cap W_{j})=0$, hence $\displaystyle\mathop{\sum}_{i\in I} \rho (W_{i})\leq 1$. It follows that, for any $\varepsilon>0$, the set $\{W_{j} ; j\in I$ and $\rho(W_{j})\geq \varepsilon\}$ has cardinality at most $\frac{1}{\varepsilon}$, hence $\rho(W_{i})$ reachs its maximum on a finite set $\{W_{j} \ ; \ j\in J \subset I\}$ of affine subspaces. Then the stationarity equation $\lambda * \rho=\rho$ gives on such a subspace $W_{j}$,
\[\rho (W_{j})=\int \rho (h^{-1} W_{j}) d \lambda (h).\]
Since $\rho(h^{-1} W_{j})\leq \rho (W_{j})$ we get, for any $h\in \supp \lambda$, $j\in J$,
\[\rho (h^{-1} W_{j})=\rho (W_{j}), \ i.e \ h^{-1} W_{j}=W_{i},\] for some $i\in J$.

   In other words the set $\{W_{j} \ ;\ j\in J\}$ is $\supp \lambda$-invariant. If $\textrm{dim} W_{j}>0$, one gets that the set of directions $\overline{W}_{j} (j\in J)$ is a $\supp \mu$-invariant finite set of subspaces of $V$, which contradicts condition i-p for the semigroup $T$. Hence each $W_{j} (j\in J)$ is reduced to a point $w_{j}$. Then the barycenter of the finite set $\{w_{j} \ ; \  j\in J\}$ is invariant under $\supp\lambda$, which contradicts the hypothesis. Hence $\rho(W)=0$ for any affine subspace $W$.
\end{proof}

  In order to state the main result of this section we consider the compactification $V \cup \mathbb S^{d-1}_{\infty}$,  and the natural projection of $\mathbb S_{\infty}^{d-1}$ on the unit sphere $\mathbb S^{d-1}$. We denote by $\widetilde{\Lambda}^{\infty} (T)$ (resp $\Lambda_{+}^{\infty} (T),\ \Lambda_{-}^{\infty} (T))$ the inverse image of $\widetilde{\Lambda} (T)$ (resp $\Lambda_{+} (T),\ \Lambda_{-}(T))$ in $\mathbb S_{\infty}^{d-1}$ (see section 2, paragraph 3). The closure $\overline{\Lambda_{a}(\Sigma)}$ of $\Lambda_{a}(\Sigma)$ in the compact space $V \cup \mathbb S^{d-1}_{\infty}$ is $T$-invariant hence $\overline{\Lambda_{a}(\Sigma)} \cap \mathbb S^{d-1}_{\infty}=\Lambda^{\infty}_{a} (\Sigma)$, which is non  void  if $\supp\rho=\Lambda_{a}(\Sigma)$ is unbounded and is a closed $T$-invariant subset of $\mathbb S_{\infty}^{d-1}$.

  Then Proposition \ref{pro:2.15} applied to $\Lambda_{a}^{\infty}(\Sigma) \subset \mathbb S^{d-1}_{\infty}$ gives the following trichotomy, since condition i-p is satisfied by $T$
\begin{itemize}
\item
case I: $T$ has no invariant proper convex cone and $\Lambda_{a}^{\infty} (\Sigma) \supset \widetilde{\Lambda}^{\infty} (T)$,
\item
case II': $T$ has  an invariant proper convex cone and $\Lambda_{a}^{\infty}(\Sigma)\supset \widetilde{\Lambda}^{\infty} (T)$,
\item
case II'': $T$ has an invariant proper convex cone and $\Lambda_{a}^{\infty}(\Sigma)$ contains only one of the sets $\Lambda^{\infty}_{+} (T), \Lambda^{\infty}_{-} (T)$, say $\Lambda^{\infty}_{+} (T)$, hence $\Lambda_{a}^{\infty} (\Sigma) \cap \Lambda_{-}^{\infty} (T)=\emptyset$.
\end{itemize}

  We assume $\alpha \in  ]0, s_{\infty}[$ exists with $k(\alpha)=1$ (see Proposition \ref{pro:4.8} for a sufficient condition).
As in Theorem \ref{the:4.7'}, we consider the $\widetilde{P}^{\alpha}$-invariant measures $\widetilde{\nu}^{\alpha}, \nu^{\alpha}_{+}, \nu^{\alpha}_{-}$.

   The following implies   Theorem C of section I  and describes the asymptotics of the probability measure  $t\cdot \rho$ when $t\rightarrow 0_{+}$:  $\rho$ has a Pareto distribution of index $\alpha$ (\cite{47}, p.74).

\begin{theo}
\label{the:5.2}
  With the above notation  assume $L_{\mu}<0$, $\Sigma$ has no fixed point in $V$, $T$ satisfies condition i-p, there exists $\alpha \in ]0, s_{\infty}[$ such that $k(\alpha)=1$ and $\mathbb E(|B|^{\alpha+\varepsilon})<\infty$, $\mathbb E(|A|^{\alpha} \gamma^{\varepsilon} (A))<\infty$ for some $\varepsilon>0$. If $d=1$ assume also that $\mu$ is non arithmetic.

  Then $\supp \rho$ is unbounded and  we have the following vague convergence on $V\setminus \{0\}$
\[\displaystyle\mathop{\lim}_{t\rightarrow 0_{+}} t^{-\alpha} (t \cdot  \rho)= \Lambda= C (\sigma^{\alpha}\otimes \ell^{\alpha})\]
  where $C>0$, $\sigma^{\alpha} \in M^1 (\widetilde{\Lambda} (T))$ are defined by the above formula and the measure $\Lambda=C(\sigma^{\alpha} \otimes \ell^{\alpha})$ satisfies $\mu * \Lambda=\Lambda$.

  In case I,  we have $\sigma^{\alpha}=\widetilde{\nu}^{\alpha}$.

  In case II', there exist $C_{+}, C_{-}>0$ with $C\sigma^{\alpha}=C_{+} \nu_{+}^{\alpha}+C_{-} \nu_{-}^{\alpha}$,

  In case II'', $\sigma^{\alpha}=\nu_{+}^{\alpha}$ .

  In case I the Radon measure $\widetilde{\nu}^{\alpha}\otimes \ell^{\alpha}$ on $V\setminus \{0\}$ is a minimal $\mu$-harmonic measure, and $\Lambda$ is symmetric.

  In cases II, $\nu^{\alpha}_{+}\otimes \ell^{\alpha}$ and $\nu^{\alpha}_{-}\otimes \ell^{\alpha}$ are minimal $\mu$-harmonic measures on $V\setminus\{0\}$.
\end{theo}
  Theorem 5.2 is proved in several steps, using the function $\widehat\rho$ on $V$,
\[\widehat \rho(v)=\widehat{\mathbb{P}}\{\langle R,v \rangle >1\}.\]

  A first step, based on the renewal theorems of section 4, shows the existence of the directional tails (see Corollary 5.8) i.e. existence of the limit $\displaystyle\mathop{\lim}_{t\rightarrow \infty} t^{\alpha} \widehat{\rho}(t^{-1}u)$. A second step is to study the positivity of these tails (see Proposition 5.9). It is based on Kac's recurrence theorem (see \cite{50}) for an associated random walk on a dual $H$-homogeneous space of $V$ (see Lemma 5.13); this recurrence property allow us to express $\widehat{\rho}$ as a potential of a non negative function  on $V\times \mathbb R$ (see Lemma 5.12), to which a weak renewal theorem can be applied (see Proposition 5.14). The action of $H$ on affine hyperplanes of $V$ leads us, via Radon transforms, to consider the natural linear representation of $H$ in the vector space $V\times \mathbb R$. The corresponding linear $\lambda$-random walk is studied in paragraph 3 below. Finally the homogeneity at infinity of $\rho$ follows from a Choquet-Deny type result (see Theorem 5.17).

\subsection{Asymptotics of directional tails}

   We apply Theorem \ref{the:4.7'} to $\mu^{*}$-potentials of suitable functions; we pass, using the map $\eta\rightarrow \widehat{\eta}$, from the convolution equation $\lambda * \rho=\rho$ to a Poisson type equation on $V\setminus \{0\}$ which involves $\mu^{*}$ and $\widehat{\rho}$ and we note that $t^{\alpha}(t^{-1} \cdot  \rho) (H^+_{u})=t^{\alpha} \widehat{\rho}(u,t)$. The corresponding convergences will play an essential role in the proof of Theorem 5.2. We denote by $\rho_{1}$ the law of $R-B$ and we consider the signed measure $\rho_{0}=\rho-\rho_{1},$ hence $\rho_{0}(V)=0$. Also we show that $\rho_{0}$ is "small at infinity", and we define $C,C_{+}, C_{-}, \sigma^{\alpha}$. With the hypothesis of Theorem 5.2, we denote by $^{*}\widetilde{\nu}^{\alpha}_{u}$ the positive kernel on $\mathbb S^{d-1}$ given by Theorem \ref{the:4.7'} and associated with $\mu^{*}$.

\begin{prop}
\label{prop: 5.3}

   One has the equations  on $V\setminus \{0\}$,
\[\rho=\displaystyle\mathop{\Sigma}_{0}^{\infty} \mu^k * (\rho-\rho_{1}), \ \ \widehat{\rho}(v)=\displaystyle\mathop{\Sigma}_{0}^{\infty} ((\mu^{*})^k *\delta_{v}) (\widehat{\rho}-\widehat{\rho}_{1}).\]

   For $u \in \mathbb S^{d-1}$, the function $t\rightarrow t^{\alpha-1} \widehat{\rho}_{0} (u,t)$ is Riemann-integrable in generalized sense  on $]0,\infty[$ and,  one has with $r_{\alpha} (u)=\int_{0}^{\infty} t^{\alpha-1} \widehat{\rho}_{0} (u,t) dt$, $p(\alpha)=\int |\langle x,  y \rangle|^{\alpha} d\nu^{\alpha}(x) d^{\alpha}\nu^{\alpha}(y)$
\[\lim_{t\rightarrow \infty} t^{\alpha} \widehat{\rho} (u,t) = \frac{^{*}e^{\alpha} (u)}{L_{\mu}(\alpha)} {^{*}\widetilde{\nu}^{\alpha}_{u}}(r_{\alpha})= C (\sigma^{\alpha} \otimes \ell^{\alpha}) (H_{u}^{+}),\]
  where $C=2 \frac{^{*}\widetilde{\nu}^{\alpha}(r_{\alpha}) \alpha}{L_{\mu}(\alpha)p (\alpha)} \geq 0$ and $\sigma^{\alpha} \in M^1 (\widetilde{\Lambda} (T))$ are defined by the above convergence and $\sigma^{\alpha}\otimes \ell^{\alpha}$ satisfies $\mu * (\sigma^{\alpha} \otimes \ell^{\alpha})=\sigma^{\alpha} \otimes \ell^{\alpha}$.

  Furthermore $\supp\rho$ is unbounded and,
\begin{itemize}
\item
  In case I: $\sigma^{\alpha}=\widetilde{\nu}^{\alpha}$,
\item
  In case II: ${^{*}\widetilde{\nu}^{\alpha}_{u}} (r_{\alpha}) \sigma^{\alpha}=\frac{1}{2} (^{*}\nu^{\alpha}_{+}(r_{\alpha}) \nu_{+}^{\alpha}+{^{*}\nu}^{\alpha}_{-} (r_{\alpha}) \nu_{-}^{\alpha})$ where ${^{*}\nu_{+}^{\alpha}} (r_{\alpha})\geq 0$, ${^{*}\nu}^{\alpha}_{-} (r_{\alpha})\geq 0$.
  \end{itemize}
\end{prop}

  The proof will follow from a series of lemmas.

  We start with the following  simple Tauberian lemma.

\begin{lemm}
 \label{lemm:5.4}
  For a non negative and non increasing function $f$ on $\mathbb R^{*}_{+}$ and $s\geq 0$, we denote, $f^s (t)=\frac{1}{t} \int^{t}_{0} x^s f(x) dx$. Then the condition $\displaystyle\mathop{\lim}_{t\rightarrow \infty} f^s (t)=c$ implies $\displaystyle\mathop{\lim}_{t\rightarrow \infty} t^s f(t)=c$.
\end{lemm}
\begin{proof}

  Let $b$ be a positive real number with $b>1$ and let us observe that, since $f$ is non increasing
\[\frac{1}{t} \int^{bt}_{t} x^s f(x) dx \leq f(t) \frac{1}{t} \int^{bt}_{t} x^s dx=\frac{t^s}{s+1} (b^{s+1}-1) f(t).\]
  It follows
\[\frac{b^{s+1}-1}{s+1} t^s f(t) \geq b \ f^s (bt)-f^s (t).\]

  Then the hypothesis gives :
\[\liminf_{t\rightarrow \infty} \frac{b^{s+1}-1}{s+1} t^s f(t) \geq (b-1) c.\]

  Using the relation $\displaystyle\mathop{\lim}_{b\rightarrow 1} \frac{b^{s+1}-1}{(s+1) (b-1)}=1$ we get
  $\displaystyle\mathop{\liminf}_{t\rightarrow \infty} (t^s f(t))\geq c$. An analogous argument gives
  $\displaystyle\mathop{\limsup}_{t\rightarrow \infty} (t^s f(t)) \leq c$. It follows
  $\displaystyle\mathop{\lim}_{t\rightarrow \infty} t^s f(t)=c$.

\end{proof}
  We will use below the multiplicative structure of the group $\mathbb R^{*}_{+}=]0, \infty[$, and we recall that Haar measure $\ell$ on the multiplicative group $\mathbb R^{*}_{+}$ is given by $\frac{dt}{t}$.

\begin{lemm}
\label{lemm:5.5}
  Assume that the $V$-valued random variable $R$ satisfies equation $(S)$,  and $\mathbb E (|B|^{\alpha+\delta})<\infty$, with $\delta>0$.
 For $u\in \mathbb S^{d-1}$ and $t,x>0$ we write
$r^{\alpha}(u,t)=\frac{1}{t} \int ^t_{0} x^{\alpha}  \widehat{\rho}_{0} (u,x) dx$.

  Then $|r^{\alpha} (u,t)|\leq 2\frac{t^{\alpha}}{\alpha+1}$. For  $\delta'$ small there exists $C(\delta')>0$ such that if $t\geq 1$, $|r^{\alpha} (u,t)| \leq C(\delta') t^{-\delta'}$. In particular the function $r^{\alpha}(u,t)$ is $b.R.i$ on $\mathbb S^{d-1} \times \mathbb R^{*}_{+}$.
\end{lemm}
\begin{proof}

  Since $\rho_{0}$ is the difference of the laws of $R$ and $R-B$ we have by definition of $\widehat{\rho}_{0}, \ |\widehat{\rho}_{0} (u,t)|\leq 2$, hence $|r^{\alpha} (u,t)| \leq 2 \frac{t^{\alpha}}{\alpha+1}$. Also
$\widehat{\rho}_{0} (u,x)=r_{1}(u,x)-r_{2}(u,x)$  where,
\begin{align*}
 r_{1} (u,x)&=\widehat{\mathbb P}\{x-\langle B,u \rangle \ <\langle R-B,u \rangle\leq x\},\\
 r_{2} (u,x)&=\widehat{\mathbb P}\{x <\langle R-B,u\rangle \leq x-\langle B,u \rangle\}.
\end{align*}
 Furthermore $r^\alpha=r_1^\alpha-r_2^\alpha$  with $r_{1}^{\alpha} (u,t)=\frac{1}{t} \int_{0}^{t} x^{\alpha}  r_{1}(u, x) dx$ and $ r_{2}^{\alpha} (u,t)=\frac{1}{t} \int_{0}^{t} x^{\alpha} r_{2} (u, x) dx$.

  In order to estimate $r^{\alpha}_{1}$, we choose $\varepsilon \in ]0,1[$ with $\varepsilon > \frac{\alpha}{\alpha+\delta}$ and write, for $t\geq 2$,
$r^{\alpha}_{1}(u,t)\leq \frac{1}{t}\int^t_{2} x^{\alpha} \widehat{\mathbb P} \{<B,u>\geq x^{\varepsilon}\} dx + \frac{1}{t}\int_{2}^{t} x^{\alpha} \widehat{\mathbb P} \{ x-x^{\varepsilon} < \langle  R-B, u\rangle  \leq x\} dx +\frac{2^{\alpha+1}}{(\alpha+1)t}$.

  Then Markov's inequality gives  $\widehat{\mathbb P}\{\langle B,u \rangle\geq x^{\varepsilon}\} \leq x^{-(\alpha+\delta)\varepsilon} \mathbb E(|B|^{\alpha+\delta})$.

  Hence the first term $I^{\varepsilon}_{1} (t)$ in the above inequality satisfies
\[I^{\varepsilon}_{1} (t) \leq \mathbb E (|B|^{\alpha+\delta}) \frac{1}{t} \int_{2}^t x^{\alpha-\varepsilon(\alpha+\delta)} dx \leq \mathbb E(|B|^{\alpha+\delta}) t^{\alpha-\varepsilon(\alpha+\delta)}.\]

  For $t-t^{\varepsilon} \geq 2$, the second term $I^{\varepsilon}_{2}(t)$ satisfies
\[I^{\varepsilon}_{2} (t) \leq \frac{1}{t} \int_{2}^t x^{\alpha} \widehat{\mathbb P} \{<R-B,u>>x-x^{\varepsilon}\} dx - \frac{1}{t} \int_{2}^{t-t^{\varepsilon}} x^{\alpha} \widehat{\mathbb P }\{\langle R-B,u \rangle > x\} dx.\]

  In the second integral above we use the change of variables $x\rightarrow x-x^{\varepsilon}$ and we get
\[I^{\varepsilon}_{2} (t) \leq \frac{1}{t} \int^t_{2} [x^{\alpha}-(x-x^{\varepsilon})^{\alpha} (1-\varepsilon \ x^{\varepsilon-1})] \widehat{\mathbb P}\{<R-B,u>> x-x^{\varepsilon}\} dx + \frac{k_{0}(\varepsilon)}{t} \]
with $0<k_{0}(\varepsilon)<\infty$.

  We observe that there exists $k_{1}(\varepsilon)<\infty$ such that for any $x\geq 2$,
\[x^{\alpha}-(x-x^{\varepsilon})^{\alpha} (1-\varepsilon \ x^{\varepsilon-1}) \leq k_{1} (\varepsilon) x^{\alpha+\varepsilon-1}.\]
  For any $\beta \in ]0, \alpha[$, Proposition 5.1 implies that $\mathbb E(|R|^{\beta})<\infty$. Also $R$ satisfies equation $(S)$ and $A$, $R\circ \widehat{\theta}$ are independent. Hence Markov's inequality gives
  \[\widehat{\mathbb P}\{\langle R-B,u \rangle x\} \leq x^{-\beta} \mathbb E(|A|^{\beta}) \mathbb E (|R|^{\beta})\leq k_{2}(\beta) x^{-\beta} \textrm{ with } k_{2} (\beta)<\infty.\]
  It follows that for any $t$ with $t-t^{\varepsilon}>2$
\[I^{\varepsilon}_{2} (t) \leq \frac{k_{0}(\varepsilon)}{t}+k_{1} (\varepsilon) k_{2} (\beta) \frac{1}{t} \int^t_{2} \frac{x^{\alpha+\varepsilon-1}}{(x-x^{\varepsilon})^{\beta}}   dx \leq k_{3} (\varepsilon,\beta) t^{\alpha-\beta+\varepsilon-1}.\]
  It remains to choose $\delta, \varepsilon, \beta$ in order to obtain $\alpha+\varepsilon-1-\beta<0$.
  We take $\delta$ so small that $\delta(\alpha+\delta)^{-1}<\alpha$ and $\varepsilon=(\alpha+\delta/p) (\alpha+\delta)^{-1}$ with $p\in \mathbb N$, $p\geq 2$, hence $\varepsilon \in ]\alpha (\alpha+\delta)^{-1}, 1[$. Also $\alpha+\varepsilon-1=\alpha-\delta(1-1/p) (\alpha+\delta)^{-1}>0$.

  We take $\gamma \in ]0,1[$ and $\beta=\alpha-\gamma \delta(1-1/p) (\alpha+\delta)^{-1}$ so that
\[\alpha+ \varepsilon-1-\beta=(\gamma-1) \delta (1-1/p) (\alpha+\delta)^{-1}.\]

  With $p=2$, $\gamma=1/2$ we get $\alpha+\varepsilon-1-\beta=-\delta/4 (\alpha+\delta)^{-1}$. We write $\delta'$=inf $(\delta/4 (\alpha+\delta)^{-1}, \delta/2)$.

  Hence, there exists $k_{3}<\infty$ and $\delta'>0$ such that for $t\geq 1$, $r^{\alpha}_{1}(u,t)\leq k_{3} t^{-\delta'}$.

  The same argument is valid for $r^{\alpha}_{2}$, hence for some $\delta'>0$ and $t\geq 1$, we have
$r^{\alpha} (u,t)\leq C(\delta') t^{-\delta'},$
with $C(\delta')<\infty$.
 Furthermore,  for $t\in ]0,1]$ we have $|r^{\alpha} (u,t)| \leq 2\frac{t^{\alpha}}{\alpha+1}$, hence the function $r^{\alpha} (u,t)$ is $b.R.i.$ on $\mathbb S^{d-1} \times \mathbb R_{+}^{*}$.
\end{proof}

\begin{lemm}
\label{lemm:5.6}
  We denote by $r$ the finite measure on $\mathbb R^{*}_{+}$ defined by $r(d x)=1_{]0,1[} (x) x^{\alpha} dx$ and we write $\rho_{0}=\rho-\rho_{1}$. Then the function $h_{\alpha}$ on $V\setminus \{0\}$ defined by

\[h_{\alpha} (v)=|v|^{-\alpha} (r*\widehat{\rho}_{0}) (v)=\frac{1}{t}\int_{0}^t x^{\alpha} \widehat{\rho}_{0} (\tilde v,x) dx\],

   is $b.R.i$ and one has  $(\delta_{u}\otimes \ell^{\alpha}) (r*\widehat{\rho}_{0})=\int_{0}^{\infty} t^{\alpha-1} \widehat{\rho}_{0} (u,t) dt=r_{\alpha}(u)$ where $t\rightarrow t^{\alpha-1} \widehat{\rho}_{0} (u,t)$ is Riemann-integrable on $]0,\infty[$ in generalised sense.
\end{lemm}
\begin{proof}

  By definition
\begin{align*}
h_{\alpha} (v)&=|v|^{-\alpha} (r*\widehat{\rho}_{0}) (v)=t^{\alpha} (r*\widehat{\rho}_{0}) (\frac{u}{t})=\frac{1}{t} \int_{0}^t y^{\alpha} \widehat{\rho}_{0} (\frac{u}{y}) dy=r^{\alpha} (u,t),\\
(\delta_{u}\otimes \ell^{\alpha}) (r* \widehat{\rho}_{0})&=
\int_{0}^{\infty} t^{\alpha-1} \frac{dt}{t^{\alpha+1}} \int_{0}^t y^{\alpha} \widehat{\rho}_{0} (\frac{u}{y}) dy=\displaystyle\mathop{\lim}_{T\rightarrow \infty} \int_{0}^T \frac{dt}{t^2} \int_{0}^t y^{\alpha} \widehat{\rho}_{0} (\frac{u}{y}) dy.
\end{align*}

  Lemma 5.5 implies that $h_{\alpha} (v)$ is b.R.i and has limit 0 at $|v|=\infty$. Integration by parts in the above formula gives

  $\int_{0}^T \frac{dt}{t^2} \int_{0}^t y^{\alpha} \widehat{\rho}_{0} (\frac{u}{y}) dy=-\frac{1}{T} \int_{0}^T y^{\alpha} \widehat{\rho}_{0} (\frac{u}{y}) dy+\int_{0}^T t^{\alpha-1}\widehat{\rho}_{0} (\frac{u}{t}) dt= -r^{\alpha} (u,T)+\int_{0}^T t^{\alpha-1} \widehat{\rho_{0}} (\frac{u}{t}) dt$.

  Since, using Lemma 5.5, $\displaystyle\mathop{\lim}_{T\rightarrow \infty} r^{\alpha} (u,T)=0$, it follows that $\int_{0}^T t^{\alpha-1} \widehat{\rho}_{0} (\frac{u}{t}) dt$ has a finite limit at $T=\infty$, hence $t\rightarrow t^{\alpha-1}\widehat{\rho}_{0} (\frac{u}{t})$ is Riemann-integrable on $\mathbb R_{+}^{*}$ in generalised sense and, for $u\in \mathbb S^{d-1}$,
\[(\delta_{u}\otimes \ell^{\alpha}) (r*\widehat{\rho}_{0})=\int_{0}^{\infty} t^{\alpha-1}\widehat{\rho}_{0} (\frac{u}{t}) dt=r_{\alpha}(u).\]
\end{proof}

\begin{lemm}
\label{lemm:5.7}
  If $\displaystyle\mathop{\lim}_{s\rightarrow s_{\infty}} k(s)>1$, then $\supp\rho$ is unbounded.
  \end{lemm}
\begin{proof}

  In order to show that $\supp\rho$ is unbounded, in view of Proposition 5.1, it suffices to show that there exists $g\in T$ with $\displaystyle\mathop{\lim}_{n\rightarrow \infty} |g^n|^{1/n}>1$. If not, then the trace $Trg$ of $g$  is bounded by $d$ on $T$. On the other hand, condition i-p implies the irreducibility of the action of $T$ on $V\otimes \mathbb C$, as shown now. Let $W\subset V \otimes \mathbb C$ be a proper $T$-invariant subspace of $V\otimes \mathbb C$ and $\overline{W}$ its complex conjugate. Then $W\cap\overline{W}$ and $W+\overline{W}$ are complexified subspaces of subpaces of $V$ which are also $T$-invariant. Using irreducibility of $T$ we get $W\cap \overline{W}=\{0\}$, $W+\overline{W}=V\otimes \mathbb C$ hence $V\otimes \mathbb C=W\oplus \overline{W}$. Let $g\in T^{\textrm{prox}} $ and $v\in V\setminus \{0\}$ with $g v=\lambda_{g} v$ and write $v=w+\bar w$ with $w\in W$, hence since $W, \overline{W}$ are $g$-invariant $gw=\lambda_{g} w$, $g\bar w=\lambda_{g}\bar w$. Since $\lambda_{g}$ is a simple eigenvalue we get $w=\bar w$ i.e $W\cap \overline{W}\neq \{0\}$ which gives a contradiction with condition i-p, hence $T$ acts irreducibly on $V\otimes \mathbb C$. Then Burnside's density theorem implies that $\textrm{EndV}\otimes \mathbb C)$ is generated as an algebra by $T$, i.e. there exists a base $g_{i}\in T$ $(i=1, ., d^2)$ of $\textrm{End V}$. Then the linear forms $u\rightarrow Tr(ug_{i})\ \ (i=1,.,d^2)$ form a basis of the dual space of $\textrm{ End V}$. In particular, for some constant $c>0$ we have for any $g\in T$, $|g|\leq c \displaystyle\mathop{\Sigma}_{i=1}^{d^2} |Tr (g g_{i})|\leq c d^3$. Then for any $n\in \mathbb N$, and $s>0$,
$\int|g|^s d\mu^n (g)\leq c d^{3s},\ \ k(s)\leq 1$.
\end{proof}

\textit{Proof of Proposition 5.3}
   If $u= t v$ with $u\in \mathbb S^{d-1}$, $t>0$ the function $\widehat{\rho} (v)=\widehat{\mathbb P}\{\langle R,u \rangle> t\}$ is bounded,  right continuous and non increasing.
 Since equation $(S)$ can be written as $R-B=A R\circ \widehat{\theta}$ and $A, R\circ \widehat{\theta}$ are independant we have
$\rho_{1}=\mu * \rho, \ \ \rho-\mu * \rho=\rho-\rho_{1}=\rho_{0}.$ Furthermore,

\[\rho=\displaystyle\mathop{\Sigma}_{0}^{n} \mu^k *\rho_{0}+\mu^{n+1}*\rho, \ \ \widehat{\rho}(v)=\displaystyle\mathop{\Sigma}_{0}^{n} ((\mu^{*})^k *\delta_{v}) (\widehat{\rho}_{0})+((\mu^{*})^{n+1}*\delta_{v}) (\widehat{\rho}).\]

  Also if $r^{-}$ denote the push-forward of $r$ by $x\rightarrow x^{-1}$
then
$r^-*\rho=\displaystyle\mathop{\Sigma}_{0}^{n} \mu^k * (r^- *\rho_{0})+\mu^{n+1}*(r^-*\rho)$.

  Since $L_{\mu}<0$ the subadditive ergodic theorem applied to $\log |S_{n} (\omega)|$ gives the convergence of $S_{n}(\omega) v$ to 0. In particular, for $\xi \in M^1 (V)$, the sequence $\mu^n *\xi$ converges in law to $\delta_{0}$, hence  $\displaystyle\mathop{\lim}_{n\rightarrow \infty} \widehat{\mu^n *\xi} (v)=(\mu^{*})^n * \delta_{v} (\widehat{\xi})=0$.

  From the above  convergence on $V$,we  have
  \[\rho-\delta_{0}=\displaystyle\mathop{\Sigma}_{0}^{\infty} \mu^k *\rho_{0},\ \ \widehat{\rho} (v)=\displaystyle\mathop{\Sigma}_{0}^{\infty} ((\mu^{*})^{k}*\delta_{v}) (\widehat{\rho_{0}}).\]

  But, by Proposition 5.1,  $\rho(\{0\})=0$, hence the stated vague convergence  of $\displaystyle\mathop{\Sigma}_{0}^{\infty} \mu^k *\rho_{0}$.
Since the sequence $(\mu^{n+1}*(r^-*\rho)) (\psi)$ converges to zero for any bounded Borel function $\psi$ on $V$ such that $\displaystyle\mathop{\lim}_{v\rightarrow 0} \psi (v)=0$, we have on such functions :
$r^-*\rho=\displaystyle\mathop{\Sigma}^{\infty}_{0} (\mu^{*})^k * (r^- * \rho_{0}).$
  We observe that, for any bounded measure $\xi$, we have $(\mu^k *\widehat{\xi}) (v)=((\mu^*)^k * \delta_{v}) (\widehat{\xi})$  and  $\widehat{r^-*\xi}=r* \widehat{\xi}$. It follows from the above equality that the potential $\displaystyle\mathop{\Sigma}^{\infty}_{0} ((\mu^*)^k * \delta_{v}) (r * \widehat{\rho}_{0})$ is finite and equal to $(r * \widehat{\rho}) (v)$.

   We have observed in Lemma 5.6 that the function $v\rightarrow |v|^{-\alpha} (r*\widehat{\rho}_{0}) (v)$ is b.R.i, hence the renewal Theorem \ref{the:4.7'} applied to $\mu^*$ and to the function $r*\widehat{\rho}_{0}$ gives for $u\in \mathbb S^{d-1}$,
\[\lim_{t\rightarrow \infty} t^{\alpha} (r*\widehat{\rho}) (u,t) = \frac{^{*}e^{\alpha}(u)}{L_{\mu}(\alpha)} (^{*}\widetilde{\nu}^{\alpha}_{u} \otimes \ell^{\alpha}) (r*\widehat{\rho}_{0}).\]
  Since  for fixed $u$, $\widehat{\rho} (u,x)=\widehat{\mathbb P}\{\langle R, u \rangle  >x\}$ is non increasing, Lemma 5.4 gives
\[\lim_{t\rightarrow \infty} t^{\alpha} \widehat{\rho} (u,t)=\frac{^{*}e^{\alpha}(u)}{L_{\mu} (\alpha)} {^{*}\widetilde{\nu}^{\alpha}_{u} }(r_{\alpha}).\]
  In particular, we have ${^{*}\widetilde{\nu}^{\alpha}_{u}} (r_{\alpha})\geq 0$. In case I this gives $^{*}\widetilde{\nu}^{\alpha} (r_{\alpha})\geq 0$ since ${^{*}\widetilde{\nu}^{\alpha}_{u}}={^{*}\widetilde{\nu}^{\alpha}}$.

  In case II, taking $u\in \Lambda_{+}(T^{*})$ and using $^{*}p^{\alpha}_{+} (u)=1$, this gives $^{*}\nu^{\alpha}_{+} (r_{\alpha})\geq 0$. Also, in the same way $^{*}\nu^{\alpha}_{-} (r_{\alpha})\geq 0$. Furthermore, in case II, using Theorem \ref{the:4.7'},
\[^{*}\widetilde{\nu}^{\alpha}_{u} (r_{\alpha})={^{*}p^{\alpha}_{+}} (u) {^{*}\nu^{\alpha}_{+}} (r_{\alpha})+{^{*}p^{\alpha}_{-}} (u){ ^{*} \nu_{-}^{\alpha}}(r_{\alpha}).\]

   If ${^{*}\widetilde{\nu}^{\alpha}} (r_{\alpha})>0$, in case II we can define a probability measure $\sigma^{\alpha}$ on $\widetilde{\Lambda} (T)$ by
\[^{*}\widetilde{\nu}^{\alpha} (r_{\alpha}) \sigma^{\alpha}=\frac{1}{2} (^{*}\nu^{\alpha}_{+} (r_{\alpha}) \nu_{+}^{\alpha}+{^{*}\nu^{\alpha}_{-}} (r_{\alpha}) \nu_{-}^{\alpha}),\]
   while in case I , $\sigma^{\alpha}=\widetilde{\nu}^{\alpha}$. If $^{*}\widetilde{\nu}^{\alpha} (r_{\alpha})=0$, we have also $^{*}\nu^{\alpha}_{+} (r_{\alpha})={^{*}\nu^{\alpha}_{-}} (r_{\alpha})=0$, hence we can leave $\sigma^{\alpha}$  with projection $\nu^{\alpha}$ on $\mathbb P^{d-1}$ undefined in the above formulae. In any case $\sigma^{\alpha}\otimes \ell^{\alpha}$ is $\mu$-harmonic.

  We get another expression for the above limit, by using the formulae for $^{*}e^{\alpha} (u)$, $^{*}p^{\alpha}_{+} (u)$, $^{*}p^{\alpha}_{-} (u), p(\alpha)$ of section 2 paragraph 3, for case  II as follows (see Theorem \ref{thm:2.17}), with
\begin{align*}
p(\alpha)&=\int |\langle x,  y \rangle|^{\alpha} d\nu^{\alpha}(x) d{^{*}\nu}^{\alpha} (y),\\
^{*}e^{\alpha} (u) ^{*}p^{\alpha} _{+}(u) p(\alpha)&=\int \langle u,u' \rangle^{\alpha}_{+} d \nu_{+}^{\alpha} (u')=\alpha (\nu_{+}^{\alpha} \otimes \ell^{\alpha}) (H_{u}^{+}).
\end{align*}

  From above, we get
\[^{*}e^{\alpha} (u) ^{*}\widetilde{\nu}^{\alpha}_{u} (r_{\alpha})=\frac{\alpha}{ p(\alpha)} ((^{*}\nu^{\alpha}_{+} (r_{\alpha}) \nu^{\alpha}_{+} + ^{*}\nu^{\alpha}_{-} (r_{\alpha}) \nu^{\alpha}_{-} ) \otimes \ell^{\alpha}) (H^{+}_{u}).\]

  Hence, with $C=2  \frac{^{*}\widetilde{\nu}^{\alpha} (r_{\alpha})\alpha}{L_{\mu}(\alpha) p(\alpha)}$,  $\sigma^{\alpha}$ as above and $C_{+} =  \frac{{^{*}\widetilde{\nu}^{\alpha}_{+}} (r_{\alpha})\alpha}{L_{\mu}(\alpha) p(\alpha)},\ C_{-}= \frac{{^{*}\widetilde{\nu}^{\alpha}_{-}} (r_{\alpha})\alpha}{L_{\mu}(\alpha) p(\alpha)}$,
\[\mathop{\lim}_{t\rightarrow \infty} t^{\alpha} \widehat{\rho} (u,t) =C (\sigma^{\alpha} \otimes \ell^{\alpha}) (H^{+}_{u}), C \sigma^{\alpha}= C_{+} \nu^{\alpha}_{+} + C_{-} \nu^{\alpha}_{-}.\]

  In case I we get the corresponding formula.  The fact that $\supp\rho$ is unbounded follows from Lemma 5.7, since $k'(\alpha)>0$ and $k(\alpha)=1$.
\hfill$\square$

\begin{coro}
\label{cor:5.8}
  For any $v\in V\setminus \{0\}$, we have
\[\lim_{t\rightarrow \infty} t^{\alpha} \widehat{\mathbb P} \{|\langle R,v\rangle|>t\}=C \frac{p(\alpha)}{\alpha} ({^{*}e^{\alpha}}\otimes h^{\alpha}) (v),\]
  with $p(\alpha)=\int |\langle x,  y \rangle|^{\alpha} d\nu^{\alpha} (x) d^{*}\nu^{\alpha}(y)$ and $C=2 \frac{{^{*}\widetilde{\nu}^{\alpha}}(r_{\alpha})\alpha}{L_{\mu}(\alpha) p(\alpha)}\geq0$.
In particular, there exists $b>0$ such that $\widehat{\mathbb P}\{|R|> t\}\leq b t^{-\alpha}$.
\end{coro}
\begin{proof}

  By definition of ${^{*}\widetilde{\nu}^{\alpha}_{u}}$ and since ${^{*}p^{\alpha}_{+}} (u)={^{*}p^{\alpha}_{-}} (-u)$ we have, $\frac{1}{2}({^{*}\nu^{\alpha}_{u}}+^{*}\nu^{\alpha}_{-u})= {^{*}\widetilde{\nu}^{\alpha}}$. Hence, using  Proposition 5.3,

\[ \lim_{t\rightarrow \infty} t^{\alpha} \widehat{\mathbb P}\{ |\langle R,u \rangle> t\}=2 \frac{^{*}e^{\alpha}(u)}{L_{\mu}(\alpha)} {^{*}}\widetilde{\nu}^{\alpha} (r_{\alpha})=C \frac{p(\alpha)}{\alpha} \ {^{*}e^{\alpha}} (u).\]

  The formula in the corollary follows by $\alpha$-homogeneity, since it is valid for $|v|=1$.

  We take a base $u_{i}\in V$ $(1\leq i\leq d)$ and write $|R|\leq \displaystyle\mathop{\Sigma}_{i=1}^{d} |\langle R, u_{i}\rangle|$ . For $t$ large:
$\widehat{\mathbb P} \{|R|>t\}\leq \displaystyle\mathop{\Sigma}_{i=1}^{d} \widehat{\mathbb P} \{ |\langle R, u_{i}\rangle| >t \leq (C'+\varepsilon) t^{-\alpha}  b'^{\alpha} \displaystyle\mathop{\Sigma}_{i=1}^{d}  {^{*}e^{\alpha}}  (u_{i})$,  with $\varepsilon>0$, $C'=C \frac{p(\alpha)}{\alpha}$, hence the result.
\end{proof}

\subsection{A dual Markov walk and the positivity of directional tails}

  The following proposition will play an essential role in the discussion of positivity for $ C_{+}, C_{-}$ and $C$ as defined in Proposition 5.3. We denote by $\Lambda_{a}^{*}(\Sigma)$ the set of elements $u$  in $\mathbb S^{d-1}$ such that the projection of $\Lambda_{a}(\Sigma)$ on the half line $\mathbb R_{+} u$ is unbounded. Here, instead of the vector space $V$ used in paragraph 2, duality in the context of Radon transforms lead us to consider a $\lambda$-random walk on the larger vector space $V\times \mathbb R$ and to use ideas of \cite{36}, for the analysis of corresponding measurable fibered  kernels (see subsection 4.1).  However, the continuity hypotheses of \cite{36} are not in general satisfied by these kernels.
\begin{prop}
\label{pro:5.9}
  With the hypothesis of Theorem 5.2  if $M\subset \mathbb S^{d-1}$ is $T^{*}$-minimal and $\Lambda^{*}_{a} (\Sigma)\supset M$, then for any $u\in M$,
\[C (u)=\displaystyle\mathop{\lim}_{t\rightarrow \infty} t^{\alpha}  \widehat{\mathbb P}\{\langle R,u\rangle >t\}=\displaystyle\mathop{\lim}_{t\rightarrow \infty} t^{\alpha}\widehat{\rho} (u,t)>0.\]
  In particular with the above notations, we have $C>0$.
\end{prop}
   We observe that $R_{n}=\displaystyle\mathop{\Sigma}_{0}^{n-1} A_{1}\cdots A_{k} B_{k+1}$ satisfies the relation $\langle R_{n+1},v\rangle =\langle R_{n},v\rangle+ \langle B_{n+1}, S'_{n} v\rangle$ where $S'_{n}=(A_{1}\cdots A_{n})^{*}$. Also $h=(g,b)\in H$ acts on $E=(V\setminus \{0\}) \times \mathbb R$ according to the formula $h(v,r)=(g^{*} v, r+\langle b,v\rangle)$, hence the pair,
\[(S'_{n}v,r+ \langle R_{n},v\rangle)=(v_{n}, r_{n}),\]
  is a random walk on the  right homogeneous $H$-space $E$. Actually, $V\times \mathbb R$ is a vector space and the above formula for $h(v,r)$ defines a right linear representation of $H$ in $V\times \mathbb R$ which leaves invariant $E\subset V\times \mathbb R$.  In particular, using the radial $\mathbb R^{*}_{+}-$fibration of this vector space, we see that the radial projection $(v,r)\rightarrow (u,p)$ with $v=|v| u$, $r=p |v|$ defines an $H$-equivariant projection from $E$ to $\mathbb S^{d-1} \times \mathbb R$, hence $(S'_{n}v, r+ <R_{n},v>)$ is  also a $\mathbb R_{+}^{*}$-fibered Markov chain over $\mathbb S^{d-1} \times \mathbb R$. Then we can write
\[E=(\mathbb S^{d-1}\times \mathbb R) \times \mathbb R_{+}^{*}\subset V\times \mathbb R.\]
  The action of $h=(g,b)$ on $\mathbb S^{d-1}\times \mathbb R$ is given by $h(u,p)=(g^{*}\cdot u, h^{u}p)$ with
\[h^{u}p=\frac{1}{|g^{*}u|} (p+\langle b,u\rangle).\]
  The proof of Proposition 5.9 is based on the relation $\langle R,v\rangle =\displaystyle\mathop{\lim}_{n\rightarrow \infty}$ $<R_{n},v>$ and on the dynamics of the random walk $(v_{n}, r_{n})=(S'_{n} v, r+ <R_{n},v>)$ on  $V\times \mathbb R$. We denote by ${^{*} \widehat P}$ its corresponding fibered Markov kernel   and we study an associated  ladder process $(x_{\tau_{n}}, W_{\tau_{n}})$, defined below. In terms of this process we can give a  new expression for $\widehat{\rho}(v)= \widehat{\mathbb P}\{\langle R, u \rangle  >t\}$ as a potential of a non negative function on $E$. Then we can use a weak renewal theorem for $(x_{\tau_{n}}, W_{\tau_{n}})$ and obtain Proposition 5.9.

   We will consider successively the two components $u_{n}=\frac{v_{n}}{|v_{n}|},  p_{n}=\frac{r_{n}}{|v_{n}|}$ and finally the fibered Markov chain $(v_{n}, r_{n})$ above $\mathbb S^{d-1}\times \mathbb R$.

  Let $M$ be a $T^{*}$-minimal subset of $\mathbb S^{d-1}$,  hence (see section 2, paragraphe 3), $M=\Lambda_{+}(T^{*})$ (or $\Lambda_{-} (T^{*}))$ in cases II, or $M=\Lambda(T)$ in case I.

  The following says that $\Lambda_{a}^{*} (\Sigma)$ is "large".
\begin{lemm}
\label{lem:5.10}
  In cases I or II' : $\Lambda_{a}^{*} (\Sigma)=\mathbb S^{d-1}$.
In case II'' : $\Lambda_{a}^{*} (\Sigma) \supset \Lambda_{+} (T^{*})$.
\end{lemm}
\begin{proof}
  Let $\Lambda_{a}^{\infty} (\Sigma)=\overline{\Lambda_{a}(\Sigma)} \cap \mathbb S^{d-1}_{\infty}$ and $u\in \mathbb S^{d-1}$,  $u'_{\infty} \in \Lambda_{a}^{\infty} (\Sigma)$ corresponds to $u'\in \mathbb S^{d-1}$. If $\langle u',u\rangle  >0$, then $u\in \Lambda_{a}^{*}(\Sigma)$. Hence the complement of $\Lambda_{a}^{*}(\Sigma)$ in $\mathbb S^{d-1}$ is contained in the set $\{u \in \mathbb S^{d-1}$; $\langle u,u'\rangle  \leq 0\ \forall u'_{\infty}\in \Lambda_{a}^{\infty} (\Sigma)\}$. From the discussion at the beginning of this section we know that $\Lambda_{a}^{\infty} (\Sigma) \neq \phi$ is $T$-invariant and closed, hence contains $\widetilde{\Lambda}^{\infty} (T)$ in cases I, II' or only $\Lambda_{+}^{\infty} (T)$ in case II'' with $\Lambda_{-}^{\infty} (T) \cap \Lambda_{a}^{\infty} (\Sigma)=\phi$.

  Since $\widetilde{\Lambda}^{\infty} (T)$ is symmetric and condition i-p is valid, it follows $\Lambda_{a}^{*} (\Sigma)=\mathbb S^{d-1}$ in cases I, II'. In case II'', we know from the end of proof of Theorem \ref{thm:2.17} that the complement of $\Lambda_{a}^{*}(\Sigma)$ is contained in $\widehat{\Lambda}_{+} (T^{*})=\{u\in \mathbb S^{d-1}\ ;\ \forall u' \in \Lambda_{+}(T), \langle u,u'\rangle >0\}$. Since $\Lambda_{+} (T^{*}) \cap-\widehat{\Lambda}_{+} (T^{*})=\phi$, we get $\Lambda_{+} (T^{*})\subset \Lambda_{a}^{*}(\Sigma)$.
  \end{proof}

    The random walk $ (v_{n}, r_{n})=(S'_{n} v, r+ \langle R_{n}, v\rangle )$ has $H$-equivariant projection $S'_{n} v$ on $V\setminus\{0\}$, the kernel ${^{*}\widehat{P}}$ has projection ${^{*} P}$ already defined in section 2, hence the positive homogeneous function ${^{*}e^{\alpha}} \otimes h^{\alpha}$, lifted to $E$, satisfies ${^{*}\widehat P} ({^{*}e}^{\alpha} \otimes h^{\alpha})= {^{*}e}^{\alpha} \otimes h^{\alpha}$, and we can consider the new relativized fibered Markov kernel ${^{*}\widehat P}_{\alpha}$ on $E$. If $(u,p)\in M\times \mathbb R$, the projection $x_{n}=(u_{n},p_{n})$ of $(v_{n}, r_{n})$  on $M\times \mathbb R$ depends on the kernel ${^{*}\widehat Q^{\alpha}}$ given by
    \[{^{*}\widehat Q^{\alpha}} \varphi (u,p)=\int \varphi (g^{*}
    \cdot u, h^u p)^{*} q^{\alpha} (u,g) d\lambda (h)\]
    where ${^{*}q^{\alpha}}$ corresponds to $q^{\alpha}$ as in section 3.  The important fact for the proof of Proposition 5.9 is that ${^{*}\widehat{Q}}^{\alpha}$ has a stationary probability $\kappa$ with $\kappa(M\times[t,\infty[)>0$ for any $t>0$, such that $p$ is not "too large" with respect to $\kappa$.

    For the analysis of ${^{*}\widehat{Q}^{\alpha}}$ we consider on $\widehat{\Omega}$ the projective limit ${^{*}\widehat {\mathbb Q}^{\alpha}_{u}}$ of the system ${^{*} q}^{\alpha}_{n} (u,\cdot) \lambda^{\otimes n}$ $(n\in \mathbb N)$ and, by abuse of notation, the corresponding expectation will be written $\mathbb E_{u}^{\alpha}$. Given a ${^{*}\widetilde Q^{\alpha}}$-stationary measure $\widetilde{\pi}_{M}^{\alpha}$, we  write  ${^{*}\widehat{\mathbb Q}^{\alpha}}=\int \delta_{u}\otimes{^{*}\widehat{\mathbb Q}^{\alpha}_{u}} d \widetilde{\pi}_{M}^{\alpha} (u)$ and we denote by $\mathbb E^{\alpha}$  the corresponding expectation symbol. We denote by $\widehat{\theta}^{\#}$ the map of $M\times H^{\mathbb Z}$ into itself defined by $\widehat{\theta}^{\#} (u,\widehat{\omega})=(g^{*}_{1}\cdot u, \widehat{\theta} \widehat{\omega})$ where $\widehat{\theta}$ is the bilateral shift on $H^{\mathbb Z}$, and ${^{*}\widehat{\mathbb Q}^{\alpha}}$ will again denote the natural $\widehat{\theta}^{\#}$-invariant measure on $M\times H^{\mathbb Z}$. Also we extend $S'_{n}(\omega)$ as a $G$-valued $\mathbb Z$-cocycle. If $\eta$ is a  probability measure on $X=M\times \mathbb R$, the associated Markov measure on ${^{a}\widehat {\Omega}}=X\times \widehat{\Omega}$, is denoted by ${^{*}\widehat {\mathbb Q}^{\alpha}_{\eta}}$,  and the extended shift by ${^{a}\widehat{ \theta}}$  where ${^{a}\widehat {\theta}} (x,\widehat{\omega})=(h_{1}x,\widehat{\theta} \widehat{\omega})$. Also if $\eta$ is ${^{*}\widehat{Q}^{\alpha}}$-stationary we will consider the bilateral associated system $(\Omega^{\#}, {^{a}\widehat{\theta}}, \eta^{\#})$ where $\Omega^{\#}=X\times H^{\mathbb Z}, {^{a}\widehat{\theta}}$ is the bilateral shift and $\eta^{\#}$ is the unique ${^{a}\widehat{\theta}}$-invariant measure with projection ${^{*}\widehat{\mathbb Q}^{\alpha}_{\eta}}$ on $X\times \widehat{\Omega}$.

\begin{lemm}
\label{lem:5.11}
  Let $M$ be a $T^{*}$-minimal subset of $\mathbb S^{d-1}$, $\pi_{M}^{\alpha}$ the unique ${^{*}\widetilde{Q}^{\alpha}}$-stationary measure on $M$. With the above notations, we consider the Markov chain $x_{n}=(u_{n},p_{n})$, on $X=M\times \mathbb R$ given by
\[u_{n+1}=g_{n+1}^{*} . u_{n},\ \ p_{n+1}=\frac{p_{n}+\langle b_{n+1},u_{n}\rangle }{|g_{n+1}^{*} u_{n}|},\ \  p_{0}=p,\ u_0=u,\]
  where $(g_{n}, b_{n})$ is distributed according to ${^{*}\widehat{\mathbb Q}_{u}^{\alpha}}$. Then, for any $p\in \mathbb R$, $x_{n}$ converges in ${^{*}\widehat{\mathbb Q}^{\alpha}}$-law to the  unique ${^{*}\widehat Q^{\alpha}}$-stationary measure $\kappa$, the projection of $\kappa$ on $M$ is $\pi_{M}^{\alpha},\  \kappa (M\times \{p\})=0$  and $\int |p|^{\varepsilon} d\kappa (u,p)<\infty$ for $\varepsilon$ small. We have ${\kappa^{\#}-\textrm{a.e}}$,
\[\limsup_{n\rightarrow \infty} |S'_{n} u| |p_{n}|=\infty, \ \ \lim_{n\rightarrow \infty} |S'_{-n}u| |p_{-n}|=0.\]
  If $\Lambda_{a}^{*} (\Sigma) \supset M$, then $\kappa(M\times ]t,\infty[)>0$ for any $t>0$ and $\displaystyle\mathop{\lim\sup}_{n\rightarrow \infty} |S'_{n} u| p_{n}=\infty$,  ${^{*}\widehat{\mathbb Q}_{\kappa}^{\alpha}}-\textrm{a.e}$.
\end{lemm}

\begin{proof}
  If will be convenient to use the functions $a(g,u), b(h,u)$ defined by $h^u p=a(g,u) p+b(h,u)$, and the random variables $a_{k}, b_{k}$ $(k\in \mathbb Z)$ defined by $a_{k}(\widehat{\omega},u)=a(g_{k}, S'_{k-1}\cdot u)$, $b_{k}(\widehat{\omega}, u)=b(h_{k}, S'_{k-1}\cdot u)$. Then we can express the action of $h_{n}\cdots h_{1}\in  H$ on $X$ as,
$u_{n}=S'_{n}\cdot u, \ \ y_{n}^p (u)=(h_{n}\cdots h_{1})^u p$, where
\[y_{n}^p (u)=a (S'_{n},u) p+ y_{n}^{\circ}(u) \textrm{ and } y_{n}^{\circ} (u)= \Sigma^{n}_{1} a_{k+1}^n (u) b_{k} (u),\]
  with $a_{k}^n (u)=a(g_{n}\cdots g_{k}, S'_{k-1}\cdot u)$.
The random variables $a_{k}, b_{k}$ are ${^{*}\widehat{\mathbb Q}^{\alpha}}$-stationary and $y_{n}^{\circ}$ has the same law as $p_{n}^{\circ} (\widehat{\omega},u)=\displaystyle\mathop{\Sigma}_{0}^{n-1} a_{-1} \cdots a_{-k} b_{-k-1}$.

  We estimate $\mathbb E^{\alpha} (|p_{n}^{\circ}|^{\varepsilon})$ for $0<\varepsilon<\tau$ and $\varepsilon$ small, where $p_{n}^{\circ}(\widehat{\omega},u)=\displaystyle\mathop{\Sigma}^{n-1}_{0} a_{-1}\cdots a_{-k}b_{-k-1}$. Since $(a_{k}\cdots a_{1}) (\omega,u)=a(S'_{k} (\omega), u)$ we get $\mathbb E^{\alpha} (|a_{-1}\cdots a_{-k}|^{\varepsilon})=\mathbb E^{\alpha} (|S'_{k} u|^{-\varepsilon})$. Hence Corollary 3.18 gives $\displaystyle\mathop{\lim}_{k\rightarrow \infty} (\mathbb E^{\alpha}(|a_{-1}\cdots a_{-k}|^{\varepsilon}))^{1/k}<1$ since $k'(\alpha)>0$ and $\mathbb E (|A|^{\alpha+\tau}) < \infty$. Also for $\varepsilon$ small,
\[\mathbb E^{\alpha} (|b_{k}|^{\varepsilon})=\mathbb E^{\alpha}(|\frac{\langle B_{1},u\rangle }{|g^{*}u|}|^{\varepsilon})\leq \mathbb E^{\alpha}(|B_{1}|^{\varepsilon} \gamma^{\varepsilon} (A))<\infty,\]
  using ${^{*}\widehat{\mathbb Q}^{\alpha}}$-stationarity, H\"older inequality and the condition $\mathbb E(|B_{1}|^{\alpha+\tau})+\mathbb E(|A_{1}|^{\alpha+\tau})<\infty$.

  Since for $0<\varepsilon<1$  $|p_{n}^{\circ}(\widehat{\omega},u)|^{\varepsilon}\leq \displaystyle\mathop{\Sigma}^{n-1}_{\circ}|a_{-1}\cdots a_{-k}|^{\varepsilon} |b_{-k-1}|^{\varepsilon}$, we get that $\mathbb E^{\alpha} (|p_{n}^{\circ}|^{\varepsilon})$ is bounded. The ${^{*}\widehat{\mathbb Q}^{\alpha}}$-\textrm{a.e}.\ convergence of  the partial sum $p_{n}^{\circ}(\widehat{\omega},u)$ to $p(\widehat{\omega},u)=\displaystyle\mathop{\Sigma}^{\infty}_{0} a_{-1}\cdots a_{-k} b_{-k-1}$ and the finiteness of $\mathbb E^{\alpha} (|p|^{\varepsilon})$ follows. By definition, $p(\widehat{\omega}, u)$ satisfies the functional equation  $p\circ{^{a}\widehat{\theta}}=a p+b$ where $p$ and $(a,b)$ are independent. It follows that the probability measure $\kappa$ on $M\times \mathbb R$ given by the formula
 \[\kappa=\int \delta_{u} \otimes \delta_{p(\widehat{\omega},u)} d\  {^{*}\widehat{\mathbb Q}^{\alpha}} (\widehat{\omega},u)\]
  is ${^{*}\widehat{Q}^{\alpha}}$-invariant.

  As observed above, the ${^{*}\widehat{\mathbb Q}^{\alpha}}$-laws of $y_{n}^{\circ}$ and $p_{n}^{\circ}$ are the same. Since the product of $\widetilde{\pi}^{\alpha}_{M}$ with the law of $y_{n}^{\circ}$ is $({^{*}\widehat Q^{\alpha}})^n$ $(\pi_{M}^{\alpha}\otimes \delta_{0})$ we have in weak topology : $\displaystyle\mathop{\lim}_{n\rightarrow \infty}$  ${(^{*}\widehat Q^{\alpha}})^n$ $(\pi_{M}^{\alpha}\otimes \delta_{0})=\kappa$.

  Since $|y_{n}^p (\widehat{\omega},u)-y_{n}^{p'} (\widehat{\omega},u)|=a(S'_{n}(\omega), u)| p-p'|$ \ \ and \ \ $a(S'_{n}(\omega),u)=|S'_{n} u|^{-1}$ converges ${^{*}\widehat{\mathbb\ Q}^{\alpha}_{u}}$-\textrm{a.e} to zero, we get the convergence of $({^{*}\widehat Q^{\alpha}})^n$ $(\pi_{M}^{\alpha}\otimes \delta_{p})$ to $\kappa$, for any $p$. On the other hand, if $\eta'$ is a ${^{*}\widehat Q^{\alpha}}$-stationary measure on $M\times \mathbb R$, its projection on $M$ is ${^{*}\widetilde Q^{\alpha}}$-stationary, hence equal to $\pi_{M}^{\alpha}$, since $M$ is $T^{*}$-minimal. Then, from above $({^{|*}\widehat Q^{\alpha}})^n \eta'$ converges to $\kappa$, hence $\eta'=\kappa$. The ${^{*}\widehat Q^{\alpha}}$-ergodicity of $\kappa$ implies the ${^{a}\widehat{\theta}}$-ergodicity of ${^{*}\widehat{\mathbb Q}^{\alpha}_{\kappa}}$ and $\kappa^{\#}$.

  If $\Lambda_{a}^{*}(\Sigma)\supset M$ assume $\kappa$ $(M\times ]t,\infty[)=0$, for some $t>0$,  i.e. the ${^{*}\widehat Q^{\alpha}}$-invariant set $\supp \kappa$ is contained in $M\times ]-\infty,t]$. Then, for any $(u,p)\in \supp \kappa$ we have $p+\langle R_{n},u\rangle  \leq t |S_{n}'u|$ ${^{*}\widehat{\mathbb Q}_{u}^{\alpha}}$-\textrm{a.e.}, i.e. $p+\langle R_{n}, u\rangle \leq t |S_{n}'u|$ ${^{*}q_{n}^{\alpha}} (u,\cdot) \lambda^{\otimes n}$-\textrm{a.e.} for any $n\in \mathbb N$ and for some $p\in\mathbb R$, $t$, $u\in M$. It follows $p+\langle R_{n},u\rangle \leq |S_{n}' u|$, $\lambda^{\otimes n}$-\textrm{a.e.} and since $\displaystyle\mathop{\lim}_{n\rightarrow \infty} |S_{n}' u|=0$, $\displaystyle\mathop{\lim}_{n\rightarrow \infty} R_{n}=R$, ${\mathbb P}$-\textrm{a.e.}, we have $\langle R,u\rangle  \leq -p$, $\widehat{\mathbb P}$-\textrm{a.e.} . This implies that the support of the projection of $\rho$ on $\mathbb R u$ is bounded in direction $u$; since by Proposition 5.1 we have $\supp \rho=\Lambda_{a}(\Sigma)$ this  contradicts the condition $\Lambda_{a}^{*}(\Sigma) \supset M$. Hence $\kappa(M\times ]t,\infty[)>0$ for any $t>0$.

  Furthermore, arguments as in the proof of Proposition 5.1, using that $\supp \lambda$ has no fixed point in $V$, show $\kappa (M\times \{p\})=0$ for any $p\in \mathbb R$. From Theorem 3.10 we know that $\displaystyle\mathop{\lim}_{n\rightarrow \infty} \frac{1}{n} \log |S'_{n}u|=L_{\mu}(\alpha)>0$, ${^{*}\mathbb Q^{\alpha}_{u}}$-\textrm{a.e.}\ \  .Furthermore, since $\kappa$ is ${^{*}\widehat Q^{\alpha}}$-ergodic and $\kappa (M\times \{0\})=0$ we have $\displaystyle\mathop{\lim\sup}_{n\rightarrow \infty} |p_{n}|>0$,
  $\widehat{\mathbb Q}_{\kappa}^{\alpha}$-\textrm{a.e.}. Then we get
  $\displaystyle\mathop{\lim\sup}_{n\rightarrow \infty}  |S'_{n}u|\  |p_{n}|=\infty$, ${^{*}\widehat{\mathbb Q}_{\kappa}^{\alpha}}$-\textrm{a.e}. If $\Lambda_{a}^{*} (\Sigma) \supset M$, from above we have
  $\kappa (M\times ]0,\infty[>0$, and again using ergodicity, $\displaystyle\mathop{\lim\sup}_{n\rightarrow \infty} p_{n}>0$.
  Since  $\displaystyle\mathop{\lim}_{n\rightarrow \infty}  |S'_{n} u|=\infty$ ${^{*}\widehat{\mathbb Q}^{\alpha}_{u}}-\textrm{a.e}$, it follows  $\displaystyle\mathop{\lim\sup}_{n\rightarrow \infty}  |S'_{n}u| p_{n}=\infty$, ${^{*}\widehat{\mathbb Q}_{\kappa}^{\alpha}}$-\textrm{a.e}.
  Using Theorem 3.2, we have also, for any $u \in M$ and ${^{*}\widehat{\mathbb Q}^{\alpha}_{u}}$-\textrm{a.e}.: $\displaystyle\mathop{\lim}_{n\rightarrow \infty} \frac{1}{n} \log |S'_{-n} u|=-L_{\mu} (\alpha) <0$.
  The condition $\int |p|^{\varepsilon} d\kappa (u,p)<\infty$ implies $\displaystyle\mathop{\lim\sup}_{n\rightarrow \infty} \frac{\log| p_{n}|}{|n|}\leq 0,  \kappa^ \#$-\textrm{a.e}.
  Then we get, $\displaystyle\mathop{\lim}_{n\rightarrow \infty} |S'_{-n} u| |p_{-n}|=0,  \kappa^{\#}$-\textrm{a.e}.
  \end{proof}

  For the analysis of ${^{*}\widehat{P}_{\alpha}}, \widehat{\rho}$ we consider the optional time $\tau$ on $X\times \widehat{\Omega}$  given for $p\neq 0$ by:
$\tau=\inf  \{n>0\ ;\ p^{-1} \langle R_{n},u\rangle >0\}$, $\tau=\infty$ if $p^{-1}\langle R_{n},u\rangle \leq 0$ for every $n$.

  We observe that $\tau$ is independent  on $p$ as long as $p>0$ or $p<0$.
By definition of $p_{n}$:
\[p+\langle R_{n},u\rangle=p_{n} |S'_{n} u|, p^{-1} p_{n}|S'_{n} u|=1+p^{-1}\langle R_{n},u\rangle.\]

  We note that $(u_{n}, p_{n})$ is the radial projection of $(v_{n}, r_{n})$ on
  $\mathbb S^{d-1} \times \mathbb R \subset V\times \mathbb R$ while the projection $(0, r_{n})$ of $(v_{n}, r_{n})$ on $\{0\}\times \mathbb R$ satisfies $r_{n}=r (p^{-1}p_{n}|S'_{n}u|)$.
  In $(u,p,r)$ coordinates on the set $\{r\neq 0\}$, the process
  $(v_{n}, r_{n})$ can be written as $(x_{n},r_{n})$, and for any $t>0$   the dilation $(v,r)\rightarrow t(v,r)$ reduces to $(x,r)\rightarrow (x,tr)$. Since $\kappa(M\times\{0\})=0$, $^*\hat{P}$ is also a measurable Markov fibered kernel above $\mathbb{S}^{d-1}\times \mathbb{R}^*$ (see subsection 4.1). This measurable setting will be useful below.Also,
\[\tau=\inf  \{n>0\ ;\ p^{-1} p_{n} |S'_{n} u|>1\},\  \tau=\infty \textrm{ if } p^{-1} p_{n} |S'_{n}u|\leq 1, \textrm{ for every } n.\]
  In particular $p^{-1} p_{\tau} >0$ where the notation $p_{\tau}$ is used if $\tau$ is finite. Also we define $\sigma_{n}=\tau \circ ({^{a}\widehat{\theta}})^{\sigma_{n-1}}$ and $\sigma_{0}=0$, $\tau_{n}=\displaystyle\mathop{\Sigma}_{1}^{n} \sigma_{k}$, so that $\tau$ can be seen as the first ladder index and $p^{-1} p_{\tau} |S'_{\tau} u|$ as the first ladder height of the $\mathbb R^{*}$-valued measurable $\mathbb Z$-cocycle ,
\[W_{n} (u,p,\widehat{\omega})=p^{-1} p_{n} |S'_{n}u|\]
   over the dynamical system $(\Omega^{\#}, {^{a}\widehat{\theta}},\kappa^{\#})$, which is well defined  since $\kappa (M\times \{0\})=0$. The random times $\tau_{n}$ can be seen as the successive times of increase for $r^{-1} r_{k}=p^{-1} p_{k} |S'_{k} u|$ along the random walk $(v_{n}, r_{n})$.  On the other hand, by Poincar\'e recurrence theorem we have $\kappa^{\#}$-\textrm{a.e}, $\displaystyle\mathop{\lim\sup}_{n\rightarrow \infty}  p^{-1} p_{n}\geq 1$. Since $\displaystyle\mathop{\lim}
_{n\rightarrow \infty} |S'_{n}u|=\infty$,then  $\tau$, $\tau_{n}$ are finite ${^{*}\widehat{\mathbb Q}^{\alpha}_{\kappa}}$-\textrm{a.e} and we have $p^{-1} p _{\tau_{n}}>0$.
\vskip 3mm
   With the above notation the $\lambda$-random walk $(v_{n}, r_{n})$ on $E\subset V\times \mathbb R$ can be written as
\[u_{n+1}=g^{*}_{n+1} . u_{n},\ \ p_{n+1}=\frac{p_{n}+\langle b_{n+1}, u_{n}\rangle}{|g^{*}_{n+1} u_{n}|},\ \ r_{n+1}=r_{n} |g^{*}_{n+1} u_{n}| p_{n+1} p^{-1}_{n}.\]
 We denote
$X=M\times \mathbb R$ and  $X_{+}=M\times ]0,\infty[ \subset X,$
  hence $x\in X_{+}$ implies $x_{\tau}\in X_{+}$ if $\tau<\infty$.
 We consider also the submarkovian stopped operator ${^{*}\widehat{P}^{\tau}}$ on $V\setminus \{0\} \times \mathbb R$. In $(u,p,r)$ coordinates on the set $\{r \neq 0\}$ with $x=(u,p)$ the associated process at time $n$  starting at $(x,r)\in X_{+}\times \mathbb R^{*}_{+}$ is $(x_{\tau_{n}}, r |S'_{\tau_{n}}u| p_{\tau_{n}} p^{-1})$, hence the restriction of ${^{*}\widehat P^{\tau}}$ to $X_{+}\times \mathbb R_{+}^{*}$ is well defined and  is a  measurable fibered kernel on $X_{+}\times \mathbb R^{*}_{+} \subset E$. This restriction will be again denoted by ${^{*}P^{\tau}}$. Since ${^{*}\widehat{P}} ({^{*}e^{\alpha}} \otimes h^{\alpha} )={^{*}e^{\alpha}} \otimes h^{\alpha}$ and $\tau$ is finite ${^{*}\mathbb Q^{\alpha}_{\kappa}}$-\textrm{a.e.} , the  kernel ${^{*}\widehat{P}^{\tau}_{\alpha}}$ given by :
\[{^{*}\widehat{P}^{\tau}_{\alpha}} \varphi=({^{*}e^{\alpha}} \otimes h^{\alpha})^{-1} \widehat P^{\tau} ( {^{*}e^{\alpha}} \otimes h^{\alpha} \varphi),\]
  is  also a measurable Markov fibered  kernel on $X\times \mathbb R^{*}_{+}$ which satisfies ${^{*}\widehat{P}_{\alpha}^{\tau}}1=1$, $\kappa \otimes \ell$ \textrm{a.e}.

  The following  lemma expresses the function $\widehat{\rho}(p^{-1}v)=\widehat{\mathbb P}\{p^{-1}\langle R,u\rangle >t\}=\psi (v,p)$ on $E$ as a ${^{*}\widehat{P}^{\tau}}$-potential of a non negative function on $E$.

\begin{lemm}
\label{lem:5.12}
  With $t v=u \in \mathbb S^{d-1}$, $t>0$  and $p\neq 0$, we write
\[\tau=\inf  \{n>0 ; p^{-1}\langle R_{n},u\rangle  >0\}, \tau=\infty \textrm{ if } p^{-1}\langle R, u \rangle\leq 0 \textrm{ for any } n\in \mathbb N,\]
$\psi(v,p)=\widehat{\mathbb P}\{p^{-1}\langle R,u\rangle >t\}$,
$\psi_{\tau} (v,p)=\widehat{\mathbb P}\{t<p^{-1}\langle R,u\rangle  \leq t+p^{-1} \langle R_{\tau}, u\rangle  ; \tau <\infty\}$, where $R_{\tau}=\displaystyle\mathop{\Sigma}_{0}^{\tau-1} A_{1}\cdots A_{k} B_{k+1}$, $\psi^{\alpha}=({^{*}e^{\alpha}}\otimes h^{\alpha})^{-1} \psi$, $\psi_{\tau}^{\alpha}=({^{*}e^{\alpha}}\otimes h^{\alpha})^{-1}  \psi_{\tau}.$
Then,
\[\psi=\sum_{0}^{\infty} ({^{*}\widehat{P}^{\tau}})^k \psi_{\tau}, \ \ \psi^{\alpha}=\sum_0^\infty
(^*\widehat{P}^{\tau}_{\alpha})^k \psi_{\tau}^{\alpha}.\]
\end{lemm}
\begin{proof}
We write $\langle R-R_{n},v \rangle =\langle R^{n}, S'_{n}v\rangle $, where $R^{n}=\displaystyle\mathop{\Sigma}_{n}^{\infty} A_{n}\cdots A_{k} B_{k+1}$ hence if $\tau < \infty$ : $\langle R-R_{\tau}, v\rangle =\langle R^{\tau}, S'_{\tau} v\rangle $. By definition of $\tau$, since $p^{-1}\langle R,u\rangle >t>0$ and the convergence of $R_{n}$ to $R$  imply $\tau<\infty$,  we have
\[\psi_{\tau} (v,p)=\psi (v,p)-\widehat{\mathbb P}\{\langle R-R_{\tau},u\rangle p^{-1}>t\ ;\ \tau<\infty\}.\]
  On the other hand, since $p^{-1} p_{\tau}>0$,
  \[\widehat{\mathbb P}\{\langle R-R_{\tau}, u\rangle p^{-1} > t\ ;\ \tau <\infty \}=\widehat{\mathbb P}\{\langle R^{\tau}, u_{\tau}\rangle p^{-1} >\frac{t}{|S'_{\tau} u|}\ ;\ \tau <\infty\}={^{*}\widehat{P}^{\tau}} \psi (v,p).\]
It follows $\psi_{\tau}=\psi - {^{*}\widehat{P}^{\tau}} \psi$, \ \ $\psi=\displaystyle\mathop{\Sigma}_{0}^{n-1} ({^{*}\widehat{P}^{\tau}})^k \psi_{\tau}+ {^{*}\widehat{P}^{\tau_{n}}} \psi$ with
\[{^{*}\widehat{P}^{\tau_{n}}} \psi (v,p)=\widehat{\mathbb P}\{t |S'_{\tau_{n}}u|^{-1} < p^{-1}\langle R^{\tau_{n}}, u_{\tau_{n}}\rangle\ ;\ \tau_{n}<\infty\}.\]

  For $(x, \widehat{\omega})\in M\times \widehat{\Omega}$ we have either $\tau_{n} (x,\widehat{\omega})=\infty$ for some $n$ hence ${^{*}P}^{\tau_{n}}\psi=({^{*}P}^{\tau})^{k} \psi=0$ or $\displaystyle\mathop{\lim}_{n\rightarrow \infty} \tau_{n}(x,\widehat{\omega})=\infty$. In the second case, $\displaystyle\mathop{\lim}_{n\rightarrow \infty} |S'_{\tau_{n}}u|^{-1}=\infty$, $\widehat{\mathbb P}-\textrm{a.e.}$ and, since $R_{\tau_{n}}$ converges  $\widehat{\mathbb P}-\textrm{a.e}$ to  $R$, we have   $\displaystyle\mathop{\lim}_{n\rightarrow \infty} {^{*}\widehat{P}^{\tau_{n}}} \psi=0$, $\psi=\displaystyle\mathop{\Sigma}_{0}^{\infty} ({^{*}\widehat{P}^{\tau}})^k \psi_{\tau}$.
 The last relation follows from the definitions of $\psi_{\tau}^{\alpha}$ and ${^{*}\widehat{P}^{\tau}_{\alpha}}$, since $\psi_{\tau}^{\alpha}$ is non negative.
 \end{proof}

\begin{rema}
  The function $\psi$ satisfies the basic property of potentials for ${^{*}\widehat{P}^{\tau}}$,
\[{^{*}\widehat{P}}^{\tau} \psi \leq \psi, \displaystyle\mathop{\lim}_{n\rightarrow \infty} ({^{*}\widehat{P}}^{\tau})^{n}\  \psi=0.\]
  This property is a key for understanding non triviality of the  Cramer-type estimate for $\psi$. It is also  valid for other natural functions such as  $\overline{\psi}$ defined by
\[\overline{\psi} (v,p)=\widehat{\mathbb P}\{\displaystyle\sup_{n\geq 1} |W_{n}|>t\}.\]
\end{rema}
  Let $\tau$ be as above and $\Lambda_{a}^{*} (\Sigma) \supset M$, hence using Lemma 5.11 we have  $\kappa(X_{+})>0$.   Let ${^{*}\widehat{Q}^{\alpha,\tau}} (x,.)$ be the law of $x_{\tau}$ under ${^{*}\widehat{\mathbb Q}^{\alpha}_{x}}$.  Since $\tau<\infty, {^{*}\widehat{\mathbb Q}}_{\kappa}^{\alpha}$-\textrm{a.e.} ,
 this kernel is  a measurable Markov kernel with respect to $\kappa$ on $X_{+}$, hence ${^{*}\widehat{Q}}^{\alpha,\tau}$ is not $\kappa$-ergodic in general and it is natural to consider the first return time to $X_{+}$ as well as the corresponding induced Markov operator ${^{*}\widehat Q^{\alpha}_{+}}$ on $X_{+}$.
 Following the idea of (\cite{36}, Lemma 2) we interpret the ladder index $\tau$ as a first return time to a subset of $\Omega^{\#}$, we construct a stationary measure for ${^{*}\widehat{Q}^{\alpha,\tau}}$ on $X_{+}$, and we show the finiteness of the corresponding expectation of $\tau$.

\begin{lemm}
\label{lem:5.13}
  Assume $\Lambda_{a}^{*} (\Sigma)\supset M$ and write
\begin{align*}
\tau&=\inf  \{n>0\ \ ;\ \ p^{-1} p_{n} |S'_{n} u|>1\}=\inf  \{n>0\ ; p^{-1}\langle  R_{n},u\rangle >0\},\\
\tau&=\infty \textrm{   if } p^{-1} p_{n}|S'_{n}u|\leq 1\textrm{ for any } n\in \mathbb N.
\end{align*}
Then the stopped operator ${^{*}\widehat{Q}^{\alpha,\tau}}$  preserves $X_{+}$ and admits a  stationary ergodic  probability $\kappa^{\tau}$ on $X_{+}$ which is absolutely continuous with respect to $\kappa$.
The integral $\mathbb E_{0}^{\alpha} (\tau)=\int \mathbb E_{u}^{\alpha} (\tau) d\kappa^{\tau} (u,p)$ is finite.

With $\gamma_{\tau}^{\alpha}=L_{\mu} (\alpha) \mathbb E_{0}^{\alpha} (\tau)$ we have $\displaystyle\mathop{\lim}_{n\rightarrow \infty} \frac{1}{n} \log (|S'_{\tau_{n}} u| \frac{p_{\tau_{n}}}{p})=\gamma_{\tau}^{\alpha} \in ]0,\infty[, {^{*}\widehat{\mathbb Q}^{\alpha}_{\kappa^{\tau}}}$-\textrm{a.e.} .
 In particular $W_{\tau}= \log (p_{\tau} p^{-1}  |S'_{\tau} u|)$ has  finite expectation with respect to ${^{*}\widehat{\mathbb Q}^{\alpha}_{\kappa^{\tau}}}$.
\end{lemm}

\begin{proof}
Since $\Lambda_{a}^{*} (\Sigma)\supset M$, Lemma 5.11 gives $\kappa(X_{+})>0$. In order to deal only with positive values of the $\mathbb R^{*}$-valued $\mathbb Z$-cocycle $W_{n}(u,p,\widehat{\omega})=|S'_{n} u| p_{n} p^{-1}$  under ${^{*}\widehat{\mathbb Q}}_{\kappa}^{\alpha}$, it is convenient to consider the two sided Markov chain $x_{n_{k}} (k\in \mathbb Z,\ x_{\circ} \in X_{+})$ induced on $X_{+}$ by $x_{n}$, hence $n_{1}$ is the first return time of $x_{n}$ to $X_{+}$ and $n_{1}\leq \tau$ since $p^{-1} p_{\tau}>0$. We note that the normalized restriction $\kappa_{+}$ of $\kappa$ is a stationary ergodic measure for $x_{n_{k}}$. Also the relativized Markov kernel ${^{*}\widehat{P}_{\alpha}}$ on $E=X\times \mathbb R_{+}^{*}$ induces a fibered  measurable Markov kernel ${^{*}\widehat{P}_{\alpha,+}}$ on $X_{+}\times \mathbb R^{*}_{+}$ with projection ${^{*}\widehat Q_{+}^{\alpha}}$ on $X_{+}$, which satisfies ${^{*}P}_{\alpha,+} (\kappa_{+}\otimes \ell)=\kappa_{+}\otimes \ell$. Since $p^{-1}p_{n_{k}}>0$,in(u,p,r) coordinates, the corresponding bilateral Markov chain can be written as $(x_{n_{k}}, r_{n_{k}})$ with $r_{n_{k}}=r p^{-1}\ p_{n_{k}} |S'_{n_{k}} u|$.  We denote by $\Omega_{+}^{\#}$ the subset of $\Omega^{\#}$ defined by the conditions $x\in X_{+}$, $x_{n}\in X_{+}$ infinitely often for $n>0$ and $n<0$,  by $\kappa^{\#}_{+}$ the normalized restriction of $\kappa^{\#}$ to $\Omega_{+}^{\#}$ and by ${^{a}\widehat{\theta}_{+}}$ the induced shift. Also let $\Omega_{0}^{\#}$ be the subset of $\Omega^{\#}_{+}$ defined by the conditions $x\in X_{+}$, $\displaystyle\sup_{k>0} (p_{n_{-k}}\ p^{-1}\  |S'_{n_{-k}} u|)<1$. From Lemma 5.11, we know that $\kappa_{+}^{\#}$-\textrm{a.e.} ,
$\lim_{n\rightarrow \infty} |S'_{-n} u|\ p_{-n}\ p^{-1}=0$, hence $\lim_{k\rightarrow \infty} |S'_{n_{-k}} u|\ p_{n_{-k}}\ p^{-1} =0$ and $|S'_{n_{-k}} u| p_{n_{-k}}  p^{-1}>0$.

  Then the index $-\nu_{0}\leq0$ of the strict last maximum of the sequence $p_{n_{-k}}\ p^{-1}\ |S'_{n_{-k}} u|=V_{-k}$  $(k\geq0)$ is finite $\kappa_{+}^{\#}$-\textrm{a.e.} . We have $\Omega^{\#}_{0}=\{\nu_{0}=0\}$ and, using Lemma 5.11,$\displaystyle\mathop{\lim}_{k\rightarrow \infty} V_{-k}=0$. It follows,using the $(^{a}\widehat{\theta}_{+}) $-invariance of $\kappa_{+}^{\#}$,
\[1=\displaystyle\mathop{\Sigma}^{ \infty}_{k=0} \kappa^{\#}_{+} \{\nu_{0}=k\}\leq \displaystyle\mathop{\Sigma}^{ \infty}_{k=0} \kappa^{\#}_{+}\{V_{-k}>\displaystyle\sup_{j<-k} V_{j}\}=\displaystyle\mathop{\Sigma}^{ \infty}_{j=0} \kappa^{\#}_{+} \{V_{0}>\displaystyle\sup_{j<0} V_{j}\}=\displaystyle\mathop{\Sigma}^{ \infty}_{j=0} q,\]
  where $q=\kappa^{\#}_{+}(\Omega_{0}^{\#})$, hence $\kappa^{\#}(\Omega_{0}^{\#})>0$.
Furthermore if $\omega^{\#}\in \Omega_{0}^{\#}$,using the positivity of $V_{k}$ and an observation of  \cite{36},we see that  $\tau (\omega^{\#})$ is the first return time of $({^{a}\widehat{\theta}})^n (\omega^{\#})$ to $\Omega^{\#}_{0}$; hence ${^{a}\widehat{\theta}^{\tau}}$ is the transformation on $\Omega^{\#}_{0}$ induced by ${^{a}\widehat{\theta}}$ or ${^{a}\widehat{\theta}_{+}}$ on $\Omega_{0}^{\#}$ and $\tau_{n} (\omega^{\#})$ is the sequence of return times   to $\Omega_{0}^{\#}$. This allow us to proceed as in (\cite{36} Lemma 2) with the $\mathbb R^{*}$-valued $\mathbb Z$-cocycle $|S'_{n} u|\ p_{n}\ p^{-1}$. Since $\kappa^{\#}$ is ${^{a}\widehat{\theta}}$-invariant  and $\kappa^{\#} (\Omega_{0}^{\#})>0$  we can apply Kac's recurrence theorem to $\Omega_{0}^{\#}$, ${^{a}\widehat{\theta}^{\tau}}$ and $\Omega^{\#}$ (see \cite{50}), hence  the normalized restriction $\kappa_{0}^{\#}$ of $\kappa^{\#}$ to $\Omega_{0}^{\#}$  is ${^{a}\widehat{\theta}^{\tau}}$-ergodic and stationary, the return time $\tau$ has finite expectation $\mathbb E_{0}^{\alpha} (\tau)$ and $\displaystyle\mathop{\lim}_{n\rightarrow \infty} \frac{\tau_{n}}{n}=\mathbb E_{0}^{\alpha} (\tau),\ \kappa_{0}^{\#}$-\textrm{a.e.} . Since $\kappa_{0}^{\#}$ is absolutely continuous with respect to $\kappa^{\#}$, Theorem 3.10 gives
\[\lim_{n\rightarrow \infty} \frac{1}{n} \log |S'_{\tau_{n}} u|=(\displaystyle\mathop{\lim}_{n\rightarrow \infty} \frac{1}{\tau_{n}} \log |S'_{\tau_{n}}u|) (\displaystyle\mathop{\lim}_{n\rightarrow \infty} \frac{\tau_{n}}{n})=\mathbb E_{0}^{\alpha} (\tau) L_{\mu} (\alpha)=\gamma_{\tau}^{\alpha}>0, \kappa_{0}^{\#}-\textrm{a.e.}  .\]

  Using Birkhoff's theorem for the non-negative increments of $W_{\tau_{n}}^{a}=\log (\frac{p_{\tau_{n}}}{p}\ |S'_{\tau_{n}}u|)$ and ${^{a}\widehat{\theta}^{\tau}}$, we get the $\kappa_{0}^{\#}$-\textrm{a.e.} convergence of $\frac{1}{n}$ $W_{\tau_{n}}^{a}$. Since $\kappa_{0}^{\#}$ is ${^{a}\widehat{\theta}^{\tau}}$-invariant,the sequence $\frac{1}{n} \log \frac{p_{\tau_{n}}}{p}$ converges to zero in $\kappa_{0}^{\#}$-measure, hence using the $\kappa_{0}^{\#}-\textrm{a.e.}$ convergence of $\frac{1}{n} \log |S'_{\tau_{n}}u|$, we get the $\kappa_{0}^{\#}$-\textrm{a.e.} convergence of $\frac{1}{n} W_{\tau_{n}}^{a}$ to $\gamma_{\tau}^{\alpha}$. In particular $\mathbb E^{\alpha}_{0}(W_{\tau}^{a})= \gamma_{\tau}^{\alpha}\in ]0,\infty[$.

  In order to relate $\kappa_{0}^{\#}$ and the kernel ${^{*}\widehat{Q}^{\alpha,\tau}}$ we consider the Markov kernel adjoint to ${^{*}\widehat{Q}^{\alpha}_{+}} (x,\cdot)$ with respect to $\kappa_{+}$, and we denote by ${^{*}\widehat{\mathbb Q}^{\alpha}_{-}} \otimes \delta_{x}$ the corresponding Markov measure on $H^{\mathbb Z_{-}} \times X_{+}$ with $\mathbb Z_{-}=-\mathbb N \cup \{0\}$. Also  we write ${^{*}\widehat{\mathbb Q}^{\alpha}_{x}}=\delta_{x}\otimes {^{*}\widehat{\mathbb Q}^{\alpha}_{+}}$ where ${^{*}\widehat{\mathbb Q}^{\alpha}_{+}}$ is supported on $H^{\mathbb N}$ and $\kappa^{\#}_{+}=\int {^{*}\widehat{\mathbb Q}^{\alpha}_{-}} \otimes \delta_{x} \otimes {^{*}\widehat{\mathbb Q}^{\alpha}_{+}} d\kappa_{+} (x)$, in particular \
   \[\kappa^{\#}(\Omega^{\#}_{0}) \kappa^{\#}_{0}=\int (1_{\Omega_{0}^{\#}})  {^{*}\widehat{\mathbb Q}^{\alpha}_{-}} \otimes \delta_{x} \otimes {^{*}\widehat{\mathbb Q}^{\alpha}_{+}} d\kappa_{+} (x)\]
with $\Omega_{0}^{\#}=\Omega_{0}^{-} \times H^{\mathbb N}$ and $\Omega_{0}^{-} \subset H^{\mathbb Z_{-}}\times X_{+}$. We denote by $\kappa^{\tau}$ the projection of $\kappa_{0}^{\#}$ on $X_{+}$, hence $\kappa^{\tau}$ has density  $u(x)$ given by $\kappa^{\#} (\Omega^{\#}_{0}) u(x)=({^{*}\widehat{\mathbb Q}^{\alpha}_{-}} \otimes \delta_{x}) (\Omega^{-}_{0})$ with respect to $\kappa_{+}$. It follows that the projection of $\kappa_{0}^{\#}$ on $X_{+}\times \widehat{\Omega}$ can be expressed as
\[\int u(x) \delta_{x}\otimes {^{*}\widehat{\mathbb Q}^{\alpha}_{+}} d\kappa_{+}(x)=\int \delta_{x} \otimes {^{*}\widehat{\mathbb Q}^{\alpha}_{+}} d\kappa^{\tau}(x)={^{*}\widehat{\mathbb Q}^{\alpha}_{\kappa^{\tau}}}.\]

  Since $\kappa_{0}^{\#}$ is invariant and ergodic with respect to the bilateral shift ${^{a}\widehat{\theta}^{\tau}}$, the same is valid for ${^{*}\widehat{\mathbb Q}^{\alpha}_{\kappa^{\tau}}}$  with respect to the associated unilateral shift ${^{a}\widehat{\theta}^{\tau}}$. Since the kernel $x\rightarrow {^{*}\widehat{\mathbb Q}^{\alpha}_{x}}$ commutes with ${^{a}\widehat{\theta}^{\tau}}$ and ${^{*}\widehat{\mathbb Q}^{\alpha,\tau}}$, the ${^{*}\widehat{Q}^{\alpha,\tau}}$-invariance and ergodicity of $\kappa^{\tau}$ follows. Also we have $\mathbb E^{\alpha}_{0}(\tau)=\int \mathbb E_{x}^{\alpha} (\tau) d\kappa^{\tau} (x)$ and the above convergences are valid ${^{*}\widehat{\mathbb Q}^{\alpha}_{\kappa^{\tau}}}$-\textrm{a.e.} .
\end{proof}

\begin{rema}
  If $S=X_{+}$, $P={^{*}\widehat{P}}_{\alpha,+}$, $\pi=\kappa_{+}$,and in the corresponding measurable setting, the measure $\kappa^{\tau}$ is closely connected with the measure $\overline{\chi}$ of Theorem \ref{the:4.4} .  The measure $\kappa^{\tau}$ can be caracterized as the unique ${^{*}\widehat{Q}}^{\alpha,\tau}$-stationary measure which is absolutely continuous with respect to $\kappa_{+}$. However, in contrast to \cite{36}, the function $\log |p|$ is not known to be $\kappa$-integrable, but we know that $\displaystyle\mathop{\lim}_{n\rightarrow \infty} |S'_{-n} u| p_{-n}\  p^{-1}=0$.
  \end{rema}

  The following weak renewal theorem for a general fibered Markov chain  will allow us to control potentials of the measurable fibered Markov kernel ${^{*}\widehat{P}}^{\tau}_{\alpha}$  on $X_{+}\times \mathbb R_{+}^{*}$ . We recall some  notation of section 4 as follows.

  Let $(S,\pi)$ be a complete separable metric space, where $\pi$ is a probability measure.  We consider a general measurable Markov chain on $S\times \mathbb R$ with kernel  $P$, we assume that $P$ commutes with the $\mathbb R$-translations and we denote    Lebesgue measure on $\mathbb R$ by $\ell$ . We assume that the measure $\pi \otimes \ell$ is $P$-invariant. Here, in contrast to section 4, our setting is the measurable one;in particular the symbol $\sup$ means essential supremum. We write a path of the corresponding Markov chain as $(x_{n}, V_{n})$ where $x_{n} \in S$ and $V_{n}\in \mathbb R$, we denote by  $^{a}\mathbb P_{x}$ the Markov measure on the paths starting from $x \in S$ and we write  $^{a}\mathbb P_{\pi}=\int {^{a}\mathbb P_{x}} d\pi (x)$, $^{a}\mathbb E_{x}$ for the corresponding expectation symbol. In this context the following weak analogue of the renewal Theorem \ref{the:4.4}  holds.

\begin{prop}
\label{pro:5.14}
  With the above notation and hypotheses, assume that $\psi$ is a compactly supported bounded non negative measurable function on $S \times \mathbb R$,  the potential $U \psi= \displaystyle\mathop{\Sigma}_{0}^{\infty} P^k \psi$ is essentially bounded on $S \times [-c,c]$ for any $c> 0$ and we have for any $\varepsilon>0$
\[\lim_{n\rightarrow \infty} {^{a}\mathbb P_{\pi}}\{ |\frac{V_{n}}{n}-\gamma|>\varepsilon\}=0,\]  with $\gamma>0$.
 Then we have
\[\lim_{t\rightarrow \infty}  \frac{1}{t} \int_{-t}^{0} ds \int_{S} U \psi (x,s) d\pi (x)=\frac{1}{\gamma} (\pi \otimes \ell) (\psi).\]

  Furthermore if $\psi$ is a  non negative measurable function on $S\times \mathbb R$ and $\displaystyle\mathop{\lim}_{t\rightarrow -\infty} U \psi (x,t)=0$, $\pi$-\textrm{a.e.} .\ then $\psi=0$, $\pi\otimes \ell$-\textrm{a.e.} .

  If $\psi$ is a measurable function on $S\times \mathbb R$ which satisfies:
\[|\psi|_{b}=\displaystyle\mathop{\sum}_{\ell=-\infty}^{\ell=\infty} \sup\{|\psi (x,s)|\ ;\ x\in S,\ s\in [\ell, \ell+1[\}<\infty, \]
then the above convergence is also valid.
\end{prop}

\begin{proof}
  We observe that the maximum principle implies $|U\psi|=\displaystyle\sup_{x,t} |U \psi| (x,t)<\infty$, since $U\psi$ is essentially locally bounded. For $\varepsilon>0$, $t>0$ we denote $n_{1}=n_{1} (t)=[\frac{1}{\gamma} \varepsilon t], \ n_{2}=n_{2} (t)=[\frac{1}{\gamma} (1+\varepsilon) t]$ where $[t]$ denotes integer part of $t>0$.We write:$\displaystyle\mathop{\Sigma}_{0}^{\infty} P^k \psi=U \psi$,
 $\displaystyle\mathop{\Sigma}_{0}^{n-1} P^k \psi=U_{n} \psi$, $\displaystyle\mathop{\Sigma}_{n+1}^{\infty} P^k \psi=U^n \psi$, $\displaystyle\mathop{\Sigma}_{n}^{m} P^k \psi=U_{n}^m \psi$,
 $I(t)= \frac{1}{t} \int_{-t}^{0} ds \int_{S} (U \psi) (x,s) d\pi (x)=\displaystyle\mathop{\Sigma}_{1}^{3} I_{k} (t)-I_{4}(t)$ where
\begin{align*}
I_{1}(t)&=\frac{1}{t}\int_{S} d\pi (x) \int_{-\infty}^{\infty} U_{n_{1}}^{n_{2}} \psi (x,s) ds,
I_{2}(t)=\frac{1}{t}\int_{S} d\pi (x) \int_{-t}^{0} U_{n_{1}}\psi (x,s) ds,\\
I_{3}(t)&=\frac{1}{t}\int_{S} d\pi (x) \int_{-t}^{0} U^{n_{2}} \psi (x,s) ds,
I_{4}(t)=\frac{1}{t}\int_{S} d\pi (x) \int_{\mathbb R\setminus  [-t,0]} U_{n_{1}}^{n_{2}} \psi (x,s) ds.
\end{align*}

  We estimate each term $I_{k} (t)$ separately.
  We have, since the measure $\pi \otimes \ell$ is $P$-invariant,
$I_{1}(t)=\frac{n_{2}-n_{1}+1}{t}(\pi \otimes \ell)  (\psi)$, hence $\displaystyle\mathop{\lim}_{t\rightarrow \infty} I_{1} (t)=\frac{1}{\gamma} (\pi \otimes \ell) (\psi)$.

  Furthermore,
\[|I_{4}(t)|\leq \frac{|\psi|}{t} (n_{2}-n_{1}+1) \displaystyle\sup_{n_{1}\leq n\leq n_{2}} \int_{S}  (^{a}\mathbb P_{x} \{V_{n}\leq a\}+{^{a}\mathbb P_{x}} \{V_{n} \geq t-a\}) d\pi (x).\]
  Since $n^{-1} V_{n}$ converges to $\gamma>0$ in probability, the above integral has limit zero, hence $\displaystyle\mathop{\lim}_{t\rightarrow \infty} I_{4}(t)=0$.

  We have also
$|I_{2} (t)| \leq \frac{ \varepsilon }{\gamma} \pi \otimes \ell (\psi).$

  In order to estimate $I_{3} (t)$ we denote, for $n\in \mathbb N$ and  $s>0$, $ \rho_{n}^s=\inf \{k\geq n\ ;\ -a \leq V_{n}-s \leq a\}$, where $\psi$ is supported on $[-a,a]$, and we use the interpretation of $U^n
\psi$ as the expected number of visits to $\psi$ after time $n$ :
\[U^n \psi (x,s) \leq |U \psi| {^{a}\mathbb P_{x}} \{\rho_{n}^s <\infty\}.\]
Taking $n=[\frac{(1+\varepsilon) t}{\gamma}]=n_{2}$ we get
\[I_{3} (t)\leq |U\psi| \int_{S} {^{a}\mathbb P_{x}}\{V_{k}-t \leq a, \textrm{  for some } k\geq [\frac{(1+\varepsilon) t}{\gamma}]\} d\pi (x).\]
Since $ \frac{V_{n}}{n}$ converges to $\gamma>0$ in probability, we get $\displaystyle\mathop{\lim}_{t\rightarrow \infty} I_{3}(t)=0$.

  Since $\varepsilon$ is arbitrary we get finally,$\displaystyle\mathop{\lim}_{t\rightarrow \infty} I(t)=\frac{1}{\gamma} (\pi \otimes \ell) (\psi)$.

  The second conclusion follows by restriction and truncation of $\psi$ on $S\times [-a,a]$.

  For the proof of the last assertion we observe that for any $\ell \in \mathbb Z$,
$\Delta=|U1_{S \times [0,1]}|=|U 1_{S\times [\ell,\ell+1[}|<\infty$. Writing $\psi=\displaystyle\mathop{\Sigma}_{\ell=-\infty}^{\ell=\infty} \psi 1_{S\times [\ell, \ell+1[}$ we get $|U\psi|\leq \Delta |\psi|_{b}$.

  The quantity $\psi\rightarrow |\psi|_{b}$ is a norm on the space  $\mathcal H$ of measurable functions on $S\times \mathbb R$ such that $|\psi|_{b}<\infty$. Since the set of  essentially bounded functions supported on $S\times [-c,c]$ for some $c>0$ is dense in $\mathcal H$ and $\psi\rightarrow (\pi \otimes \ell) (\psi)$ is a continuous functional on $\mathcal H$, the above relation extends by density to any $\psi \in \mathcal H$.
  \end{proof}

\textit{Proof of Proposition 5.9.}
Assume that for some $u\in M$, $C (u)=\displaystyle\mathop{\lim}_{t\rightarrow \infty} t^{\alpha} \widehat{\mathbb P}\{\langle R, u \rangle  >t\}=0$.
For $p>0$, with the notations of Lemma 5.12, this means $\displaystyle\mathop{\lim}_{t\rightarrow \infty} \psi^{\alpha}(v,p)=0$. Using Proposition 5.3 we know that this implies $\displaystyle\mathop{\lim}_{t\rightarrow \infty}\psi^{\alpha}(u,p)=0$ for any $u=t v\in M\ (t>0)$. Also, using Lemma 5.13, we have, since $\Lambda_{a}^{*} (\Sigma)\supset M$,
\[\lim_{n\rightarrow \infty}\ \frac{1}{n} \log (|S'_{\tau_{n}} u|\ \frac{p_{\tau_{n}}}{p})=\gamma_{\tau}^{\alpha} >0, {^{*}\widehat{\mathbb Q}_{\kappa^{\tau}}^{\alpha}}-\textrm{\textrm{a.e}.}\]
Since the canonical Markov measure associated with $\kappa^{\tau}$ and ${^{*}\widehat{Q}^{\alpha,\tau}}$ is a push-forward of ${^{*}\widehat{\mathbb Q}^{\alpha}_{\kappa^{\tau}}}$, this convergence is also valid with respect to this canonical Markov measure.
Then, using Lemma 5.12 and $\psi_{\tau}^{\alpha}\geq 0$, we can apply Proposition 5.14 with $V_{n}=\log (p^{-1} p_{\tau_{n}}) |S'_{\tau_{n}}u|$, $\gamma=\gamma_{\tau}^{\alpha}>0$, $S=X_{+}$, to the measurable Markov kernel $P={^{*}\widehat{P}_{\alpha}^{\tau}}$ on $X_{+}\times \mathbb R_{+}^{*}$, to the potential $\displaystyle\mathop{\Sigma}_{0}^{\infty} ({^{*}\widehat{P}_{\alpha}^{\tau}})^k \psi_{\tau}^{\alpha}$ of the non negative function $\psi^{\alpha}_{\tau}\leq ({^{*}e^{\alpha}}\otimes h^{\alpha})^{-1}$ and to the ${^{*}\widehat{Q}^{\alpha,\tau}}$-stationary measure $\pi=\kappa^{\tau}$; we get $\psi_{\tau}^{\alpha}=0$, $\kappa^{\tau}\otimes \ell$-\textrm{a.e}., hence
\[\widehat{\mathbb P}\{t < p^{-1} \langle R, u \rangle   \leq t+ p^{-1}\langle R_{\tau},u\rangle, \tau<\infty\}=0.\]

  Since $p^{-1}\langle R_{\tau}, u \rangle  >0$, this gives $p^{-1}\langle R, u \rangle  \leq 0$ $\kappa^{\tau} \otimes \widehat{\mathbb P}-\textrm{a.e.}$ on $\{\tau<\infty\}$, in particular for some $(u,p)\in X_{+}$,  we have  $p^{-1} \langle R, u \rangle   \leq 0$  i.e. $\langle R, u \rangle   \leq 0$, $\widehat{\mathbb P}$-\textrm{a.e.} \  on $\{\tau<\infty\}$. But, since $\Lambda^{*}(\Sigma)\supset M$,  for any $u\in M$ the set $\{\langle R, u \rangle  >0\ ; \ \tau<\infty\}=\{\langle R, u \rangle  >0\}$ is not $\widehat{\mathbb P}$-negligible, hence the required contradiction. Since, using Proposition 5.3, we have $C(u)=C(\sigma^{\alpha}\otimes \ell^{\alpha}) (H_{u}^{+})$ it follows $C>0$.
  \hfill$\square$
\begin{rema}
  The proof of Proposition 5.9 given above uses the $\mathbb R^{*}$-valued multiplicative cocycle $W_{n}$. The interpretation of the ladder index $\tau$ as a first return time to a subset of $\Omega^{\#}$ depends on the reduction of $W_{n}$ to a positively valued cocycle; hence the use of the seemingly artificial inducing procedure on $X_{+}\times \mathbb R^*$, as was done above.
  \end{rema}

\subsection{A Choquet-Deny type property}

  Here, as in section 4, we consider a fibered Markov chain on $S\times \mathbb R$, but we reinforce the hypothesis on the Markov kernel $P$, using spectral gap properties instead of equicontinuity properties. Hence $S$ is a compact metric space, $P$ commutes with $\mathbb R$-translations and acts continuously on $C_{b} (S\times \mathbb R)$. We define for $t\in \mathbb R$, the Fourier operator $P^{it}$ on $C (S)$ by
\[P^{it} \varphi(x)=P(\varphi \otimes e^{it\cdot}) (x,0).\]

  For $t=0$ the operator  $P^{it}=P^0$ is equal to $\overline{P}$, the factor operator on $S$ defined by $P$. We assume that for each $t\in \mathbb R$, $P^{it}$ preserves the space $H_{\varepsilon}(S)$ of $\varepsilon$-H\"older functions and is a bounded operator therein. Moreover we assume that $P^{\cdot it} (t\in \mathbb R)$, and $P$ satisfies the following condition D.

  A weaker Choquet-Deny type result under a similar condition was shown in (\cite{30}, Proposition 3.5). The stronger form given here allow us to deal with polynomially bounded $P$-harmonic measures and is needed for the proof of the homogeneity at infinity of $\rho$ in Theorem 5.2.Condition D is as follows
\begin{enumerate}
\item
For any $t\in \mathbb R$, one can find $n_{0}\in \mathbb N$, $\rho(t) \in [0,1[$ and $C(t)>0$ for which
\[[(P^{it})^{n_{0}}\varphi]_{\varepsilon} \leq \rho (t) [\varphi]_{\varepsilon}+C(t) |\varphi|.\]
\item
For any $t\in \mathbb R$, the equation $P^{it}\varphi=e^{i\theta} \varphi$, $\varphi \in H_{\varepsilon} (S)$, $\varphi \neq 0$ has only the trivial solution $e^{i\theta}=1$, $t=0$, $\varphi=$constant.
\item
 For some $\tau>1\:\sup \{\int |a|^{\tau} P(x,0), d(y,a)\ ;\ x\in S\}<\infty$.
\end{enumerate}
  Conditions 1,2 above imply that $\overline{P}$ has a unique stationary measure $\pi$ and the spectrum of $\overline{P}$ in $H_{\varepsilon} (S)$  is of the form $\{1\} \cup\Delta$ where $\Delta$ is a compact subset of the open unit disk (see \cite{33} ). They imply also that, for any $t\neq 0$, the spectral radius of $P^{it}$ is less than one.

  With the notations of section 4, condition 3  above will allow us to estimate $^{a}\mathbb E (|V_{n}|^p)$ for $p\leq \delta$ and to show the continuity of $t\rightarrow |P^{it}|$.

  The following is a simple consequence of conditions 1,2,3 above.

\begin{lemm}
\label{lem:5.15}
With the above notation, let $I\subset \mathbb R$ be a compact subset of $\mathbb R\setminus \{0\}$. Then there exists $D>0$ and $\sigma \in [0,1[$ such that for any $n\in \mathbb N$,
$\displaystyle\sup_{t\in I}|(P^{it})^n|\leq D\sigma^n$.
\end{lemm}

\begin{proof}
  Conditions 1 and 2 for $P^{it} (t\neq 0)$ imply that the spectral radius $r_{t}$ of $P^{it}$ satisfies $r_{t}<1$. Hence there exists $C_{t}>0$ such that for any $n\in \mathbb N$,$|(P^{it})^n |\leq C_{t} (\frac{1}{2}+\frac{r_{t}}{2})^n$. On the other hand $t\rightarrow |P^{it}|$ is continuous as the following calculation shows. For $a, t, t' \in \mathbb R, \delta' \in [0,1]$, we have $|e^{i a t}-e^{i a t'}| \leq 2|a|^{\delta'} \ |t-t'|^{\delta'}$ hence we have
\begin{align*}
|P^{it} \varphi (x)-P^{it'} \varphi (x)|&\leq 2 |\varphi|\  |t-t'|^{\delta'} \int |a|^{\delta'} P((x,0), d(y,a)),\\
 |P^{it}-P^{it'}|&\leq 2 M_{\delta'} |t-t'|^{\delta'}.
 \end{align*}

  For each $t\in I$ we fix $n_{t}\in \mathbb N$ such  that $|(P^{it})^{n_{t}}|\leq \frac{1}{3}$. Then the above continuity of $P^{it}$, hence of $(P^{it})^{n_{t}}$, gives that for $t'$ sufficiently chose to $t,|(P^{it'})|\leq \frac{1}{2}$. Using compactness of I we find $n_{1},\cdots n_{k}$ such that  one of the inequalities $|(P^{it})^{n_{j}}| \leq \frac{1}{2}\ \  (1\leq j \leq k)$ is valid at any given point of $I$. Then, since $|P^{it}|\leq 1$ with $n_{0}=n_{1}\cdots n_{k}$ we get $|P^{it})^{n_{0}}|\leq \frac{1}{2}$. Using Euclidean division of $n$ by $n_{0}$, we get the required inequality. \qquad $\square$
\vskip 2mm
  We are interested  in the action of $P^n$ on functions on $S\times \mathbb R$ which are of the form $u\otimes f$ where $u\in C (S)$, $f\in \mathbb L^1  (\mathbb R)$ and we are interested  also in $P$-harmonic Radon measures which satisfy boundedness conditions. In some proofs, since $H_{\varepsilon} (S)$ is dense in $C(S)$, it will be convenient to assume $u\in H_{\varepsilon} (S)$.
  \end{proof}

\begin{defi}{}
\label{def:5.16}
  We say  that the Radon measure $\theta$ on $S\times \mathbb R$ is translation-bounded if for any compact subset $K$ of $S\times \mathbb R$, any $a\in \mathbb R$, there exists $C(K)>0$ such that
$|\theta(a+K)|\leq C(K)$
  where $a+K$ is the compact subset of $S\times \mathbb R$ obtained from $K$ using translation by $a\in \mathbb R$.
  \end{defi}

  We are led to consider a positive function $\omega$ on $\mathbb R^d$ which satisfies $\omega(x+y)\leq \omega (x) \omega(y)$. For example, if $p\geq 0$, such a function $\omega_{p}$ is defined by $\omega_{p} (a)=(1+|a|)^p$.
 We denote $\mathbb L^1_{\omega} (\mathbb R)=\{ f\in \mathbb L^1 (\mathbb R)$; $\omega f\in \mathbb L^1 (\mathbb R)\}$ and we observe that $f\rightarrow \|\omega f\|_{1}=\|f\|_{1,\omega}$ is a norm under which $\mathbb L^1_{\omega} (\mathbb R)$ is a Banach algebra. The dual space of $\mathbb L^1_{\omega} (\mathbb R)$ is the space $\mathbb L^{\infty}_{\omega} (\mathbb R)$ of measurable functions $g$ such that $g \omega^{-1}\in \mathbb L^{\infty} (\mathbb R)$ and the duality is given by $\langle g,f\rangle=\int g(a) f (a) da$. The Fourier transform $\widehat f$ of $f\in \mathbb L^1_{\omega} (\mathbb R)$ is well defined by $\widehat f (t)=\int f(a) e^{ita}da$. We denote by  $J^c$ the ideal of $\mathbb L^1(\mathbb R)$ which consists of functions $f\in \mathbb L^1 (\mathbb R)$ such that $\widehat f$ has a compact support not containing 0 and we  write $J^c_{\omega}= \mathbb L^1_{\omega} (\mathbb R) \cap J^c$. Also we denote by $\mathbb L^1_{0}(\mathbb R)$ the ideal of $\mathbb L^1 (\mathbb R)$ defined by the condition $\widehat f(0)=0$. It is well known that $J^c$ is dense in $\mathbb L_{0}^1 (\mathbb R)$, hence for $\omega=\omega_{p}$ the ideal  $J^c_{\omega^2}$ is dense in $\mathbb L_{\omega}^1 (\mathbb R)$  (see \cite{32} p. 187).

\begin{theo}
\label{thm:5.17}
 With the above notation, assume that the family $P^{it} (t\in \mathbb R)$ satisfies conditions D and let $\omega=\omega_{p}$ with $p<\delta$. Then for any $f\in \mathbb L_{\omega}^1 (\mathbb R) \cap J^c$, $u\in C(S)$, we have the convergence
$$\displaystyle\mathop{\lim}_{n\rightarrow \infty} \displaystyle\sup_{x\in S}\  \|P^n (u\otimes f) (x,\cdot)\|_{1,\omega}=0.$$
  If $\theta$ is a $P$-harmonic Radon measure which is translation-bounded, then $\theta$ is proportional to $\pi\otimes \ell$. In particular $\pi\otimes \ell$ is a minimal $P$-harmonic Radon measure.
  \end{theo}

  The proof follows from the above considerations and the following lemmas.
\begin{lemm}
\label{lem:5.18}
  Assume $\pi_{n}$ is a sequence of bounded measures on $\mathbb R$, $\omega$ is a positive Borel function on $\mathbb R$ such that for any $x,y\in \mathbb R$, $\omega(x+y)\leq \omega (x) \omega(y)$ and assume that the total variation measures  $|\pi_{n}|$ of $\pi_{n}$  satisfy $\sup \{|\pi_{n}| (\omega)\  ; \ n\in \mathbb N\}<\infty$. Let $f\in \mathbb L_{\omega}^1 (\mathbb R) \cap \mathbb L^2 (\mathbb R)$ and assume $A_{n}, B_{n}$ are sequences of Borel subsets of $\mathbb R$ such that, with $A'_{n}=\mathbb R \setminus A_{n} \ ,\ B'_{n}=\mathbb R\setminus B_{n}$,
\begin{enumerate}
\item
 $\displaystyle\mathop{\lim}_{n\rightarrow \infty} |\pi_{n}| (\omega) \|f 1_{B'_{n}}\|_{1,\omega}=0$,
\item
$\displaystyle\mathop{\lim}_{n\rightarrow \infty} |\pi_{n}| (\omega 1_{A'_{n}})=0$,
\item
 $\displaystyle\mathop{\lim}_{n\rightarrow \infty} \|\pi_{n} * f\|_{2} \ \|\omega^2 1_{A_{n}+B_{n}}\|_{1}^{1/2}=0$.
\end{enumerate}
  Then we have $\displaystyle\mathop{\lim}_{n\rightarrow \infty} \|\pi_{n} * f\|_{1,\omega}=0$. Furthermore, if the measures $\pi_{n}$ depend of a parameter $\lambda$ and if the convergences in 1--3 are uniform in $\lambda$, then the convergence of $ \|\pi_{n}* f\|_{1,\omega}$ is also uniform in $\lambda$.
\end{lemm}

\begin{proof}
Let $\eta, \eta'$ be two bounded measures on $\mathbb R$ and let $A, B$ be Borel subsets with complements $A', B'$ in $\mathbb R$. Observe that, since $0\leq \omega (x+y) \leq \omega (x) \omega(y)$,
\[\omega(x+y)\leq \omega(x+y)1_{A+ B} (x+y)+(\omega 1_{A}) (x) (\omega 1_{B'}) (y)+ (\omega 1_{A'}) (x) \omega (y).\]
It follows
\[|(\eta *\eta')| (\omega) \leq |\eta*\eta'| (\omega 1_{A+B})+|\eta| (\omega)|\ \  |\eta'| (\omega 1_{B'})+|\eta| (\omega 1_{A'}) |\eta'| (\omega).\]
Then we take $\eta=\pi_{n}$, $\eta'=f (a) da$, $A=A_{n}$, $B=B_{n}$ and we get
\[\|\pi_{n}*f\|_{1,\omega} \leq \int |\pi_{n}*f| (a) |\omega 1_{A_{n}+B_{n}}| (a) da+|\pi_{n}|(\omega 1_{A_{n}}) \|f 1_{B'_{n}}\|_{1,\omega}+|\pi_{n}|(\omega 1_{A'_{n}}) \|f\|_{1,\omega}.\]

  Conditions 1, 2 imply that the two last terms in the above inequality have limits zero. Using condition 3 and Schwarz inequality we see that the first term has also limit zero. If $\pi_{n}$ depends of a parameter $\lambda$, the uniformity of the convergence of $\|\pi_{n}*f\|_{1,\omega}$ follows directly from the bound for $\|\pi_{n}*f\|_{1,\omega}$ given above. \end{proof}
  The following lemma is an easy consequence of condition 3 on the Markov kernel $P$ and of H\"older inequality.

\begin{lemm}
\label{lem:5.19}
  For any $p\in [1,\delta]$,
there exists $C_{p}>0$ such that  ${\displaystyle\sup_{x}} ^{a}\mathbb E_{x} ((|V_{n}|)^p)\leq C_{p} n^p$.
In particular, for any $L>0$, ${\displaystyle\sup_{x,n}} ^{a}\mathbb P_{x} \{|V_{n}|>nL\}\leq \frac{C_{p}}{L^p}$.
  \end{lemm}

  We leave to the reader the proof of the well known first inequality. The second one follows from Markov's inequality.

  For a Radon measure $\theta$ on $S\times \mathbb R$ and $b\in \mathbb R$,
  we denote by $\theta *\delta_{b}$ the Radon measure defined by $(\theta *\delta_{b}) (\varphi)=\int \varphi (x,a+b) d\theta (x,a)$; for $\varphi\in C_{c}(S\times \mathbb R)$,
  $\theta$  translation-bounded we write
  $|\theta|_{\varphi}=\displaystyle\sup \{ |\theta* \delta_{b} (\varphi)|; \  b\in \mathbb R\}$. For such measures and any bounded measure $r$ on $\mathbb R$, $\theta*r$ is well defined by
  $(\theta* r) (\varphi)=\int (\theta*\delta_{b}) (\varphi) dr(b)$ and we have $|\theta*r|_{\varphi}\leq|r| \theta_{\varphi}$ where $|r|$ is the total variation of $r$. In particular, $f\in \mathbb L^{1}(\mathbb R)$ can be identified with the measure $r_{f}=f(a) da$, we can define $\theta*f=\theta*r_{f}$ and if $f_{n}\in \mathbb L^1 (\mathbb R)$ converges in $\mathbb L^1$-norm to $f\in \mathbb L^1 (\mathbb R)$, then $\theta*f_{n}$ converges to $\theta*f$ in the vague topology. On the other hand, if $r$ has compact support and $\theta$ is a Radon measure on $\mathbb R$, $\theta*r$  is well defined as a Radon measure.

\begin{lemm}
\label{lem:5.20}
  With the above notation, assume that $\theta$ is a translation-bounded non negative Radon measure on $S\times \mathbb R$. Let $r$ be a non negative continuous function on $\mathbb R$ with compact support containing 0. Then for $p>1$, there exists a non negative bounded measure $\overline{\theta}$ on $S$ such that $\theta*r\leq (1_{S} \otimes \omega_{p}) (\overline{\theta} \otimes \ell)$.
\end{lemm}
\begin{proof}
 For simplicity of notation, assume $r>0$ on $[0,1]$. We denote by $\theta_{k} *\delta_{k}$ the restriction of $\theta$ to $S\times [k,k+1[$ $(k\in \mathbb Z)$ and we write $\theta=\displaystyle\mathop{\Sigma}_{k\in \mathbb Z} \theta_{k} * \delta_{k}$ with $\supp \ \theta_{k} \subset S\times [0,1]$.  We observe that, since $\theta$ is translation-bounded, the mass of $\theta_{k}$ is bounded for $k\in \mathbb Z$, hence $\displaystyle\mathop{\theta}^{.}=\displaystyle\mathop{\Sigma}_{k\in \mathbb Z} (1+|k|)^{-p} \theta_{k}$ is a bounded measure supported on $S\times [0,1]$. We have clearly
\[\theta_{k} \leq (1+|k|)^p \displaystyle\mathop{\theta}^{.},
\theta\leq \displaystyle\mathop{\theta}^{.}*\displaystyle\mathop{\Sigma}_{k\in \mathbb Z} (1+|k|)^{p} \delta_{k},
\theta*r\leq \displaystyle\mathop{\theta}^{.}*(r*\displaystyle\mathop{\Sigma}_{k\in \mathbb Z} (1+|k|)^{p} \delta_{k}).\]

  But, by definition of $\omega_{p}$ and since $\supp r$ is compact, we have
$r*\displaystyle\mathop{\Sigma}_{k\in \mathbb Z} (1+|k|)^{p} \delta_{k}\leq c \omega_{p}$ for some $c>0$ and it follows, $\theta*r\leq c\displaystyle\mathop{\theta}^{.} * \omega_{p}$. We desintegrate the bounded measure $\displaystyle\mathop{\theta}^{.}$ as $\displaystyle\mathop{\theta}^{.}=\int \delta_{x} \otimes \displaystyle\mathop{{\theta^x}}^{.} d\widetilde{\theta} (x)$ where $\widetilde{\theta}$ is the projection of $\displaystyle\mathop{\theta}^{.}$ on $S$ and ${\displaystyle\mathop{\theta^x}^{.}}$ is a probability measure. Hence $\displaystyle\mathop{\theta}^{.}*\omega_{p}=\int \delta_{x}\otimes (\displaystyle\mathop{{\theta^x}}^{.} *\omega_{p}) d\widetilde{\theta} (x)$. But, since $\displaystyle\mathop{{\theta^x}}^{.}$ is
supported on
$[0,1]$, we have $\displaystyle\mathop{{\theta^x}}^{.} *\omega_{p} \leq 2^p \displaystyle\mathop{{\theta^x}}^{.} ([0,1]) \omega_{p}$. Hence $\displaystyle\mathop{\theta}^{.}\leq 2^p(1_{S} \otimes \omega_{p}) (\widetilde{\theta}\otimes \ell)$ and finally $\theta*r\leq 2^p c(1_{S} \otimes \omega_{p}) (\widetilde{\theta}\otimes \ell)=(1_{S} \otimes \omega_{p}) (\overline{\theta} \otimes \ell)$ with $\overline{\theta}=2^p c \widetilde{\theta}$.
\end{proof}

\textit{Proof of  Theorem 5.17.}
We fix $p\in [1,\delta[$, $\omega=\omega_{p} , u\in H_{\varepsilon} (S)$, $u\geq 0$ and for $x\in S$, we define the positive measure $\pi_{n}^x$ on $\mathbb R$ by $\pi_{n}^x (\varphi)=P^n (u\otimes \varphi) (x,0)$ where $\varphi$ is a non negative Borel function on $\mathbb R$. We observe that $\pi_{n}^x (1)\leq |u|$, and  for $f\in C_{b} (\mathbb R) \cap \mathbb L^1 (\mathbb R)$, we have
\[P^n (u\otimes f) (x,a)=\pi_{n}^x (f*\delta_{a})=(\pi_{n}^x * f^{*})(a),\]
where $f^{*} (a)=\overline{f(-a)}$. It follows $\|P^n (u\otimes f) (x,\cdot)\|_{1}\leq |u| \  \|f\|_{1}$ and also for $f\in \mathbb L^1_{\omega} (\mathbb R)$, $\|P^n (u\otimes f) (x,\cdot)\|_{1,\omega}\leq \|f\|_{1,\omega} \sup\{ \pi_{n}^x (\omega)\  ;\ x\in S\}$. From condition 3, since $p<\delta$ we have
\[pi_{n}^x(\omega)\leq |u| \ ^{a}\mathbb E_{x} (1+|V_{n}|)^p\leq 2^p C_{p}n^p\]
where we have used Lemma 5.19  in order to bound $^{a}\mathbb E_{x} (1+|V_{n}|)^p$.

  We fix $\delta>1$, we denote $B_{n}=\{a \in \mathbb R\ ;\ |a| \leq n^{1+\delta}\}$, $A_{n}=\{a\in \mathbb R\ ;\ |a|\leq n^{c+1}\}$ with $c>0$ to be defined later and we verify the  conditions 1--3 of Lemma 5.18 for $\pi_{n}=\pi_{n}^x$, uniformly in $x\in S$ for $f\in \mathbb L_{\omega^2}^{1}(\mathbb R)$.
Since  $f\omega^2 \in \mathbb L^1 (\mathbb R)$, Markov's inequality gives $\|f 1_{B'_{n}}\|_{1,\omega} \leq \|f\|_{1,\omega^2} n^{-p (1+\delta)}$. Then, using the bound of $\pi_{n}^x (\omega)$ given above,
\[\pi_{n}^x (\omega) \|f 1_{B'_{n}} \|_{1,\omega} \leq  2^p C_{p} n^{-p\delta} \|f\|_{1,\omega^2}.\]
Hence condition 1 of lemma 5.18 is satisfied.

  We write $\pi_{n}^x (\omega 1_{A'_{n}}) \leq |u| ^{a}\mathbb E_{x} (\omega (V_{n}) 1_{A'_{n}} (V_{n}))$ and use H\"older inequality for $p'>1, \frac{1}{q'}=1-\frac{1}{p'}$ and $p'p<\delta$:
\[\pi_{n}^x (\omega 1_{A'_{n}})\leq C' 2^p n^p \mathbb P\{|V_{n}| \geq n^{c+1}\})^{1/q'}\leq C''\  n^{p-cp/q'} ,\]
 where we have used the fact that ${\displaystyle\sup_{x}} ^{a}\mathbb E_{x} (|U_{1}|^{pp'})<\infty$ for $pp'<\delta$ and Lemma 5.19. If we take $c>q'=\frac{p'}{p'-1}$, we see that $\displaystyle\mathop{\lim}_{n\rightarrow \infty} \pi_{n}^x$ $(\omega 1_{A'_{n}})=0$ uniformly, hence condition 2 is satisfied.

  In order to verify condition 3 we observe that $\widehat{\pi}_{n}^x (t)=(P^{it})^n u(x)$. For $f\in J^c\subset \mathbb L^2 (\mathbb R)$
 we denote $Y=\supp \widehat{f}\subset \mathbb R\setminus \{0\}$ and we know from Lemma 5.17, that there exist $D>0$ and $\sigma \in [0,1[$ such that for any $t\in Y$, $n\in \mathbb N$: $|(P^{it})^n|\leq D \sigma^n$.
From Plancherel formula, we get
\[\|\pi_{n}^x * f\|_{2}= (\int |\widehat{\pi}_{n}^x (t)|^2 |\widehat{f} (t)|^2 dt)^{1/2} \leq \|f\|_{2} |u| \displaystyle\sup_{t\in Y} |(P^{it})^n|\leq D|u|\sigma^n \|f\|_{2}.\]

  On the other hand, $\|\omega^2 1_{A_{n}+B_{n}}\|_{1}$ is bounded by a polynomial in $n$. Since $\sigma<1$, condition 3 is satisfied. Hence, Lemma 5.18 gives:
  $\displaystyle\mathop{\lim}_{n\rightarrow \infty} \|\pi_{n}^x * f\|_{1,\omega}=0$ uniformly. By density, the same relation is valid for all $f\in \mathbb L_{\omega}^1 (\mathbb R)$.

  Now, let us choose $p\in ]1,\delta[$ and $\omega=\omega_{p}$. Since $\theta$ is translation-bounded we can assume $\theta$ to be non negative and translation-bounded. Then Lemma 5.20 gives for any $r$ as in the lemma,$\theta*r\leq (1_{S} \otimes \omega) (\overline{\theta} \otimes \ell)$. Taking $r=r_{n}$ as an approximate identity we have $\theta= \displaystyle\mathop{\lim}_{n\rightarrow \infty} \theta * r_{n}$ in the weak sense. Hence we can assume $\theta\leq (1_{S} \otimes \omega) (\overline{\theta}\otimes \ell)$ where $\overline{\theta}$ is a bounded measure on $S$. Let $f\in \mathbb L_{\omega}^1 (\mathbb R) \cap J^c(\mathbb R)$, $u\in H_{\varepsilon} (S)$ be as above,  hence $f,u$ satisfy  for every $n\in \mathbb N$ the following relations
\begin{align*}
\theta (u\otimes f)&=(P^n \theta) (u\otimes f)=\theta (P^n (u\otimes f))=\int P^n (u\otimes f) (x,a) d\theta (x,a),\\
|\theta (u\otimes f)|&\leq \int d\overline{\theta} (x) \int |P^n (u \otimes f) (x,a)| \omega (a) da\leq \|\overline{\theta}\| \displaystyle\sup_{x} \|P^n (u\otimes f) (x,\cdot)\|_{1,\omega}.
\end{align*}

  From the first part of the proof we get  $\theta (u\otimes f)=0$ for any $u\in H_{\varepsilon} (S)$, $f\in J^c_{\omega}=\mathbb L^1_{\omega} (\mathbb R) \cap J^c$. This relation remains valid for $f$ in the ideal $I_{\omega}^c$ of $\mathbb L_{\omega}^1 (\mathbb R)$ generated by $J_{\omega}^c$. Using regularisation on Fourier transforms we see that the closure in $\mathbb L^1 (\mathbb R)$ of $J_{\omega}^c$ contains $J^c$, hence the unique Fourier exponential which vanishes on $I_{\omega}^c$ is 1. Then using classical Fourier Analysis (see \cite{32} p.187) we get that $I_{\omega}^c$ is dense  in $\mathbb L_{0}^1 (\mathbb R)$. As observed above, since $\theta$ is translation-bounded, this implies $\theta (u\otimes f)=0$ for any $f\in \mathbb L_{0}^1 (\mathbb R)$. Since $H_{\varepsilon} (S)$ is dense in $C (S)$, we get that $\theta$ is invariant by $\mathbb R$-translation. Then we have $\theta=\overline{\theta} \otimes \ell$ where $\overline{\theta}$ is a positive measure on $S$ which satisfies $\overline{P}\  \overline{\theta}=\overline{\theta}$. Using  parts 1,2 of condition D, this implies that $\overline{\theta}$ is proportional to $\pi$, hence $\theta$ is proportional to $\pi\otimes \ell$.

  For the final assertion we observe that, if $\theta$ is a $P$-harmonic positive Radon measure with $\theta\leq c$ $\pi \otimes \ell$, for some $c>0$, then $\theta$ is translation bounded. Hence as above, $\theta$ is proportional to $\pi\otimes \ell$.
  \hfill $\square$

\subsection{Homogeneity at infinity  of  the stationary measure}

  For the proof of Theorem 5.2 we  prepare the following  propositions and lemmas. If $\alpha\notin \mathbb N$, it follows from \cite{5} that Theorem 5.2 is a consequence of Propositions 5.3 and 5.9. If $\alpha \in \mathbb N$, as follows from \cite{51}, the situation is different in general. More precisely, as shown in (\cite{5}, p 706), if $\eta \in M^{1}(V)$ is suitably choosen and $\alpha\in \mathbb N$, convergence of $t^{-\alpha} (t\cdot\eta)$ $(t\rightarrow 0_{+})$ on the sets $H_{v}^{+}$ for every $v\in V$ does not imply vague convergence.

  Here,  we will need to use the Choquet-Deny type results of paragraph 4. We start with an improvment of Corollary 5.8.
\begin{prop}
\label{pro:5.21}
For any $u\in \mathbb S^{d-1}$,
$\displaystyle\mathop{\lim}_{t\rightarrow \infty} t^{\alpha} \  \widehat{\mathbb P} \{|\langle R, u \rangle  |>t\}=C\frac{p(\alpha)}{\alpha} \ {^{*}e^{\alpha}} (u)>0$ .

  In cases I, for any $u\in \mathbb S^{d-1}$ :
$\displaystyle\mathop{\lim}_{t\rightarrow \infty} t^{\alpha} \  \widehat{\mathbb P} \{\langle R, u \rangle   >t\}=\frac{1}{2} C \frac{p(\alpha)}{\alpha}\  {^{*}e^{\alpha}} (u)>0$ .

 In case II, for any $u\in \Lambda_{+} (T^{*})$, if $\overline{\Lambda_{a}(\Sigma)} \supset \Lambda_{+}^{\infty} (T)$:
$\displaystyle\mathop{\lim}_{t\rightarrow \infty} t^{\alpha} \   \widehat{\mathbb P} \{\langle R, u \rangle  >t\}= \frac{p(\alpha)}{\alpha} C_{+} {^{*}e^{\alpha}}(u)>0$.
\end{prop}

\begin{proof}
This a trivial consequence of Proposition 5.3, Corollary 5.8, Proposition 5.9 and Lemma 5.10.
\end{proof}

  The following is a corollary of the proof of Proposition 5.9 and of Proposition 5.21.
\begin{coro}
\label{cor:5.22}
With the above notation we write
\[\gamma_{\tau}^{\alpha}=L_{\mu} (\alpha) \mathbb E^{\alpha}_{0} (\tau), \psi_{\tau}^{\alpha} (v,p)=\widehat{\mathbb P}
\{t<p^{-1}\langle R, u \rangle   \leq t+p^{-1}\langle R_{\tau}, u\rangle \ ,\ \tau<\infty\} {^{*}e^{\alpha}} (u)^{-1} t^\alpha\]
 and we denote by $\kappa^{\tau}$ the  ${^{*}\widehat{Q}}^{\alpha,\tau}$-stationary measure on $X_{+}$, given by Lemma 5.12. Then, if $\overline{\Lambda_{a}(\Sigma)} \supset \Lambda_{+}^{\infty} (T)$ we have
\[C_{+}\geq \frac{\alpha}{p(\alpha) \gamma^{\alpha}_{\tau}} \int_{]0,\infty[\times X_{+}} \psi_{\tau}^{\alpha} (v,p) t^{-1} d\kappa^{\tau} (u,p) dt>0.\]
\end{coro}

\begin{proof}
 With the notation of Lemma 5.12, we have $\psi^{\alpha}=\displaystyle\mathop{\Sigma}_{0}^{\infty}
({^{*}\widehat P^{\tau}_{\alpha}})^k \psi_{\tau}^{\alpha}$ where ${^{*}\widehat P^{\tau}_{\alpha}}$ is a fibered Markov kernel on $X_{+}\times \mathbb R_{+}^{*}$ which satisfies the conditions of Proposition 5.14 and $\psi_{\tau}^{\alpha}$ is bounded by $({^{*}e^{\alpha}} \otimes h^{\alpha})^{-1}$. Hence, as in the proof of the proposition, if $\overline{\psi}_{\tau}^{\alpha}$ is a Borel function with compact support, bounded by $\psi_{\tau}^{\alpha}$,then
$\displaystyle\mathop{\lim}_{t\rightarrow \infty} \psi^{\alpha} (v,p) \geq \displaystyle\mathop{\lim\sup}_{t\rightarrow \infty} \displaystyle\mathop{\Sigma}_{0}^{\infty} \ {^{*}\widehat P^{\tau}_{\alpha}} \overline{\psi}_{\tau}^{\alpha} (v,p)$,
and, using Proposition 5.21, since $\displaystyle\mathop{\lim}_{t\rightarrow \infty} \psi^{\alpha} (v,p)$ is constant on $X_{+}$,we get
\[
\displaystyle\mathop{\lim}_{t\rightarrow \infty} \psi^{\alpha} (v,p)\geq \frac{1}{\gamma^{\alpha}_{\tau}} \int_{]0,\infty[\times X_{+}} \overline{\psi}_{\tau}^{\alpha} (v',p) d\kappa^{\tau} (u',p) \frac{dt}{t}.\]
 Hence, approximating from below $\psi_{\tau}^{\alpha}$ by $\overline{\psi}_{\tau}^{\alpha}$, we have
\[\displaystyle\mathop{\lim}_{t\rightarrow \infty} \psi^{\alpha} (v,p) \geq \frac{1}{\gamma^{\alpha}_{\tau}} \int_{]0, \infty[\times X_{+}} \psi_{\tau}^{\alpha} (v',p) \frac{dt}{t} d\kappa^{\tau} (u',p),\]
\[C (u) = \displaystyle\mathop{\lim}_{t\rightarrow \infty} t^{\alpha} \widehat{\mathbb P}\{\langle R, u \rangle  >t\}= {^{*}e^{\alpha}} (u) \displaystyle\mathop{\lim}_{t\rightarrow \infty} \psi^{\alpha} (v,p),\]
\[
C (u) \geq \frac{{^{*}e^{\alpha}}(u)}{\gamma_{\tau}^{\alpha}} \int_{]0, \infty[\times X_{+}} \psi_{\tau}^{\alpha} (v',p) t^{-1} dt d\kappa^{\tau} (u',p).\]

  The final formula follows from Proposition 5.21.
  \end{proof}

\begin{rema}

  We observe that, if $d=1$, and $A, B$ are positive, a formula of this type for $C=C_{+}$, with equality, is given in \cite{14}.  We don't know if such an equality is valid in our setting.
\end{rema}
\begin{lemm}

  For any compact subset $K$ of $V\setminus \{0\}$, there exists a constant $C(K)>0$ such that $\displaystyle\sup_{t>0} t^{-\alpha} (t \cdot  \rho) (K)\leq C(K)$. In particular the family $\rho_{t}=t^{-\alpha} (t \cdot  \rho)$ is relatively compact for the topology of vague convergence and any cluster value $\eta$ of the family $\rho_{t}$ satisfies $\displaystyle\sup_{t>0}  t^{-\alpha} (t\cdot\eta) (K) \leq C(K)$. Hence $\displaystyle\sup_{t>0} t\cdot((e^{\alpha}\otimes h^{\alpha})\eta) (K)\leq C'(K)$ with $C'(K) >0$.
\end{lemm}
\begin{proof}

  For some $\delta>0$ we have $K\subset \{x\in V\ ;\ |x| > \delta\}$, hence using Corollary 5.8, $t^{-\alpha} \mathbb P\{|R|>
\frac{\delta}{t}\}\leq \frac{b}{\delta^{\alpha}}=C(K)$. The relative  compactness of the family $\rho_{t}$ follows. Also,

\[(t t_{n})^{-\alpha} ( t t_{n} \cdot  \rho) (K) \leq C(K),\ \  t^{-\alpha} (t\cdot\eta) (K)= \displaystyle\mathop{\lim}_{n\rightarrow \infty} (t t_{n})^{-\alpha} (t t_{n} \cdot  \rho) (K)\leq C(K).\]

  Hence $\sup_{t>0} t^{-\alpha} (t\cdot\eta) (K) \leq C(K)$.
  Since $e^{\alpha} \otimes h^{\alpha}$ is $\alpha$-homogeneous we have $t\cdot((e^{\alpha}\otimes h^{\alpha})\eta)=(e^{\alpha}\otimes h^{\alpha}) (t^{-\alpha} (t\cdot\eta))$.
 With $C_{K}=\displaystyle\sup_{v\in K} (e^{\alpha} \otimes h^{\alpha}) (v)$, we get

\[t\cdot((e^{\alpha} \otimes h^{\alpha}) \eta) (K) \leq C_{K} \  t^{-\alpha} (t\cdot\eta) (K)\leq C_{K} C(K)=C'(K).\]

\begin{lemm}
\label{lemm:5.24}
  Assume $\eta$ is the vague limit of $t^{-\alpha}_{n} (t_{n} \cdot  \rho)$ $(t_{n}\rightarrow 0_{+})$. Then $\eta$ is $\mu$-harmonic, i.e $\mu*\eta=\eta$.
\end{lemm}
\begin{proof}

 Let $\varphi$ be $\varepsilon$-H\"older continuous on $V$ with compact support contained in the set $\{x\in V\ ;\ |x|\geq \delta\}$ with $\delta>0$, and let us show $\displaystyle\mathop{\lim}_{t\rightarrow 0_{+}} t^{-\alpha} I_{t} (\varphi)=0$ where $I_{t} (\varphi) =(t \cdot  \rho-t. (\mu*\rho)) (\varphi)$. By  definition $I_{t} (\varphi)=\mathbb E (\varphi (t R)-\varphi (tA_{1} Ro \widehat{\theta}))$ with $\varphi (t R)=0$ if $|t R|<\delta$ and $\varphi (t A_{1} Ro\widehat{\theta})=0$ if $|t A_{1} Ro \widehat{\theta}|<\delta$. Hence, $I_{t} (\varphi) \leq [\varphi]_{\varepsilon} t^{\varepsilon} \mathbb E (|B_{1}|^{\varepsilon} 1_{\{|t R|>\delta\}}+ |B_{1}|^{\varepsilon} 1_{\{t |A_{1} Ro \widehat{\theta}|>\delta\}})$.

We write
\[I_{t}^1=t^{\varepsilon-\alpha} \mathbb E(|B_{1}|^{\varepsilon} 1_{\{|t R|>\delta}\]
$I_{t}^2= t^{\varepsilon-\alpha} \mathbb E (|B_{1}|^{\varepsilon} 1_{\{t |A_{1} Ro \widehat{\theta}|>\delta\}})$
and we estimate $I_{t}^1, I_{t}^2$ as follows. We have $I_{t}^1\leq \delta^{\varepsilon-\alpha} \mathbb E (|B_{1}|^{\varepsilon}|R|^{\alpha-\varepsilon} 1_{\{|t R|>\delta\}}|$. Since $|R|^{\alpha-\varepsilon}\leq c(|A_{1}Ro \widehat{\theta}|^{\alpha-\varepsilon}+|B_{1}|^{\alpha-\varepsilon})$,  using independence of $Ro\widehat{\theta}$ and $|B_{1}|^{\varepsilon} |A_{1}|^{\alpha-\varepsilon}$, we get
\[ E(|B_{1}|^{\varepsilon} |R|^{\alpha-\varepsilon})\leq c \ E (|B_{1}|^{\alpha})+c \ E (|A_{1}|^{\alpha-\varepsilon} |B_{1}|^\varepsilon) \ E (|R|^{\alpha-\varepsilon}.\]
Using H\"older inequality we get $\mathbb E(|A_{1}|^{\alpha-\varepsilon} |B_{1}|^{\varepsilon})<\infty$. Also using Proposition 5.1, we get $\mathbb E (|R|^{\alpha-\varepsilon})<\infty$. It follows that $|B_{1}|^{\varepsilon} |R|^{\alpha-\varepsilon} 1_{\{|t R|>\delta\}}$ is bounded by the integrable function $|B_{1}|^{\varepsilon} |R|^{\alpha-\varepsilon}$. Then  by dominated convergence, $\displaystyle\mathop{\lim}_{t\rightarrow 0} I_{t}^1=0$.
 In the same way we have

\[I_{t}^2 \leq \delta^{\varepsilon-\alpha} \mathbb E (|B_{1}|^{\varepsilon} |A_{1} Ro \widehat{\theta}|^{\alpha-\varepsilon} 1_{\{ t|A_{1} Ro \widehat{\theta}|>\delta\}}.\]

  Also, using independence and H\"older inequality, we have

\[E(|B_{1}|^{\varepsilon} |A_{1} Ro \widehat{\theta}|^{\alpha-\varepsilon}) \leq \mathbb E (|B_{1}|^{\varepsilon} |A_{1}|^{\alpha-\varepsilon}) \mathbb E (|R|^{\alpha-\varepsilon})<\infty.\]

  Then by dominated convergence $\displaystyle\mathop{\lim}_{t\rightarrow 0_{+}}  I_{t}^2 =0$. Hence $\displaystyle\mathop{\lim}_{t\rightarrow 0_{+}} t^{-\alpha} I_{t} (\varphi)=0$.
  By definition of $\eta$ we have, for any $g\in G , {\{\lim}_{t_{n}\rightarrow 0_{+}} t_{n}^{-\alpha} (t_{n}\cdot {(g\rho) (\varphi)=(g\eta) (\varphi)}$.
 Furthermore we have $|\varphi (x)|\leq |\varphi| 1_\{x \in V ; |x|\geq \delta\}$ and,
\[|(g\eta) (\varphi)|\leq |\varphi| \eta \{x\in V ; |g x|\geq \delta\}\leq |\varphi| \displaystyle\mathop{\lim}_{n\rightarrow \infty} t_{n}^{-\alpha} \mathbb P \{|R| > \frac{\delta}{|g| t_{n}}\}.\]

   Using Corollary  5.8 we get $|(g\eta) (\varphi)|\leq  \frac{b}{\delta^{\alpha}} |\varphi | |g|^{\alpha}$.

  Since $\int |g|^{\alpha} d\mu (g) <\infty$ and for any $g\in G$, $\displaystyle\mathop{\lim}_{t_{n}\rightarrow 0_{+}} t_{n}^{-\alpha} (t_{n} \cdot (g\rho) (\varphi)=g  \eta (\varphi)$, we have by dominated convergence  $\displaystyle\mathop{\lim}_{t_{n}\rightarrow 0_{+}} t_{n}^{-\alpha} (t_{n}\cdot (\mu*\rho) (\varphi)=(\mu*\eta) (\varphi)$. Then the property $\displaystyle\mathop{\lim}_{t\rightarrow 0_{+}} t^{-\alpha} I_{t}=0$ implies $(\mu*\eta) (\varphi)=\eta (\varphi)$, hence $\mu*\eta=\eta$.
\end{proof}
\begin{lemm}
\label{lemm:5.25}
  Assume $\eta$ and $\sigma\otimes \ell^{\alpha}$ are $\mu$-harmonic Radon measures on $V\setminus \{0\}$ with $\sigma\in M^{1} (\mathbb S^{d-1})$. Assume also that for any $v\in V \setminus \{0\},\ \eta (H_{v}^{+})=(\sigma \otimes \ell^{\alpha}) (H_{v}^{+})$. Then we have $\eta=\sigma \otimes \ell^{\alpha}$.
\end{lemm}

\begin{proof}

  As in the proof of Corollary 5.8, we observe that the condition $\eta(H_{v}^{+})=(\sigma\otimes \ell^{\alpha}) (H_{v}^{+})$ implies for any $\delta>0$,
\[\sup_{t>0} t^{-\alpha}  (t\cdot\eta) \{ x\in V ; |x| >\delta\}<\infty .\]
  Hence, as in Lemma 5.25, for any compact $K\subset V\setminus \{0\}$, with $\eta^{\alpha}=(e^{\alpha} \otimes h^{\alpha}) \eta$, we have
\[\sup_{t>0} t^{-\alpha} (t\cdot\eta) (K) \leq C(K),\ \    \sup_{t>0} (t\cdot\eta^{\alpha}) (K) \leq C'(K).\]

  It follows that $\eta^{\alpha}$ is dilation-bounded.

  We recall that $P($resp $\breve{P})$ is the convolution operator by $\mu$ on $V\setminus\{0\}$ (resp $\breve{V}$), hence $\breve{P}(\overline{\sigma} \otimes \ell^{\alpha})=\overline{\sigma}\otimes \ell^{\alpha}$ and $\breve{P}(e^{\alpha}\otimes h^{\alpha})=e^{\alpha}\otimes h^{\alpha}$. We denote by $\breve{Q}_{\alpha}$ the Markov operator on $\breve{V}$ deduced from $\breve{P}$ by Doob's relativisation with respect to $e^{\alpha}\otimes h^{\alpha}$.
  On the other hand, the projection $\overline{\sigma} \otimes \ell^{\alpha}$ (resp $\breve{\eta}$) of $\sigma \otimes \ell^{\alpha}$ (resp $\eta$) on $\breve{V}$ satisfies

   $\mu*(\overline{\sigma} \otimes \ell^{\alpha})=\overline{\sigma}\otimes \ell^{\alpha}$ (resp $\mu*\breve{\eta}=\breve{\eta}$), hence $\breve{Q}_{\alpha} (e^{\alpha} \overline{\sigma} \otimes \ell)=e^{\alpha} \overline{\sigma} \otimes \ell$, $\breve{Q}_{\alpha} ( \breve{\eta}^{\alpha})=\breve{\eta}^{\alpha}$.

   We observe that the fibered Markov operator $\breve{Q}_{\alpha}$ satisfies condition $D$ of subsection 5.4, in view of Corollary 3.20 and of the moment condition on $A$.

  Then Theorem 2.6 implies $\overline{\sigma}=\nu^{\alpha}$. Also, in view of Corollary 3.20 for $s=\alpha$ and the above observations, we can apply the second part of Theorem 5.17  to $\breve{\eta}^{\alpha}$ with $P= \breve{Q}_{\alpha}$, hence $\breve{\eta}^{\alpha}$ is proportional to $\pi^{\alpha}\otimes \ell$, i.e. $\breve{\eta}$ is proportional to $\nu^{\alpha} \otimes \ell^{\alpha}$. Since $\overline{\sigma}=\nu^{\alpha}$ and $\eta(H_{v}^{+})=(\sigma \otimes \ell^{\alpha}) (H_{v}^{+})$ we get $\breve{\eta}=\nu^{\alpha}\otimes \ell^{\alpha}$.

  We denote for $v\in\mathbb S^{d-1}$, $\lambda_{v}=|\langle v, \cdot\rangle^{\alpha} \sigma \otimes \ell^{\alpha}$, $\eta_{v}=|\langle v,\cdot\rangle|^{\alpha} \eta$.
  Since $\sigma\otimes \ell^{\alpha}$ and $\eta$ are $\mu$-harmonic, we have $\int g \lambda_{g^{*}v}   d\mu(g)=\lambda_{v}$, $\int g \eta_{g^{*}v} d\mu (g)=\eta_{v}$.
   The projections $\breve{\lambda}_{v}$ and $\breve{\eta}_{v}$ on $\breve{V}$ satisfy the same equation, hence are equal.
  As in section 3, we get that the sequences of Radon measures $g_{1}\cdots g_{n} \lambda_{g^*_{n}\cdots g_{1}^* v}$ and $g_{1}\cdots g_{n} \eta_{g_{n}^*\cdots g_{1}^* v}$ are vaguely bounded $^{*}\mathbb Q^{\alpha}_{v}$-martingales.
  On $\breve{V}$ we get, using Theorem 3.2, for some $z(\omega)\in \mathbb P^{d-1}$ and $^{*}\mathbb Q^{\alpha}_{v} - \textrm{a.e.}$
\[\lim_{n\rightarrow \infty} g_{1}\cdots g_{n} \breve{\lambda}_{g^*_{n}\cdots g^*_{1}v}=\lim_{n\rightarrow \infty} g_{1}\cdots g_{n} \breve{\eta}_{g^*_{n}\cdots g_{1}^* v}=\delta_{z(\omega)}\otimes \ell.\]

  Let $z_{+}(\omega)$ and $z_{-}(\omega)$ be opposite points on $\mathbb S^{d-1}$ with projection $z(\omega)$ on $\mathbb P^{d-1}$. The martingale convergence on $V\setminus \{0\}$ gives that $g_{1}\cdots g_{n} \lambda_{g^*_{n}\cdots g^*_{1}v}$ (resp $g_{1}\cdots g_{n} \eta_{g^*_{n}\cdots g^*_{1}v}$) converges vaguely to $p(\omega) \delta_{z_{+}(\omega)}+q(\omega) \delta_{z_{-}(\omega)}$ (resp $p'(\omega) \delta_{z_{+}(\omega)}+q'(\omega) \delta_{z_{-}(\omega)})$ with $p(\omega)+ q(\omega)=p'(\omega)+q'(\omega)=1$. The condition $\eta (H_{v}^{+})=(\sigma \otimes \ell^{\alpha})  (H_{v}^{+})$ implies in the limit:
\[p(\omega \delta_{z_{+}(\omega)}+q(\omega) \delta_{z_{-}(\omega)}=p'(\omega) \delta_{z_{+}(\omega)}+q'(\omega) \delta_{z_{-}(\omega)}.\]

   Hence, taking expectations we get $\eta=\sigma \otimes \ell^{\alpha}$.
   \end{proof}
   \textit{Proof of theorem 5.2.}
  The convergence of $\rho_{t}=t^{-\alpha}(t \cdot  \rho)$ to $C(\sigma^{\alpha}\otimes \ell^{\alpha})$ on the sets $H_{v}^+$ and the positivity properties of $C,C_{+}, C_{-}$ follows from Corollary 5.21. For the vague convergence of $\rho_{t}$ we observe that Lemma 5.25 gives the vague compactness of $\rho_{t}$. If $\eta=\displaystyle\mathop{\lim}_{t_{n}\rightarrow 0_{+}} t_{n}^{-\alpha} (t_{n} \cdot  \rho)$, Lemma 5.24 gives the $\mu$-harmonicity of $\eta$. Since $\eta (H_{v}^+)=C(\sigma^{\alpha}\otimes \ell^{\alpha}) (H_{v}^+)$, Lemma 5.25 gives $\Lambda=C(\sigma^{\alpha} \otimes \ell^{\alpha})=\eta$, hence the vague convergence of $\rho_{t}$ to $\Lambda$. The detailed form of $\Lambda$ follows from Proposition 5.3.

  For the final minimality assertions one uses the second part of Theorem 5.17 and we replace $\eta$ by $(e^{\alpha}\otimes h^{\alpha}) \eta$. We verify condition $D$ for $\widetilde{Q}^{\alpha}$ or operators associated with $\widetilde{Q}^{\alpha}$ as follows. Part 1 of condition $D$ follows directly from Corollary 3.21.  In case I, $\widetilde{Q}^{\alpha}$ satisfies part 2 of condition $D$ by Corollary 3.21, 1 is a simple eigenvalue of $\widetilde{Q}^{\alpha}$ and, if $\theta$ is positive Radon measure with $\mu * \theta=\theta$, $\theta\leq \widetilde{\nu}^{\alpha} \otimes \ell^{\alpha}$,  hence $\theta$ is proportional to $\widetilde{\nu}^{\alpha} \otimes \ell^{\alpha}$. In case II, one restricts $\widetilde{Q}^{\alpha}$ to the convex cone generated to $\Lambda_{+} (T)$, so as to achieve the simplicity of 1 as an eigenvalue of $\widetilde{Q}^{\alpha}$ and the absence of other unimodular eigenvalue. Then Corollary 3.21 shows that part 2 of condition $D$ is satisfied for the corresponding operator, hence the above argument is also valid for $\nu_{+}^{\alpha}\otimes \ell^{\alpha}$. Part 3 of condition $D$ follows from the moment conditions assumed on $\lambda$.
  \end{proof}

\begin{rema}
  In general $\overline{\supp (\sigma^{\alpha}\otimes \ell^{\alpha})} \cap \mathbb S^{d-1}_{\infty}$ is smaller than $\Lambda^{\infty}_{a} (\Sigma)$  and $\supp \sigma^{\alpha}$ has a  fractal structure.

  In the context of extreme value theory for the process $X_{n}$, the convergence stated in the theorem plays a basic role and implies that $\rho$ has multivariate regular variation.

   Actually, using  the properties of $\Lambda$, \cite{5} gives also the weak convergence for any $\alpha$ (resp  $\alpha\notin 2 \mathbb N)$ in case II'' (resp case I).

  This is valid too for $\alpha\notin 2 \mathbb N$ if $C_{+}=C_{-}$ in case II', for example if the law of $B_{1}$ is symmetric (see \cite{38}).
\end{rema}
  For the last assertion in Theorem C we need the following.

\begin{prop}
\label{prop:5.26}
  Let $\mathcal{B}_{\varepsilon,\alpha}$ be the set of locally bounded Borel functions on $V\setminus \{0\}$ such that the set of discontinuities of $f$ is $\Lambda$-negligible and for $\varepsilon>0$,

\[K_{f}(\varepsilon)=\sup \{|v|^{-\alpha} |\log v|^{1+\varepsilon} |f(v)|\ ;\ v\neq 0\}<\infty .\]

  Then for any $f\in \mathcal{B}_{\varepsilon,\alpha}$ \qquad
$\displaystyle\mathop{\lim}_{t\rightarrow 0_{+}} t^{-\alpha} (t \cdot \rho) (f)=\Lambda (f)$.
\end{prop}

  The proof depends of two lemmas in which we will use the norm
  $\|v\|=\displaystyle\sup_{1\leq i\leq d} |\langle x,e_{i}
  \rangle|$ instead of $|v|$  where $e_{i} (1\leq i\leq d)$ is a basis of $V$. Also for $\delta>0$ and $0<\delta_{1}<\delta_{2}$ we write

$B_{\delta}=\{v\in V\ ;\ \|v\| \leq \delta\},  B_{\delta_{1},\delta_{2}} = B_{\delta_{2}}\setminus B_{\delta_{1}}, B'_{\delta}=V\setminus B_{\delta}$.

\begin{lemm}
\label{lemm:5.27}
  For any $f\in \mathcal{B}_{\varepsilon,\alpha}$, $0<\delta_{1}<\delta_{2}$,
\[{\lim}_{t\rightarrow 0_{+}} t^{-\alpha} (t \cdot  \rho) (f 1_{B_{\delta_{1}, \delta_{2}}})=\Lambda (f 1_{B_{\delta_{1}, \delta_{2}}}).\]
\end{lemm}
\begin{proof}
  From the fact that $\nu^\alpha$ gives measure zero to any projective subspace and the homogeneity of $\Lambda=\sigma^{\alpha} \otimes \ell^{\alpha}$, we know that $\Lambda$ gives measure zero to any affine hyperplane, hence the boundary of $B_{\eta_{1},\eta_{2}}$ is $\Lambda$-negligible. Then the proof follows from the vague convergence of $t^{-\alpha} (t \cdot  \rho)$ to $\Lambda$ and the hypothesis of $\Lambda$-negligibility of the discontinuity set of $f$.
\end{proof}
\begin{lemm}
\label{lemm:5.28}
   There exists $C>0$ such that for any $f\in \mathcal{B}_{\varepsilon,\alpha}$, $t>0$ and $\delta_{2}>e$,
\[|t^{-\alpha} (t \cdot  \rho) (f 1_{B_{\delta_{2}}})|\leq C K_{f} (\varepsilon) |\log \delta_{2}|^{-\varepsilon}.\]
 Furthermore  there exists $C(\varepsilon)>0$ such that for any $f\in \mathcal{B}_{\varepsilon,\alpha}, t>0,\ \delta_{1}<e^{-1}$
\[|t^{-\alpha} (t \cdot  \rho) (f 1_{B_{\delta_{1}}})|\leq C(\varepsilon) K_{f} (\varepsilon) |\log \delta_{1}|^{-\varepsilon}.\]
\end{lemm}
\begin{proof}
   Let $\varphi_{\varepsilon} (x)$ be the function on $\mathbb R_{+}\setminus \{1\}$ given by $\varphi_{\varepsilon} (x)=x^{\alpha} |\log x|^{-1-\varepsilon}$. For $x\geq e$ we have $\varphi'_{\varepsilon} (x)\leq \alpha \ x^{\alpha-1} |\log x|^{-1-\varepsilon}$. We denote $F_{t}(x)=\widehat{\mathbb P}(|t R|>x)$ and we observe that, using Proposition 5.1, the non increasing function  $F_{t}$ is continuous. We have
\[|t^{-\alpha} (t \cdot  \rho) (f 1_{B'_{\delta_{2}}})|\leq t^{-\alpha} K_{f} (\varepsilon) \int_{\delta_{2}}^{\infty} \varphi_{\varepsilon} (x) dF_{t} (x).\]
  Integrating by parts, we get
\[|t^{-\alpha} (t \cdot  \rho) (f 1_{B'_{\delta_{2}}})|\leq t^{-\alpha} K_{f} (\varepsilon)  [\varphi_{\varepsilon} (x) F_{t} (x)]_{\delta_{2}}^{\infty}+t^{-\alpha} K_{f} (\varepsilon) \int_{\delta_{2}}^{\infty} \varphi'_{\varepsilon} (x) F_{t} (x) dx.\]

  From Corollary 5.7 we know that, for some $C>0$, $F_{t} (x) \leq C\ t^{\alpha} x^{-\alpha}$. Then, using the above estimation of $\varphi'_{\varepsilon} (x)$, we get
\[|t^{-\alpha} (t \cdot  \rho) (f 1_{B'_{\delta_{2}}})|\leq C K_{\varepsilon} (f) \int_{\delta_{2}}^{\infty}\frac{\alpha}{|\log x|^{1+\varepsilon}}\ \frac{dx}{x} \leq C K_{\varepsilon} (f) |\log \delta_{2}|^{-\varepsilon}.\]

   The proof of the second assertion follows the same lines and uses the estimation of $|\varphi'_{\varepsilon} (x)|$ by
   $(\alpha+1+\varepsilon) x^{\alpha-1} |\log x|^{-1-\varepsilon}$ for $x\leq e^{-1}$.
  \end{proof}

\textit {Proof of proposition 5.26.}
  For $\delta\geq e$ and with $D>0$ we have $\int_{B'_{\delta}}\varphi_{\varepsilon} (\|v\|) d \Lambda (v) \leq D \int_{\delta}^{\infty} \frac{x^{\alpha}}{(\log x)^{1+\varepsilon}} \frac{dx}{x^{\alpha+1}}=\frac{D}{|\log \delta|^{\varepsilon}}$ hence $\displaystyle\mathop{\lim}_{\delta\rightarrow \infty} \int_{B'_{\delta}} \varphi_{\varepsilon} (\|v\|) d\Lambda (v)=0$.
 Also for $\delta<1$,
\[\int_{B'_{\delta}} \varphi_{\varepsilon} (\|v\|) d\Lambda (v) \leq D \int_{0}^{\delta} \frac{x^{\alpha}}{|\log x|^{1+\varepsilon}} \frac{dx}{x^{\alpha+1}}=\frac{D}{|\log \delta|^{\varepsilon}},\]
  hence $\displaystyle\mathop{\lim}_{\delta\rightarrow \infty} \int_{B_{\delta}} \varphi_{\varepsilon} (\|v\|) d\Lambda (v)=0$.

  Then the Proposition  follows from the lemmas.
$\hfill\square$

\textit {Proof of  Theorem C.}
  Except for the last assertion, Theorem C is a direct consequence of Theorem 5.2. The last assertion is the content of Proposition 5.26. \hfill $\square$

\appendix
\section{An analytic approach to tail-homogeneity}

  Under the hypothesis of compact support for $\lambda$ and density for $\mu$, an analytic proof of tail-homogeneity of $\rho$ is given below. The full hypothesis is only used in the study of positivity properties of $C, C_{+}, C_{-}$, hence we split the presentation into two parts according to the hypothesis at hand on $\lambda$. Also the argument gives analytic expressions for $C, C_{+}, C_{-}$. We  recall  Wiener-Ikehara's theorem (see\cite{52}, p.233). Assume $A(x)$ is non negative,  increasing on $[1,\infty[$, $f(s)=\int_{1}^{\infty} x^{-s-1} A(x) dx$ is finite for $s>1$, $f$ extends as a function $\tilde f$ meromorphic in an open set $D\supset\{0<Rez \leq \alpha\}$ and $\tilde {f}$ has only a possible unique simple pole at $z=\alpha$, with $\displaystyle\mathop{\lim}_{s\rightarrow \alpha_{-}} (\alpha-s) f(s)=A$; then one has $\displaystyle\mathop{\lim}_{x\rightarrow \infty}
 x^{-1}  A (x)=A$. The use of this result  (see Lemma A.1 below)will give the tail-homogeneity of $\rho$. On the other hand, if $\beta$ denotes the convergence abcissa    of the Mellin transform $f(s)=\int_{0}^{\infty} x^{s} d\nu (x)$, a lemma of E. Landau (see \cite{52} p. 58) says that $f$ cannot be extended holomorphically to a neighbourhood of $\beta$. This will allow us to show $A>0$.

 \begin{lemm}
 \label{lemm:A.1}
   Let $\nu$ be a probability on $[1,\infty[$, $\alpha>0$ such that $f(x)=\int_{1}^{\infty} x^{s} d\nu (x)$ is finite for $s<\alpha$, $f(s)$ extends to an open set $D\supset\{0<$Rez$\leq \alpha\}$ as  a meromorphic function $\tilde{f}$ which has a simple pole at $z=\alpha$ with residue $A>0$. Then one has $\displaystyle\mathop{\lim}_{x\rightarrow \infty} x^{\alpha} \nu (x,\infty)=\alpha^{-1} A$.
 \end{lemm}
 \begin{proof}

   We write $A(x)=\int_{1}^{x}  t^{\alpha} \nu (t, \infty) dt$ for $x\geq 1$ and we observe that the finiteness of $f(s)$ for $s<\alpha$ implies $\displaystyle\mathop{\lim}_{x\rightarrow \infty} x^{\alpha-\varepsilon} \nu (x,\infty)=0$ for $\varepsilon >0$.

   Integrating by parts we have for $s>1$

 \[\int_{1}^{\infty} x^{-s-1} A(x) dx =s^{-1} \int_{1}^{\infty} x^{\alpha-s} \nu (x,\infty) dx =s^{-1} (\alpha-s+1)^{-1} \int_{1}^{\infty} x^{\alpha-s+1} d\nu(x).\]

   Since $A(x)\geq 0$ is increasing, we can use Wiener-Ikehara's theorem: $\displaystyle\mathop{\lim}_{x\rightarrow \infty} x^{-1} A(x)=\alpha^{-1} A$. Then one can apply the Tauberian Lemma 5.4 to the decreasing function $\nu(x,\infty)$ and gets $\displaystyle\mathop{\lim}_{x\rightarrow \infty} x^{\alpha} \nu (x,\infty)=\alpha^{-1} A$.
 \end{proof}
 The connection with the spectral gap properties in section 3 depends on the following. The hypothesis is as in Corollary 3.21.
 \begin{lemm}
 \label{lemm:A.2}
   There exists an open set $D\subset \mathbb C$ which contains the set $\{Rez \in ]0, \alpha]\}$ such that $(I-\widetilde{P}^z)^{-1}$ is meromorphic in $D$ with a unique simple pole at $z=\alpha$. We have

 In case I, $\displaystyle\mathop{\lim}_{z\rightarrow \alpha} (\alpha-z) (I-\widetilde{P}^z)^{-1}=k'(\alpha)^{-1} (\widetilde{\nu}^{\alpha} \otimes e^{\alpha}).$

 In cases II, $\displaystyle\mathop{\lim}_{z\rightarrow \alpha} (\alpha-z) (I-\widetilde{P}^z)^{-1}=k'(\alpha)^{-1} (\nu^{\alpha}_{+} \otimes e^{\alpha}_{+}+\nu^{\alpha}_{-}\otimes e^{\alpha}_{-})$.
 \end{lemm}
 \begin{proof}

  We restrict to case I, since the proof is similar in cases II. The operator $\widetilde{P}^z$ on $H_{\varepsilon}(\mathbb S^{d-1})$ where $z=s+it$,  is defined by the formula $\widetilde{P}^z\varphi (x)=\int |g x|^z \varphi ( g \cdot x) d\mu (g)$ and  $\frac{1}{k(s)} \widetilde P^z$ is conjugate to the operator $\widetilde{Q}^z$ considered in Corollary 3.21. From this corollary we deduce that $\widetilde{Q}^{\alpha}$ satisfies a Doeblin-Fortet condition and $r(\widetilde P^{\alpha+it})=r(\widetilde{Q}^{\alpha+it})<1$ if $t\neq 0$. On the other hand the function $z\rightarrow \widetilde{P}^z$ is holomorphic in the set $\{0<Rez <s_{\infty}\}$ since for any loop $\gamma$ in this set we have $\int_{\gamma} \widetilde{P}^z dz=\int \varphi ( g \cdot x) d\mu (g) \int _{\gamma} |g x|^z dz=0$. It follows that there exists $\varepsilon>0$ such that for $|z-\alpha|<\varepsilon$ there exists a holomorphic function $k(z)$ such that $k(z)$ is a simple dominant eigenvalue of $\widetilde{P}^z$ with $k(z)=1+k'(\alpha) (z-\alpha)+\circ (z-\alpha)$. Since $k'(\alpha)\neq 0$, we have $k(z)\neq 1$ for $z\neq \alpha$ and $|z-\alpha|$ small. Also for $|z-\alpha|$ small, we have in case I  the decomposition $\widetilde{P}^z=k(z) \widetilde{\nu}^z \otimes e^z+U(z)$ where $\widetilde{\nu}^z\otimes e^z$ is a projector on the line $\mathbb C e^z$, $U(z)$ satisfies $U(z) (\widetilde{\nu}^z \otimes e^z)=(\widetilde{\nu}^z\otimes e^z) U(z)=0$, $r(U(z))<1$, and $\widetilde{\nu}^z \otimes e^z$, $U(z)$ depend holomorphically on $z$. We consider also the projection $p^z=I- \widetilde{\nu}^z \otimes e^z$ and we write $I-\widetilde{P}^z=(1-k (z)) (\widetilde{\nu}^z \otimes e^z)+ p^z (I-U(z))$. Hence, for $|z-\alpha|$ small
    \[(I-\widetilde{P}^z)^{-1}=(1-k(z))^{-1}(\widetilde{\nu}^z \otimes e^z)+p^z (I-U(z))^{-1}.\]

   In particular, $\displaystyle\mathop{\lim}_{z\rightarrow \alpha} (\alpha-z) (I-\widetilde{P}^z)^{-1}=k'(\alpha)^{-1} (\widetilde{\nu}^{\alpha}\otimes e^{\alpha})$, and $(I-\widetilde{P}^z)^{-1}$ is meromorphic in a disk $B_{0}$ centered at $\alpha$ with radius $\varepsilon' \leq \varepsilon$, with unique pole at $z=\alpha$. For $z=\alpha+it$ with $|t|\geq \varepsilon'$, we get from above that there exists a disk $B_{t}$ centered at $\alpha+it$ such that $r(\widetilde{P}^z)<1$ for $z\in B_{t}$, hence $(I-\widetilde{P}^z
)^{-1}$ is a bounded operator depending holomorphically on $z$ for $z\in B_{t}$. If $Rez \in ]0,\alpha[$, then $r(\widetilde P^z)\leq r (\widetilde P^s) <1$ hence $(I-\widetilde P^z)$ is invertible and the function $z\rightarrow (I-\widetilde P^z)^{-1}$ is holomorphic in the domain $\{0< Rez <\alpha\}$. Then the open set $D=(\displaystyle\mathop{\cup}_{t} B_{t}) U\{Rez \in ]0,\alpha[\}$, where $t=0$ or $|t| \geq \varepsilon'$ satisfies the conditions of the lemma, hence the formula for $\displaystyle\mathop{\lim}_{z\rightarrow \alpha} (\alpha-z) (I-\widetilde{P}^z)^{-1}$ is valid.
\end{proof}

  We denote, for $u\in \mathbb S^{d-1}$ and $Rez=s<\alpha$,
\[\bar f_{z} (u)=\mathbb E (\langle R, u \rangle  _{+}^{z})\quad ,\quad \bar d_{z} (u)=\mathbb E(\langle R, u \rangle  _{+}^{z}-\langle R-B,u \rangle_{+}^{z}).\]

\begin{prop}
\label{prop:A.3}
  With the notation and hypothesis of Theorem 5.2, we have the  convergence :
\[\lim_{t\rightarrow 0_{+}} t^{-\alpha} (t\cdot \rho) (H_{u}^{+})=C(\sigma^{\alpha}\otimes \ell^{\alpha}) (H^{+}_{u})=C(u) = \alpha^{-1} \displaystyle\mathop{\lim}_{s\rightarrow \alpha_{-}} (\alpha-s) \bar f_{s}(u),\]
  where $C\geq 0$ and $\sigma^{\alpha}\in M^{1} (\widetilde{\Lambda}(T))$ satisfies $\mu * (\sigma^{\alpha}\otimes \ell^{\alpha})=\sigma^{\alpha} \otimes \ell^{\alpha}$.

  In case $I$, $\sigma^{\alpha}$ is symmetric, $\supp\sigma^{\alpha}=\widetilde{\Lambda}(T)$ and $C(u)=(\alpha k' (\alpha))^{-1} {^{*}\widetilde{\nu}}^{\alpha} (\bar d_{\alpha}) {^{*}e}^{\alpha}(u)$.

  In cases II, $C(u)=(\alpha k' (\alpha))^{-1} [{^{*}\nu}_{+}^{\alpha} (\bar d_{\alpha}) {^{*}e}_{+}^{\alpha} (u) + {^{*}\nu}_{-}^{\alpha} (\bar d_{\alpha}) {^{*}e}_{-}^{\alpha} (u)]$.

\end{prop}
\begin{proof}
We write equation $(S)$ of section 5 in the form: $R-B=A R\circ\widehat{\theta}$.

  For any $v\in V\setminus \{0\}$, $Rez=s\in [0,\alpha[$ we define
\[f_{z}(v)=\mathbb E(\langle R,v\rangle_{+}^{z}),\ \ f_{z}^{1}(v)=\mathbb E(\langle R-B, v \rangle_{+}^{z}).\]

  Then equation $(S)$ implies: $\langle R-B, v \rangle_{+}=\langle R\circ \widehat{\theta},\ \ A^{*} v \rangle_{+}$, hence ${^{*}P} f_{z}=f_{z}^{1}$, $(I-{^{*}P}) f_{z}=f_{z}-f_{z}^{1}=d_{z}$ with $d_{z}(v)=\mathbb E(\langle R, v \rangle_{+}^{z} -
  \langle R-B, v \rangle_{+}^{z})$. We write a continuous $z$-homogeneous function $f$ on $V\setminus\{0\}$ as $f=\bar{f}\otimes h^{z}$ with $\bar{f}\in C(\mathbb S^{d-1})$, and we recall that, as in section 2, ${^{*}P}f={^{*}\widetilde{P}}^{z} \bar{f}\otimes h^{z}$. Then, since $f_{z}$ and $d_{z}$ are $z$-homogeneous and continuous, equation $(S)$ gives, $(I-{^{*}\widetilde{P}}^{z}) \bar{f}_{z}=\bar{d}_{z}$.

  For $u\in \mathbb S^{d-1}$ and $\varepsilon(z)=|z| 1_{[1,s_{\infty}[} (s)$, $\varepsilon'(z)=1_{[0,1]} (s)$, $\bar d_{z} (u)$ is dominated by
\[\mathbb E(\langle R, u \rangle  _{+}^{z}-\langle R-B , u \rangle _{+}^{z})|\leq \varepsilon'(z) \mathbb E (|B|^{s})+\varepsilon (z) \mathbb E (|B| (|B|+\langle R, u \rangle  _{+} )^{s-1}).\]
  Hence using H\"older inequality and the moment hypothesis we get that for $u$ fixed, $\bar{d}_{z} (u)$ is a holomorphic function in the domain $Rez \in ]0,\alpha+\delta[$. On the other hand, Lemma A.2 for $\mu^{*}$ implies that the operator-valued function $(I-{^{*}\widetilde{P}}^{z})^{-1}$ is meromorphic in an open set $D$ which contains the set $\{Rez \in ]0,\alpha]\}$, with unique simple pole at $\alpha\in D$. The above estimation of $\bar{d}_{z} (u)$ shows that the same meromorphy property is valid for $\bar f_{z}=(I-{^{*}\widetilde{P}}^{z})^{-1} \bar{d}_{z}$. If we denote by $\rho_{u}$ the law of $\langle R, u \rangle  _{+}$, we have $\bar{f}_{s} (u)=\int x^{s} d\rho_{u}(x)$, hence $\bar{f}_{s}(u)$ is the Mellin transform of the positive measure $\rho_{u}$. Then we can apply Lemma A.1   to $\bar{f}_{s} (u)$, $\rho_{u}$ and obtain the tail of $\rho_{u}$ in the form: $\displaystyle\mathop{\lim}_{t\rightarrow \infty} t^{\alpha} \rho_{u} (t,\infty)=\displaystyle\mathop{\lim}_{s\rightarrow \alpha_{-}} \alpha^{-1}(\alpha-s) \bar f_{s} (u)$. Hence using Lemma A.2 we have

In case I, $\displaystyle\mathop{\lim}_{s\rightarrow \alpha_{-}} (\alpha-s) \bar{f}_{s} (u)=\frac{1}{k'(\alpha)} {^{*}\widetilde{\nu}^{\alpha}} (\bar{d}_{\alpha})^{*} e^{\alpha}(u)$.

In case II, $\displaystyle\mathop{\lim}_{s\rightarrow \alpha_{-}} (\alpha-s) \bar{f}_{s} (u)=\frac{1}{k'(\alpha)} [{^{*}\nu^{\alpha}_{+}} (\bar{d}_{\alpha}){^{*} e^{\alpha}_{+}}(u)+{^{*}\nu^{\alpha}_{-}} (\bar{d}_{\alpha}){^{*}e^{\alpha}_{-}}(u)]$.

  Using the expressions of ${^{*}e^{\alpha}} (u), {^{*}e^{\alpha}_{+}} (u), {^{*}e^{\alpha}_{-}} (u)$ given by Theorem \ref{thm:2.17} we obtain in the two cases,
\[\lim_{s\rightarrow \alpha_{-}} (\alpha-s) \bar{f}_{s} (u)=\alpha C (\sigma^{\alpha}\otimes \ell^{\alpha}) (H_{u}^{+})=\alpha C(u),\]
  with  certain constants $C\geq 0$, $C(u)\geq0$,  a certain $\sigma^{\alpha}\in M^{1}(\mathbb S^{d-1})$ which satisfies
\[\mu * (\sigma^{\alpha} \otimes \ell^{\alpha})=\sigma^{\alpha}\otimes \ell^{\alpha}.\]

  In case I, we see that $\sigma^{\alpha}$ is symmetric and $\supp \sigma^{\alpha}=\widetilde{\Lambda} (T)$.

  In cases II, the detailed expression of $\sigma^{\alpha}$ shows that $C\sigma^{\alpha}= C_{+} \sigma_{+}^{\alpha}+C_{-} \sigma_{-}^{\alpha}$ where $\supp \sigma_{+}^{\alpha}=\Lambda_{+} (T)$, $\supp \sigma_{-}^{\alpha}=\Lambda_{-} (T)$. As a result we have $\displaystyle\mathop{\lim}_{t\rightarrow \infty} t^{\alpha} \rho_{u} (t,\infty)=C(\sigma^{\alpha}\otimes \ell^{\alpha}) (H_{u}^{+})=C(u)$.
\end{proof}
\begin{prop}
\label{prop:A.4}
  With the notation and hypothesis of Proposition A.3, assume futhermore that $\supp\lambda$ is compact and $\mu$ has a density on $G$. Then $C>0$. In cases II, we have $C(u)>0$ if $u\in \Lambda_{+}(T^{*})$.
\end{prop}

\begin{proof}
  Assume $C=0$ and observe that Proposition A.3 gives $C(u)=0$ for any $u \in \widetilde{\Lambda} (T^{*})$. The equation $(I-{^{*}\widetilde{P}}^{s}) \bar{f}_{s}=\bar{d}_{s}$ which occurs in the proof of Proposition A.3 gives by integration against ${^{*}\widetilde{\nu}^{s}}$, since ${^{*}\widetilde{P}}^{s} {(^{*}\widetilde{\nu}^{s})}=k(s) {(^{*}\widetilde{\nu}^{s}})$,
\[(1-k(s)) {^{*}\widetilde{\nu}^{s}} (\bar{f}_{s})={^{*}\widetilde{\nu}^{s}} (\bar{d}_{s})\quad (s<\alpha).\]
  Here,  the estimation  of $\bar{d}_{z}$  in the proof of Proposition A.3 gives the analyticity of ${^{*}\widetilde{\nu}^{z}} (\bar{d}_{z})$ in an open set containing $]0, \alpha+\delta[$ where $\tilde{\nu}^{z}$ is defined in the proof of Lemma A.2 by perturbation.
Since $C(u)=0$, the function $\tilde f_{z}$ considered in the proof of Proposition A.3 is meromorphic with no pole at $\alpha$, hence   Landau's lemma mentioned above gives  that the convergence abcissa of $\int x^{s} d \rho_{u} (x)$ is larger than $\alpha$ and $\bar{f}_{s} (u)<\infty$ if $0<s<\alpha+\delta$.
\vskip 2mm
  Now we consider the case $s\geq \alpha+\delta$. We observe that in case I, the density hypothesis on $\mu$ implies that $\supp {^{*}\widetilde{\nu}^s}=\supp {^{*}\widetilde{\nu}}=\mathbb S^{d-1}$ and the compactness hypothesis of $\supp \lambda$ implies that ${^{*}\widetilde P^z}$ defines a compact operator on $C(\mathbb S^{d-1})$. Let $\Delta \subset \{Rez >0\}$ be the set where $I-{^{*}\widetilde{P}^z}$ is not invertible and observe that, since $r({^{*}\widetilde{P}^z})<1$ if $Rez \in ]0,\alpha[$, we have $\Delta \cap \{0<Rez <\alpha\}=\phi$. Since the function $z\rightarrow {^{*}\widetilde{P}^z}$ is holomorphic, the extension of Riesz-Schauder theory given in \cite{34} implies that $\Delta$ is discrete without any accumulation point and $(I-{^{*}\widetilde{P}^z})^{-1}$  is meromorphic in the domain $\{0<Rez\}$, with possible poles in $\Delta$ only. The same property is valid for the function $\tilde{f}_{z}=(I- {^{*}\widetilde{P}^z})^{-1} (\tilde{d}_{z})$ in any domain $D\subset \{Rez >0\}$ where $\tilde{d}_{z}$ extends $\bar{d}_{s}$ holomorphically. We define
\[\beta=\sup\{s>0\ ;\ \mathbb E(\langle R, u \rangle  _{+}^{s})<\infty \textrm{ for } u\in \mathbb S^{d-1}\},\]
  hence from above, $\beta\geq \alpha+\delta$, and we show below $\beta=\infty$. Assume $\beta<\infty$ and let us show $\displaystyle\mathop{\lim}_{s\rightarrow \beta_{-}} \ {^{*}\widetilde{\nu}^s} (\bar{f}_{s})<\infty$. We observe that $\bar f_{z}$ (resp $\bar{d}_{z})$ is well defined and holomorphic in the domain $\{0< Rez <\beta\}$ (resp $D=\{0< Rez <\beta+1\}$, because $\supp \lambda$ is compact and the estimation of  $\bar{d}_{z}$ given in the proof of Proposition A.3. Hence the function $z\rightarrow \tilde{f}_{z}= (I-{^{*}\widetilde{P}^z})^{-1}$ $(\bar{d}_{z})$ is a meromorphic extension of $\bar{f}_{z}$ to $D$. It follows that ${^{*}\widetilde{\nu}^z} (\tilde{f}_{z})$ extends meromorphically ${^{*}\widetilde{\nu}^s} (\bar{f}_{s})$ to a possibly smaller domain $D' \subset D$ which contains $]0,\beta+1[$.  Also the equation $(1-k(s))$ ${^{*}\widetilde{\nu}^s} (\bar{f}_{s})={^{*}\widetilde{\nu}^s} (\bar{d}_{s})$ extends meromorphically to $D'$, with $k(z), {^{*}\widetilde{\nu}^z}$ defined as in the proof of Lemma A.2.  Here, since $\supp \mu$ is compact we have $s_{\infty}=\infty$. Since $k(s)>1$ for $s>\alpha$ and $\tilde{d}_{z}$ is holomorphic in $D'$, the function ${^{*}\widetilde{\nu}^z} (\tilde{f}_{z})=(1-k(z))^{-1} {^{*}\widetilde{\nu}^z} (\tilde{d}_{z})$ is holomorphic in a domain which contains the interval $[\alpha+\delta, \beta+1[$, hence $\displaystyle\mathop{\lim}_{s\rightarrow \beta_{-}} {^{*}\widetilde{\nu}^s}  (\bar{f}_{s})<\infty$. On the other hand, the meromorphy of $\tilde{f}_{z}=(I-{^{*}\widetilde{P}^z})^{-1} (\tilde{d}_{z})$ in $D$ implies
\[\tilde{f}_{z}=\displaystyle\mathop{\Sigma}_{1}^{m} \varphi_{j} (\beta-z)^{-j}+\psi_{z} \textrm{  with } \psi_{z}, \varphi_{j}\in C (\mathbb S^{d-1}), \varphi_{m}=\displaystyle\mathop{\lim}_{s\rightarrow \beta_{-}} (\beta-s)^{m} \bar{f}_{s}.\]
 Since $f_{s}(u)\geq 0$, we have $\varphi_{m}(u)\geq 0$. Also,
\[{^{*}\widetilde{\nu}^{\beta}} (\varphi_{m})=\displaystyle\mathop{\lim}_{s\rightarrow \beta_{-}} {^{*}\widetilde{\nu}^s} (\varphi_{m})=\displaystyle\mathop{\lim}_{s\rightarrow \beta_{-}} (\beta-s)^m {^{*}\widetilde{\nu}^s} (\bar{f}_{s}).\]
  From above we know that this limit is zero, hence ${^{*}\widetilde{\nu}^{\beta}} (\varphi_{m})=0$. Since $\varphi_{m}$ is non negative we get  $\varphi_{m} (u)=0$ ${^{*}\widetilde{\nu}^{\beta}}-\textrm{a.e}$ and the continuity of $\varphi_{m}$ implies $\varphi_{m}(u)=0$ for $u\in \supp {^{*}\widetilde{\nu}^{\beta}}=\mathbb S^{d-1}$. By induction we get $\varphi_{j}=0$ for any $j\geq 1$, hence the function $z\rightarrow \tilde{f}_{z}$ is holomorphic in a domain $D_{\varepsilon}$ which contains $]0, \beta+\varepsilon[$ for some $\varepsilon >0$ depending on the possible poles of $(I-{^{*}\widetilde{P}^z})^{-1} (\tilde{d}_{z})$, in $\mathbb R_{+}^{*}$.

  Then, as above, Landau's lemma gives that $\bar{f}_{s} (u)$ is finite for $s<\beta+\varepsilon$ and $u\in \mathbb S^{d-1}$. Hence $\mathbb E(\langle R, u \rangle  _{+}^{s})$ is finite for $u\in \mathbb S^{d-1}$, $s<\beta+\varepsilon$ and this gives the required contradiction. Then we have for $s>0$: $(1-k(s)) {^{*}\widetilde{\nu}^s} (\bar{f}_{s})={^{*}\widetilde{\nu}^s} (\bar{d}_{s})$. We observe also that $(I-{^{*}\widetilde{P}^z})^{-1} (\bar{d}_{z})$ is well defined and holomorphic in a domain which contains $]0,\infty[$.

  It follows for $s>0$: $(k(s) {^{*}\widetilde{\nu}^s} (\bar f_{s}))^{1/s}=({^{*}\widetilde{\nu}^s} (\bar f^1_{s}))^{1/s}$. Since $\bar f^1_{s}(u)=\mathbb E(\langle R-B,u \rangle^s_{+})$ we have for $s>1$,
\[{^{*}\widetilde{\nu}^s} (\bar f^1_{s}))^{1/s} \leq \mathbb E (|B|^s)^{1/s}+({^{*}\widetilde{\nu}^s} (\bar f_{s}))^{1/s},\ (k(s)^{1/s}-1) ({^{*}\widetilde{\nu}^s} (\bar f _{s})^{1/s} \leq \mathbb E (|B|^s)^{1/s}.\]
  Since $\lambda$ has compact support we have $\displaystyle\mathop{\lim}_{s\rightarrow \infty} \mathbb E (|B|^s)^{1/s}=d <\infty$; also Proposition \ref{pro:4.8} gives $\displaystyle\mathop{\lim}_{s\rightarrow \infty} k(s)^{1/s}=c>1$. Hence $\displaystyle\mathop{\lim}_{s\rightarrow \infty} ({^{*}\widetilde{\nu}^s} (\bar f_{s}))^{1/s}\leq (c-1)^{-1} d<\infty$.

  By definition of $\bar f_{s}$ it follows $\sup \{\langle R , u \rangle _{+};  (\omega,u)\in \widehat{\Omega}\times \widetilde{\Lambda}(T^{*})\} <\infty$; using Lemma 5.10, this contradicts the fact that $\langle R, u \rangle  _{+}$ is unbounded on $\widehat{\Omega} \times \widehat{\Lambda} (T^{*})$.

  In cases II, the above argument can easily be modified,  ${^{*}\widetilde{\nu}^{s}}$ replaced by ${^{*}\nu^{s}_{+}}$ and $\mathbb S^{d-1}$ by $\supp {^{*}\nu_{+}}$. Then we get that $\langle R, u \rangle  _{+}$ is bounded on $\widehat{\Omega} \times \Lambda_{+} (T^{*})$, which contradicts Lemma 5.10.
\end{proof}
\begin{rema}
  If $\supp\lambda$ is compact and $\mu$ has a density, Corollary 3.21  and the use of \cite{34}, allow us to avoid the use of the renewal theorem of section 4 and of Kac's formula in the proof of Theorem 5.2. Furthermore, if $\alpha\notin \mathbb N$, Theorem 5.2 then follows from the properties of Radon transforms of positive Radon measures (see \cite{5}). However, in the general case one needs to use Lemma 5.13 and Proposition 5.9. On the other hand, the density assumption on $\mu$ is not necessary for the validity of  Proposition A.4 as follows from Theorem 5.2.
\end{rema}








\end{document}